\documentclass[11pt, letterpaper]{amsart}
\usepackage{amsmath}
\usepackage{amssymb}
\usepackage[left=1in,right=1in,bottom=0.8in,top=0.8in]{geometry}
\usepackage{amsfonts}
\usepackage{graphicx}
\usepackage[font=small,labelfont=bf]{caption}
\usepackage{epstopdf}
\usepackage{hyperref}
\usepackage[dvipsnames]{xcolor}
\usepackage{amsthm}
\usepackage{float}
\usepackage{pgfplots}
\usepackage{listings}
\usepackage{longtable}
\usepackage{mathrsfs}
\usepackage{bbm}
\usepackage{euscript}
\usepackage{mathtools}
\usepackage{tikz}
\usepackage{comment}
\usepackage{tikz-cd}
\usepackage{amsxtra}
\usepackage{pifont}
\usepackage{amsbsy}
\usepackage{wasysym}
\usepackage{calligra}
\usepackage{url}
\usepackage{enumerate}
\usepackage[
backend=biber,
style=alphabetic,
]{biblatex}

\makeatletter
\renewcommand{\tocsection}[3]{%
  \indentlabel{\@ifnotempty{#2}{\bfseries\ignorespaces#1 #2.\quad}}\bfseries#3}
  \renewcommand{\tocsubsection}[3]{%
  \indentlabel{\hspace{2.5em}\@ifnotempty{#2}{\ignorespaces#1 #2.\quad}}#3}
\makeatother

\newtheorem{theorem}{Theorem}[section]
\newtheorem{corollary}[theorem]{Corollary}
\newtheorem{proposition}[theorem]{Proposition}
\newtheorem{lemma}[theorem]{Lemma}

\theoremstyle{definition}
\newtheorem{definition}[theorem]{Definition}
\newtheorem{construction}[theorem]{Construction}
\newtheorem{setup}[theorem]{Setup}

\newtheorem{notation}[theorem]{Notation}

\theoremstyle{remark}
\newtheorem{remark}[theorem]{Remark}

\definecolor{energy}{RGB}{114,0,172}
\definecolor{refblue}{RGB}{30,90,180}
\definecolor{freq}{RGB}{45,177,93}
\definecolor{spin}{RGB}{251,0,29}
\definecolor{signal}{RGB}{203,23,206}
\definecolor{circle}{RGB}{217,86,16}
\definecolor{average}{RGB}{203,23,206}

\colorlet{shadecolor}{gray!20}
\pgfplotsset{compat=1.9}

\usepgflibrary{fpu}

\newcommand{\RR}{\mathbb{R}}
\newcommand{\ZZ}{\mathbb{Z}}

\newcommand{\QQ}{\mathbb{Q}}

\newcommand{\PP}{\mathbb{P}}

\newcommand{\bA}{\mathbb{A}}

\newcommand{\FF}{\mathbb{F}}

\newcommand{\cB}{\mathcal{B}}
\newcommand{\cC}{\mathcal{C}}
\newcommand{\cD}{\mathcal{D}}

\newcommand{\cG}{\mathcal{G}}

\newcommand{\cK}{\mathcal{K}}

\newcommand{\cO}{\mathcal{O}}
\newcommand{\cP}{\mathcal{P}}
\newcommand{\cQ}{\mathcal{Q}}
\newcommand{\cR}{\mathcal{R}}

\newcommand{\cT}{\mathcal{T}}
\newcommand{\cV}{\mathcal{V}}
\newcommand{\cU}{\mathcal{U}}
\newcommand{\cW}{\mathcal{W}}

\newcommand{\fm}{\mathfrak{m}}

\newcommand{\bP}{\boldsymbol p}
\newcommand{\bx}{\boldsymbol x}

\newcommand{\Br}{\operatorname{Br}}

\newcommand{\cha}{\operatorname{char}}

\newcommand{\cores}{\operatorname{cores}}

\newcommand{\divs}{\operatorname{div}}

\newcommand{\cEnd}{\mathscr{E}\text{\kern -3pt {\calligra\large nd}}\,}
\newcommand{\cBun}{\mathscr{B}\text{\kern -3pt {\calligra\large un}}\,}
\newcommand{\cCoh}{\mathscr{C}\text{\kern -3pt {\calligra\large oh}}\,}

\newcommand{\et}{\operatorname{\acute{e}t}}

\newcommand{\Frac}{\operatorname{Frac}}

\newcommand{\cHom}{\mathscr{H}\text{\kern -3pt {\calligra\large om}}\,}
\newcommand{\id}{\operatorname{id}}

\newcommand{\im}{\operatorname{im}}

\newcommand{\nr}{\operatorname{nr}}

\newcommand{\ord}{\operatorname{ord}}

\newcommand{\pr}{\operatorname{pr}}

\newcommand{\res}{\operatorname{res}}

\newcommand{\ev}{\operatorname{ev}}

\newcommand{\Spec}{\operatorname{Spec}}

\newcommand{\Tor}{\operatorname{Tor}}

\newcommand{\one}{\mathbf{1}}
\newcommand{\bd}{\mathbf{d}}

\newcommand{\CH}{\operatorname{CH}}

\newcommand{\Supp}{\operatorname{Supp}}

\hypersetup{colorlinks=true,linkcolor=refblue,citecolor=refblue,urlcolor=refblue,
  pdftitle={Decomposition of the diagonal and complexes of unramified cohomology}}
\numberwithin{equation}{section}
\usepackage[T1]{fontenc}
\usepackage{lmodern,booktabs,microtype}

\newcommand{\floor}[1]{\left\lfloor#1\right\rfloor}
\newcommand{\ceil}[1]{\left\lceil#1\right\rceil}

\providecommand{\Tortwo}{\Tor^{\ZZ/2}}
\providecommand{\Torthree}{\Tor^{\ZZ/3}}
\hypersetup{
  pdftitle={Decompositions of the diagonal for hypersurfaces of small slope},
  pdfauthor={Jinting Liang, Zhiyu Tian, Shitan Xu}}

\makeatletter
\@ifundefined{U+wasy}{\input{uwasy.fd}}{}
\makeatother
\DeclareFontShape{U}{wasy}{b}{n}{<-> wasyb10}{}
\providecommand{\Ftwo}{\mathbb F_2}
\providecommand{\Nm}{\operatorname{Nm}}

\hypersetup{pdftitle={Decompositions of the diagonal of complete intersections in products of projective spaces}}
\title[Decompositions of the diagonal]{Decompositions of the diagonal of complete intersections in products of projective spaces}
\author[Jinting Liang]{Jinting Liang}
\address{Department of Mathematics\\University of British Columbia\\Vancouver, BC V6T 1Z2, Canada}
\email{liangj@math.ubc.ca}
\author[Shitan Xu]{Shitan Xu}
\address{Beijing International Center for Mathematical Research\\Peking University\\100871, Beijing, China}
\email{shitanxu@bicmr.pku.edu.cn}
\date{October 2026}

\begin{document}
\begin{abstract}
We study rationality of very general complete intersections in products of projective spaces. By constructing specific reference varieties and apply specialization method in the form of relative torsion orders, we extend \cite[Theorem~1.1]{Hotchkiss} and \cite[Theorem~1.1]{BohningBothmerPirutkaZhang2026} to very
general complete intersections in products in every characteristic.
In particular, we prove that very general quartic hypersurfaces of positive dimension admit no decomposition of the diagonal in every characteristic.
\end{abstract}
\maketitle
\section{Introduction}

Let $k$ be a field and let $X$ be a proper geometrically integral
$k$-variety of dimension $n$. Rationality problems ask how similar $X$
is to projective space. The following notions distinguish several
degrees of similarity:
\enlargethispage{3pt}
\begin{center}
\begin{tabular}{@{}l@{\;}c@{\;}p{0.62\textwidth}@{}}
$X$ is \emph{$k$-rational}
& : &
if $X$ is $k$-birational to $\PP^n_k$.
\\[7pt]

\multicolumn{1}{c}{\scalebox{1.5}{$\Downarrow$}}
& &
\\[7pt]

$X$ is \emph{stably $k$-rational}
& : &
if $X\times_k\PP^m_k$ is $k$-rational for some integer $m\geq0$.
\\[7pt]

\multicolumn{1}{c}{\scalebox{1.5}{$\Downarrow$}}
& &
\\[7pt]

$X$ is \emph{retract $k$-rational}
& : &
if there exist nonempty open subsets $U\subset X$ and
$V\subset\PP^N_k$, for some integer $N\geq n$, together with
$k$-morphisms
$U\xrightarrow{i}V\xrightarrow{r}U$, such that
$r\circ i=\id_U$.
\\[7pt]

\multicolumn{1}{c}{\scalebox{1.5}{$\Downarrow$}}
& &
\\[7pt]

$X$ is \emph{$k$-unirational}
& : &
if there exists a dominant rational map
$\PP^N_k\dashrightarrow X$ defined over $k$ for some integer $N\geq n$.
\\[7pt]

\multicolumn{1}{c}{\scalebox{1.5}{$\Downarrow$}}
& &
\\[7pt]

$X$ is \emph{rationally connected}
& : &
if there exists a dense open subset $W\subset X\times_k X$ such that,
for every algebraically closed extension $L/k$ and every pair
$(x_1,x_2)\in W(L)$, there is an $L$-morphism
$f:\PP^1_L\longrightarrow X_L$ such that
$f(0)=x_1$ and $f(\infty)=x_2$.
\end{tabular}
\end{center}

A main question in rationality asks whether the reverse implications
hold. We refer to \cite[Sections~1--2]{SchreiederSurvey2026} for a recent survey of
this question. Suppose now that $X$ is smooth and projective.
An \emph{(integral) decomposition of the diagonal} is a relation in
$\CH_n(X\times_kX)$ of the form
$$
 [\Delta_X]=z\times X+\Gamma,
$$
where $z$ is a zero-cycle of degree $1$ on $X$ and $\Gamma$ is supported on
$X\times E$ for a proper closed subset $E\subsetneq X$.
This is equivalent to \emph{universal $\CH_0$-triviality}: the degree
map $\deg:\CH_0(X_L)\to\ZZ$ is an isomorphism for every field
extension $L/k$, see \cite[Proposition~2.11]{VoisinSurvey} and
\cite[Section~7.1]{Schreieder}. The \emph{torsion order} $\Tor(X)$ is
the least positive integer $e$ for which
$e[\Delta_X]=z\times X+\Gamma$ with $\deg z=e$ and the same support
condition, and is $\infty$ if no such $e$ exists
\cite[Section~7.2]{Schreieder}. Thus $\Tor(X)=1$ exactly when $X$
admits an integral decomposition of the diagonal. Retract rationality
implies this property, so $\Tor(X)\ne1$ obstructs both retract and
stable rationality \cite[Lemma~7.5]{Schreieder}. Admitting no
decomposition of the diagonal also obstructs $\bA^1$-connectedness by
\cite[Theorem~2.4.3]{AM11} and \cite[Lemma~7.3]{Schreieder}.
We recall the generic-point criterion in Section~\ref{sec:preliminaries}.

We are particularly interested in rationality problems for hypersurfaces
and complete intersections in products of projective spaces. For a field $k$, let $k_0$ be its prime field. Fix positive integers $r$, $c$
and $N_j$ for $1\leq j\leq r$, and let $n:=\sum_jN_j-c>0$. Consider $P=\prod_{j=1}^r\PP^{N_j}_k$.
Fix homogeneous coordinates $x_{j0},\ldots,x_{jN_j}$ on each factor.
Let $\bd=(d_{ij})_{1\leq i\leq c,\,1\leq j\leq r}$ be a
$c\times r$ matrix of nonnegative integers.

\begin{definition}\label{def:intersection}
A \emph{complete intersection} of degree matrix $\bd$ and
codimension $c$ in $P$ is a nonempty closed subscheme $X\subset P$
of pure codimension $c$ which is the scheme-theoretic common zero
locus of sections $\{f_i\}_{1\leq i\leq c}$ for some
$f_i\in H^0\bigl(P,\cO_P(d_{i1},\ldots,d_{ir})\bigr)$, where $1\leq i\leq c$.
\end{definition}

Note that the nonemptiness and codimension requirements exclude rows whose
entries are all zero. Its defining hypersurfaces are all ample precisely when every $d_{ij}$ is positive. We use the following
projective version of the independent-coefficient convention of
\cite[Definition~2.4 and Lemma~2.5]{LS}.

\begin{definition}\label{def:general}
A complete intersection $X\subset P$ of degree matrix $\bd$
is \emph{very general} if it admits defining sections $f_1,\ldots,f_c$
as in Definition~\ref{def:intersection} such that, after
multiplying each $f_i$ by a nonzero scalar to make one fixed
monomial coefficient equal to $1$, all other coefficients of
all these sections are jointly algebraically independent over $k_0$.
\end{definition}

\begin{remark}\label{rmk:families}
Definition~\ref{def:general} describes the generic point of
$S_0=\prod_i\PP H^0(P_0,\cO_{P_0}(d_{i1},\ldots,d_{ir}))$, where
$P_0=\prod_j\PP^{N_j}_{k_0}$, as in
\cite[Definition~2.11 and Example~2.14]{LangeZhang2026}.
It is independent of the coefficients normalized to $1$.
We assume that $k$ has enough algebraically independent elements for
the very general members under consideration to exist.
Every nonempty very general complete intersection of the expected
codimension is smooth in every characteristic, as explained in
Section~\ref{subsec:specialization}.
For a geometrically integral member, the geometric relative torsion orders
with a fixed coordinate boundary are independent of the field of definition
by \cite[Lemma~3.6(e)]{LS}.
\end{remark}

The \emph{decomposition method} goes back to Bloch and Srinivas
\cite[Proposition~1]{BS}. Its specialization form, developed by Voisin, and later by
Colliot-Th\'el\`ene and Pirutka, turns a suitable singular special fiber into an obstruction for a very general
smooth member \cite[Theorem~2.1]{Voisin}, \cite[Theorem~1.14]{CTP}. Koll\'ar's degeneration to positive
characteristic gave linear degree bounds for nonrationality
\cite[Theorem~2]{Kollar}. Totaro combined inseparable covers with specialization
of the diagonal \cite[Theorem~2.1]{Totaro}. Schreieder obtained logarithmic bounds using unramified classes whose
pullbacks vanish over the singular locus of a special fiber after an
alteration \cite[Theorem~6]{SchQuadric},
\cite[Theorem~1.1 and Proposition~3.1]{SchSmall}.
His later torsion-order bounds use universal relations in Milnor
$K$-theory and apply in arbitrary characteristic
\cite[Theorem~1.1]{SchTorsion}.

In characteristic zero, motivic and tropical degenerations give another
approach through subdivisions of Newton polytopes
\cite[Theorem~3.15]{NO}, \cite[Theorem~3.16]{Moe}. Lange and Schreieder
introduced \emph{relative torsion orders}, which permit
specialization in nonproper families \cite[Section~3 and Theorem~4.3]{LS}.
There have been further breakthroughs recently. Ottem proved that a very general complex quartic sixfold
admits no decomposition of the diagonal \cite[Theorem~1.1]{Ottem2026}.
His proof uses an explicit quartic reference variety birational to a
quadric bundle of Colliot-Th\'el\`ene--Ojanguren type \cite{CTO}.
Very recently, Hotchkiss proved that a very general hypersurface
of degree at least five and any positive dimension admits no
decomposition of the diagonal over an uncountable algebraically
closed field of characteristic different from two
\cite[Theorem~1.1]{Hotchkiss}. B\"ohning, Graf von Bothmer, Pirutka and Zhang proved the corresponding
quartic theorem in every positive dimension in characteristic different
from two \cite[Theorem~1.1]{BohningBothmerPirutkaZhang2026}.

For complete intersections in projective space, Chatzistamatiou and
Levine refined Roitman's argument to obtain a factorial upper bound
for the torsion order in the Fano range, and used inseparable degenerations and
Hodge-Witt cohomology to obtain prime-power lower bounds in
characteristic zero \cite[Proposition~5.2 and Theorem~8.2]{CL}.
Lange and Zhang used relative torsion orders and iterated affine
degenerations of Schreieder's reference hypersurfaces to prove
logarithmic torsion-order bounds for very general complete intersections
and hypersurfaces in products of projective spaces
\cite[Theorems~1.1--1.3 and Section~3]{LangeZhang2026}.
Starting from the quadric surface bundle of Hassett, Pirutka and
Tschinkel \cite[Example~8 and Proposition~11]{HPTQuadric}, Fiammengo and L\"uders used the
cycle-theoretic degeneration methods of \cite[Theorem~1.2]{PS}, \cite[Theorem~4.3]{LS} to raise degree
and dimension and obtain relative diagonal obstructions in low
bidegree outside characteristic two
\cite[Theorems~1.1--1.3]{FiammengoLuders}.

In this paper, we consider construction to complete
intersections in products in every characteristic. Our first result focuses on products of projective spaces of equal dimension:

\begin{theorem}\label{thm:cubic}
Let $(X,\bd)\subset(\PP^N_k)^r$ be a positive-dimensional geometrically
integral very general complete intersection of ample hypersurfaces,
where $r\geq2$. If one row of $\bd$ is at least $(3,2,\ldots,2)$,
up to permutation, then $X_{\bar k}$ admits no decomposition of the
diagonal in every characteristic.
\end{theorem}

\begin{corollary}\label{cor:bidegrees}
Let $(X,\bd)$ be a very general hypersurface of bidegree $(d,e)$ in
$\PP^N_k\times\PP^N_k$ with $d,e\geq2$ and $d+e\geq5$. Then $X_{\bar{k}}$ admits no decomposition of the diagonal.
\end{corollary}

If there are at least three equal-dimensional factors in the product, we have the following stronger result in every characteristic.
\begin{theorem}\label{thm:quadratic}
Let $(X,\bd)\subset(\PP^N_k)^r$ be a positive-dimensional geometrically
integral very general complete intersection of ample hypersurfaces,
where $r\geq3$. If one row of $\bd$ is at least $(2,\ldots,2)$,
then $X_{\bar k}$ admits no decomposition of the diagonal in every
characteristic.
\end{theorem}

In particular, we have a criterion of rationality of a very general hypersurface in $(\PP^N_k)^r$ when $r\geq 3$.

\begin{corollary}\label{cor:rationality}
For a very general ample hypersurface $(X,\bd=(d_1,d_2,\cdots,d_r))$ in $(\PP^N_k)^r$, with $r\geq3$, the following are equivalent:
\begin{itemize}
    \item $X_{\bar{k}}$ is rational. 
    \item $X_{\bar{k}}$ is stably rational. 
    \item $X_{\bar{k}}$ admits a decomposition of the diagonal.
    \item $d_j=1$ for some $1\leq j\leq r.$ 
\end{itemize}
\end{corollary}

The next theorem focuses on products of projective spaces of unequal dimensions. 
\begin{theorem}\label{thm:uniform}
Let $r\geq1$ and $t\geq0$. Let $Q=\prod_{j=1}^r\PP^{a_j}_k$ and
$R=\prod_{\ell=1}^{t}\PP^{N_\ell}_k$ where $a_j$ and $N_\ell$ are positive integers.
%with $a_j,N_\ell\geq 1$. 
Let $X\subset Q\times R$ be a
positive-dimensional smooth geometrically integral very general complete
intersection of hypersurfaces of nonnegative multidegrees.
Suppose one defining row is at least
$(a_1+1,a_2+1,\ldots,a_r+1;3,\ldots,3)$,
and every other hypersurface is relatively ample over $Q$. Then $X_{\bar k}$ admits no integral decomposition of the diagonal if 
$$\dim X\leq \sum_ja_j+B_d-1,$$ where $d$ is the minimum of the degrees of $X$ in the factors of $R$ and $B_d$ is defined in Definition~\ref{def:norm-counting}. In particular, $X_{\bar k}$ is neither stably rational nor retract rational.
\end{theorem}

The characteristic-two construction also gives the following extension of the quartic theorem of
B\"ohning, Graf von Bothmer, Pirutka and Zhang
\cite[Theorem~1.1]{BohningBothmerPirutkaZhang2026}.
\begin{corollary}\label{cor:hypersurfaces}
A very general hypersurface $X\subset\PP^{n+1}$ of positive dimension $n$ and degree
$d\geq4$ admits no geometric decomposition of the diagonal in every
characteristic. 
\end{corollary}

The same affine constructions give a result for arbitrary products.
\begin{corollary}\label{cor:quartic-products}
Let $X\subset\prod_{j=1}^r\PP^{N_j}_k$ be a positive-dimensional
geometrically integral very general complete intersection of ample
hypersurfaces. If one defining row is at least $(4,\ldots,4)$, then
$X_{\bar k}$ admits no integral decomposition of the diagonal in every
characteristic.
\end{corollary}

The structure of this paper is as follows: 
Section~\ref{sec:preliminaries} recalls relative torsion orders and
establishes the specialization and comparison tools.
Section~\ref{sec:varieties} constructs the reference varieties.
Section~\ref{sec:char2-quartics} proves the quartic theorem in characteristic two.
Section~\ref{sec:bounds} derives coordinate bounds and their consequences. The proof of the results mentioned in the introduction is contained in Section~\ref{sec:proofs}.

\section*{Acknowledgements}\label{sec:ack}
 We would like to thank Zhiyu Tian for his invaluable discussions and support. SX thanks Lu Qi, Ting Gong, Jie Xu, Zheng Xu, and Qizheng Yin for useful discussions and encouragement during the preparation of this work. JL was supported in part by the Natural Sciences and Engineering Research Council of Canada.

\section*{AI disclosure}
ChatGPT and Codex assisted with the mathematical development and
preparation of this work. These AI tools suggested arguments and constructions, did calculations and reference checks, and assisted in auditing the proofs. The authors rewrote the paper to make it more accessible. The authors take full responsibility for the mathematical correctness of the paper.

\section{Relative torsion order and specialization}\label{sec:preliminaries}
We recall torsion orders and the relative torsion orders of
\cite[Section~3]{LS}. We then explain how an affine generic-point
obstruction specializes to very general complete intersections.

\subsection{Torsion orders and relative torsion orders}\label{subsec:orders}
We express decomposition of the diagonal in terms of generic-point
classes and recall the specialization properties used below.
A variety is an integral separated scheme of finite type over a field.
For an algebraic scheme $X$ and a commutative ring $\Lambda$ with identity,
let $\CH_i(X,\Lambda)=\CH_i(X)\otimes_{\ZZ}\Lambda$.
Let $\delta_X\in\CH_0(X_{k(X)})$ be the class of the diagonal $k(X)$-point.
We use the same notation for its image with coefficients $\Lambda$.

\begin{definition}[{\cite[Definition~7.9]{Schreieder}}]\label{def:torsion}
Let $X$ be a proper variety over $k$. Its \emph{torsion order} $\Tor(X)$
is the least integer $e\geq1$ such that $e\delta_X=z_{k(X)}$ in
$\CH_0(X_{k(X)})$ for some $z\in\CH_0(X)$.
If no such integer exists, let $\Tor(X)=\infty$.
\end{definition}

\begin{definition}[{\cite[Definition~3.1]{LS}}]\label{def:decomposition}
Let $X$ be a variety over $k$ and let $W\subsetneq X$ be closed.
We say that $X$ admits a \emph{$\Lambda$-decomposition of the diagonal
relative to $W$} if
$$
\delta_X\in\im\bigl(\CH_0(W_{k(X)},\Lambda)
\longrightarrow\CH_0(X_{k(X)},\Lambda)\bigr).
$$
\end{definition}

Let $U=X\setminus W$. Localization \cite[Section~1.8]{Fulton} gives
$$
\CH_0(W_{k(X)},\Lambda)\longrightarrow\CH_0(X_{k(X)},\Lambda)
\longrightarrow\CH_0(U_{k(X)},\Lambda)\longrightarrow0.
$$
Thus Definition~\ref{def:decomposition} is equivalent to
$\delta_U=0$ in $\CH_0(U_{k(X)},\Lambda)$.

\begin{definition}[{\cite[Definition~3.3]{LS}}]\label{def:relative}
The \emph{relative torsion order} is
$$
\Tor^\Lambda(X,W)=\ord\bigl(\delta_U\in\CH_0(U_{k(X)},\Lambda)\bigr),
$$
with value $\infty$ if this class has infinite order.
\end{definition}

\begin{remark}\label{rmk:notation}
For proper $X$, localization and restriction to the generic point of the
second factor express the relation in Definition~\ref{def:torsion} as
$$
e[\Delta_X]=z\times X+\Gamma
$$
in $\CH_{\dim X}(X\times_kX)$, where
$\Supp\Gamma\subset X\times D$ for a proper closed subset $D\subsetneq X$.
Thus $\Tor(X)=1$ precisely when $X$ admits a decomposition of the diagonal.

Since $k(U)=k(X)$, we write
$\Tor^\Lambda(U)=\Tor^\Lambda(U,\varnothing)=\Tor^\Lambda(X,W)$.
This depends only on $U$. The notation also makes sense for an algebraic
scheme $X$ whose complement $U$ is nonempty and integral
\cite[Remarks~3.4--3.5]{LS}.
\end{remark}

Enlarging the boundary gives
$\Tor^\Lambda(X,W)\mid\Tor^\Lambda(X,W')$ when $W'\subset W$.
Over an algebraically closed ground field, relative torsion orders are
invariant under every field extension \cite[Lemma~3.6(b), (e)]{LS}.
We also use \cite[Lemma~3.8]{LS}: in a separated flat family of finite type
over a DVR with a closed boundary, if both complementary fibers are
nonempty and geometrically integral, the geometric relative torsion order
of the special fiber divides that of the generic fiber. These facts hold for
every commutative coefficient ring, including coefficients of the same
characteristic as the ground field.

\subsection{Specialization}\label{subsec:specialization}
We first recall the prime-to-characteristic pairing argument of
\cite[Proposition~3.1]{SchSmall}. We then construct the affine
specializations used in the main proofs. See also
\cite[Theorem~8.6 and Remark~8.9]{Schreieder} and
\cite[proof of Theorem~6.1]{LS}.

\setcounter{theorem}{5}
\begin{lemma}\label{lem:pairing}
Let $k$ be algebraically closed and let $m$ be a prime different from
$\cha k$. Let $Z$ be a proper integral variety with a smooth open subset
$U=Z\setminus W$. Let $\pi:Y\to Z$ be a proper dominant generically finite
morphism, where $Y$ is smooth, proper and integral, and
$[k(Y):k(Z)]$ is coprime to $m$.
Suppose that $0\ne\gamma\in H^q_{\nr}(k(Z)/k,\mu_m^{\otimes q})$ for some
$q>0$ and that $(\pi^*\gamma)|_{\kappa(y)}=0$ for every
$y\in\pi^{-1}(W)$. Then $\Tor^{\ZZ/m}(U)=m$.
\end{lemma}
\begin{proof}
Let $F=k(Z)$, $V=\pi^{-1}(U)$ and $D=Y\setminus V$.
The restrictions of $\pi^*\gamma$ to scheme points of $Y$ are defined by
\cite[Proposition~4.8 and Corollary~4.9]{Schreieder}.
Suppose that $\delta_U=0$ in $\CH_0(U_F,\ZZ/m)$.
The morphism $f=\pi|_V:V\to U$ between smooth varieties is a local
complete intersection morphism of virtual relative dimension zero, so
its Gysin pullback $f_F^!$ is defined \cite[Section~6.6]{Fulton}.
Near the generic point of $U$, the morphism is finite flat.
Thus $f_F^!(\delta_U)$ is the class of
$\xi=\Spec k(Y)\subset V_F$ with multiplicity one
\cite[Lemma~42.59.5 (\href{https://stacks.math.columbia.edu/tag/0FF4}{Tag~0FF4})]{stacks-project}.
Localization then gives $[\xi]=i_*z$ in $\CH_0(Y_F,\ZZ/m)$ for a zero-cycle
$z$ on $D_F$, where $i:D_F\hookrightarrow Y_F$.

The Merkurjev pairing \cite[Corollary~5.2]{Schreieder} gives
$$
\Phi:\CH_0(Y_F,\ZZ/m)\longrightarrow H^q(F,\mu_m^{\otimes q}),
\qquad [P]\longmapsto
\cores_{\kappa(P)/F}\bigl(((\pi^*\gamma)_F)|_{\kappa(P)}\bigr).
$$
Every closed point of $D_F$ maps to a scheme point of $D$, so the
hypothesis and compatibility with residue-field extensions give
$\Phi(i_*z)=0$. On the other hand,
$$
\Phi([\xi])=\cores_{k(Y)/F}\res_{k(Y)/F}(\gamma)
=[k(Y):F]\gamma\ne0.
$$
This contradiction proves $\delta_U\ne0$, whose order is $m$.
\end{proof}

\begin{remark}\label{rmk:boundary}
In Lemma~\ref{lem:pairing}, suppose that
$(\pi^{-1}(W))_{\mathrm{red}}$ has smooth irreducible components $B_i$.
It suffices that $(\pi^*\gamma)|_{k(B_i)}=0$ for every $i$.
Then \cite[Proposition~4.8 and Corollary~4.9]{Schreieder} give vanishing
at every scheme point of the boundary.
\end{remark}

\begin{lemma}\label{lem:comparison}
Let $X$ be a smooth projective variety of positive dimension over an
algebraically closed field, and let $\varnothing\ne W\subsetneq X$ be closed.
If $\Tor(X)<\infty$, then $\CH_0(X)\simeq\ZZ$ and, for every $m\geq2$,
$$
\Tor^{\ZZ/m}(X,W)\mid\Tor^\ZZ(X,W)\mid\Tor(X).
$$
In particular, a decomposition of the diagonal implies
$\Tor^\ZZ(X,W)=1$.
\end{lemma}
\begin{proof}
Let $e=\Tor(X)$. Acting with a diagonal relation on degree-zero cycles
and moving their support away from its proper closed subset gives
$e\CH_0(X)^0=0$.
The group $\CH_0(X)^0$ is generated by pushforwards of degree-zero
classes on smooth projective curves. Their Jacobians have divisible
groups of $k$-points in every characteristic
\cite[Proposition~39.9.11(4) (\href{https://stacks.math.columbia.edu/tag/03RP}{Tag~03RP})]{stacks-project}.
Thus $\CH_0(X)^0$ is divisible and killed by $e$, so it is zero.
Choose $p\in W(k)$.
By \cite[Lemma~3.6(a), (b), (d)]{LS},
$$
\Tor^{\ZZ/m}(X,W)\mid\Tor^\ZZ(X,W)
\mid\Tor^\ZZ(X,\{p\})=\Tor(X).
$$
\end{proof}

We now specialize the ordinary generic-point class of an affine model.
For $P=\prod_{j=1}^r\PP^{N_j}_k$ and $E\subseteq\{1,\ldots,r\}$, let
$D_E=\bigcup_{j\in E}\{x_{j0}=0\}$. Let
$P_E^\circ=P\setminus D_E=\prod_{j\in E}\bA^{N_j}\times
\prod_{j\notin E}\PP^{N_j}$.

\begin{theorem}\label{thm:transfer}
Let $(X,\bd)\subset P$ be a geometrically integral very general complete
intersection of codimension $c$ and degree matrix $\bd=(d_{ij})$ over $k$.
Fix $E\subseteq\{1,\ldots,r\}$ such that $U=X\cap P_E^\circ$ is nonempty.
Let $\Omega$ be algebraically closed of the same characteristic as $k$.
Suppose that sections
$h_i\in H^0(P_\Omega,\cO_{P_\Omega}(d_{i1},\ldots,d_{ir}))$, for
$1\leq i\leq c$, have a nonempty integral scheme-theoretic common zero
locus $V=Z(h_1,\ldots,h_c)\cap(P_E^\circ)_\Omega$ of codimension $c$.
Then, for every commutative coefficient ring $\Lambda$ with identity,
$$
\Tor^\Lambda(V)\mid
\Tor^\Lambda\bigl(X_{\bar k},(X\cap D_E)_{\bar k}\bigr)
\mid\Tor^\Lambda(X,X\cap D_E).
$$
\end{theorem}
\begin{proof}
At every point of $V$, the equations $h_1,\ldots,h_c$ generate an ideal
of height $c$ in a regular local ring. They form a regular sequence
\cite[Lemma~10.104.2 (\href{https://stacks.math.columbia.edu/tag/02JN}{Tag~02JN})]{stacks-project}.
Let $\{M_{i\alpha}\}_\alpha$ be the full monomial basis in row $i$.
Write $h_i=\sum_\alpha h_{i\alpha}M_{i\alpha}$ and let
$G_i=\sum_\alpha a_{i\alpha}M_{i\alpha}$, with all $a_{i\alpha}$ jointly
algebraically independent over $\Omega$.
Let $L=\Omega(\{a_{i\alpha}\})$ and let $t$ be transcendental over $L$.
Let $R=L[t]_{(t)}$ and $K=L(t)$. Consider
$\cV=Z(h_1+tG_1,\ldots,h_c+tG_c)\subset(P_E^\circ)_R$.
Its special fiber is $V_L$, which is geometrically integral because
$\Omega$ is algebraically closed.

At a point of the special fiber, let $A$ be the ambient local ring and
write $F_i=h_i+tG_i$ after trivializing the line bundles.
The sequence $t,F_1,\ldots,F_c$ is regular. Permuting it shows that $t$
is a nonzerodivisor in $A/(F_1,\ldots,F_c)$
\cite[Lemma~10.68.4 (\href{https://stacks.math.columbia.edu/tag/00LJ}{Tag~00LJ})]{stacks-project}. This quotient is torsion-free over
$R$, hence flat \cite[Lemma~15.22.11(1) (\href{https://stacks.math.columbia.edu/tag/0AUW}{Tag~0AUW})]{stacks-project}.
Flatness on the generic fiber is immediate, so $\cV$ is flat over $R$.

Let $b_{i\alpha}=h_{i\alpha}+ta_{i\alpha}$.
Since $a_{i\alpha}=(b_{i\alpha}-h_{i\alpha})/t$, these coefficients are
jointly algebraically independent over $\Omega(t)$.
Normalize each row by the coefficient $b_{i0}$ of the same monomial
used to normalize the equations of $X$. The ratios
$\beta_{i\alpha}=b_{i\alpha}/b_{i0}$, for $\alpha\ne0$, are also jointly
algebraically independent over $\Omega(t)$.

Let $F=k_0(z_{i\alpha}:1\leq i\leq c,\ \alpha\ne0)$.
Let $P_0=\prod_j\PP^{N_j}_{k_0}$. Over $F$, let
$U_{\mathrm{gen}}\subset((P_0)_E^\circ)_F$ be defined by
$M_{i0}+\sum_{\alpha\ne0}z_{i\alpha}M_{i\alpha}=0$ for $1\leq i\leq c$.
The normalized coefficients of $X$ give an embedding $F\hookrightarrow k$,
and $z_{i\alpha}\mapsto\beta_{i\alpha}$ gives an embedding $F\hookrightarrow K$.
Thus $U$ and $\cV_K$ are base changes of $U_{\mathrm{gen}}$.
A common algebraically closed extension of $\bar k$ and $\bar K$ over $F$
compares the two geometric fibers. In particular, $\cV_K$ is nonempty
and geometrically integral. Applying \cite[Lemmas~3.6(e) and~3.8]{LS}
to this separated flat family with empty boundary gives
$$
\Tor^\Lambda(V)=\Tor^\Lambda(V_{\bar L})
\mid\Tor^\Lambda(\cV_{\bar K})=\Tor^\Lambda(U_{\bar k}).
$$
Finally, field extension sends $\delta_U$ to $\delta_{U_{\bar k}}$, so
$\Tor^\Lambda(U_{\bar k})\mid\Tor^\Lambda(U)$.
\end{proof}

We first justify the smoothness of the very general complete intersections
used below. For a nonempty set $I$ of rows, let $r(I)$ be the dimension
of the product of factors on which at least one row in $I$ has positive
degree. In a local trivialization, the value and first derivatives of
each section in these factors are independent coefficient parameters
where that row has positive degree. The constant and coordinate
monomials on the affine chart give these first jets in every characteristic.
For a fixed point of this product and a projective relation
$\sum_{i\in I}\lambda_i df_i=0$ with all $\lambda_i\ne0$, the equations
$f_i=0$ and this relation give $|I|+r(I)$ independent linear conditions
on the coefficients. The point and the relation vary in dimension
$r(I)+|I|-1$. Their incidence image lies in a proper closed subset of
the coefficient space. There are finitely many sets $I$, so very general
coefficients admit no dependence among the differentials at any common
zero. The intersection is then smooth of codimension $c$ wherever nonempty.

We next place a model of codimension $q\leq c$ in the full affine chart
by adding $c-q$ coordinate equations. For a monomial $M$, let $\deg_jM$
be its total degree in the $j$th coordinate group.

\begin{corollary}\label{cor:embedding}
Let $X\subset P=\prod_{j=1}^r\PP^{N_j}_k$ be a geometrically integral
very general complete intersection of positive dimension
$n=\sum_jN_j-c$, with degree matrix $\bd=(d_{ij})$.
Let $D=D_{\{1,\ldots,r\}}$ and suppose that $X\cap(P\setminus D)\ne\varnothing$.
Let $\Omega$ be algebraically closed of the same characteristic as $k$.
Let $V=Z(g_1,\ldots,g_q)\subset\bA^{n+q}_\Omega$ be a nonempty integral
affine complete intersection of codimension $q$, where $0\leq q\leq c$.
Suppose the following placement is possible.
\begin{enumerate}
\item Assign the variables of $V$ injectively to the slots
$u_{j\ell}=x_{j\ell}/x_{j0}$, for $1\leq j\leq r$ and $1\leq\ell\leq N_j$.
\item Match the equations $g_s$ to distinct rows $i_s$ so that
$\deg_jM\leq d_{i_sj}$ for every monomial $M$ of $g_s$ and every factor $j$.
\item Match each other row $i$ to a distinct unused slot
$u_{j(i),\ell(i)}$ with $d_{i,j(i)}>0$.
\end{enumerate}
Then, for every commutative coefficient ring $\Lambda$ with identity,
$$
\Tor^\Lambda(V)\mid\Tor^\Lambda(X_{\bar k},(X\cap D)_{\bar k})
\mid\Tor^\Lambda(X,X\cap D).
$$
If $\Tor^{\ZZ/m}(V)=m$ for some prime $m$, then $X_{\bar k}$ admits no
decomposition of the diagonal.
\end{corollary}
\begin{proof}
There are $n+c$ affine slots, of which $n+q$ are used by $V$.
Multihomogenize each $g_s$ to row $i_s$ by multiplying each monomial $M$
by $\prod_jx_{j0}^{d_{i_sj}-\deg_jM}$.
For each other row $i$, use the section
$x_{j(i),\ell(i)}x_{j(i),0}^{d_{i,j(i)}-1}\prod_{h\ne j(i)}x_{h0}^{d_{ih}}$.
All exponents are nonnegative. On the affine chart the equations become
$g_1,\ldots,g_q$ and the $c-q$ distinct unused coordinates, so their
scheme-theoretic zero locus is isomorphic to $V$ and has codimension $c$.
Theorem~\ref{thm:transfer} gives the divisibilities.

The divisor $D$ belongs to the ample system $\cO_P(1,\ldots,1)$.
Thus $X\cap D$ is nonempty and proper because $\dim X>0$ and
$X\cap(P\setminus D)\ne\varnothing$.
The very general complete intersection is smooth by the argument above.
Lemma~\ref{lem:comparison} then gives the conclusion.
\end{proof}

The next lemma gives affine generic-point obstructions from top forms.
\begin{lemma}\label{lem:forms}
Let $Y$ be a smooth projective variety of dimension $n>0$ over an
algebraically closed field $k$. Let $D\subset Y$ be an effective divisor
whose irreducible components are smooth.
If $H^0(Y,\omega_Y)\ne0$, then $\Tor^\QQ(Y,D)=\infty$ in characteristic
zero, and $\Tor^{\ZZ/p}(Y,D)=p$ in characteristic $p>0$.
For every $e\geq n+2$, a smooth degree-$e$ hypersurface in $\PP^{n+1}_k$
has an affine hyperplane complement with this property.
\end{lemma}
\begin{proof}
Let $U=Y\setminus D$. Use coefficients $\Lambda=\QQ$ in characteristic
zero and $\Lambda=\ZZ/p$ in characteristic $p>0$.
If $\delta_U=0$, localization and spreading out, followed by transposition,
give
$$
[\Delta_Y]=A+B
$$
in $\CH_n(Y\times Y,\Lambda)$, where $A$ is supported on $Y\times D$
and $B$ on $E\times Y$ for a proper closed subset $E\subsetneq Y$.
Let $D_i$ be the irreducible components of $D$.
By the action of correspondences on forms in \cite[Lemma~2.2]{Totaro},
the action of $A$ on $H^0(Y,\omega_Y)$ factors through
$\bigoplus_iH^0(D_i,\Omega^n_{D_i})=0$.
The action of $B$ restricts to zero on $Y\setminus E$, and restriction
of global top forms to this open subset is injective. Both actions
vanish, contradicting the identity action of the diagonal.
Thus $\delta_U\ne0$, of infinite order over $\QQ$ and order $p$ over
$\ZZ/p$. The last statement follows from
$\omega_Y\simeq\cO_Y(e-n-2)$ and Bertini's theorem
\cite[Lemma~33.47.3 (\href{https://stacks.math.columbia.edu/tag/0FD6}{Tag~0FD6})]{stacks-project}.
\end{proof}

\setcounter{theorem}{12}
\begin{remark}\label{rmk:graphs}
Adjoin coordinates $T_1,\ldots,T_h$. For $1\leq\nu\leq h$, use the monic
equation $T_\nu=P_\nu(z_1,\ldots,z_K,\allowbreak T_1,\ldots,T_{\nu-1})$.
Elimination gives an isomorphism of the new coordinate ring with the
original coordinate ring, so the affine scheme is the same.
If the original scheme is an integral complete intersection of codimension
$q$, the new one is an integral complete intersection of codimension $q+h$.
Indeed, its ideal has height $q+h$ and has $q+h$ generators in the regular
ambient ring \cite[Lemma~10.104.2 (\href{https://stacks.math.columbia.edu/tag/02JN}{Tag~02JN})]{stacks-project}. Its generic-point class and
torsion orders are preserved by the isomorphism.

We call each $T_\nu$ a \emph{graph coordinate}. A \emph{coordinate copy}
is the case $T_\nu=z_i$. Different occurrences of $z_i$ may then belong
to different coordinate groups. Every graph equation and every modified
equation must satisfy the group-degree bounds of Corollary~\ref{cor:embedding}.
\end{remark}

We end with a specialization criterion for forms after extending constants.
A \emph{form} is a homogeneous polynomial of positive degree. It is
\emph{isotropic} if it vanishes at a nonzero vector, and \emph{anisotropic} otherwise
\cite[Definition~1.1]{KrashenMatzri}.
\begin{lemma}\label{lem:constant}
Let $\Omega$ be algebraically closed and let
$F\in\Omega[t,x_1,\ldots,x_N]$.
Suppose that $F_0=F|_{t=0}$ is irreducible and
$\deg_{\bx}F=\deg_{\bx}F_0$.
Let $K=\Omega(t)$, $V=Z(F)\subset\bA^N_K$ and
$V_0=Z(F_0)\subset\bA^N_\Omega$.
Then $V$ is geometrically integral. If a form $q$ over $K(V)$ has
coefficients regular at the generic point of $V_0$ and anisotropic
reduction over $\Omega(V_0)$, then $q$ remains anisotropic over
$\bar K(V_{\bar K})$.
\end{lemma}
\begin{proof}
Suppose that $F$ factors over a finite extension $K'/K$.
Let $R'$ be a DVR of $K'$ above $t=0$. Its residue field is $\Omega$.
Gauss's lemma gives primitive factors over $R'$ of positive degrees.
Their reductions have the same degrees, because
$\deg_{\bx}F=\deg_{\bx}F_0$. This contradicts the irreducibility of $F_0$,
so $V$ is geometrically integral.

Let $\cV=Z(F)\subset\bA^N_{\Omega[t]_{(t)}}$.
After base change to any such $R'$, the local ring at the generic point
of the special fiber is a one-dimensional local domain whose maximal
ideal is generated by the uniformizer of $R'$. It is a DVR with fraction
field $K'(V_{K'})$ and residue field $\Omega(V_0)$.
The coefficients of $q$ are regular there.
A nonzero isotropic vector can be scaled so that its coordinates are
integral and at least one is a unit. Reduction then contradicts the
anisotropy of the reduced form. Every isotropic vector over
$\bar K(V_{\bar K})$ is defined over a finite extension $K'/K$, which
proves the claim.
\end{proof}

\section{Reference varieties}\label{sec:varieties}
We construct affine varieties with nonzero generic-point classes for
specialization. Throughout this section, $k$ is algebraically closed and
valuations are trivial on the indicated ground field unless another
convention is given. We use $\Omega$
when another algebraically closed constant field is needed. In
Sections~\ref{sec:bounds} and~\ref{sec:proofs}, $k$ again denotes the field
of the very general complete intersection.

\begin{definition}\label{def:reference}
Let $k$ be algebraically closed and let $\Lambda$ be a nonzero commutative
ring with identity. In this paper, a \emph{reference variety} with
coefficients $\Lambda$ is a positive-dimensional integral affine
$k$-variety $V$ such that $\delta_V\ne0$ in
$\CH_0(V_{k(V)},\Lambda)$.
\end{definition}

For a proper compactification $Z\supset V$ with boundary $W=Z\setminus V$,
the boundary is nonempty and localization gives
$\Tor^\Lambda(V)=\Tor^\Lambda(Z,W)$.
Thus this definition concerns ordinary Chow groups of the affine variety.
When $Z$ is smooth and projective, Lemma~\ref{lem:comparison} shows that
$\delta_V\ne0$ excludes a decomposition of the diagonal on $Z$.
We will use the same implication after specialization to a smooth
projective member.

\setcounter{theorem}{4}
\begin{definition}\label{def:profile}
A \emph{profile} $(C_1,\ldots,C_r;d_1,\ldots,d_r)$ is a partition of
the affine coordinates into groups of sizes $C_j$, with total degree
at most $d_j$ in the $j$th group for each defining equation.
\end{definition}

For an equation with base coordinates $\bx$ and other coordinates
$\boldsymbol y$, clear denominators in $k[\bx]$ and write the nonzero
polynomial as $G=\sum_\alpha a_\alpha(\bx)\boldsymbol y^\alpha$.
Let $c=\gcd_\alpha a_\alpha$ in $k[\bx]$.
Its \emph{primitive numerator} is $G/c$, defined up to a nonzero constant.
Dividing by $c$ preserves the generic fiber over $k(\bx)$ and does not
increase any total or group degree. All statements below about closed points or
boundary points concern the affine variety defined by this primitive
numerator.

\subsection{Norm forms and valuation criteria}\label{subsec:norms}

In this subsection, we first develop the symbol, form, and evaluation tools used
throughout the constructions. Then we prove our main criterion
Theorem~\ref{thm:norm-criterion} using convex geometry. Applying the criterion to Pfister-form substitutions will provide the degree bounds in Section~\ref{subsec:degree-bounds}. We start with the notations for coordinates:
\begin{notation}\label{not:coordinates}

\begin{itemize}
\item For $s\geq1$, write $[s]=\{1,\ldots,s\}$.
For a tuple $\bx=(x_1,\ldots,x_s)$ and $I\subseteq[s]$, we use the lowercase letter to denote the product $x_I=\prod_{i\in I}x_i,  x_\varnothing=1.$ 
\item An explicitly declared family $(X_I)_{I\subseteq[s]}$ denotes $2^s$ independent
coordinates, one for each subset $I$. We also use $Z_I,P_I,Q_I,R_I$
and explicitly declared lowercase indexed families in this way.
For a named tuple of base variables, such as $(z_1,\ldots,z_s)$,
the notation $z_I$ instead denotes $\prod_{i\in I}z_i$.
\item We write $e_I$ for the corresponding standard basis vector, so that $\sum_I X_Ie_I$ is a vector in $F^{2^s}.$
\item When a subset occurs as a subscript, we simply omit the notation of a set, e.g. $e_0$ means $e_\varnothing$, $e_{i,j}$ means $e_{\{i,j\}}$, and $z_{i,j,k}$ denotes $z_{\{i,j,k\}}$.
\item For $\nu\in\ZZ^s$ denote $\bx^\nu=\prod_i x_i^{\nu_i}$.
\item The \emph{support} of a nonzero Laurent polynomial
$P=\sum_\nu c_\nu\bx^\nu$ is $\operatorname{supp}(P)=\{\nu:c_\nu\ne0\}$.
\item For a condition $\cC$, let $\delta_{\cC}$ be $1$ if
$\cC$ holds and $0$ otherwise.
\end{itemize}
\end{notation}
\begin{definition}\label{def:kummer-symbols}
Let $m\geq2$ be invertible in a field $F$.
Write $\mu_m$ for the group of $m$th roots of unity and
$H^q(F,\mu_m^{\otimes r})$ for Galois cohomology \cite[Section~2]{Schreieder}.
The \emph{Kummer class} of $a\in F^\times$ is its image under $F^\times/F^{\times m}\xrightarrow{\ \sim\ }H^1(F,\mu_m)$.
We denote it by $(a)_m$ \cite[Proposition~4.3.6]{GilleSzamuely}.
The \emph{symbol} \cite[Equation~(2.2)]{Schreieder} associated with $a_1,\ldots,a_s\in F^\times$ is $(a_1,\ldots,a_s)_m =(a_1)_m\cup\cdots\cup(a_s)_m \in H^s(F,\mu_m^{\otimes s})$.
We omit the subscript when $m=2$.
The entries $a_i$ are also called the \emph{slots} of the symbol.

\end{definition}
\begin{remark}\label{rmk:unramified}
    For a field extension $F'/F$ we write $\alpha_{F'}$ or
$\alpha|_{F'}$ for restriction \cite[Section~3.1]{Schreieder}. An extension \emph{splits} a symbol if its restriction there is zero. For a discrete valuation $\mu$ with $m$ invertible in its residue
field, the residue map has the form \cite[Equation~(3.2)]{Schreieder} $\partial_\mu:H^q(F,\mu_m^{\otimes r}) \longrightarrow H^{q-1}(\kappa(\mu),\mu_m^{\otimes(r-1)})$.
A class is \emph{unramified at $\mu$} if its residue is zero \cite[Definition~4.3 and Remark~4.4]{Schreieder}.
The \emph{specialization} $s_\mu(\alpha)$ of an unramified class is the restriction to the
residue field of its unramified extension over the henselization of $\mu$ \cite[Construction~6.8.6 and Corollary~6.8.8]{GilleSzamuely} and \cite[Expos\'e~X, Lemma~2.2.1]{SGA4III}.
\end{remark}
We use restriction, corestriction, residues and specialization with
the conventions of \cite[Sections~2--4]{Schreieder}.
Corestriction is understood for every finite field extension $F'/F$,
including inseparable extensions, with $\cores_{F'/F}\circ\res_{F'/F}=[F':F]\id$
\cite[Remark~1.11 and rule~R2d]{Rost1996}.
After fixing a primitive $m$th root of unity, we also use the
induced identification $H^2(F,\mu_m^{\otimes2})\simeq\Br(F)[m]$
when writing a cyclic symbol as a Brauer class \cite[Corollary~4.4.9]{GilleSzamuely}.

\begin{definition}\label{def:kato-symbols}
Suppose $\cha F=2$ and let $\wp(b)=b^2+b$.
The \emph{Artin--Schreier class} of $a\in F$ is
$[a]\in F/\wp(F)=H^1_{\et}(F,\ZZ/2)$ \cite[Section~2.1]{HuSun}.
For $z_1,\ldots,z_s\in F^\times$, we use the \emph{Kato symbol} \cite[Section~3.1]{HuSun} $[a;z_1,\ldots,z_s]\in H^{s+1,s}(F) :=H^1_{\et}(F,\Omega^s_{F,\log})$,
represented by $a\,d\log z_1\wedge\cdots\wedge d\log z_s$ in
Kato's differential presentation \cite[Section~4, before Lemma~4.1]{Otabe}. Here $\Omega^s_{F,\log}$ is the
\'etale subsheaf of differential forms generated by wedges
of logarithmic differentials $d\log u=du/u$ of units \cite[Section~2]{Otabe}.
\end{definition}
\begin{remark}\label{rmk:milnor}
We also use \emph{Milnor $K$-theory modulo two}, $\mathrm k_r^M(F)=K_r^M(F)/2$.
Milnor $K$-theory is the graded tensor algebra on the multiplicative
abelian group $F^\times$, modulo the Steinberg relations
$\{b,1-b\}=0$ for $b\ne0,1$ \cite[Section~1]{Rost1996}. We use
braces $\{b_1,\ldots,b_r\}$ to denote its symbols modulo two.
In characteristic two, the logarithmic differential map sends $\{b_1,\ldots,b_r\}\longmapsto d\log b_1\wedge\cdots\wedge d\log b_r$ in $\Omega_F^r$
\cite[Section~2]{Otabe}. 
\end{remark}

The next definition collects the notions we need for forms:
\begin{definition}\label{def:forms}
\begin{itemize}
\item The notation $q\perp q'$ denotes the sum on separate coordinates,
and $q\subset q'$ means that an injective linear map realizes $q$ as a restriction of $q'$.
 \item We say that $q$ \emph{represents} $b$ if $q(v)=b$ for some vector $v$.
\item An \emph{isometry} is an invertible linear map preserving a form. 
We write $q\simeq q'$ for isometric forms $q$ and $q'$.
\item A form $q$ is \emph{round} if every nonzero represented value
$b$ satisfies $bq\simeq q$ \cite[Section~9.A]{EKM}.
\item For a quadratic form $q$, its \emph{polar form} is
$b_q(u,v):=q(u+v)-q(u)-q(v)$ \cite[Definition~7.1]{EKM}. A nonsingular quadratic form has nondegenerate polar form. The \emph{dimension} of a quadratic form is the dimension of its underlying vector space.
\item A \emph{hyperbolic plane} is a quadratic space isometric to
$(F^2,h)$, where $h:F^2\longrightarrow F$ sends $(P,Q)$ to $PQ$.
Here $P,Q$ are the two independent coordinates on $F^2$. A quadratic form is \emph{hyperbolic} if it is isometric to an orthogonal sum of hyperbolic planes. Equivalently, in a suitable
basis of its underlying vector space $F^{2r}$, it has the expression $$
 q(P_1,Q_1,\ldots,P_r,Q_r)=\sum_{i=1}^r P_iQ_i. $$
These definitions apply in every characteristic
\cite[Section~7]{EKM}.
\end{itemize}
\end{definition}

Next, we recall the definition of Pfister forms.
\begin{definition}\label{def:pfister}
Let $F/k$ be a field extension. For $s\geq1$ and
$a_1,\ldots,a_s\in F^\times$, if $\cha F\ne2$, an \emph{$s$-fold Pfister form} \cite[Section~9.B]{EKM} is $$
 \langle\!\langle a_1,\ldots,a_s\rangle\!\rangle
  =\bigotimes_{i=1}^s\langle1,-a_i\rangle
  =\mathop{\perp}_{I\subseteq[s]}\langle(-1)^{|I|}a_I\rangle. $$
Here $\langle c_1,\ldots,c_n\rangle$ denotes
$\sum_jc_jZ_j^2$. Thus an $s$-fold Pfister form has dimension
$2^s$ and is written in the independent coordinates $Z_I$ of
Notation~\ref{not:coordinates}.
Because $-1$ is a square in $k$, over extensions of $k$
we may replace its signed coefficients by $a_I$ after rescaling
coordinates.

If $\cha F=2$, write $[1,a](P,Q)=P^2+PQ+aQ^2$ \cite[Notation~7.3(2)]{EKM}.
The \emph{bilinear Pfister form}
$\bigotimes_{i=1}^s\langle1,z_i\rangle_b$ is the diagonal symmetric
bilinear form with diagonal entries $z_I$ \cite[Section~4.B]{EKM}.
The associated \emph{quadratic $(s+1)$-fold Pfister form} \cite[Section~9.B]{EKM} is $$
 \Pi(a;\boldsymbol z)
   =\left(\bigotimes_{i=1}^s\langle1,z_i\rangle_b\right)\otimes[1,a]
   =\sum_{I\subseteq[s]}z_I(P_I^2+P_IQ_I+aQ_I^2). $$
The pairs $(P_I,Q_I)$ are independent, so this quadratic form has
dimension $2^{s+1}$.
Note that the purely diagonal polynomial $\sum_I z_IZ_I^2$ in characteristic
two is instead a \emph{quasilinear quadratic form} \cite[Remark~7.24 and Section~10]{EKM}.
\end{definition}

The symbol $(a_1,\ldots,a_s)_2$ (resp.\ $[a;\boldsymbol z]$)
vanishes exactly when the corresponding quadratic Pfister form $ \langle\!\langle a_1,\ldots,a_s\rangle\!\rangle$ (resp.\ $\Pi(a;\boldsymbol z)$) is
isotropic \cite[Corollary~9.10, Fact~16.2 and Theorem~23.7]{EKM}. We then recall the meaning of symbols associated to a form.
\begin{definition}\label{def:association}
We use the convention of \cite[Definition~1.2]{KrashenMatzri}.
Let $F/k$ be a field extension, let $m$ be a prime invertible in
$F$, and let $\Phi\in F[Z_1,\ldots,Z_n]$ be a homogeneous form
of degree $m$, with $n\geq1$. Fix a symbol $\beta=(a_1,\ldots,a_s)_m\in H^s(F,\mu_m^{\otimes s})$, with $a_i\in F^\times$ and $s\geq1$.
For a field extension $L/F$, write $\Phi_L$ for the form obtained
by extending scalars and $\beta_L=\res_{L/F}(\beta)$.

We say that $\Phi$ is \emph{associated with $\beta$} if it
\emph{splits} and \emph{neutralizes} $\beta$. Explicitly, for every
field extension $L/F$, both of the following hold:
\begin{enumerate}
\item \emph{Splitting.}
If $\Phi_L(\boldsymbol v)=0$ for some $\boldsymbol v\in L^n\setminus\{0\}$, then $\beta_L=0$.
\item \emph{Neutralizing.}
For every $\boldsymbol v\in L^n$ with
$b=\Phi_L(\boldsymbol v)\ne0$, one has $\beta_L\cup(b)_m=0$ in $H^{s+1}(L,\mu_m^{\otimes(s+1)})$,
where $(b)_m\in H^1(L,\mu_m)$ is the Kummer class of $b$.
\end{enumerate}
\end{definition}

Note that in the above definition, association is a property of the fixed pair $(\Phi,\beta)$.
The symbol need not be uniquely determined by the form. In the next definition, we specify a few norms and their associated symbols used later.
\begin{definition}\label{def:norm-forms}
Let $F/k$ be a field extension, let $m$ be a prime invertible
in $F$, and choose a primitive $m$th root $\zeta\in k$.
For a finite \'etale $F$-algebra $E$, the \emph{norm}
$N_{E/F}(u)$ is the determinant of multiplication by $u$ on $E$. 
\begin{itemize}
    \item For $a\in F^\times$ the \emph{Kummer norm form} for $E=F[i]/(i^m-a)$ is defined by  $$
 N_a(Z_0,\ldots,Z_{m-1})
   =N_{E/F}\left(\sum_{r=0}^{m-1}Z_ri^r\right). $$
Its associated symbol is $(a)_m$.
\item For $a,b\in F^\times$, let $\cD$ be the cyclic algebra with
basis $(i^rj^t)_{0\leq r,t<m}$ with the relations $i^m=a$, $j^m=b$ and $ij=\zeta ji$.
The \emph{reduced norm} $\operatorname{Nrd}$ and \emph{reduced trace} $\operatorname{Trd}$ on $\cD$
are the determinant and trace after extension to a splitting field
and identification with $m\times m$ matrices. These expressions
descend to $F$ \cite[Construction~2.6.1]{GilleSzamuely}.
The \emph{cyclic reduced norm form} is defined to be  $$
 N_{a,b}((Z_{rt}))
  =\operatorname{Nrd}_{\cD/F}
       \left(\sum_{r,t=0}^{m-1}Z_{rt}i^rj^t\right), $$
with associated symbol $(a,b)_m$, using the root-of-unity convention
of \cite[Section~1]{KrashenMatzri}.
\item For $m=3$ and $a,b,c\in F^\times$, the \emph{first Tits norm form} is $$
 \begin{aligned}
 N_c(E_0,E_1,E_2)={}&\operatorname{Nrd}(E_0)
     +c\operatorname{Nrd}(E_1)+c^2\operatorname{Nrd}(E_2)-c\operatorname{Trd}(E_0E_1E_2).
 \end{aligned} $$
Here $E_q=\sum_{r,t=0}^2 Z_{rtq}i^rj^t$ for $0\leq q\leq2$,
where the $27$ coefficients $Z_{rtq}$ are independent
coordinates. This cubic norm is the first Tits construction
in rescaled coordinates \cite[Section~2.4, Equation~(2.6)]{PeterssonThakur2026},
with associated symbol
$(a,b,c)_3$ \cite[Example~1.7]{KrashenMatzri}.
\end{itemize}
\end{definition}

\begin{definition}\label{def:suspensions}
Let $\Phi$ be a degree-$m$ form associated with a symbol $\beta$.
Its \emph{binary suspension} by $d\in F^\times$ is the form defined by  $\Phi(\boldsymbol Z)-d\Phi(\boldsymbol Z')$.
The coordinate tuples $\boldsymbol Z$ and $\boldsymbol Z'$
are independent. By \cite[Proposition~1.9]{KrashenMatzri},
this form is associated with $\beta\cup(d)_m$.
A \emph{coordinate subform} is the restriction obtained by letting
some coordinates be equal to zero.
\end{definition}
\begin{remark}\label{rmk:specified-forms}
We use the Kummer, cyclic and first Tits norm forms, their iterated
binary suspensions, and their coordinate subforms.
These forms split and neutralize their associated symbols in the
sense of Definition~\ref{def:association}, see \cite[Lemma~1.4 and Examples~1.5--1.7]{KrashenMatzri} and \cite[Proposition~1.9]{KrashenMatzri}. Note that both properties pass to coordinate subforms, since their zeros and
represented values are also zeros and represented values of the
original form. In particular, $\beta_L\ne0$ implies that $\Phi_L$
is anisotropic. Finally it is clear that replacing the associated symbol by a nonzero scalar
multiple in $\ZZ/m$ does not change either property.
\end{remark}

\begin{definition}\label{def:weights}
For the forms in Remark~\ref{rmk:specified-forms}, we use their coordinate
presentations, in which the coefficient of every pure power $Z_\alpha^m$ is a single nonzero monomial $c_\alpha\bx^{\omega_\alpha}, c_\alpha\in k^\times.$
We define the \emph{weight} of $Z_\alpha$ to be $\omega_\alpha$. The finite set of weights is denoted by $\cW$.
\end{definition}

\begin{remark}\label{rmk:degree-control}
The Kummer, cyclic and Tits bases have weights $r$, $(r,t)$ and
$(r,t,q)$, respectively. Suspension by a new parameter adds an
entry $0$ or $1$.

Let $\Phi_0(\bx;\boldsymbol Z)$ be one of these degree-$m$
forms in its fixed polynomial presentation, with independent parameter
variables $\bx=(x_1,\ldots,x_s)$ and coordinates
$(Z_\omega)_{\omega\in\cW}$. For independent invertible parameters
$\boldsymbol t=(t_1,\ldots,t_s)$, its norm presentation gives $$
 \Phi_0\bigl((t_i^m x_i)_i;
              (\boldsymbol t^{-\omega}Z_\omega)_\omega\bigr)
       =\Phi_0(\bx;\boldsymbol Z). $$
Indeed, rescaling a Kummer or cyclic generator rescales its coefficient
inversely and preserves the represented element. For the first Tits norm,
the additional change is $c\longmapsto v^3c$ and $(E_0,E_1,E_2)\longmapsto(E_0,v^{-1}E_1,v^{-2}E_2)$.
It preserves both the cubic norm terms and the term
$c\operatorname{Trd}(E_0E_1E_2)$. A suspension by $x_j^a$ preserves
the identity by rescaling its second block by $t_j^{-a}$, and
subforms also inherit it. Thus every nonzero monomial
$c\bx^\alpha\prod_\omega Z_\omega^{b_\omega}$ of $\Phi_0$ satisfies $$
 \sum_\omega b_\omega=m,\qquad
 m\alpha=\sum_\omega b_\omega\omega. $$
Here the second equality follows by comparing powers of the $t_i$.

If we do the substitution $Z_\omega=\xi_\omega$ in the final affine
variables, for every set of coordinates $\cB$, every nonzero substituted
monomial thus satisfies $$
\begin{aligned}
\deg_{\cB}\left(\bx^\alpha\prod_\omega\xi_\omega^{b_\omega}\right)
&=\deg_{\cB}(\bx^\alpha)
  +\sum_{\substack{\omega\in\cW\\b_\omega>0}}
       b_\omega\deg_{\cB}(\xi_\omega)\\
&=\frac1m
  \sum_{\substack{\omega\in\cW\\b_\omega>0}}
       b_\omega\deg_{\cB}(\bx^\omega)
  +\sum_{\substack{\omega\in\cW\\b_\omega>0}}
       b_\omega\deg_{\cB}(\xi_\omega)\\
&=\sum_{\substack{\omega\in\cW\\b_\omega>0}}
       \frac{b_\omega}{m}
       \deg_{\cB}(\bx^\omega\xi_\omega^m).
\end{aligned} $$
Here $\deg_{\cB}$ is total degree in that group. Since the coefficients
$b_\omega/m$ sum to one, the degrees of pure terms give an upper bound for all mixed
terms.
\end{remark}

Fix a field extension $E/k$ and an integer $m\geq2$. For field extensions $F/E$, let $A(F)$ be abelian groups annihilated
by $m$, with functorial restriction maps. For every finite extension
$F'/F$, suppose that there is a corestriction map
$\cores_{F'/F}:A(F')\to A(F)$. We require
compatibility with composition in finite towers and the degree identity: $\cores_{F''/F}=\cores_{F'/F}\circ\cores_{F''/F'}$ for $F\subset F'\subset F''$, and $\cores_{F'/F}\circ\res_{F'/F}=[F':F]\id_{A(F)}$.
Let $U$ be a separated scheme of finite type over $E$, with fixed
classes $a_p\in A(\kappa(p))$ at its closed points. We propose the following condition:

\theoremstyle{definition}
\newtheorem*{curvecondition}{Condition (C)}
\theoremstyle{remark}
\begin{curvecondition}
For every integral curve $C$ closed in $U$, with normalization
$\nu:\widetilde C\to C$, and every $g\in E(C)^\times$, one has $$
 \sum_{w\in\widetilde C^{(1)}}\ord_w(g)
 \cores_{\kappa(w)/E}
 \left(\res_{\kappa(w)/\kappa(\nu(w))}
          a_{\nu(w)}\right)=0. $$
\end{curvecondition}

Under Condition~\textup{(C)}, evaluation induces a homomorphism
from $\CH_0(U,\ZZ/m)$ to $A(E)$.
We denote it by $\ev_U$, which sends $[p]$ to
$\cores_{\kappa(p)/E}(a_p)$.
To see this, first extend the rule linearly to $Z_0(U)$.
For the normalization $\nu:\widetilde C\to C$ of an integral closed
curve, let $w\in\widetilde C$ be closed and let $p=\nu(w)$.
The degree identity and transitivity of corestriction give
$$
 \ev_U(\nu_*[w])
 =[\kappa(w):\kappa(p)]\cores_{\kappa(p)/E}(a_p)
 =\cores_{\kappa(w)/E}
       \left(\res_{\kappa(w)/\kappa(p)}a_p\right).
$$
Condition~\textup{(C)} then gives
$\ev_U(\nu_*\divs_{\widetilde C}(g))=0$ for every
$g\in E(C)^\times$. These cycles generate rational equivalence,
so evaluation descends to $\CH_0(U)$. Since $mA(E)=0$, it factors
through $\CH_0(U,\ZZ/m)$.

\setcounter{theorem}{19}
\begin{lemma}\label{lem:reciprocity}
Let $E/k$ be a field extension, let $m\geq2$ be invertible in $E$,
and fix integers $q,r\geq0$. Let $U$ be a separated scheme of
finite type over $E$, with fixed classes $a_p\in H^q(\kappa(p),\mu_m^{\otimes r})$ at its closed points. For every integral curve $C$ closed in $U$, let
$\nu:\widetilde C\to C$ be its normalization, and choose a normal
projective curve $\overline C$ containing $\widetilde C$ as a
dense open subcurve. At every closed point $w\in\overline C$, let $b_w=\res_{\kappa(w)/\kappa(\nu(w))}(a_{\nu(w)})$ if $w\in\widetilde C$, and let $b_w=0$ if $w\in\overline C\setminus\widetilde C$.
Suppose that, for each such $C$, there exists a class $ \alpha_C\in H^q(E(C),\mu_m^{\otimes r})$ unramified at every closed point of $\overline C$ and satisfying
$s_w(\alpha_C)=b_w$ (Remark~\ref{rmk:unramified}) at every such point. Then the fixed classes satisfy Condition~\textup{(C)}.
Hence, they induce a homomorphism $$
 \ev_U:\CH_0(U,\ZZ/m)
       \longrightarrow H^q(E,\mu_m^{\otimes r}),
 \qquad
 [p]\longmapsto\cores_{\kappa(p)/E}(a_p). $$
\end{lemma}

\begin{proof}
Fix an integral curve $C$ closed in $U$ and a function
$g\in E(C)^\times$. Denote $K=E(C)=E(\overline C)$. We first compute a residue at a closed point
$w\in\overline C$. Since $\overline C$ is a normal curve,
$R_w=\cO_{\overline C,w}$ is a discrete valuation ring.
Choose a uniformizer $\pi_w$ and let
 $g=\pi_w^{n_w}u_w,$ $n_w=\ord_w(g)\in\ZZ,$ and $u_w\in R_w^\times.$
The Kummer homomorphism is additive with respect to
multiplication in $K^\times$, so $$
 (g)_m=n_w(\pi_w)_m+(u_w)_m
 \in H^1(K,\mu_m). $$

Because $\alpha_C$ is unramified at $w$, its restriction to the
fraction field of the henselization $R_w^h$ extends to a class
in $H^q(R_w^h,\mu_m^{\otimes r})$. The reduction of this class
is $s_w(\alpha_C)=b_w$. Residues are compatible with passage to
the henselization, and the residue--cup-product formula
\cite[Lemma~3.1]{Schreieder} thus gives $$
 \partial_w\bigl((f)_m\cup\alpha_C\bigr)
   =\partial_w((f)_m)\cup b_w
 \qquad(f\in K^\times). $$
By \cite[Lemma~3.2]{Schreieder}, the residue of a Kummer class
is its valuation modulo $m$. It follows that $$
 \partial_w\bigl((\pi_w)_m\cup\alpha_C\bigr)=b_w,
 \qquad
 \partial_w\bigl((u_w)_m\cup\alpha_C\bigr)=0. $$
Combining these identities, we obtain
$$\partial_w\bigl((g)_m\cup\alpha_C\bigr) = n_w\,\partial_w\bigl((\pi_w)_m\cup\alpha_C\bigr)+\partial_w\bigl((u_w)_m\cup\alpha_C\bigr)=n_w b_w =\ord_w(g)b_w
  \in H^q(\kappa(w),\mu_m^{\otimes r}). $$

Now apply \cite[Remark~2.5 and Proposition~2.2~\textup{(RC)}]{Rost1996} to  $(g)_m\cup\alpha_C
       \in H^{q+1}(K,\mu_m^{\otimes(r+1)}).$ We have  $$
 \sum_{w\in\overline C^{(1)}}
       \cores_{\kappa(w)/E}
       \bigl(\partial_w\bigl((g)_m\cup\alpha_C\bigr) )=0
 \quad\in H^q(E,\mu_m^{\otimes r}). $$
 Hence $$
\begin{aligned}
 0
 &=
 \sum_{w\in\overline C^{(1)}}
       \cores_{\kappa(w)/E}
       \bigl(\partial_w\bigl((g)_m\cup\alpha_C\bigr) )\\
 &=
 \sum_{w\in\overline C^{(1)}}\ord_w(g)\,
       \cores_{\kappa(w)/E}(b_w)\\
 &=
 \sum_{w\in\widetilde C^{(1)}}\ord_w(g)\,
       \cores_{\kappa(w)/E}(b_w)\\
 &=
 \sum_{w\in\widetilde C^{(1)}}\ord_w(g)\,
       \cores_{\kappa(w)/E}
       \left(
        \res_{\kappa(w)/\kappa(\nu(w))}
        (a_{\nu(w)})
       \right).
\end{aligned} $$
The third equality uses $b_w=0$ at every point of
$\overline C\setminus\widetilde C$, and the fourth uses its
definition on $\widetilde C$. All sums are finite, since
$\ord_w(g)=0$ except at finitely many points. The last expression is exactly the sum in
Condition~\textup{(C)}. Since $C$ and $g$ were arbitrary,
that condition holds. The argument following Condition~\textup{(C)} gives the
homomorphism on $\CH_0(U,\ZZ/m)$.
\end{proof}
The following corollary uses these evaluation maps on Chow groups to detect the generic-point class.

\begin{corollary}\label{cor:affine-detection}
Let $m\geq2$ be an integer invertible in $k$, and let
$V\subset\bA^{s+n}_k$ be an integral closed subvariety, with
coordinate functions $x_1,\ldots,x_s,y_1,\ldots,y_n$, where
$s\geq1$, such that  $0\ne\beta=(x_1,\ldots,x_s)_m \in H^s(k(V),\mu_m^{\otimes s})$.
Suppose that, for every henselian discretely valued extension
$(F,\mu)/k$ with $\mu|_k=0$, every $F$-point of $V$ at which
all $x_i$ are nonzero and $(x_1,\ldots,x_s)_m\ne0$ has $\mu(x_i)=0$ for all $i$ and $\mu(y_j)\geq0$ for all $j$.
For a field extension $L/k$ and a closed point $p\in V_L$,
let $\kappa(p)$ be its residue field. Define $a_p=(x_1(p),\ldots,x_s(p))_m$ if $\prod_i x_i(p)\ne0$, and let $a_p=0$ otherwise, as a class in $H^s(\kappa(p),\mu_m^{\otimes s})$.
Then these classes induce, for every $L/k$, a homomorphism $$
 \ev_{V_L}:\CH_0(V_L,\ZZ/m)
       \longrightarrow H^s(L,\mu_m^{\otimes s}),\qquad
 \sum_p n_p[p]\longmapsto
       \sum_p n_p\cores_{\kappa(p)/L}(a_p), $$
such that $\ev_{V_{k(V)}}(\delta_V)=\beta\neq0 $. In particular, $\delta_V\neq 0$. If furthermore $m$ is prime, $\Tor^{\ZZ/m}(V)=m$.
\end{corollary}

\begin{proof}
Fix a field extension $L/k$ and let $D_L=V_L\cap\{x_1\cdots x_s=0\}.$ We apply Lemma~\ref{lem:reciprocity} with $U=V_L$, $E=L$ and $q=r=s$. Fix an integral closed curve $C\subset V_L$.
Let $\nu:\widetilde C\to C$ be its normalization, and choose a normal
projective compactification $\overline C$.

If $C\subset D_L$, then $a_p=0$ at every closed point
$p\in C$. Taking $\alpha_C=0$ satisfies all the requirements of
Lemma~\ref{lem:reciprocity} for this curve. Now suppose $C\not\subset D_L$. Since $C$ is integral, none of
$x_1|_C,\ldots,x_s|_C$ is identically zero. Hence they belong
to $L(C)^\times$, and we may define $\alpha_C=(x_1|_C,\ldots,x_s|_C)_m \in H^s(L(C),\mu_m^{\otimes s})$.
We verify its unramifiedness and its specializations at all
closed points of $\overline C$.

For a closed point $w\in\overline C$, let $R_w=\cO_{\overline C,w}$,
 $F_w=\Frac(R_w^h),$ and $\alpha_{C,w}
   =\res_{F_w/L(C)}(\alpha_C),$ where $R_w^h$ is the henselization of the discrete valuation
ring $R_w$. Let $\mu_w$ be the valuation on $F_w$ extending $\ord_w$. It is trivial on $L$, hence on $k$,
and its residue field is $\kappa(w)$. The composite $$
 \Spec{}F_w
 \longrightarrow\Spec{}L(C)
 \longrightarrow C\longrightarrow V_L\longrightarrow V $$
defines an $F_w$-point of $V$. By assumption either  $\alpha_{C,w}=0$, or
all $x_i$ are units. Thus $\alpha_C$ is unramified at $w$ in both cases.  Let $\widetilde\alpha_{C,w}\in H^s(R_w^h,\mu_m^{\otimes s})$
be the unique class whose restriction to $F_w$ is
$\alpha_{C,w}$. For the closed immersion
$i_w:\Spec\kappa(w)\hookrightarrow
\Spec{}R_w^h$, we have $s_w(\alpha_C)=i_w^*(\widetilde\alpha_{C,w}),$ see Remark~\ref{rmk:unramified}. 

Now we check the specialization $s_w(\alpha_C)$ for all closed points $w\in \overline C$. 

For a closed point $w\in\widetilde C$, let $p=\nu(w)$. If every $x_i(p)$ is nonzero, then every $x_i|_C$ is a unit
in $R_w$. Functoriality of Kummer classes and cup products gives $s_w(\alpha_C)=\res_{\kappa(w)/\kappa(p)}(a_p).$ If some $x_i(p)=0$, then $x_i|_C$ is a nonzero element of
the maximal ideal of $R_w$, so
$\ord_w(x_i|_C)>0$. The valuation hypothesis
forces $\alpha_{C,w}=0$, and hence $s_w(\alpha_C)=0$.
Since $a_p=0$ by definition, we again have $s_w(\alpha_C)=\res_{\kappa(w)/\kappa(p)}(a_p).$

For a closed point $w\in\overline C\setminus\widetilde C$, we show that $\alpha_{C,w}=0$. Suppose instead that $\alpha_{C,w}\ne0$.
The valuation hypothesis makes all the $x_i|_C$ units and
all the $y_j|_C$ integral in $R_w^h$.
These functions belong to $L(C)$, and henselization
preserves their valuations, so they already belong
to $R_w$. Since $C$ is closed in $\bA^{s+n}_L$,
its coordinate ring $A=\Gamma(C,\cO_C)$ is generated
over $L$ by the restrictions of the $x_i,y_j$.
Their integrality gives an inclusion $A\hookrightarrow R_w.$ Equivalently, the coordinate values define a morphism
$\Spec{}R_w\to C$: every defining equation of $C$ vanishes in $L(C)$ and hence in its subring $R_w$. Let $\widetilde A$ be the integral closure of $A$ in $L(C)$,
so $\widetilde C=\Spec\widetilde A$. Since $R_w$ is integrally
closed, $\widetilde A\subset R_w$.
Thus the above morphism lifts to $h:\Spec{}R_w\longrightarrow\widetilde C$.

Finally, the composite of $h$ with the open immersion
$\widetilde C\hookrightarrow\overline C$ and the canonical
morphism $\Spec\cO_{\overline C,w}\to\overline C$
coincide on the generic point as they both induce the identity on
$L(C)$. Since $\overline C$ is separated, they coincide
everywhere. Hence their images of the closed point must coincide. This shows that $w\in\widetilde C$, a contradiction.

We have checked the hypotheses of
Lemma~\ref{lem:reciprocity} for every integral closed curve
in $V_L$. Applying it gives the induced homomorphism
$\ev_{V_L}$. For $L=k(V)$, the diagonal point evaluates to
$\beta\ne0$, so $\delta_V\ne0$. If $m$ is prime, this nonzero
class in a $\ZZ/m$-module has order $m$.
\end{proof}

\begin{remark}\label{rmk:localization-criterion}
Let $V/k$ be integral,
$K=k(V)$ and $D\subsetneq V$ closed. If there is a homomorphism $\varphi:\CH_0(V_K,\Lambda)\to G$ such that $\varphi(\delta_V)\neq 0$ and $\varphi$
annihilates the image of $\CH_0(D_K,\Lambda)$ inside $\CH_0(V_K,\Lambda)$, then the localization
sequence $$
 \CH_0(D_K,\Lambda)\longrightarrow\CH_0(V_K,\Lambda)
 \longrightarrow\CH_0((V\setminus D)_K,\Lambda)\longrightarrow0 $$
shows that $\delta_{V\setminus D}\ne0$ and hence
$\Tor^\Lambda(V\setminus D)\ne1$. This also gives $\Tor^\Lambda(V')\ne1$ for every integral variety
$V'$ containing $V\setminus D$ as an open subvariety because restriction
sends $\delta_{V'}$ to $\delta_{V\setminus D}$.
\end{remark}

The following lemma calculates valuations of the forms above, see Remark~\ref{rmk:specified-forms}.

\begin{lemma}\label{lem:norm-values}
Let $k$ be algebraically closed, let $m\in\{2,3\}$ be invertible in
$k$, and let $(F,\mu)$ be a henselian discretely valued extension
of $k$, with $\mu|_{k^\times}=0$. Fix
$\bx=(x_1,\ldots,x_s)\in(F^\times)^s$ such that $0\ne\beta=(x_1,\ldots,x_s)_m \in H^s(F,\mu_m^{\otimes s})$.
Let $\Phi(\boldsymbol Z)$ be a form from
Remark~\ref{rmk:specified-forms}, in its fixed coordinate presentation
and associated with a nonzero scalar multiple of $\beta$. Here $\boldsymbol Z=(Z_\omega)_{\omega\in\cW}$ is the tuple of independent form variables, indexed by the 
weights of Definition~\ref{def:weights}. Let $\bP=(\mu(x_1),\ldots,\mu(x_s))$ and $\mu(0)=\infty$.
Then, for every  $\boldsymbol z=(z_\omega)_{\omega\in\cW}\in F^{\cW}$,
$$
 \mu\bigl(\Phi(\boldsymbol z)\bigr)
  =\min_{\omega\in\cW}
       \{\omega\cdot\bP+m\mu(z_\omega)\}. $$
 The formula includes $\boldsymbol z=0$, when both
sides are $\infty$.
\end{lemma}

\begin{proof}
By Remark~\ref{rmk:specified-forms},  $\beta\neq 0$ implies that $\Phi$ is anisotropic. The lemma is immediate for $\boldsymbol z=0$, so we may assume $\boldsymbol z\ne0$. By anisotropy,  this implies $\Phi(\boldsymbol z)\neq0$.

For every $\omega\in\cW$, let $e_\omega$ denote
the corresponding coordinate basis vector. By the fixed
coordinate presentation of Definition~\ref{def:weights}, $\Phi(z_\omega e_\omega)
   =c_\omega\bx^\omega z_\omega^m,$ with $c_\omega\in k^\times.$ Since $\mu$ is trivial on $k^\times$, we have $$
 \mu\bigl(\Phi(z_\omega e_\omega)\bigr)
   =\omega\cdot\bP+m\mu(z_\omega). $$
We prove that the valuation of the full form is the minimum of these coordinate valuations. For $m=2$, the proof of \cite[Lemma~3.12]{Hotchkiss}
uses only Hensel's lemma, so it also works over the henselian field $F$.
Indeed, if the valuation of a diagonal anisotropic form exceeded that
of a least-valued term $az^2$, the sum of the other terms divided by
$-az^2$ would lie in $1+\fm_\mu$ and would be a square. Rescaling $z$
would then give a nonzero isotropic vector.

Suppose $m=3$. Let $N$ be one of the three types of norms in Definition~\ref{def:norm-forms} and let $E$ denote
its underlying vector space. Define $u_N(e)=\mu(N(e))/3$ for $e\ne0$, and let $u_N(0)=\infty$.
Since the associated symbol is nonzero, the
corresponding Kummer algebra is a field, the cyclic
algebra is a division algebra, and the first Tits
algebra is a Jordan division algebra. For the Kummer and cyclic norms, valuation theory gives
$$ u_N(e+e')\geq\min\{u_N(e),u_N(e')\},\qquad u_N(fe)=\mu(f)+u_N(e)
       \qquad(f\in F^\times).$$
See \cite[Corollary~2.2 and Equation~(2.7)]
{WadsworthValuation}. Here the first inequality is called 
the \emph{ultrametric} inequality. For the first Tits norm, we pass from $F$ to its completion
$\widehat F_\mu$. Completion induces an isomorphism of absolute Galois groups for a
henselian discretely valued field
\cite[Expos\'e~X, Lemma~2.2.1]{SGA4III}.
Hence restriction on Galois cohomology is an isomorphism, and the
associated symbol remains nonzero over $\widehat F_\mu$.
Thus the first Tits algebra remains division.
Then
\cite[Section~7.2, Equation~(7.1)]{PeterssonThakur2026} and \cite[Satz~5.1]{Petersson1973} show that the ultrametric inequality and $u_N(fe)=\mu(f)+u_N(e)$ hold for the first Tits norm over
the completion and hence, by restriction, over $F$.

We next show that the ultrametric inequality suffices to prove the lemma. Let $e=\sum_{\omega\in\cW}z_\omega e_\omega\in E$. The following linear operators act on their corresponding underlying vector spaces and preserve the corresponding norm form:
\begin{itemize}
\item For the Kummer norm form, the automorphism $i\mapsto\zeta i$.
\item For the cyclic reduced norm form, conjugation by $i$ and by $j$.
\item For the first Tits norm form, simultaneous conjugation
by $i$ and $j$ in all three blocks, together with $(E_0,E_1,E_2) \longmapsto(E_0,\zeta E_1,\zeta^2E_2)$.
\end{itemize}
Here $\zeta\in k$ is a primitive cube root of unity. These operators generate a finite abelian group
$\cG$ of order a power of $3$. Every element
of $\cG$ preserves $N$, so $u_N(ge)=u_N(e)$ for any $g\in\cG.$
In each case, the fixed coordinate lines $Fe_\omega$
are precisely the distinct character spaces of
$\cG$. Indeed, the Kummer automorphism distinguishes
the powers of $i$, the two cyclic conjugations distinguish
the monomials $i^rj^t$, and the additional block-scaling
operator distinguishes the three Tits blocks.

Let $\chi_\omega:\cG\to k^\times$ be the
character of $Fe_\omega$. By character orthogonality,
the projection onto this coordinate line is given by $$
 \pr_\omega
   =\frac{1}{|\cG|}
      \sum_{g\in\cG}\chi_\omega(g)^{-1}g.$$
Since $3$ is invertible in $k$ and $\mu|_{k^\times}=0$,
every nonzero coefficient in this projection expression
has valuation zero. Using the ultrametric inequality, we obtain
$$ u_N(\pr_\omega e)
 =
 u_N\left(
   \frac{1}{|\cG|}
   \sum_{g\in\cG}\chi_\omega(g)^{-1}ge
      \right)
 \geq
 \min_{g\in\cG}
 u_N\left(
   \frac{\chi_\omega(g)^{-1}}{|\cG|}ge
      \right)
 =u_N(e). $$
Since $\pr_\omega e=z_\omega e_\omega$, the above inequality implies $u_N(e)\leq
 \min
     u_N(z_\omega e_\omega).$ Conversely, applying ultrametricity to $e=\sum_\omega z_\omega e_\omega$ gives $u_N(e)\geq
 \min
     u_N(z_\omega e_\omega).$
Combining the two inequalities shows that  $$
 u_N(e)=
 \min_{\omega\in\cW}
     u_N(z_\omega e_\omega). $$
Then the discussion at the beginning of the proof shows that 
$$\mu(N(e)) =\min_{\omega\in\cW}
       \mu(N(z_\omega e_\omega))
 =\min_{\omega\in\cW}
       \{\omega\cdot\bP+3\mu(z_\omega)\}. $$
Thus the formula holds for all three forms.
It also holds for their coordinate subforms, since
letting omitted coordinates be equal to zero simply removes
the corresponding entries from the minimum. It remains to prove that the formula is preserved under
binary suspension (Definition \ref{def:suspensions}). 

Suppose the lemma holds for one of these degree-$m$ forms $\Psi$, associated with a symbol
$\gamma\ne0$. Let $\Phi(\boldsymbol Z,\boldsymbol Z')
   =\Psi(\boldsymbol Z)-d\Psi(\boldsymbol Z')$ be its binary suspension by $d$ associated with $\gamma\cup(d)_m\ne0$. For substituted vectors $\boldsymbol z,\boldsymbol z'$, denote $ A=\Psi(\boldsymbol z),B=\Psi(\boldsymbol z').$
We claim that $$
 \mu(A-dB)=\min\{\mu(A),\mu(d)+\mu(B)\}. $$
If $A=0$ or $B=0$, this is immediate. Suppose both
are nonzero. The valuation inequality always gives $\mu(A-dB)\geq\min\{\mu(A),\mu(dB)\}.$ If this inequality were strict, we must have $\mu(A)=\mu(dB)$ and $\mu\left(1-\frac{dB}{A}\right)>0.$
Thus $\frac{dB}{A}\in1+\fm_\mu.$
Since $m$ is invertible in the residue field,
Hensel's lemma applied to $T^m-dB/A$ shows that $dB/A$ is an $m$th power in $F$.

On the other hand, the neutralizing property in
Definition~\ref{def:association} gives $\gamma\cup(A)_m=0$ and $\gamma\cup(B)_m=0.$ Because $dB/A$ is an $m$th power, we would have $$
0=\gamma\cup(\frac{dB}{A})_m=\gamma\cup(d)_m+\gamma\cup(B)_m-\gamma\cup(A)_m=\gamma\cup(d)_m, $$
contrary to the hypothesis. So $\mu(A-dB)=\min\{\mu(A),\mu(dB)\},$ and this case is reduced to the first half of the proof.\qedhere

\end{proof}

We also need the valuations of Laurent polynomials proved in the next lemma. The following lemma also holds for valuations of arbitrary rank.
\setcounter{theorem}{24}
\begin{lemma}\label{lem:laurent-values}
Let $k$ be algebraically closed, let $m\in\{2,3\}$ be invertible in
$k$, and let $(F,\mu)$ be a henselian valued extension of $k$, with
value group $\Gamma$ of arbitrary rank and $\mu|_{k^\times}=0$.
Fix $\bx=(x_1,\ldots,x_s)\in(F^\times)^s$ with
 $0\ne\beta=(x_1,\ldots,x_s)_m\in H^s(F,\mu_m^{\otimes s}),$ and denote $\bP=(\mu(x_1),\ldots,\mu(x_s)).$ 
 
 Then for every nonzero Laurent polynomial $P(\bx)=\sum_{\nu\in S}c_\nu\bx^\nu
       \in k[X_1^{\pm1},\ldots,X_s^{\pm1}],$
where $S\subset\ZZ^s$ is finite and nonempty and each
$c_\nu\in k^\times$, one has $$
 P(\bx)\ne0,\qquad
 \mu\bigl(P(\bx)\bigr)
       =\min_{\nu\in S}\nu\cdot\bP. $$
\end{lemma}
\begin{proof}
The subgroup of $\Gamma$ generated by the $\mu(x_i)$ is a finitely generated free abelian group. Make an invertible integral monomial change from
$\bx$ to
$(a_1,\ldots,a_r,b_1,\ldots,b_{s-r})$ such that the
$\mu(a_i)$ form a basis of the subgroup generated by the
$\mu(x_i)$ and all $b_j$ are units. Since $k$ is algebraically closed, alternation gives $$
 \beta=\pm(a_1,\ldots,a_r,b_1,\ldots,b_{s-r})_m. $$
Let $A=\cO_\mu$, let $\kappa=\kappa(\mu)$, and let $d=s-r$.
The units $b_1,\ldots,b_d$ define a class $\alpha\in H^d(A,\mu_m^{\otimes d})$.
Henselian invariance
\cite[Expos\'e~VIII, Corollary~8.6]{SGA4II} gives $$
 H^d(A,\mu_m^{\otimes d})
 \xrightarrow{\sim}
 H^d(\kappa,\mu_m^{\otimes d}). $$
If $(\bar b_1,\ldots,\bar b_d)_m$ were zero, then
$\alpha=0$, and restriction to $F$ followed by cup product
with $(a_1,\ldots,a_r)_m$ would give $\beta=0$.
Thus $(\bar b_1,\ldots,\bar b_d)_m\ne0$, with the empty
symbol interpreted as $1$.
The cohomological-dimension bound
\cite[Equation~(2.3)]{Schreieder} then makes the $\bar b_j$ algebraically independent over $k$.
Express the Laurent polynomial after the monomial change as $P=\sum_{e\in\ZZ^r}a^eP_e(b_1,\ldots,b_{s-r})$,
where only finitely many $P_e$ are nonzero. Algebraic independence
of the residues makes each nonzero $P_e(b)$ a unit. Distinct exponent
vectors $e$ give distinct valuations $\sum_i e_i\mu(a_i)$, so the
summands have a unique least valuation. This proves both
$P(\bx)\ne0$ and the minimum formula.
\end{proof}

For a discrete valuation on a nonhenselian field,
Lemmas~\ref{lem:norm-values} and~\ref{lem:laurent-values} can be applied
to $F_\mu^h$ by verifying $\beta|_{F_\mu^h}\ne0$. We now perform our first construction, which is modified based on the twisted Pfister equations
used in \cite[Section~7]{SchTorsion} and
\cite[Section~6]{LS}.

\begin{construction}\label{constr:norm-data}
Let $\Omega$ be an algebraically closed field, let $m\in\{2,3\}$
be invertible in $\Omega$, and let $\lambda$ be
transcendental over $\Omega$. Let $k=\overline{\Omega(\lambda)}.$ Let $\bx=(x_1,\ldots,x_s)$, with $s\geq1$,
be independent base variables. Choose one of these forms
as in Remark~\ref{rmk:specified-forms}, written as $\Phi(\bx;\boldsymbol Z)\in\Omega[\bx][Z_\omega:\omega\in\cW]$, where $\boldsymbol Z=(Z_\omega)_{\omega\in\cW}$. Its associated symbol over $\Omega(\bx)$ is
a nonzero scalar multiple of $(x_1,\ldots,x_s)_m.$
Here the $Z_\omega$ are independent form variables,
and $\cW$ is their weight set as in
Definition~\ref{def:weights}. Choose a distinguished weight $\eta\in\cW$ and introduce further independent variables
$\boldsymbol z=(z_1,\ldots,z_n)$ with an extra variable $W$, all independent
of $\bx$. Let $\xi_\eta=W$, let $\xi_\omega\in\Omega[\bx,\boldsymbol z]$ for $\omega\in\cW\setminus\{\eta\}$, and write $\boldsymbol\xi=(\xi_\omega)_{\omega\in\cW}$.
Choose an exponent vector $\gamma=(\gamma_1,\ldots,\gamma_s)\in\ZZ_{\geq0}^s$
and define $$
 G(\bx,\boldsymbol z,W)
  =\Phi(\bx;\boldsymbol\xi)
      -\lambda\bx^\gamma W
      \in k[\bx,\boldsymbol z,W]. $$
Since no $\xi_\omega$ with
$\omega\ne\eta$ involves $W$, the polynomial $G$ has
degree $m$ in $W$, with leading coefficient
$c_\eta\bx^\eta$. Let $\tau=\frac{m\gamma-\eta}{m-1}\in\QQ^s\subset\RR^s$.
Let $G_{\mathrm{prim}}=G/c(\bx)$, see Definition~\ref{def:profile}. Define $V=Z(G_{\mathrm{prim}})\subset\bA_k^{s+n+1}$. 

Finally, to apply Lemma~\ref{lem:norm-values}, for every henselian discretely valued extension $(H,\mu)/k$ with $\mu|_{k^\times}=0$, consider every zero $(\bx,\boldsymbol z,W)
 \in(H^\times)^s\times H^n\times H$ of $G=0$ such that $\beta:=(x_1,\ldots,x_s)_m\ne0$. We refer to these zeros as \emph{admissible zeros}. Again denote $\bP=(\mu(x_1),\ldots,\mu(x_s)),$ and $\mu(0)=\infty.$

\end{construction}
The next result shows that Construction~\ref{constr:norm-data} provides a reference variety with carefully fixed valuations. For a finite nonempty set $T\subset\ZZ^s$, let $\operatorname{conv}(T)$ denote its convex hull.
We use the notation $\operatorname{int}(\operatorname{conv}(T))$
to denote its interior in the ambient space $\RR^s$. We then have:

\begin{theorem}\label{thm:norm-criterion}
Use Construction~\ref{constr:norm-data}.
Suppose that there exists a finite nonempty set $T\subset\ZZ^s$ that satisfies

\begin{itemize}
    \item $\tau
 \in\operatorname{int}(\operatorname{conv}(T)).$
 \item at every admissible zero, 
 $\min_{\omega\in\cW} \bigl\{\omega\cdot\bP+m\mu(\xi_\omega)\bigr\} \leq\min_{t\in T}t\cdot\bP$.
\item if at an admissible zero we have $\mu(x_i)=0$ and $\mu(\xi_\omega)\geq0$ for all $1\leq i\leq s$ and $\omega\in \cW $, then $\mu(z_j)\geq0\quad(1\leq j\leq n).$ 

\end{itemize}

Then the affine scheme $V$ is geometrically integral and satisfies $\Tor^{\ZZ/m}(V)=m.$

\end{theorem}

\begin{proof}
We would like to apply Corollary~\ref{cor:affine-detection}. Let us start by checking the valuation conditions needed using Lemma~\ref{lem:norm-values}.   

Fix a henselian discretely valued extension $(H,\mu)/k$
and an admissible solution $(\bx,\boldsymbol z,W)$.
Since $T$ is finite and nonempty, the number
$\min_{t\in T}t\cdot\bP$ is finite.
Then the inequality in our second assumption implies that at least
one $\xi_\omega$ is nonzero.
By Remark~\ref{rmk:specified-forms}, the nonzero symbol
$\beta=(x_1,\ldots,x_s)_m$ makes the form anisotropic, so we have  $\Phi(\bx;\boldsymbol\xi)\ne0.$ The equation $ \Phi(\bx;\boldsymbol\xi)
   =\lambda\bx^\gamma W$ then gives $W\ne0$, since $\lambda$ and all $x_i$ are
nonzero. Let $w=\mu(W)$.

Because $\lambda\in k^\times$ and $\mu$ is trivial on
$k^\times$, Lemma~\ref{lem:norm-values} gives $$
 \gamma\cdot\bP+w
 =\mu\bigl(\lambda\bx^\gamma W\bigr)
 =\mu\bigl(\Phi(\bx;\boldsymbol\xi)\bigr)
 =\min_{\omega\in\cW}
       \bigl\{\omega\cdot\bP
                      +m\mu(\xi_\omega)\bigr\}. $$
In particular, the distinguished coordinate $\xi_\eta=W$ gives $\gamma\cdot\bP+w
 \leq\eta\cdot\bP+mw.$ On the other hand the assumed upper bound implies that $\gamma\cdot\bP+w
 \leq t\cdot\bP.$ Thus $$
 \frac{(\gamma-\eta)\cdot\bP}{m-1}
 \leq w
 \leq(t-\gamma)\cdot\bP
 ,
$$ for every $t\in T$.
The above inequalities also imply that  $$
 \tau\cdot\bP
 =
 \frac{(m\gamma-\eta)\cdot\bP}{m-1}
 \leq t\cdot\bP
 ,
$$ for every $t\in T$.
By linearity, the same inequality holds for every
$t\in\operatorname{conv}(T)$. 

We claim that this implies that $\bP=0$. In fact, if $\bP\neq0$, because $\tau$ is an interior point in
$\RR^s$, there is a sufficiently small $\varepsilon>0$ such that $\tau-\varepsilon\bP\in\operatorname{conv}(T)$.
But $$
 (\tau-\varepsilon\bP)\cdot\bP
 =
 \tau\cdot\bP
 -\varepsilon\sum_{i=1}^s\mu(x_i)^2
 <
 \tau\cdot\bP, $$
gives a contradiction. Hence $\mu(x_i)=0$ for every $i$. The bounds on $w$ now force $w=\mu(W)=0$, and the consequence of Lemma~\ref{lem:norm-values} becomes
$0=\min_{\omega\in\cW}m\mu(\xi_\omega).$ This means every $\xi_\omega$ is integral. By our third assumption here, this further implies that every $z_j$ is integral as well. Hence we have proved that  $\mu(x_i)=0$, $\mu(W)=0$ and $\mu(z_j)\geq0$ for all $i,j$,
as needed by Corollary~\ref{cor:affine-detection}. To show $V=Z(G_{\mathrm{prim}})\subset\bA_k^{s+n+1}$ is geometrically integral, let $R=k[\bx,\boldsymbol z], F_0=k(\bx,\boldsymbol z).$ By Construction~\ref{constr:norm-data},
$G\in R[W]$ has degree $m$ in $W$ and leading coefficient $c_\eta\bx^\eta$.

We first prove that $G$ is irreducible in $F_0[W]$.
Let $\widehat F_0=k(x_2,\ldots,x_s,\boldsymbol z)((x_1))$ be the completion of $F_0$ for the valuation $\ord_{x_1}$. The coordinate symbol $\beta_0=(x_1,\ldots,x_s)_m \in H^s(\widehat F_0,\mu_m^{\otimes s})$
is nonzero. Indeed, taking the residue first at $x_1$
and then successively at $x_2,\ldots,x_s$ gives $$
 \partial_{x_s}\cdots\partial_{x_1}(\beta_0)
 =
 1\in H^0(k(\boldsymbol z),\ZZ/m)=\ZZ/m, $$
by \cite[Lemmas~3.1--3.2]{Schreieder}. Suppose that $G$ were reducible over $F_0$. Since $m\in\{2,3\}$, $G$ has a root in $F_0$.
Let $\mu=\ord_{x_1}$ on $\widehat F_0$.
Now the images of $\bx,\boldsymbol z$, together
with this root, form a solution of $G=0$ over $\widehat F_0$.
Its coordinate symbol is nonzero, so the valuation
condition proved above holds and gives $\mu(x_1)=0$,
whereas $\mu(x_1)=1$.
This contradiction proves that $G$ is irreducible in $F_0[W]$.

We now go back to $R[W]$ and consider the polynomial
$c(\bx)$ used in Construction~\ref{constr:norm-data}.
Let $d\in R$ be the gcd of the coefficients of $G$
when $G$ is regarded as a polynomial in $W$.
Since $d$ divides its leading coefficient
$c_\eta\bx^\eta$, every irreducible factor
of $d$ is one of the $x_i$.
Thus, up to a nonzero scalar, $d$ is a monomial in
$\bx$ and belongs to $k[\bx]$. But a polynomial in $k[\bx]$ divides all the
coefficients of $G$ as a polynomial in $W$ over $R$
if and only if it divides all the coefficients
$A_\alpha(\bx)$ in the expansion $G=\sum_\alpha A_\alpha(\bx) (\boldsymbol z,W)^\alpha$.
So $d=c(\bx)$ up to a nonzero
scalar. Hence $G_{\mathrm{prim}}=G/c(\bx)$
is primitive in $R[W]$.
Gauss's lemma and the irreducibility over $F_0$ show
that $G_{\mathrm{prim}}$ is irreducible in $k[\bx,\boldsymbol z,W]$.
Thus $V$ is integral. It is also geometrically
integral because $k$ is algebraically closed.

Next, we show the generic symbol $\beta\neq 0$. It suffices to use $G$ instead of $G_{\mathrm{prim}}$ for this purpose. In fact, the $x_i$ are nonzero
in $k(V)$, and $V^\circ=V\cap\{x_1\cdots x_s\ne0\}$ is a nonempty dense open subset. As we have shown above, $c(\bx)$ is a monomial in $\bx$ up to a
nonzero scalar, so it is invertible on $V^\circ$. Thus $G=0$ and $G_{\mathrm{prim}}=0$ define the same
scheme there.

Let $t=\lambda^{-1}$ and work first over $\Omega(t)$,
before extending constants to $k=\overline{\Omega(t)}$. Note that the coefficients of $\Phi$ and all the assignments
belong to $\Omega[\bx,\boldsymbol z,W]$, independent of $t$.

Let $E/\Omega(t)$ be any finite extension contained
in $k$. Choose a discrete valuation ring $R_E$ of $E$
lying above the $t$-adic valuation ring $\Omega[t]_{(t)}$.
Its residue field is $\Omega$, since it is a finite
extension of the algebraically closed field $\Omega$.
Consider the model over the DVR $R_E$: $$
 \cV_E
 =
 \Spec
 \frac{
 R_E[x_1^{\pm1},\ldots,x_s^{\pm1},
                \boldsymbol z,W]
 }{
 \bigl(t\Phi(\bx;\boldsymbol\xi)
             -\bx^\gamma W\bigr)
 }. $$
Write $U_E$ for its generic fiber.
Since $t$ is invertible in $E$, its defining equation
is equivalent to $G=0$.
Its base change to $k$ is $V^\circ$, so $U_E$ is
integral. Denote its function field by $K_E=E(U_E).$ The special fiber of $\cV_E$ is $$
 \Spec
 \frac{
 \Omega[x_1^{\pm1},\ldots,x_s^{\pm1},
                   \boldsymbol z,W]
 }{
 (\bx^\gamma W)
 }
 \simeq
 \Spec
 \Omega[x_1^{\pm1},\ldots,x_s^{\pm1},\boldsymbol z], $$
because $\bx^\gamma$ is invertible.
In particular, this special fiber is integral, with
function field $\Omega(\bx,\boldsymbol z)$.

Let $\zeta_E$ be the generic point of the special fiber. The derivative of the defining equation with respect
to $W$ is $t\frac{\partial}{\partial W}
       \Phi(\bx;\boldsymbol\xi)
 -\bx^\gamma.$ Its reduction at the special fiber of $R_E$ is
$-\bx^\gamma$, which is a unit.
Thus $\cV_E$ is smooth over $R_E$ at every point of its special fiber. In particular, $A_E=\cO_{\cV_E,\zeta_E}$ is a DVR with fraction field $K_E$
and residue field $\Omega(\bx,\boldsymbol z)$. Since each $x_i$ is a unit in $A_E$, the symbol
$\beta_E=(x_1,\ldots,x_s)_m
 \in H^s(K_E,\mu_m^{\otimes s})$ is thus unramified at the valuation given by $A_E$, and the
specialization of
Remark~\ref{rmk:unramified} gives $$
 s_{A_E}(\beta_E)
 =
 (x_1,\ldots,x_s)_m
 \in H^s(\Omega(\bx,\boldsymbol z),
                        \mu_m^{\otimes s}). $$
The class on the right hand side is nonzero, which can be checked by taking iterated
coordinate residues as above.
Since the specialization of the zero class is zero,
we obtain $\beta_E\ne0$.

This argument works for every finite constant extension $E/\Omega(t)$ contained in $k$.
The fields $K_E$ form a filtered union equal to $k(V)$.
By compatibility of \'etale cohomology with filtered
limits \cite[Lemma~59.51.5 (\href{https://stacks.math.columbia.edu/tag/03Q6}{Tag~03Q6})]{stacks-project}, if $\beta=0$ over $k(V)$, it would already be $0$ over some $K_E$. This is a contradiction. Hence $\beta\ne0$.

We have checked all the hypotheses of
Corollary~\ref{cor:affine-detection} for $V$. It gives $\Tor^{\ZZ/m}(V)=m$.
\end{proof}

We provide a sufficient way to verify the second
assumption in Theorem~\ref{thm:norm-criterion}.
Let $J\subseteq\cW\setminus\{\eta\}$ be nonempty,
and suppose that $\xi_\omega=P_\omega(\bx)
 \in\Omega[\bx]\setminus\{0\}$ are polynomials that depend only on $\bx$ for all $\omega\in J$.
Write $P_\omega(\bx)=\sum_{\nu\in\operatorname{supp}(P_\omega)}a_{\omega,\nu}\bx^\nu$, with $a_{\omega,\nu}\in\Omega^\times$.
At every admissible zero, Lemma~\ref{lem:laurent-values} gives $$
 P_\omega(\bx)\ne0,
 \qquad
 \mu(P_\omega(\bx))
 =\min_{\nu\in\operatorname{supp}(P_\omega)}\nu\cdot\bP. $$
Hence
$$ \min_{\omega\in\cW}
       \{\omega\cdot\bP+m\mu(\xi_\omega)\}
 \leq
   \min_{\omega\in J}
       \{\omega\cdot\bP+m\mu(P_\omega(\bx))\}
=
   \min_{\substack{\omega\in J\\
             \nu\in\operatorname{supp}(P_\omega)}}
             (\omega+m\nu)\cdot\bP
   =\min_{t\in T}t\cdot\bP, $$
where the finite nonempty set $T$ is fixed as  $T=\bigcup_{\omega\in J}(\omega+m\operatorname{supp}(P_\omega))$, or equivalently $T=\{\omega+m\nu:\omega\in J,\ \nu\in\operatorname{supp}(P_\omega)\}$.
Note that the set $T$ does not depend on the
admissible zero. Thus these assignments verify the second hypothesis of
Theorem~\ref{thm:norm-criterion}. The rest of the subsection is devoted to finding appropriate constructions with the other hypotheses of Theorem~\ref{thm:norm-criterion}.

We now verify the convexity hypothesis for the binary construction.

\setcounter{theorem}{28}
\begin{definition}\label{def:convex-notation}
Let $a_1,\ldots,a_r$ be positive integers, let $a=\sum_j a_j$, and let
$[a]=A_1\sqcup\cdots\sqcup A_r$ with $|A_j|=a_j$.
For $N\geq1$ and $I\subseteq[N]$, the \emph{indicator vector}
$\one_I\in\{0,1\}^N$ has $i$th coordinate $1$ if $i\in I$ and $0$ otherwise.
Write $\one$ for the all-ones vector and $\one_\varnothing=0$.
For nonnegative integers $k_1,\ldots,k_r$, let $K(k_1,\ldots,k_r)$ be the
convex hull of the vectors $\one_I$ with $I\subseteq[a]$ and
$|I\cap A_j|\leq k_j$ for every $j$. Equivalently,
$$
 K(k_1,\ldots,k_r)
 =\prod_{j=1}^r
 \left\{
 (v_i)_{i\in A_j}\in[0,1]^{a_j}:
 \sum_{i\in A_j}v_i\leq k_j
 \right\}.
$$
For convex sets $B,B'$ and $\alpha\geq0$, write
$B+B'=\{v+w:v\in B,\ w\in B'\}$ and
$\alpha B=\{\alpha v:v\in B\}$ for their \emph{Minkowski sum} and
\emph{scalar dilate}. For $S\subseteq[a]$, write
$[0,\one_S]=\{\alpha\one_S:0\leq\alpha\leq1\}$.
\end{definition}

\begin{construction}\label{constr:binary-data}
Use the partition in Definition~\ref{def:convex-notation}.
Let $\Omega$ be algebraically closed of characteristic different from two,
let $\lambda$ be transcendental over $\Omega$, and let
$k=\overline{\Omega(\lambda)}$.
Let $d=3$, let $g\in\{0,1\}$, and let $h=1+g$.
Let $\bx=(x_1,\ldots,x_a)$ and $\boldsymbol u=(u_1,\ldots,u_h)$ be
independent base variables. Let $Y_{\theta,J}$ be independent coordinates
for $\theta\subseteq[a]$ and $J\subseteq[h]$.
Let $C=x_{[a]}u_{[h]}$ and $W=Y_{[a],[h]}$.
The Pfister polynomial
$$
 \Pi=\Pi(\bx,\boldsymbol u;\boldsymbol Y)
 =\sum_{\theta\subseteq[a]}\sum_{J\subseteq[h]}
      x_\theta u_JY_{\theta,J}^2
$$
has associated symbol $(x_1,\ldots,x_a,u_1,\ldots,u_h)_2$.
Let $k_j=\floor{(a_j+1)/2}$ for each $j$.
Let $Q$ be the product over $1\leq j\leq r$ of the sums
$\sum_{I\subseteq A_j,\ |I|\leq k_j}x_I$.
Let $P_0=Q\prod_{i=1}^h(1+u_i)$.
Choose one index from every even-sized block $A_j$ and let $S$ be the set
of these indices. If all $a_j$ are odd, let $S=\varnothing$.
If $S=\varnothing$, let $Y_{\varnothing,\varnothing}=P_0$.
Otherwise introduce an independent coordinate $Y$ and let
$Y_{\varnothing,\varnothing}=P_0-x_SY$ and $Y_{S,\varnothing}=P_0+Y$.
Then
$$
 Y_{\varnothing,\varnothing}+x_SY_{S,\varnothing}=(1+x_S)P_0.
$$
All other $Y_{\theta,J}$ are free, and any free coordinate except $W$
may be set equal to zero. If $S\ne\varnothing$, the coordinate $Y$ may
also be set equal to zero.
Let $\Phi$ be the polynomial after these substitutions, let
$G=\Phi-\lambda CW$, and let $V=Z(G_{\mathrm{prim}})$ as in
Construction~\ref{constr:norm-data}.
Here $\eta=\gamma=\one\in\ZZ^{a+h}$ and $\tau=\one$.
\end{construction}

\begin{lemma}\label{lem:binary-detection}
The primitive hypersurfaces $V$ in Construction~\ref{constr:binary-data}
are geometrically integral and have $\Tortwo=2$.
\end{lemma}

\begin{proof}
We verify the hypotheses of Theorem~\ref{thm:norm-criterion} for every
permitted omission of free coordinates.
Write $\xi_{\theta,J}$ for the polynomial substituted for $Y_{\theta,J}$.
Let $(H,\mu)/k$ be a henselian discretely valued extension with
$\mu|_{k^\times}=0$, and consider an admissible zero.
Let $\bP=(\mu(x_1),\ldots,\mu(x_a),\mu(u_1),\ldots,\mu(u_h))$
and $\mu(0)=\infty$. Let
$M=\min_{\theta,J}\{\mu(x_\theta u_J)+2\mu(\xi_{\theta,J})\}$.
By Lemma~\ref{lem:norm-values}, $\mu(\Phi)=M$.
Let $K=K(k_1,\ldots,k_r)$ and $\cK=2K+[0,\one_S]$.
For each $j$, we have $2k_j\geq a_j$, and $|S\cap A_j|$ is $1$ if
$2k_j=a_j$ and $0$ otherwise.
Choose $0<\varepsilon<1$ and let $v=\one-\varepsilon\one_S$.
Every coordinate of $v$ lies in $(0,2)$, and
$$
 \sum_{i\in A_j}v_i=a_j-\varepsilon|S\cap A_j|<2k_j.
$$
The description of $K$ in Definition~\ref{def:convex-notation} gives
$v\in\operatorname{int}(2K)$.
Since $\varepsilon\one_S\in[0,\one_S]$, translation of a neighborhood of
$v$ gives $\one\in\operatorname{int}(\cK)$.
This also holds for $S=\varnothing$.

The exponent convex hull of $P_0$ is $K\times[0,1]^h$.
By Lemma~\ref{lem:laurent-values}, $P_0\ne0$ and
$\ell:=\mu(P_0)=\min_{v\in K\times[0,1]^h}v\cdot\bP$.
If $S=\varnothing$, then $M\leq2\ell$.
Suppose $S\ne\varnothing$ and let $\sigma=\mu(x_S)$.
Lemma~\ref{lem:laurent-values} gives $\mu(1+x_S)=\min\{0,\sigma\}$.
The identity in Construction~\ref{constr:binary-data} gives
$$
 \min\{\mu(\xi_{\varnothing,\varnothing}),
          \sigma+\mu(\xi_{S,\varnothing})\}
 \leq\ell+\min\{0,\sigma\}.
$$
If $\sigma\geq0$, either $\mu(\xi_{\varnothing,\varnothing})\leq\ell$
or $\mu(\xi_{S,\varnothing})\leq\ell-\sigma$. Hence
$M\leq2\ell$.
If $\sigma\leq0$, either $\mu(\xi_{\varnothing,\varnothing})\leq\ell+\sigma$
or $\mu(\xi_{S,\varnothing})\leq\ell$. Hence $M\leq2\ell+\sigma$.
Let $\cP_3=\cK\times[0,2]^h$. In both cases,
$$
 M\leq2\mu(P_0)+\min\{0,\mu(x_S)\}
   =\min_{z\in\cP_3}z\cdot\bP.
$$
The equality follows because the two summands in
$2(K\times[0,1]^h)+([0,\one_S]\times\{0\})$ vary independently.
Moreover, $\one\in\operatorname{int}(\cP_3)$ by the interior statement for
$\cK$. The vertices of $\cP_3$ form a finite set of integral vectors,
so they give the set $T$ in Theorem~\ref{thm:norm-criterion}.

Finally, suppose all base coordinates are units and all
$\xi_{\theta,J}$ are integral. Every free $Y_{\theta,J}$ is then integral.
The polynomial $P_0$ is integral because its coefficients lie in $\Omega$.
If the additional coordinate $Y$ is present, then
$Y=\xi_{S,\varnothing}-P_0$ is integral.
This proves the third hypothesis of Theorem~\ref{thm:norm-criterion},
which gives the conclusion.
\end{proof}

\subsection{Quadric fibrations}\label{subsec:quadric-models}
In this subsection, we construct affine reference varieties birational to a quadric fibration over rational bases. Following the construction in \cite[Theorems~3.4 and~3.14]{Hotchkiss}, we extend these
constructions to the affine models and coordinate profiles
needed for complete intersections in products. 

A field is $C_d$ if every homogeneous form of degree $e\geq1$ in more than
$e^d$ variables has a nontrivial zero. A field of transcendence degree at
most $d$ over the algebraically closed field $k$ is $C_d$
\cite[Theorem~6]{Lang}, \cite[Theorem~2a]{Nagata1957}.
For a field that is not finitely generated, the coefficients descend to a
finitely generated subfield of the same transcendence bound. In particular,
the residue field of a nontrivial valuation of a rational $d$-dimensional
base is $C_{d-1}$.

\begin{lemma}\label{lem:common-symbol}
Assume $\cha k\ne2$, and let $F/k$ be a finitely generated
field extension. Let $s\geq1$, choose
$a_1,\ldots,a_s,g_1,g_2\in F^\times$, and define
$ \rho=(a_1,\ldots,a_s)_2,$ and $\alpha_i=\rho\cup(g_i)_2$ with $i=1,2$. Suppose that $\alpha_1,\alpha_2\ne0$.
Let $\Pi$ be the Pfister form associated with $\rho$, let $q=\Pi\perp\langle-g_1g_2\rangle$, and let $Q=\{q=0\}=\{\Pi(\boldsymbol Z)-g_1g_2\zeta^2=0\}\subset\PP_F^{2^s}$.
Here $\zeta$ is the extra coordinate for homogenization. Then $Q$ is smooth and geometrically integral, and
$\alpha_1|_{F(Q)}=\alpha_2|_{F(Q)}\neq 0$. Denote it by $\alpha$. If furthermore, for every nontrivial discrete valuation $\mu$ of $F/k$,
at least one of $\alpha_1,\alpha_2$ is unramified at $\mu$,
then 
$$
 0\ne\alpha\in
 H^{s+1}_{\nr}
 \bigl(F(Q)/k,\mu_2^{\otimes(s+1)}\bigr).
$$
\end{lemma}

\begin{proof}
The form $q$ is nondegenerate of dimension $2^s+1\geq3$, so $Q$ is smooth
and geometrically integral. Let $L=F(Q)$. Its generic point shows that
$\Pi_L$ represents $g_1g_2$. Thus
$$
 \alpha_1|_L+\alpha_2|_L=\rho_L\cup(g_1g_2)_2=0
$$
by \cite[Theorem~2.2]{SchSmall}.

Suppose their common restriction is zero. The Pfister form
$\Psi=\Pi\perp(-g_1)\Pi$ is anisotropic over $F$ and isotropic over $L$.
The subform theorem gives $cq\subset\Psi$ for some $c\in F^\times$
\cite[Theorem~2.3]{SchSmall}. Since $q$ represents $1$, the form $\Psi$
represents $c$. Roundness gives $c^{-1}\Psi\simeq\Psi$, hence
$q\subset\Psi$ \cite[Corollary~9.9]{EKM}. Orthogonal decomposition and
cancellation \cite[Proposition~7.22 and Theorem~8.4]{EKM} give
$$
 \Pi\perp(-g_1)\Pi
 \simeq\Pi\perp\langle-g_1g_2\rangle\perp\psi,
 \qquad
 (-g_1)\Pi\simeq\langle-g_1g_2\rangle\perp\psi.
$$
Then $\Pi$ represents $g_2$, contradicting $\alpha_2\ne0$ by
\cite[Theorem~2.2]{SchSmall}.

For the last conclusion, let $\nu$ be a nontrivial discrete valuation of
$L/k$. If $\nu|_F=0$, all symbol entries are units. Otherwise its normalized
restriction is a discrete valuation $\mu$ of $F/k$. By hypothesis, one
$\alpha_i$ extends over $\cO_\mu$. Pulling it to $\cO_\nu$ proves that
the common class has zero residue, by the residue exact sequence
\cite[Theorem~2.3]{Schreieder}.
\end{proof}

We then apply the above lemma to quadric fibrations.

\begin{lemma}\label{lem:pfister}
Let $\cha k\neq 2$. Use the notation of Lemma~\ref{lem:common-symbol}
with $F=k(B)$ and $\alpha_1,\alpha_2\ne0$, where $B$ is a smooth integral $k$-variety.
Let $f:V\to B$ be a dominant affine morphism, with $V$
integral and
$V_F\simeq Q\setminus\{\zeta=0\}$.
Suppose that $f$ extends to a morphism
$\bar f:\overline V\to\overline B$ of proper integral
$k$-varieties containing $V$ and $B$ as open subvarieties. For a nonempty open subset $B^\circ\subset B$, define $ \cQ^\circ
 =
 \{\Pi(\boldsymbol Z)-g_1g_2\zeta^2=0\}
 \subset\PP_{B^\circ}^{2^s}.$
Assume that there exists a $B^\circ\subset B$
that satisfies the following:
\begin{itemize}

\item
$a_1,\ldots,a_s,g_1,g_2$ are regular units
on $B^\circ$.

\item
There is a
$B^\circ$-isomorphism
$
\varphi :f^{-1}(B^\circ)
\simeq
\cQ^\circ\setminus\{\zeta=0\},
$

\item
For every nontrivial discrete valuation $\mu$ of $F/k$
centered in $\overline B\setminus B^\circ$, there is
$i\in\{1,2\}$ such that
$
 \alpha_i|_{F_\mu^h}=0.
$
\end{itemize}
Then $\Tortwo(V)=2$.
\end{lemma}

\begin{proof}
Let $U=f^{-1}(B^\circ)$ and $W=\overline V\setminus U$.
The first two hypotheses make $U$ smooth. Every valuation of $F/k$ has a
center on $\overline B$. At a center in $B^\circ$, both symbols have unit
entries. At any other center, hypothesis~(3) makes one symbol zero over
the henselization, hence unramified. Lemma~\ref{lem:common-symbol} gives
$$
 0\ne\alpha\in H^{s+1}_{\nr}
 \bigl(k(\overline V)/k,\mu_2^{\otimes(s+1)}\bigr).
$$

Blow up an ideal defining $W$ and take an odd-degree alteration
$\pi:Y\to\overline V$ with $Y$ smooth and proper, such that
$D=(\pi^{-1}(W))_{\mathrm{red}}$ is a divisor with smooth irreducible components
\cite[Theorem~2.1]{IT}.
The class $\gamma=\pi^*\alpha$ is unramified and nonzero by
\cite[Proposition~4.7(a), (b)]{Schreieder} and the odd degree of $\pi$. We verify its vanishing on $D$.

For a component $E\subset D$, let $R=\cO_{Y,\eta_E}$,
$\nu=\ord_E$ and $b=(\bar f\circ\pi)(\eta_E)$.
Let $L'=k(Y)$, so the residue field of $R$ is $k(E)$.
If $b\in B^\circ$, properness of $\cQ^\circ\to B^\circ$ extends the
generic point to $h:\Spec R\to\cQ^\circ$. Its closed point lies on
$\zeta=0$. Otherwise $\varphi^{-1}\circ h$ would map $\Spec R$ to $U$.
This map and the map induced by $\pi$ coincide on $\Spec L'$, hence
everywhere by separatedness of $\overline V$, contradicting
$\pi(\eta_E)\in W$. Thus the reduced Pfister form $\overline\Pi$
has a nonzero isotropic vector, and
$$
 \gamma|_{k(E)}
 = (\bar a_1,\ldots,\bar a_s)_2\cup(\bar g_i)_2=0
$$
by \cite[Theorem~2.2]{SchSmall}.

If $b\notin B^\circ$, the local inclusion
$\cO_{\overline B,b}\hookrightarrow R$ sends a nonzero element of its
maximal ideal to an element of positive valuation. Hence $\nu|_F$ is
nontrivial, and its normalization $\mu$ has center $b$.
Hypothesis~(3) gives $\alpha_i|_{F_\mu^h}=0$ for some $i$.
The local map $\cO_\mu\to R$ extends to
$\cO_\mu^h\to R^h$. This extension is injective: a uniformizer has
nonzero image, whereas every nonzero prime ideal of the henselian DVR
contains that uniformizer. It gives $F_\mu^h\hookrightarrow L'_\nu{}^h$,
so $\gamma|_{L'_\nu{}^h}=0$. Unramified specialization then gives
$\gamma|_{k(E)}=0$.

Remark~\ref{rmk:boundary} extends the vanishing from the generic points
of the smooth components of $D$ to all of $D$.
Lemma~\ref{lem:pairing} gives $\Tortwo(U)=2$.
Restriction sends $\delta_V$ to $\delta_U\ne0$, proving
$\Tortwo(V)=2$ without an additional affineness assumption on $V$.
\end{proof}

Our main construction in this subsection is the following, based on the Gram construction as in \cite[Theorems~3.4 and~3.14]{Hotchkiss}:

\begin{construction}\label{constr:gram-cycle}
Let $\Omega$ be an algebraically closed field of
characteristic different from two, let $t$ be
transcendental over $\Omega$, and let $
k=\overline{\Omega(t)}.
$
Fix an integer $s\geq3$ and independent base variables
$\boldsymbol z=(z_1,\ldots,z_s)$ over $k$.
Denote 
$
K=k(z_1,\ldots,z_s),$ $
z_I=\prod_{i\in I}z_i$ for 
$I\subseteq[s],$ and $
z_\varnothing=1.
$ Consider the coordinate symbol
$
\beta=(z_1,\ldots,z_s)_2
\in H^s(K,\mu_2^{\otimes s})
$
and its associated Pfister form. Since $-1$ is a square in $k$, its Pfister form has the unsigned
presentation
$$
\Pi_s\left(\sum_{I\subseteq[s]}c_Ie_I\right)
=
\sum_{I\subseteq[s]}z_Ic_I^2
$$
on $K^{2^s}$, with orthogonal coordinate basis
$(e_I)_{I\subseteq[s]}$. Introduce independent fiber coordinates
$X,Y_0,Y_1,\ldots,Y_s,$ and write $\boldsymbol Y=(Y_1,\ldots,Y_s).
$ Define the linear forms $L_1=(z_1-1)Y_0+z_2Y_1$, $L_i=z_1Y_0-z_{i-1}Y_{i-1}+z_{i+1}Y_i$ for $1<i<s$, and $L_s=z_1Y_0-z_{s-1}Y_{s-1}+Y_s$ in $Y_0,\boldsymbol Y$.
Define a $K$-linear map
$
\iota_t:K^{s+1}\longrightarrow K^{2^s}
$ by
$$
\iota_t(Y_0,\boldsymbol Y)
={}2z_1Y_0e_\varnothing
   +\sum_{i=1}^sL_ie_i
+2t\left(
Y_0e_{1,2,3}
+\sum_{j=1}^{s-1}Y_je_{j,j+1}
+Y_se_{1,s}
\right).
$$
Since $t\ne0$, the coefficients recover
$Y_0,Y_1,\ldots,Y_s$, respectively. Thus $\iota_t$ is injective. Define
$\cG_s(X,Y_0,\boldsymbol Y)=4X^2+\Pi_s\bigl(\iota_t(Y_0,\boldsymbol Y)\bigr)$.
Explicitly,
$$
\cG_s
={}4X^2+4z_1^2Y_0^2+\sum_{i=1}^s z_iL_i^2
+4t^2\left(
z_1z_2z_3Y_0^2
+\sum_{j=1}^{s-1}z_jz_{j+1}Y_j^2
+z_1z_sY_s^2
\right).
$$
In particular, $\cG_s$ is homogeneous of degree two
in the fiber coordinates $X,Y_0,\boldsymbol Y$, and
the coefficients of this quadratic form are polynomials
of degree at most three in $\boldsymbol z$. The equation defines a projective
quadric family
$\cQ_s=\{\cG_s=0\}\subset\bA^s_{k,\boldsymbol z}\times\PP^{s+1}_{k,[X:Y_0:Y_1:\cdots:Y_s]}$.
Note that the projection of $\cQ_s$ to the first factor has generic fiber $
Q_s=\{\cG_s=0\}\subset\PP_K^{s+1}.
$
We will use two different affine charts of this family.

On $Y_0\ne0$, divide the fiber coordinates by $Y_0$ and use their original
names. Let $\ell_i=L_i|_{Y_0=1}$ and
$G_s^\sharp=\cG_s|_{Y_0=1}$. On $X\ne0$, divide by $X$ and let
$G_s=\cG_s|_{X=1}$. These are the two affine equations used below.

The basis supports occurring in $\iota_t$ form the set
$$
\cT_s=
\{\varnothing\}
\cup\bigl\{\{i\}:1\leq i\leq s\bigr\}
\cup\bigl\{\{1,2,3\}\bigr\}
\cup\bigl\{\{j,j+1\}:1\leq j<s\bigr\}
\cup\bigl\{\{1,s\}\bigr\}.
$$
We call the index subsets occurring in $\cT_s$ the \emph{reserved supports}. All the reserved supports are distinct for $s\geq3$,
so
$$
|\cT_s|=1+s+1+(s-1)+1=2s+2.
$$
If $I\notin\cT_s$, the basis vector $e_I$ is
orthogonal to the image of $\iota_t$.

There is also an odd-dimensional version of the above construction. Assume $s\geq4$.
Notice $\{1,3\}\notin\cT_s$, so $e_{1,3}$ is an unused
orthogonal direction. Use one more independent coordinate $Z$ and
replace $\iota_t(Y_0,\boldsymbol Y)$ by
$\iota_t(Y_0,\boldsymbol Y)+2Ze_{1,3}$.
Since $\Pi_s(e_{1,3})=z_1z_3$, the corresponding
fiber-homogeneous equation is
$\cG_s+4z_1z_3Z^2=0$.
Taking the same two charts gives $G_s^\sharp+4z_1z_3Z^2=0$ and $G_s+4z_1z_3Z^2=0$.
Each of them now has $2s+2$ affine coordinates, hence has dimension $2s+1$.

Finally, we describe the special fiber at $t=0$.
Over $\Omega(z_1,\ldots,z_s)$, the linear map becomes
$$
\iota_0(Y_0,\boldsymbol Y)
=
2z_1Y_0e_\varnothing+\sum_{i=1}^sL_ie_i.
$$
The coordinate change $(Y_0,Y_1,\ldots,Y_s)\mapsto
(2z_1Y_0,L_1,\ldots,L_s)$ has triangular matrix with diagonal
$2z_1,z_2,\ldots,z_s,1$. Its determinant is nonzero, so $\iota_0$
is an isomorphism onto the span of $e_\varnothing,e_1,\ldots,e_s$.
Hence, after also replacing $2X$ by a new scalar coordinate, $\cG_s|_{t=0}$ is isometric to
$
\langle1,1,z_1,\ldots,z_s\rangle.
$
In the odd-dimensional version it is isometric to
$
\langle1,1,z_1,\ldots,z_s,z_1z_3\rangle.
$
Both forms are nondegenerate and isotropic, since
$-1$ is a square in $\Omega$. This implies that their projective quadrics are rational
over $\Omega(z_1,\ldots,z_s)$.

\end{construction}

The next lemma proves both of the affine charts are reference varieties.

\begin{lemma}\label{lem:gram-local}
Under Construction~\ref{constr:gram-cycle}, for any $
\cU\subseteq
\{I\subseteq[s]:I\notin\cT_s\},
$ choose independent coordinates
$\boldsymbol R=(R_I)_{I\in\cU}$, and define the following affine hypersurface
$$
V_{\cU}
=
\left\{
G_s+4\sum_{I\in\cU}z_IR_I^2=0
\right\}.
$$
Then $V_{\cU}$ is geometrically integral of dimension
$2s+|\cU|$ with
$\beta=(z_1,\ldots,z_s)_2\neq 0$ over $k(V_{\cU})$ and
$\delta_{V_{\cU}}\ne0 \in\CH_0\bigl((V_{\cU})_{k(V_{\cU})},\ZZ/2\bigr).$ In particular, $\Tortwo(V_{\cU})=2$ and $V_{\cU}$ is a reference variety with
coefficients $\ZZ/2$.
The same conclusions also hold for $V^{\sharp}_{\cU}$ obtained by replacing $G_s$ by $G_s^\sharp$. 
\end{lemma}

\begin{proof}
We adapt the complementary vectors and valuation calculation in
\cite[Lemma~3.8 and Proposition~3.11]{Hotchkiss} to both affine charts
and the extra orthogonal directions.
Let $h=|\cU|$ and let $V$ denote either affine chart.
The coordinate residues show that $\beta=(z_1,\ldots,z_s)_2\ne0$ over
$K$, so $\Pi_s$ is anisotropic. Let
$\phi(Y_0,\boldsymbol Y,\boldsymbol R)=\iota_t(Y_0,\boldsymbol Y)+2\sum_{I\in\cU}R_Ie_I$
and $q_{\cU}=4X^2+\Pi_s\circ\phi$.
The extra supports are orthogonal to the image of $\iota_t$, so $\phi$
is injective and $\Pi_s\circ\phi$ is anisotropic. In characteristic
different from two it is nondegenerate. Thus $q_{\cU}$ is nondegenerate
of dimension $s+h+2\geq5$, and its projective quadric $Q_{\cU}$ is smooth
and geometrically integral. The affine equations are primitive over
$k[\boldsymbol z]$: their constant coefficient on $X=1$, or their
$X^2$ coefficient on $Y_0=1$, is $4$. Gauss's lemma proves that $V$ is
integral, hence geometrically integral over $k$.

We first verify the valuation hypothesis of
Corollary~\ref{cor:affine-detection}. Let $(F,\mu)/k$ be henselian and
discretely valued, with $\mu|_{k^\times}=0$. Consider an $F$-point of $V$
with all $z_i\ne0$ and $\beta_F\ne0$.
Lemma~\ref{lem:laurent-values} shows that evaluation embeds
$K=k(\boldsymbol z)$ into $F$, and that every nonzero Laurent polynomial
$P=\sum_\nu c_\nu\boldsymbol z^\nu$ over $k$ satisfies
$$
 \mu(P)=\min_{\nu:c_\nu\ne0}\sum_i\nu_i\mu(z_i).
$$
In particular, if every $\mu(z_i)=0$, then $\mu|_{K^\times}=0$.

Let $\Sigma_j=\sum_{i=1}^jz_i$. Define vectors in the space of $\Pi_s$ by
$$
\begin{aligned}
r_0&=e_\varnothing-\frac{1}{tz_2z_3}e_{1,2,3},\\
r_j&=\frac{(1-\Sigma_j)e_\varnothing
             +2\sum_{i=1}^j e_i-t^{-1}e_{j,j+1}}
            {1+\Sigma_j},\\
r_s&=\frac{(1-\Sigma_s)e_\varnothing
             +2\sum_{i=1}^s e_i-(tz_1)^{-1}e_{1,s}}
            {1+\Sigma_s}.
\end{aligned}
$$
The middle formula has $1\leq j<s$. Substitution in $\iota_t$ gives
$\langle r_j,\iota_t(Y_0,\boldsymbol Y)\rangle_\Pi=0$ and
$\Pi_s(r_j)=1+(t^2D_j)^{-1}$ for $0\leq j\leq s$.
Here let $D_0=z_2z_3/z_1$, $D_j=(1+\Sigma_j)^2/(z_jz_{j+1})$ for
$1\leq j<s$, and $D_s=z_1(1+\Sigma_s)^2/z_s$.

Suppose $\mu|_K\ne0$. Let $\rho_i=\mu(z_i)$ and
$m_j=\min\{0,\rho_1,\ldots,\rho_j\}$.
The Laurent minimum formula gives $\mu(1+\Sigma_j)=m_j$, so
$$
 \mu(D_0)=\rho_2+\rho_3-\rho_1,
 \qquad \mu(D_j)=2m_j-\rho_j-\rho_{j+1},
 \qquad \mu(D_s)=\rho_1+2m_s-\rho_s,
$$
where the middle formula has $1\leq j<s$.
If all these values were nonnegative, the first path inequality would give
$\rho_2\leq2\min(0,\rho_1)-\rho_1=-|\rho_1|$.
Induction along the path gives $m_j=\rho_j$ for $2\leq j\leq s$ and
$$
 \rho_s\leq\cdots\leq\rho_2\leq-|\rho_1|.
$$
The last inequality gives $\rho_s\geq-\rho_1$, hence
$\rho_1\geq0$ and $\rho_2=\cdots=\rho_s=-\rho_1$.
Then $0\leq\mu(D_0)=-3\rho_1$ forces every $\rho_i=0$, a contradiction.
Thus $\mu(D_{j_0})<0$ for some $j_0$.

Hensel's lemma gives a unit $c$ with
$c^2=1+(t^2D_{j_0})^{-1}$, since the right side is congruent to $1$
and the residue characteristic is not two. The vector
$r'_{j_0}=c^{-1}r_{j_0}$ has Pfister value $1$ and is orthogonal to
$\phi$. Thus
$$
 (X,Y_0,\boldsymbol Y,\boldsymbol R)
 \longmapsto 2Xr'_{j_0}+\phi(Y_0,\boldsymbol Y,\boldsymbol R)
$$
is an injective isometry from $q_{\cU}$ into $\Pi_s$.
The point of $V$ would make $\Pi_s$ isotropic, contrary to
$\beta_F\ne0$. Hence $\mu|_{K^\times}=0$.

The rest of the affine coordinates are integral. To see this, inclusion and
reduction give $K\hookrightarrow\cO_\mu$ and $K\hookrightarrow\kappa(\mu)$.
On the two chart boundaries,
$$
\begin{aligned}
 q_{\cU}(0,Y_0,\boldsymbol Y,\boldsymbol R)
   &=\Pi_s\bigl(\phi(Y_0,\boldsymbol Y,\boldsymbol R)\bigr),\\
 q_{\cU}(X,0,\boldsymbol Y,\boldsymbol R)
   &=\Pi_s\bigl(2Xe_\varnothing+\phi(0,\boldsymbol Y,\boldsymbol R)\bigr).
\end{aligned}
$$
The second identity uses the zero $e_\varnothing$ component of
$\phi(0,\boldsymbol Y,\boldsymbol R)$. Both boundary maps are injective
and defined over $K$, so they stay injective after reduction.
If an affine coordinate had negative valuation, multiply its homogeneous
point vector by a common power of a uniformizer until all entries are
integral and at least one is a unit. The chart coordinate, originally $1$,
then reduces to zero. Either boundary map gives a nonzero isotropic
vector $w$ for $\overline\Pi_s$ over $\kappa(\mu)$.
Choose $J$ with $w_J\ne0$ and lift its other entries to $a_I\in\cO_\mu$.
The polynomial $f(T)=z_JT^2+\sum_{I\ne J}z_Ia_I^2$ satisfies
$\bar f(w_J)=0$ and $\bar f'(w_J)=2\bar z_Jw_J\ne0$.
Hensel's lemma lifts $w$ to a nonzero isotropic vector over $F$, again
contradicting $\beta_F\ne0$. This proves the valuation condition for both
charts.

It remains to prove nonvanishing over $k(V)$.
Over $\Omega(t)$, replace each added $R_I$ by $tR_I$ and let
$\widetilde q_{\cU}=\cG_s+4t^2\sum_{I\in\cU}z_IR_I^2$.
Let $\widetilde G$ be its equation on the relevant chart, and let
$\boldsymbol y$ be the other fiber coordinates, including all $R_I$.
At $t=0$ the extra terms disappear, and the coordinate change in
Construction~\ref{constr:gram-cycle} gives
$q_0=\cG_s|_{t=0}\simeq\langle1,1,z_1,\ldots,z_s\rangle$
over $K_0=\Omega(\boldsymbol z)$.
Its projective quadric is smooth, geometrically integral and rational,
since $-1$ is a square in $\Omega$.

Let $A_0=\Omega[z_1^{\pm1},\ldots,z_s^{\pm1}]$ and
$V_0=\Spec A_0[\boldsymbol y]/(\widetilde G|_{t=0})$.
Its generic fiber over $A_0$ is the relevant affine chart of $\{q_0=0\}$
times $\bA^h$. It is rational of dimension $s+h$.
The special polynomial has the unit coefficient $4$ described above,
so it is primitive over $A_0$. Gauss's lemma proves that $V_0$ is integral.
Its function field is
$\kappa_0\simeq K_0(w_1,\ldots,w_{s+h})$, with the $w_i$ algebraically
independent. The coordinate residues give
$$
 \partial_{z_s}\cdots\partial_{z_1}(z_1,\ldots,z_s)_2
 =1\in H^0\bigl(\Omega(w_1,\ldots,w_{s+h}),\ZZ/2\bigr).
$$
Thus $\beta_0=(z_1,\ldots,z_s)_2\ne0$ over $\kappa_0$.

For any finite extension $E/\Omega(t)$ in $k$, including inseparable
extensions, choose a DVR $R_E$ above $\Omega[t]_{(t)}$.
Its residue field is $\Omega$. Let $\pi_E$ be a uniformizer and let
$\cV_E=\Spec R_E[z_1^{\pm1},\ldots,z_s^{\pm1},\boldsymbol y]/(\widetilde G)$.
Its generic fiber is integral: after extension to $k$ it is an open
subset of the integral variety $V$. The reduction of $\widetilde G$ is
the nonzero equation of $V_0$, so $\pi_E$ does not divide $\widetilde G$
in the factorial polynomial ring. The quotient is torsion-free, hence
flat over $R_E$. Its embedding into the integral generic fiber makes it
a domain. At the generic point of the special fiber, the local
ring $A_E$ satisfies $A_E/(\pi_E)=\kappa_0$. It is a one-dimensional local
domain with maximal ideal $(\pi_E)$, hence a DVR. Let $M_E$ be its
fraction field.

Every $z_i$ is a unit of $A_E$. If $\beta|_{M_E}=0$, a nonzero isotropic
Pfister vector can be scaled to have integral coordinates with at least
one unit. Its reduction contradicts $\beta_0\ne0$.
Thus $\beta|_{M_E}\ne0$ for every finite $E/\Omega(t)$ in $k$.
Their filtered union is $k(V)$, so $\beta|_{k(V)}\ne0$.
Corollary~\ref{cor:affine-detection} now gives $\delta_V\ne0$.
\end{proof}

\begin{definition}\label{def:affine-quadric}
Let $B_0=\Spec R_0$ be an integral affine $k$-variety,
let $K=\Frac(R_0)$, and write $\eta_{B_0}=\Spec K$.
Let $Q\subset\PP^N_K$ be an integral projective quadric,
and let $H\subset\PP^N_K$ be a hyperplane not containing $Q$.
An \emph{affine quadric model of $Q\setminus H$}
is a dominant affine morphism $f:V\to B_0$, with $V$ integral,
together with a fixed $K$-isomorphism
$\varphi:V\times_{B_0}\Spec K\xrightarrow{\sim}Q\setminus H$.
In particular, $V$ is affine and $k(V)$ is identified with $K(Q)$.

An \emph{$H$ coordinate change} is the restriction to $\PP^N_K\setminus H$ of a projective
linear automorphism of $\PP^N_K$ preserving $H$.
\end{definition}

Before proving the next result, we note that the statement and the proof of Lemma~\ref{lem:reciprocity} also work for Milnor $K$-theory modulo $m$, 
with no invertibility assumption on $m$. Precisely, 
under the same hypotheses, we obtain an evaluation
$\CH_0(U,\ZZ/m)\longrightarrow K_q^M(L)/m$,
where $U$ is defined over $L$.
Indeed, replace $H^q(-,\mu_m^{\otimes r})$ by
$K_q^M(-)/m$ and replace $(g)_m\cup\alpha_C$ by
$\{g\}\alpha_C$. The result follows from
\cite[rule~R3f, Proposition~2.2~\textup{(RC)}
and Remark~2.4]{Rost1996}. 

We next prove the corresponding criterion for
Artin--Schreier--Kato cohomology in characteristic
two. Recall that for a regular ring $R$ of characteristic two, we denote $
H^{r+1,r}(R)
:=
H^1_{\et}(\Spec R,\Omega^r_{R,\log}).
$

\begin{lemma}\label{lem:kato-reciprocity}
Let $E/k$ be a field extension with $\cha E=2$,
and fix an integer $r\geq0$.
Let $U$ be a separated scheme of finite type over $E$,
with fixed classes
$a_p\in H^{r+1,r}(\kappa(p))$ at its closed points.
For every integral curve $C$ closed in $U$, let
$\nu:\widetilde C\to C$ be its normalization,
and choose a normal projective curve $\overline C$
containing $\widetilde C$ as a dense open subcurve.
At every closed point $w\in\overline C$, let $b_w=\res_{\kappa(w)/\kappa(\nu(w))}\bigl(a_{\nu(w)}\bigr)$ if $w\in\widetilde C$, and let $b_w=0$ if $w\in\overline C\setminus\widetilde C$.
Suppose that, for each such $C$, there exists a class
$
\alpha_C\in H^{r+1,r}(E(C))
$
extending over $\cO_{\overline C,w}$ for every
closed point $w\in\overline C$ and satisfying
$
s_w(\alpha_C)=b_w
$
at every such point. Then the fixed classes satisfy
Condition~\textup{(C)}, and induce a homomorphism
$$
\ev_U:\CH_0(U,\ZZ/2)
\longrightarrow H^{r+1,r}(E),
\qquad
[p]\longmapsto\cores_{\kappa(p)/E}(a_p).
$$
\end{lemma}
\begin{proof}
Fix an integral closed curve $C\subset U$, let $M=E(C)$ and choose
$g\in M^\times$. At $w\in\overline C$, let
$R_w=\cO_{\overline C,w}$ and $\kappa_w=\kappa(w)$.
The extension of $\alpha_C$ over this regular DVR is unique by
\cite[Corollary~2.3(1)]{Otabe}, and its reduction is $b_w$.

The Milnor action is given on symbols by
$[a;b_1,\ldots,b_r]\cdot\{g\}=[a;b_1,\ldots,b_r,g]$,
or by exterior multiplication by $d\log g$ in the differential
presentation \cite[Section~2]{Otabe}.
Let $\pi_w$ be a uniformizer, $n_w=\ord_w(g)$ and
$u_w=g\pi_w^{-n_w}$. Then $u_w\in R_w^\times$.
The product with $u_w$ extends over $R_w$, while the uniformizer formula
has residue $b_w$. Thus $\alpha_C\cdot\{g\}$ is tame and
$$
 \partial_w\bigl(\alpha_C\cdot\{g\}\bigr)=n_wb_w
$$
by \cite[Equation~(4.2) and Lemma~4.1]{Otabe}.
Proper-curve reciprocity \cite[Remark~4.5]{Otabe} gives
$$
 0=\sum_{w\in\overline C^{(1)}}
       \ord_w(g)\,\cores_{\kappa_w/E}(b_w)
  =\sum_{w\in\widetilde C^{(1)}}\ord_w(g)\,
       \cores_{\kappa_w/E}
       \bigl(\res_{\kappa_w/\kappa(\nu(w))}(a_{\nu(w)})\bigr).
$$
This is Condition~\textup{(C)}. The degree and composition identities for
restriction and corestriction, including inseparable residue extensions,
\cite[Section~2, properties~(C1)--(C3)]{Otabe}, show that these sums are
the evaluations of the pushforwards of principal divisors from the
normalization. This gives the homomorphism on ordinary
$\CH_0(U,\ZZ/2)$.
\end{proof}

The following criterion uses the two reciprocity
lemmas to affine quadric models.
\begin{proposition}\label{prop:charts}
Let $K/k$ be a rational function field of positive transcendence
degree, and let
$
Q=\{q=0\}\subset\PP^N_K
$
be an integral projective quadric.
Let $H\subset\PP^N_K$ be a hyperplane not containing $Q$.
If $\cha k\ne2$, let
$\alpha\in H^r(K,\mu_2^{\otimes r})$ be a Kummer symbol modulo two,
with $r\geq1$, and if $\cha k=2$, let
$\alpha\in H^{r+1,r}(K)$ be an Artin--Schreier--Kato symbol,
with $r\geq0$.
Let $\Pi$ be its associated nonsingular quadratic Pfister form
as in Definition~\ref{def:pfister}.
Assume the following conditions:
\begin{enumerate}
\item $
\alpha|_{K(Q)}\ne0.
$
\item
For every nontrivial discrete valuation $\mu$ of $K/k$
and every field extension $E/\widehat K_\mu$, $Q(E)\ne\varnothing$ implies $\alpha_E=0,$ where $\widehat K_\mu$ is the completion of $K$ at $\mu$.
\item
Writing $H=\PP(W_H)$ for a vector hyperplane
$W_H\subset K^{N+1}$, there is an injective $K$-linear map
$\jmath:W_H\hookrightarrow W_\Pi$ such that
$\Pi(\jmath(v))=q(v)$ for every $v\in W_H$, where $W_\Pi$
is the vector space of $\Pi$. Equivalently,  $q|_{W_H}$ is a subform of $\Pi$.
\end{enumerate}
Then every affine quadric model $f:V\to B_0$ of $Q\setminus H$ as in
Definition~\ref{def:affine-quadric} satisfies
$$
0\neq \delta_V \in\CH_0(V_{k(V)},\ZZ/2).
$$
\end{proposition}
\begin{proof}
Write $B_0=\Spec R_0$ and fix the isomorphism
$
\varphi:
V\times_{B_0}\Spec K
\xrightarrow{\sim}Q\setminus H.
$
Let $L=k(V)=K(Q)$ and, for a field extension $F/k$, write $A(F)=H^r(F,\mu_2^{\otimes r})$ if $\cha k\ne2$, and $A(F)=H^{r+1,r}(F)$ if $\cha k=2$.
For a closed point $p\in V_L$, consider its image under $V_L\to V\to B_0$ and let $a_p=\res_{\kappa(p)/K}(\alpha)$ if $p$ maps to $\eta_{B_0}$, and $a_p=0$ otherwise.
In the first case, the restriction uses the embedding
$K\hookrightarrow\kappa(p)$ induced by $f$.
We verify the hypotheses of
Lemma~\ref{lem:reciprocity} when $\cha k\ne2$,
and Lemma~\ref{lem:kato-reciprocity} when $\cha k=2$.

Let $C\subset V_L$ be an integral closed curve,
let $\nu:\widetilde C\to C$ be its normalization,
and let $\overline C$ be its normal projective
compactification.
If $C$ does not dominate $B_0$, no point of $C$
maps to $\eta_{B_0}$.
All fixed values on $C$ are then zero,
and the zero curve class satisfies these conditions.

Now suppose that $C$ dominates $B_0$.
Consider the embedding $K\hookrightarrow M=L(C)$ induced by the generic point of $C$ and its map to $B_0$. Let $
\alpha_C=\res_{M/K}(\alpha).
$ 
For a closed point $w\in\overline C$, let $R_w=\cO_{\overline C,w}$ be the DVR with
fraction field $M$ and residue field $\kappa_w$. Denote $v_w=\ord_w,$ and $\widehat M_w=\Frac(\widehat R_w).$ Note that $v_w$ is trivial on $L$, hence on $k$.
\begin{itemize}
\item If $v_w|_{K^\times}=0$, then every nonzero element of $K$ is a unit in $R_w$,
so inclusion and reduction give
$
K\hookrightarrow R_w,$ and $
K\hookrightarrow\kappa_w.
$
Pullback from $K$ extends $\alpha_C$ over $R_w$,
and its reduction is
$b_w=\res_{\kappa_w/K}(\alpha)$.

\item If $v_w|_{K^\times}$ is nontrivial, then its value group is $e\ZZ$ for some $e>0$, so
$\mu=e^{-1}v_w|_{K^\times}$
is a normalized discrete valuation of $K/k$ and it induces an embedding
$
\widehat K_\mu\hookrightarrow\widehat M_w.
$
The generic $M$-point of $Q$ gives a
$\widehat M_w$-point, and condition~(2) implies
$
\alpha_C|_{\widehat M_w}=0.
$

If $\cha k\ne2$, compatibility of residues with
completion gives $\partial_w(\alpha_C)=0$.
For a uniformizer $\pi_w$, its unramified
specialization is
$$
b_w=s_w(\alpha_C)
=
\partial_w\bigl((\pi_w)_2\cup\alpha_C\bigr).
$$
Hence $b_w=0$ by \cite[Equation~(3.3) and Lemmas~3.1--3.2]{Schreieder}.

If $\cha k=2$ and $r\geq1$, the zero completed class makes
$\alpha_C$ tame at $w$ with zero residue. The exact sequence
$$
0\longrightarrow H^{r+1,r}(R_w)
\longrightarrow H^{r+1,r}_{\mathrm{tame},w}(M)
\xrightarrow{\partial_w}H^{r,r-1}(\kappa_w)
\longrightarrow0
$$
of \cite[Equation~(4.1)]{Otabe} gives an extension over $R_w$.
Its pullback to $\widehat R_w$ is zero by regular-local injectivity into
$H^{r+1,r}(\widehat M_w)$ \cite[Corollary~2.3(1)]{Otabe}.
Reduction then gives $b_w=0$.

For $\cha k=2$ and $r=0$, write $\alpha_C=[a]$.
Vanishing over $\widehat M_w$ gives $\widehat b$ with
$a=\widehat b^2+\widehat b$. Choose $b\in M$ with
$v_w(b-\widehat b)>0$. Then $a+b^2+b$ lies in the maximal ideal of
$R_w$, so its Artin--Schreier torsor extends $\alpha_C$ over $R_w$
and has split special fiber. Again $b_w=0$.
\end{itemize}

Hence we have shown $\alpha_C$ is unramified at any point in $\overline{C}$. Now we check the specialization compatibility. Let $w\in\widetilde C$ and let $p=\nu(w)$.
The morphism to $B_0$ induces an inclusion
$R_0\hookrightarrow R_w$, since $C$ dominates $B_0$. Our goal is to show $b_w=\res_{\kappa_w/\kappa(p)}(a_p).$

If $p$ maps to $\eta_{B_0}$, every nonzero element
of $R_0$ has nonzero image in $\kappa(p)$ and hence
in $\kappa_w$.
It is a unit in $R_w$, and taking
quotients shows that $v_w|_K=0$.
Functoriality gives
$
b_w=\res_{\kappa_w/\kappa(p)}(a_p).
$

If $p$ does not map to $\eta_{B_0}$, choose
$0\ne u\in R_0$ vanishing at its base image.
Its image in $M$ is nonzero but its reduction
in $\kappa_w$ is zero, so $v_w(u)>0$.
This gives $b_w=0$ as in the second bullet case above,
which agrees with the fixed value $a_p=0$.
Thus the compatibility holds at every point in $\widetilde C$.

Now let $w\in\overline C\setminus\widetilde C$.
If $v_w|_K$ is nontrivial, we already have $b_w=0$. Suppose that it is trivial and abbreviate $R=R_w$. The embedding $K\hookrightarrow R$ makes $\Spec R$
a $K$-scheme.
Properness of $Q/K$ then extends the generic $M$-point
to a morphism
$h:\Spec R\longrightarrow Q$.
We claim that the closed point of $\Spec R$ maps to $Q\cap H$.

If not, $h$ factors through $Q\setminus H$. Together with $L\subset R$, this gives a morphism
$\Spec R\to V_L$ extending the generic map of $C$.
It factors through $C$, because $C$ is closed
and its defining equations vanish in $M$. Write $C=\Spec B$.
Its normalization is the spectrum of the integral
closure of $B$ in $M$, and the inclusion $B\subset R$ extends to that integral
closure because $R$ is integrally closed.
Thus the morphism lifts to $\widetilde C$.
After composing with
$\widetilde C\hookrightarrow\overline C$, it agrees
generically with the canonical map
$\Spec\cO_{\overline C,w}\to\overline C$.
Separatedness makes the two maps equal and hence $w\in\widetilde C$, a contradiction.

Hence the reduction is a $\kappa_w$-point
of $Q\cap H$.
A nonzero homogeneous vector representing this
point lies in $W_H\otimes_K\kappa_w$.
Then the scalar extension of $\jmath$ remains injective
and sends it to a nonzero isotropic vector for
$\Pi_{\kappa_w}$. This implies
$
\alpha_{\kappa_w}=0,
$
so $b_w=0$.

Now apply Lemma~\ref{lem:reciprocity} or
Lemma~\ref{lem:kato-reciprocity} according to the
characteristic. We have a homomorphism
$\ev_{V_L}:\CH_0(V_L,\ZZ/2)\longrightarrow A(L)$.
The diagonal point defining $\delta_V$ is
$L$-rational and maps to $\eta_{B_0}$.
Hence 
$
\ev_{V_L}(\delta_V)
=
\alpha|_{K(Q)}
\ne0
$
by condition~(1).
Thus $\delta_V\ne0$ and $\Tortwo(V)=2$.
\end{proof}

\begin{remark}\label{rmk:quadric-changes}
Use the hypotheses of Proposition~\ref{prop:charts}.
Birational changes of the base transport $K,\alpha,Q,H$ and their
valuations along a field isomorphism. An $H$ coordinate change with
linear representative $A$ transports the quadratic form and boundary map
to $q^{\mathrm{new}}=q\circ A^{-1}$ and
$\jmath^{\mathrm{new}}=\jmath\circ A^{-1}|_{W_H}$.
The three hypotheses hold for these data. Over a polynomial base ring,
clearing denominators and taking the primitive numerator preserves the
generic affine scheme, and Gauss's lemma gives an integral affine model.
Condition~(3) must be checked for the transported quadratic form before
multiplication by a scalar in $K^\times$.

The criterion also permits extra variables. Let
$q'(\boldsymbol Z,T_1,\ldots,T_h)$ be a homogeneous quadratic form with
$q'(\boldsymbol Z,0,\ldots,0)=q(\boldsymbol Z)$, where $h\geq1$.
Let $Q'=\{q'=0\}\subset\PP_K^{N+h}$ and let $H'$ restrict to $H$
on the original $\PP_K^N$. Suppose that $Q'$ satisfies condition~(2)
and that the restriction of $q'$ to the vector hyperplane of $H'$
embeds into $\Pi$. Then $Q'$ is integral and
$\alpha|_{K(Q')}\ne0$, so every affine quadric model of $Q'\setminus H'$
has nonzero generic-point class modulo two.

Indeed, let $q_i=q'(\boldsymbol Z,T_1,\ldots,T_i,0,\ldots,0)$ and
$Q_i=\{q_i=0\}$ for $0\leq i\leq h$.
A factorization of $q_i$ would give two linear factors whose restrictions
to $T_i=0$ are nonzero, since $q_{i-1}\ne0$. It would factor $q_{i-1}$.
Thus every $Q_i$ is integral by induction. The section
$Q_{i-1}\subset Q_i$ is an integral effective Cartier divisor. Its
local ring at the generic point is a one-dimensional local domain whose
maximal ideal is generated by the section equation, hence a DVR.
Its fraction and residue fields are $K(Q_i)$ and $K(Q_{i-1})$.
Properness of the projective Pfister quadric specializes a splitting
point from the first field to the second. Thus
$$
 \alpha|_{K(Q_i)}=0\quad\Longrightarrow\quad
 \alpha|_{K(Q_{i-1})}=0,
$$
which contradicts condition~(1). This proves nonvanishing over $K(Q')$.
To check condition~(2), it suffices to embed
$q'_{\widehat K_\mu}$ into $\Pi_{\widehat K_\mu}$ whenever
$\alpha_{\widehat K_\mu}\ne0$: the embedding stays injective over every
extension of $\widehat K_\mu$, so any point of $Q'$ splits the symbol.
\end{remark}

\subsection{Hermitian constructions in characteristic two}\label{subsec:hermitian}
In this subsection, we construct analogues of the affine quadric
models in characteristic 2. Then we use Proposition~\ref{prop:charts} to show that the constructions provide reference varieties.

\begin{definition}\label{def:hermitian-data}
Assume $\cha k=2$ and let $s\geq3$. Let $K=k(a,z_1,\ldots,z_s)$,
where the coordinates are algebraically independent over $k$.
Let $E=K(\theta)$, with $\theta^2+\theta=a$ and conjugation
$\bar\theta=\theta+1$. On $\bigoplus_{I\subseteq[s]}Ee_I$ use the
Hermitian form
$$
 h\left(\sum_Iu_Ie_I,\sum_Iv_Ie_I\right)=\sum_Iz_Iu_I\bar v_I.
$$
Write $h(v)=h(v,v)$ and $\Nm(u)=u\bar u=N_{E/K}(u)$.
Its associated quadratic form is $\Pi(a;\boldsymbol z)$ from
Definition~\ref{def:pfister}, since
$\Nm(P+\theta Q)=P^2+PQ+aQ^2$.
Orthogonality refers to the Hermitian form.
\end{definition}

\begin{construction}\label{constr:hermitian-cycle}
Let $h$ be the Hermitian form of Definition~\ref{def:hermitian-data}.
Choose $t\in k^\times$ and distinct subsets $J,I\subseteq[s]$
of size at least two, neither of which is a pair
$\{j,j+1\}$, $1\leq j<s$.
Choose $\kappa_0,\kappa_s\in K^\times$ and
$\epsilon\in\{1,\theta\}$, and define the following vectors,
with $1\leq j<s$ in the formula for $w_j$:
$$
 \begin{aligned}
 w_0={}&z_1e_0+\bar\theta(z_1+1)e_1
      +z_1\bar\theta\sum_{i=2}^s e_i+t\theta\kappa_0e_J,\\
 w_j={}&z_{j+1}e_j+z_je_{j+1}+t\theta e_{j,j+1},\\
 w_s={}&e_s+t\epsilon\kappa_se_I.
 \end{aligned}
$$
Write $W_t=\langle w_0,\ldots,w_s\rangle_E$ and
$\sigma_t=h|_{W_t}$. Denote $p_0=\kappa_0e_J$ and $p_s=\kappa_se_I$. For $1\leq j\leq s$, let
$\Sigma_j=\sum_{i=1}^jz_i$, $U_j=1+\Sigma_j$ and $T_j=\Sigma_j+aU_j^2$.
Let $u_j=U_j\theta e_0+\sum_{i=1}^je_i$.
The following vectors and rational functions will be used in the
valuation calculation. In the formula for $r_j$, let $1\leq j<s$:
$$
 \begin{aligned}
 r_0&=e_0+\frac{z_1e_J}{t\bar\theta\kappa_0z_J},\\
 r_j&=\frac{u_j+(t\bar\theta)^{-1}e_{j,j+1}}{1+U_j\theta},\\
 r_s&=\frac{u_s+z_se_I/(t\bar\epsilon\kappa_sz_I)}{1+U_s\theta},\\
 \end{aligned}
$$
Let $D_0=\frac{\kappa_0^2z_J}{z_1^2}$, $D_j=\frac{T_j}{z_jz_{j+1}}$ for $1\leq j<s$, and $D_s=\frac{\Nm(\epsilon)\kappa_s^2z_IT_s}{az_s^2}$.
\end{construction}
For $0\leq j\leq s$, one checks directly that
$$
 h(r_j,W_t)=0,\qquad h(r_j)=1+(t^2aD_j)^{-1}.
$$
We choose the data so that one of the correction terms has sufficiently
large valuation. The choices below follow the method of
\cite[Sections~3.1--3.3]{Hotchkiss}.

\begin{definition}\label{def:endpoint-choices}
The \emph{admissible choices} in Construction~\ref{constr:hermitian-cycle}
are the following four families.
\begin{enumerate}
\item For the standard choice, let $J=\{1,2,3\}$, $\kappa_0=1$ and
$I=\{1,s\}$. Let $\epsilon=\theta$ and $\kappa_s=1$.
\item For the modified choice, let $J=\{1,2,3\}$, $\kappa_0=z_1/a$ and
$I=\{1,s\}$. Let $\epsilon=\theta$ and $\kappa_s=z_1^{-1}$.
\item For choice $b$, let $J=\{1,2,3\}$, $\kappa_0=z_1/a$ and
$I=\{b,s\}$. Let $\epsilon=\theta$ and $\kappa_s=a^{-1}$.
Assume either $b=1$ and $s\geq4$, or $2\leq b\leq s-2$ and $s+b\geq7$.
\item For the weighted choice, let $\kappa_0=a^{h_0}$ and $\kappa_s=a^r$,
where $h_0,r\in\ZZ$ and $s\in J\cup I$.
Let $\ell=|J\setminus\{1\}|$ and
$d_J=\sum_{j\in J\setminus\{1\}}(j-1)$. Assume $d_J+2h_0\geq2\ell$.
\end{enumerate}
For any of these choices, let
$\cT=\{\varnothing,\{1\},\ldots,\{s\},J,\{1,2\},\ldots,\{s-1,s\},I\}$
be the set of supports.
\end{definition}

\begin{definition}\label{def:hermitian-equations}
For the standard admissible choice as in Definition~\ref{def:endpoint-choices}, let $P_1,\ldots,P_s,X$ be independent
coordinates and define the following coefficients,
with $1\leq j<s$ in the formula for $\lambda_j$:
$$
 \begin{aligned}
 \lambda_0&=z_1^2+az_1+az_1^2\Sigma_s+t^2az_1z_2z_3,\\
 \lambda_j&=z_jz_{j+1}(z_j+z_{j+1}+t^2a),\\
 \lambda_s&=z_s(1+t^2az_1),
 \end{aligned}
$$
Let $H_s^1=X^2+a\lambda_0+\sum_{j=1}^s\lambda_jP_j^2$ and
$H_s^0=H_s^1|_{X=a}$.
For arbitrary admissible data, take independent coordinates $X,Q_0,R$,
the pairs $P_j,S_j$, and the pairs $U_S,V_S$. Define
$$
 q=\Nm\bigl(X+\theta(z_1Q_0+R)\bigr)
 +h\left(\theta Q_0w_0+\sum_{j=1}^s(P_j+\theta S_j)w_j
          +\sum_{S\notin\cT}(U_S+\theta V_S)e_S\right).
$$ In the affine chart $Q_0=1$, we use the same letters for the affine coordinates.
\end{definition}

The next lemma verifies condition~\textup{(2)} of
Proposition~\ref{prop:charts} for these forms.

\begin{lemma}\label{lem:endpoints}
Use Construction~\ref{constr:hermitian-cycle}, with any
admissible choice of
Definition~\ref{def:endpoint-choices}. Let
$\alpha=[a;z_1,\ldots,z_s],$ and $h=\Pi(a;\boldsymbol z).$ Let
$W_t^+=W_t\oplus\bigoplus_{S\notin\cT}Ee_S$.
Let $\mu$ be a nontrivial discrete valuation of $K/k$,
denote $F=K_\mu^h$, and suppose that $\alpha_F\ne0$.
Then there exists
$
r_\mu\in\bigoplus_{S\in\cT}(E\otimes_KF)e_S
$
such that
$
h(r_\mu)=1,$ and $
h(r_\mu,(W_t^+)_F)=0.
$
In particular, there is an injective $F$-linear
isometry

$$
\Nm_{E\otimes_KF/F}
\perp h|_{(W_t^+)_F}
\hookrightarrow h_F,
\qquad
(u,w)\longmapsto ur_\mu+w.
$$
Hence, if a quadratic form $q$ over $K$
is a subform of
$
\Nm_{E/K}\perp h|_{W_t^+},
$
then its projective quadric $Q=\{q=0\}$ satisfies
condition~\textup{(2)} of Proposition~\ref{prop:charts}. This works for the homogenization of $H_s^i$
and the form $q$ of Definition~\ref{def:hermitian-equations}.
\end{lemma}

\begin{proof}
Fix $\mu$ and $F=K_\mu^h$ with $\alpha_F\ne0$.
Then $h_F$ is anisotropic, see the paragraph after Definition~\ref{def:pfister}. In particular,
$
E_F:=E\otimes_KF=F(\theta)
$
is a field, because otherwise its norm form $[1,a]$,
which is a subform of $h_F$, would be isotropic.
We use the same notations $h$ and $\Nm$
for their scalar extensions to $F$.
Let
$A=\max\{0,-\mu(a)\},$ $\rho_i=\mu(z_i),$ and $m_j=\min\{0,\rho_1,\ldots,\rho_j\}.$ Since $\mu$ is trivial on $k$, we have $\mu(t)=0$. We prove the lemma within a few steps:
\begin{itemize}
\medskip
    \item[step 1:] In this step, we produce an estimate for orthogonal values of $h$. Notice that if $\mu(a)>0$, the polynomial $T^2+T+a$ has
the simple root $0$ modulo the maximal ideal.
Hensel's lemma would give a root in $F$, contrary
to $E_F$ being a field. Hence
$
\mu(a)=-A\leq0.
$ Let $\eta\in F$ satisfy $\mu(\eta)>A$.
Then $\mu(a\eta)>0$, so the polynomial
$
a\eta T^2+T+1
$
has the simple root $1$ after reduction.
Hensel's lemma gives $u\in\cO_\mu$ satisfying
$
u+a\eta u^2=1.
$
Hence,
$$
\Nm(1+\eta u\theta)
=
1+\eta u+a\eta^2u^2
=1+\eta(u+a\eta u^2)=
1+\eta.
$$
Thus every $1+\eta$ with $\mu(\eta)>A$ is a norm
from $E_F$.

Let $v_1,\ldots,v_n$ be nonzero mutually orthogonal vectors,
and let $f_i=h(v_i)$. Then
$$
\min_i\mu(f_i)\leq\mu\left(\sum_i f_i\right)
\leq\min_i\mu(f_i)+A.
$$
Only the upper bound needs proof. If it fails, choose $j$ with
$\mu(f_j)$ minimal and let $\eta=(\sum_i f_i)/f_j$.
There is $c\in E_F^\times$ with $\Nm(c)=1+\eta$. Orthogonality gives
$$
h\left(cv_j+\sum_{i\ne j}v_i\right)
=(1+\eta)f_j+\sum_{i\ne j}f_i=0.
$$
This vector has nonzero pairing $cf_j$ with $v_j$, contradicting
anisotropy.

\item[step 2:] In this step, we reduce the construction of $r_\mu$
to a valuation inequality. Since $T_j=\Nm(1+U_j\theta)\ne0$,
all denominators in Construction~\ref{constr:hermitian-cycle} are nonzero.
Directly, $h(u_j,w_0)=0$. Among the other pairings,
$h(u_j,w_j)=z_jz_{j+1}$ for $j<s$, $h(u_s,w_s)=z_s$, and all
others vanish. Thus
$$
h(r_j,W_t)=0,\qquad h(r_j)=1+(t^2aD_j)^{-1}.
$$
Each $r_j$ has support in $\cT$, so it is orthogonal to $(W_t^+)_F$.
If $\mu(D_j)<0$, then $\mu((t^2aD_j)^{-1})=A-\mu(D_j)>A$.
By step~1, $h(r_j)=\Nm(c)$ for some $c\in E_F^\times$.
Then $r_\mu=c^{-1}r_j$ has the desired properties.
It remains to prove that some $D_j$ has negative valuation.

\item[step 3:] In this step, we treat $A=0$.
The equality in step~1 gives $\mu(U_j)=m_j$.
If $m_j<0$, the term $aU_j^2$ has strictly smaller valuation than
$\Sigma_j$ in $T_j=\Sigma_j+aU_j^2$.
If $m_j=0$, the polynomial $X^2+X+\bar a$ is irreducible over the
residue field, so the reduction of $\Nm(1+U_j\theta)$ is nonzero.
Thus $\mu(T_j)=2m_j$ in both cases.

Suppose every $\mu(D_j)\geq0$. The path inequalities give
$\rho_j+\rho_{j+1}\leq2m_j$, hence
$\rho_2\leq2\min(0,\rho_1)-\rho_1=-|\rho_1|$.
Induction gives
$$
\rho_s\leq\cdots\leq\rho_2\leq-|\rho_1|,
\qquad m_j=\rho_j.
$$
The latter equality holds for $2\leq j\leq s$. We check the four choices.
\begin{itemize}
\item In the standard case, $\mu(D_s)=\rho_1+\rho_s\geq0$
forces $\rho_1\geq0$ and $\rho_2=\cdots=\rho_s=-\rho_1$.
Then $\mu(D_0)=-\rho_1+\rho_2+\rho_3=-3\rho_1\geq0$, so all values vanish.
\item In the modified case, $\mu(D_s)=\rho_s-\rho_1\geq0$
forces $\rho_1\leq0$ and $\rho_2=\cdots=\rho_s=\rho_1$.
Then $\mu(D_0)=\rho_1+\rho_2+\rho_3=3\rho_1\geq0$, with the same conclusion.
\item For choice $b$, one has $\mu(D_0)=\rho_1+\rho_2+\rho_3$
and $\mu(D_s)=\rho_b+\rho_s$.
If $b=1$, the latter gives $\rho_1\geq0$ and
$\rho_2=\cdots=\rho_s=-\rho_1$, while the former gives $-\rho_1\geq0$.
If $b\geq2$, the two nonpositive values $\rho_b,\rho_s$ have
nonnegative sum, so $\rho_s=0$ and all values vanish.
\item In the weighted case, powers of $a$ and $\Nm(\epsilon)$ are units, so
$$
\mu(D_0)=\sum_{j\in J}\rho_j-2\rho_1,
\qquad \mu(D_s)=\sum_{i\in I}\rho_i.
$$
If $\rho_1>0$, the first value is at most $-(\ell+2)\rho_1$ when
$1\notin J$, and at most $-(\ell+1)\rho_1$ otherwise, a contradiction.
If $\rho_1<0$, all $\rho_i$ are negative, contradicting $\mu(D_s)\geq0$.
Thus $\rho_1=0$. Since $s\in J\cup I$ and the other values are nonpositive,
one of these two inequalities forces $\rho_s=0$, hence all values vanish.
\end{itemize}

All coefficients $a,z_1,\ldots,z_s$ are now units, so the reduction
$\overline h$ is a nonsingular quadratic form in $2^{s+1}$ variables.
The residue field $\kappa$ satisfies $\operatorname{trdeg}_k\kappa\leq s$,
because $\mu$ is nontrivial on $K/k$ and henselization leaves the residue
field fixed. It is a $C_s$-field by the opening paragraph of
Subsection~\ref{subsec:quadric-models}. Since $2^{s+1}>2^s$,
$\overline h$ is isotropic. Its nonzero zero lifts by Hensel's lemma,
contradicting anisotropy of $h_F$. Thus some $\mu(D_j)<0$.

\item[step 4:] In this step, we treat $A>0$. Since $F$ is henselian, let $\nu$ be the unique extension of $\mu$ to $E_F$,
normalized by $\nu|_F=\mu$, see \cite[Chapter~II, \S6, Theorem~(6.2)]{NeukirchANT}. The same theorem also shows conjugation preserves $\nu$.
Then $\theta^2+\theta=a$ gives
$
\nu(\theta)=-A/2.
$
We make the convention that $U_0=1$ and define
$C_j=U_j+\theta^{-1}$ for $0\leq j\leq s$.
Let $x_j=\nu(C_j)$. Then $x_0=0$ and $C_j+C_{j-1}=z_j$.
Since $1+U_j\theta=\theta C_j$, we have 
$$\mu(T_j)=\mu(\Nm(1+U_j\theta))=\mu(\Nm(\theta C_j))=-A+\nu(C_j\overline{C_j})=-A+2x_j.$$
We show that $x_j\leq A$ for every $j$.
If $x_j>A$, since $A>0$, we have 
$$\nu(C_j)=x_j>A>\frac{A}{2}=\nu(\theta^{-1}).$$ Then $\mu(U_j)=\nu(U_j)=\nu(C_j+\theta^{-1})=\min\{\nu(C_j),\nu(\theta^{-1})\}=A/2$.
As a result, $\mu(aU_j^2)=0,$ and $\mu(T_j)=-A+2x_j>A.$ However,
$
T_j=\Sigma_j+aU_j^2=h(U_j\theta e_0)+\sum_{i=1}^jh(e_i)
$
is a sum of orthogonal values.
By step~1,
$
\mu(T_j)
\leq
\min\{0,\rho_1,\ldots,\rho_j\}+A
\leq A,
$
a contradiction.

Now we show $x_1<A$.
If $\rho_1\ne0$, then $C_1=C_0+z_1$ gives
$
x_1=\min(0,\rho_1)<A.
$
Suppose that $\rho_1=0$ and $x_1=A$.
Then
$$
\mu(U_1)=\mu(C_1+\theta^{-1})=\min\{A,A/2\}=A/2,
\qquad
\mu(T_1)=A.
$$
Let
$
a'=\frac{T_1}{U_1^2}.
$
Then $a'$ is a unit, and
$$
a+a'=a+\frac{z_1+aU^2_1}{U_1^2}=\frac{z_1}{U_1^2}=\frac{U_1+1}{U_1^2}=U_1^{-2}+U_1^{-1}.
$$
Hence $\theta':=\theta+U_1^{-1}$ satisfies
$(\theta')^2+\theta'=a'$ and generates the same
quadratic extension $E_F/F$. By step~1 applied with parameter $a'$, 
$
z_1=1+U_1
$
is a norm, since $\mu(U_1)=A/2>0=\max\{0,-\mu(a')\}$.
Choose $c\in E_F^\times$ such that
$\Nm(c)=z_1$.
Then
$$
h(ce_0+e_1)=z_1+z_1=0,
$$
contradicting anisotropy of $h$.
Thus $x_1<A$ in every case.
In fact, the same discussion proves that 
$$
x_1-\rho_1<A.
$$
Suppose that every $\mu(D_j)\geq0$.
For $1\leq j<s$, the inequalities are
$$
2x_j\geq A+\rho_j+\rho_{j+1}.
$$
For $j=1$, they give
$$
\rho_2\leq2x_1-\rho_1-A<x_1.
$$
Thus $C_2=C_1+z_2$ has valuation $x_2=\rho_2$.
Inductively, for $2\leq j<s$,
$$
\rho_{j+1}\leq2x_j-\rho_j-A=x_j-A<x_j,
\qquad
x_{j+1}=\rho_{j+1}.
$$
Hence, for $2\leq j\leq s$,
$$
\rho_j=x_j
\leq2x_1-\rho_1-(j-1)A
<
(3-j)A-|\rho_1|.
$$
\begin{itemize}
    \item For the standard admissible choices,
    $\mu(D_s)=x_s+\rho_1-A,$ but the bound discussed above implies 
    $x_s+\rho_1-A
<
(2-s)A+\rho_1-|\rho_1|
\leq0.$ This contradicts $\mu(D_s)\geq0$.
\item For the modified admissible choices, $\mu(D_s)=x_s-\rho_1-A,$ but the bound discussed above gives $x_s-\rho_1-A
<
(2-s)A-\rho_1-|\rho_1|
\leq0.$ Again, this contradicts $\mu(D_s)\geq0$.
\item For choice $b$,
$
\mu(D_s)=A+\rho_b+x_s.
$
If $b=1$, the assumption $s\geq4$ gives
$$
\mu(D_s)
<
(4-s)A+\rho_1-|\rho_1|
\leq0.
$$
If $b\geq2$, the assumption $s+b\geq7$ gives
$$
\mu(D_s)
<
(7-s-b)A-2|\rho_1|
\leq0.
$$
We obtain a contradiction.

\item Finally, for the weighted admissible choices,
recall that $\ell=|J\setminus\{1\}|$ and
$d_J=\sum_{j\in J\setminus\{1\}}(j-1)$.
Summing the strict bounds above gives
$$
\sum_{j\in J\setminus\{1\}}\rho_j
<(2\ell-d_J)A-\ell|\rho_1|.
$$
Using the indicator notation of Notation~\ref{not:coordinates}, we get
$$
\begin{aligned}
\mu(D_0)
&=-2h_0A-(2-\delta_{1\in J})\rho_1
  +\sum_{j\in J\setminus\{1\}}\rho_j\\
&<(2\ell-d_J-2h_0)A-(2-\delta_{1\in J})\rho_1-\ell|\rho_1|\leq0.
\end{aligned}
$$
The last inequality uses $d_J+2h_0\geq2\ell$ and
$\ell\geq2-\delta_{1\in J}$, which follows from $|J|\geq2$.
This contradicts $\mu(D_0)\geq0$.
\end{itemize}
This proves some $\mu(D_j)<0$ also when $A>0$.
\medskip
\item[step 5:] In this step, we prove injectivity of the isometry. By step~2 we have $r_{\mu}$. For $u\in E_F$ and $w\in(W_t^+)_F$, orthogonality gives
$$
h(ur_\mu+w)
=
\Nm(u)+h(w).
$$
The map $(u,w)\mapsto ur_\mu+w$ is injective. In fact, 
if $ur_\mu+w=0$, pairing with $r_\mu$ gives
$
0=h(ur_\mu+w,r_\mu)=u,
$
and then $w=0$.
This proves the above isometric embedding
of quadratic forms.
\medskip
\item[step 6:] In this step, we check the rest of the statements of the lemma.
Let $q'$ be a quadratic subform of
$\Nm_{E/K}\perp h|_{W_t^+}$, let $Q'=\{q'=0\}$,
and fix a nontrivial discrete valuation $\mu$ of $K/k$.
If $\alpha_{\widehat K_\mu}=0$, then
$\alpha_L=0$ for every extension
$L/\widehat K_\mu$.
Otherwise $K_\mu^h\hookrightarrow\widehat K_\mu$ gives
$\alpha_{K_\mu^h}\ne0$.
The embedding in step~5 gives
$$
q'_{\widehat K_\mu}\hookrightarrow h_{\widehat K_\mu}.
$$
The embedding is injective over every $L/\widehat K_\mu$.
Thus an $L$-point of $Q'$ gives a nonzero isotropic vector
of $h_L$. This is condition~\textup{(2)} of
Proposition~\ref{prop:charts}.

For the standard choice, the homogenization of $H_s^1$ is
$q_s^1=X^2+a\lambda_0Q_0^2+\sum_j\lambda_jP_j^2$.
It is the restriction of $\Nm\perp h|_{W_t^+}$ under
$$
(X,Q_0,\boldsymbol P)\longmapsto
\left(X,\theta Q_0w_0+\sum_jP_jw_j\right).
$$
This map is injective because $w_0,\ldots,w_s$ are $E$-linearly independent.
The identity follows from $\Nm(X)=X^2$ and
$h(\theta Q_0w_0+\sum_jP_jw_j)=a\lambda_0Q_0^2+\sum_j\lambda_jP_j^2$.
The homogenization of $H_s^0$ is $q_s^0=q_s^1(aQ_0,Q_0,\boldsymbol P)$,
so it too is a subform.

For arbitrary admissible data, the full-norm form is the restriction
under $(X,Q_0,R,\ldots)\mapsto(u,w)$, where
$$
\begin{aligned}
u&=X+\theta(z_1Q_0+R),\\
w&=\theta Q_0w_0+\sum_j(P_j+\theta S_j)w_j
  +\sum_{S\notin\cT}(U_S+\theta V_S)e_S.
\end{aligned}
$$
The vectors $w_0,\ldots,w_s$ and the unused $e_S$ are $E$-linearly
independent. Hence $w=0$ gives $Q_0=0$, $P_j+\theta S_j=0$
and $U_S+\theta V_S=0$. Since $1,\theta$ are independent over $K$,
all these coordinates vanish. The equality $u=0$ then gives
$X=R=0$. Thus the map is injective, as desired.\qedhere
\end{itemize}
\end{proof}

The equations of Definition~\ref{def:hermitian-equations} define the
following affine reference varieties:

\begin{proposition}\label{prop:nonvanishing}
Let $\Omega$ be algebraically closed of characteristic two,
let $t$ be transcendental over $\Omega$, and let
$k=\overline{\Omega(t)}$.
The following affine varieties $V$ from
Definition~\ref{def:hermitian-equations} are geometrically
integral and satisfy $0\neq\delta_V\in\CH_0(V_{k(V)},\ZZ/2):$
\begin{enumerate}
\item The hypersurfaces $\{H_s^i=0\}$ for $i=0,1$ and their enlargements
$$
H_s^i+\sum_{I\in\cU}z_I(U_I^2+U_IV_I+aV_I^2)=0,
$$
where $\cU\subseteq\{I\subseteq[s]:I\notin\cT\}$ for the standard choice.
Any subset of the added coordinates may be omitted.
The hypersurface $\{H_s^i=0\}$ has total degree five and dimension $2s+i$.
\item For every admissible choice, the hypersurface defined
      by the primitive numerator of $q|_{Q_0=1}$.
\end{enumerate}
Furthermore, the form in \textup{(2)} still defines a reference variety after the following modifications: letting any
      subset of $R,S_j,U_S,V_S$ be equal to zero, taking
      birational base coordinate changes, and taking $H$ coordinate
      changes, with $H=\{Q_0=0\}$, as in
      Remark~\ref{rmk:quadric-changes}\textup{(1)}.
\end{proposition}
\begin{proof}
Let $K=k(a,z_1,\ldots,z_s)$ and
$\alpha=[a;z_1,\ldots,z_s]$.
We verify the conditions of Proposition~\ref{prop:charts}.
Lemma~\ref{lem:endpoints} proves condition~\textup{(2)}.

For condition~\textup{(3)} of Proposition~\ref{prop:charts}, let $q$ be the homogeneous quadratic form of
Definition~\ref{def:hermitian-equations},  and let
$Q=\{q=0\}\subset\PP^N_K$. The generic affine fiber is $Q\setminus H$, with
equation $q|_{Q_0=1}=0$.
To verify condition~\textup{(3)}, we restrict the
homogeneous form $q$ to the vector hyperplane $Q_0=0$. Define
$$
\jmath:
(X,R,\boldsymbol P,\boldsymbol S,\boldsymbol U,\boldsymbol V)
\longmapsto{}
(X+\theta R)e_0
+\sum_{j=1}^s(P_j+\theta S_j)w_j
+\sum_{S\notin\cT}(U_S+\theta V_S)e_S.
$$
Since $h(e_0)=1$ and $Ee_0$ is orthogonal to the
other summands, we have
$$
q|_{Q_0=0}=h\circ\jmath.
$$
As $\jmath$ is injective, this proves condition~\textup{(3)} for $q$. For  $\{H_s^i=0\}$, use the standard admissible choice.
Their homogenizations are
$$
q_s^1=X^2+a\lambda_0Q_0^2+
       \sum_{j=1}^s\lambda_jP_j^2,\qquad
q_s^0=(a^2+a\lambda_0)Q_0^2+
       \sum_{j=1}^s\lambda_jP_j^2.
$$
Thus their restrictions to $Q_0=0$ are
$X^2+\sum_j\lambda_jP_j^2$ and
$\sum_j\lambda_jP_j^2$, respectively.

We check that the $K$-linear maps
$$
\jmath_1(X,\boldsymbol P)=Xe_0+\sum_{j=1}^sP_jw_j,
\qquad
\jmath_0(\boldsymbol P)=\sum_{j=1}^sP_jw_j
$$
identify these boundary forms with subforms of $h$. The calculation of the scalar forms in Lemma~\ref{lem:endpoints}
gives $h(w_j)=\lambda_j$. For adjacent indices, the off-diagonal
pairings are $z_i z_{i+1}z_{i+2}$ when $i+1<s$ and $z_{s-1}z_s$
for the final pair, and all other pairings vanish. They lie in $K$,
so their two cross terms cancel in characteristic two.
Since $h(e_0)=1$ and $h(e_0,w_j)=0$ for $j\geq1$,
we obtain $h\circ\jmath_i=q_s^i|_{Q_0=0}$ for $i=0,1$.

To prove injectivity, suppose
$Xe_0+\sum_jP_jw_j=0$.
The coefficient of $e_{j,j+1}$ is $t\theta P_j$
for $1\leq j<s$, and the coefficient of $e_{1,s}$
is $t\theta P_s$.
These supports are distinct, and $t\theta\ne0$,
so every $P_j$ vanishes.
The coefficient of $e_0$ then gives $X=0$.
Thus $\jmath_1$ is injective, and so is $\jmath_0$ because it is a restriction of 
$\jmath_1$.
This proves condition~\textup{(3)} of
Proposition~\ref{prop:charts} for both scalar models.

An added summand is $h((U_I+\theta V_I)e_I)$ with $I\notin\cT$.
It is orthogonal to $W_t$ and to every complementary vector above.
Adding these vectors extends both the local and the boundary maps
injectively, also after omitting any subset of these coordinates.
Once nonvanishing for $H_s^i$ is proved below,
Remark~\ref{rmk:quadric-changes}\textup{(2)} gives integrality and
symbol nonvanishing for the enlarged generic quadric. A common base
factor of the enlarged equation would divide the original integral
scalar equation. Thus its polynomial is primitive, and Gauss's lemma
and Proposition~\ref{prop:charts} give the affine conclusion.

We are left with proving integrality and condition~\textup{(1)} of Proposition~\ref{prop:charts}.
\begin{itemize}
\item[step 1:] In this step, we prove integrality and nonvanishing for
the scalar models. Consider $H_s^0$ over $\Omega(t)$ and its specialization at $t=0$.
Take the change of coordinates
$Z_1=a(z_1+1)+z_2P_1$, $Z_i=az_1+z_{i-1}P_{i-1}+z_{i+1}P_i$ for $1<i<s$, and $Z_s=az_1+z_{s-1}P_{s-1}+P_s$,
which has determinant $z_2\cdots z_s$.
The special fiber $H_s^0|_{t=0}$ becomes
$$
a^2+az_1^2+\sum_{i=1}^s z_iZ_i^2=0.
$$
Let $b=a/z_1$, $c_i=Z_i/z_1$, and $\tau=b+c_1^2$.
Dividing by $z_1^2$ gives
$$
\tau z_1=b^2+\sum_{i=2}^s z_ic_i^2.
$$
In $\Omega[b,c_1,\ldots,c_s,z_2^{\pm1},\ldots,z_s^{\pm1}]$,
the element $\tau$ is a prime element as it is monic and linear in the variable $b$. Besides, $\tau$ does not divide $b^2+\sum_{i=2}^s z_ic_i^2$.
Thus this linear polynomial $\tau z_1=b^2+\sum_{i=2}^s z_ic_i^2$ in $z_1$ is primitive and
irreducible. Thus the open subset $az_1\cdots z_s\ne0$ of the special fiber is integral,
with function field $M_0=\Omega(b,c_1,\ldots,c_s,z_2,\ldots,z_s)$,
obtained by solving the equation for $z_1$.

Now we calculate the restriction of the symbol $\alpha$ to $M_0$. Recall that the symbol $[u;v_1,\ldots,v_s]$ is represented by
$u\,d\log v_1\wedge\cdots\wedge d\log v_s$ modulo exact forms and Artin-Schreier relations. In particular, the symbols are additive in the first
entry and multiplicative in each logarithmic entry,
and
$$
[u^2;v_1,\ldots,v_s]=[u;v_1,\ldots,v_s],\qquad [u;u,v_2,\ldots,v_s]=0,\qquad
[vc^2;v,v_2,\ldots,v_s]=0,
$$
see \cite[Section~1.1]{ChapmanMcKinnie2018} and
\cite[Section~4, before Equation~(4.2)]{Otabe}. It follows that
$$
\begin{aligned}
[a;z_1,z_2,\ldots,z_s]
&=[bz_1;b,z_2,\ldots,z_s]\\
&=\left[b^2+\sum_{i=2}^s z_ic_i^2;
         \tau,z_2,\ldots,z_s\right]\\
&=[b;\tau,z_2,\ldots,z_s].
\end{aligned}
$$
The first equality uses $[bz_1;bz_1,z_2,\ldots,z_s]=0$.
For the second, $d\tau=db$ shows that
$[bz_1;b,z_2,\ldots,z_s]$ and
$[\tau z_1;\tau,z_2,\ldots,z_s]$ have the same differential
representative $z_1\,db\wedge d\log z_2\wedge\cdots\wedge d\log z_s$.
For the last, the terms with first entry $z_ic_i^2$ vanish
because $z_i$ is a logarithmic entry, and the Artin--Schreier
relation replaces $b^2$ by $b$.
At the discrete valuation $\tau=0$, the elements
$b=\tau+c_1^2$ and $z_2,\ldots,z_s$ are units,
with $\bar b=c_1^2$.
Thus the symbol is tame there, and the residue formula
\cite[Equation~(4.2)]{Otabe} gives
$$
\partial_\tau
\bigl([a;z_1,\ldots,z_s]\bigr)
=
[c_1^2;z_2,\ldots,z_s]
=
[c_1;z_2,\ldots,z_s].
$$
Taking iterated residues at $z_2=0,\ldots,z_s=0$ of it
then gives
$$
[c_1]\in
\Omega(c_1,\ldots,c_s)/
\wp\bigl(\Omega(c_1,\ldots,c_s)\bigr),
$$
where $\wp(u)=u^2+u.$ This class is nonzero because at the valuation $c_1=\infty$,
$c_1$ has a pole of order one. But a coboundary $u^2+u$ has either a pole of even order or no pole at this valuation. Hence the residue is nonzero, and we know
$$
[a;z_1,\ldots,z_s]|_{M_0}\ne0.
$$
The coefficients of $P_1^2$ and $P_s^2$ in
$H_s^0|_{t=0}$ are $z_1z_2(z_1+z_2)$ and $z_s$,
which are coprime. Hence the special fiber is irreducible as otherwise there would be a factor that becomes a unit after localizing at $az_1\cdots z_s$, but then it would divide both of the coefficients of $P_1^2$ and $P_s^2$, a contradiction.
Both $H_s^0$ and its special fiber at $t=0$ have total degree five.
The coefficients of $h$ are regular at the generic point of the
special fiber, and its reduction is anisotropic by the symbol calculation.
Then Lemma~\ref{lem:constant} applied to this equation and to $h$ proves
geometric integrality of $\{H_s^0=0\}$ and that the restriction of
$\alpha$ to its function field is not zero.

Finally, let $q_s^i$ be the fiber homogenizations from
Lemma~\ref{lem:endpoints}.
After the invertible change $X'=X+aQ_0$, one has
$$
q_s^1(X'+aQ_0,Q_0,\boldsymbol P)
=(X')^2+q_s^0(Q_0,\boldsymbol P).
$$
Thus $q_s^1$ is an enlargement of $q_s^0$. So Remark~\ref{rmk:quadric-changes}\textup{(2)} proves the same property for $\{H_s^1=0\}$ as well.
Here the affine equation is primitive since its coefficient of $X^2$
is one, so Gauss's lemma gives the affine conclusion.
The same identity at $t=0$, together with the integral-divisor argument
in the proof of that remark, proves integrality and symbol nonvanishing
for the special fiber of $H_s^1$.
\item[step 2:] In this step, we treat the full-norm form and its
coordinate omissions. First let
$g_t$ be the polynomial obtained by letting
$R=S_j=U_S=V_S=0$.
For every admissible choice, over $K_0=\Omega(a,z_1,\ldots,z_s)$, the same invertible
change from the $P_j$ to the $Z_j$ as in the last case gives
$$
g_0=X^2+z_1X+\sum_{j=1}^s z_jZ_j^2.
$$
All perturbation terms vanish at $t=0$ for every admissible choice.
We first prove that $g_0$ is irreducible.
After localizing at $X$, let $c_j=Z_j/X$ and
$D=1+\sum_{j=1}^s z_jc_j^2$.
Since $X$ is invertible, the equation $g_0=0$ is equivalent to
$XD+z_1=0$, and
$$
\frac{K_0[X^{\pm1},c_1,\ldots,c_s]}{(XD+z_1)}
\simeq K_0[c_1,\ldots,c_s,D^{-1}].
$$
The right side is an integral domain, with function field
$M_0=K_0(c_1,\ldots,c_s)$.
A factor of $g_0$ that becomes a unit after localizing at $X$
must be a scalar times a power of $X$.
But $X$ does not divide $g_0$, since
$\sum_jz_jZ_j^2\ne0$.
Thus $g_0$ itself is irreducible in $K_0[X,Z_1,\ldots,Z_s]$.

Again the restriction of $\alpha$ to $M_0=K_0(c_1,\ldots,c_s)$ is nonzero.
Indeed, again taking iterated residues at $z_1=0,\ldots,z_s=0$ gives
$
[a]\in\Omega(a,c_1,\ldots,c_s)/
\wp\bigl(\Omega(a,c_1,\ldots,c_s)\bigr),
$
which is nonzero by the same pole argument at $a=\infty$.
Thus the special quadratic form associated with $h$ is anisotropic.

Let $D'/\Omega(t)$ be a finite extension contained in $k$,
and choose a DVR $R_{D'}\subset D'$ above $t=0$, with
uniformizer $\pi$ and residue field $\Omega$.
Then
$A=R_{D'}[a,z_1,\ldots,z_s]_{(\pi)}$
is a DVR with fraction field $D'(a,\boldsymbol z)$ and
residue field $K_0$.
All nonzero rational functions in $a,\boldsymbol z$ with
coefficients in $\Omega$ are units of $A$.
Thus $g_t$ has coefficients in $A$. Both $g_t$ and its reduction have total degree two in
$X,P_1,\ldots,P_s$.
If $g_t$ factored over $D'(a,\boldsymbol z)$, Gauss's lemma
would give primitive factors over $A$.
Since the degree is the same under reduction, both factors would still have positive degree.
This would factor the irreducible special polynomial,
a contradiction. 

Moreover, $g_t$ is monic in $X$, so
$A[X,\boldsymbol P]/(g_t)$ is an integral, flat $A$-algebra, whose special fiber is integral.
Let
$B=\bigl(A[X,\boldsymbol P]/(g_t)\bigr)_{(\pi)}$ and $M_{D'}=\Frac\bigl(D'(a,\boldsymbol z)[X,\boldsymbol P]/(g_t)\bigr)$.
The local ring $B$ at the generic point of that fiber is a
one-dimensional local domain with maximal ideal $(\pi)$, hence a DVR.
Its fraction field is $M_{D'}$ and its residue field is $M_0$.

The coefficients of $h$ are units in $B$.
If $h$ were isotropic over $M_{D'}$, we could rescale
a nonzero isotropic vector so that its coordinates are
integral and at least one is a unit. Then
its reduction would be a nonzero isotropic vector of the
special $h$, contradicting the fact that $\alpha\neq 0$ over
$M_0$.
Hence $g_t$ is irreducible and $\alpha$ is nonzero after
every such finite extension $D'$. Since
every factorization and isotropic vector over $k$ uses only
finitely many constants, the same conclusions hold over $k$.

Finally, add any subset of the coordinates $R,S_j,U_S,V_S$.
Letting the added coordinates be zero gives $g_t$, and the local
and boundary embeddings above are still injective after any omissions.
Remark~\ref{rmk:quadric-changes}\textup{(2)} proves integrality and
condition~\textup{(1)} for the enlarged generic quadric.
Taking its primitive numerator gives an integral affine
quadric model by Gauss's lemma, hence a geometrically integral
variety because $k$ is algebraically closed.
\end{itemize}
Proposition~\ref{prop:charts} completes the proof.
The coordinate changes follow from Remark~\ref{rmk:quadric-changes}\textup{(1)}.
\end{proof}

\setcounter{theorem}{46}
\begin{proposition}\label{prop:norm-extension}

The following affine reference varieties with coefficients $\ZZ/2$ exist over an algebraically closed field $k$ of
characteristic two.
\begin{enumerate}
\item For $(A,B)=(3,3)$ and for every
$4\leq A\leq B\leq\binom{A+1}{3}-2A+3$,
there is a polynomial
$F\in k[x_1,\ldots,x_A,y_1,\ldots,y_B]$ with $\deg_{\boldsymbol x}F\leq3$ and $\deg_{\boldsymbol y}F\leq2$,
whose zero locus is geometrically integral of dimension
$A+B-1$.
\item There is such a hypersurface in
$\bA^3_k\times\bA^3_k\times\bA^3_k$ whose equation has
total degree at most $3,2,2$ in the three groups of coordinates, respectively.
\end{enumerate}
\end{proposition}
\begin{proof}
Choose $\Omega,t,k$ as in Proposition~\ref{prop:nonvanishing}.
\begin{itemize}
\item[step 1:] In this step, we construct the equations for $A\geq4$.
Let $s=A-1$ and use the standard admissible choice.
Let $\cP$ be the set of pairs in $[s]$ other than
$\{1,2\},\ldots,\{s-1,s\},\{1,s\}$, and let $\cB$
be the set of triples other than $\{1,2,3\}$.
In the defining equation of $q$ in
Definition~\ref{def:hermitian-equations}, let all $S_j=0$, let $V_J=0$ for $J\in\cB$,
and let $U_S=V_S=0$ for $S\notin\cT\cup\cP\cup\cB$.
Rename the $U_J$ by $W_J$ for $J\in\cB$. We get the following equation
$$
H_s^1+z_1X+az_1^2+XR+aR^2
+\sum_{I\in\cP}z_I\Nm(U_I+\theta V_I)
+\sum_{J\in\cB}z_JW_J^2=0.
$$
Make the birational base coordinate change
$A_0=a^{-1}$, $b_i=z_i/a$,
and do the following change of coordinates
$X=a^2Y$, $P_i=aQ_i$, $R=aZ$.
Then replace $U_{I}, V_{I}, W_{J}$ by $aU_I, aV_{I},aW_{J}$ separately. Multiplying by $A_0^5$ gives
$$
\begin{aligned}
F_s={}&A_0Y^2+A_0^2b_1(Y+1)
+b_1^2\sum_{i=1}^s b_i+t^2b_1b_2b_3\\
&+\sum_{j=1}^{s-1}b_jb_{j+1}
 (b_j+b_{j+1}+t^2)Q_j^2
+b_s(A_0^2+t^2b_1)Q_s^2\\
&+A_0^2(YZ+Z^2)
+\sum_{I\in\cP}b_I\bigl(A_0(U_I^2+U_IV_I)+V_I^2\bigr)
+\sum_{J\in\cB}b_JW_J^2.
\end{aligned}
$$
Every monomial has total degree at most three in
$(A_0,b_1,\ldots,b_s)$ and at most two in the other
coordinates.
Since the $Y^2$ coefficient is $A_0$ and the constant
coefficient is not divisible by $A_0$, the polynomial is primitive.
This is also true after any omissions of the additional coordinates,
which are covered by Proposition~\ref{prop:nonvanishing}.
Then by Remark~\ref{rmk:quadric-changes}\textup{(1)}, this equation defines
a reference variety.
\item[step 2:] In this step, we count the coordinates and construct
the degree $(3,2,2)$ hypersurface in $\bA^3_k\times\bA^3_k\times\bA^3_k$.
Note in the last construction, the group of coordinates $(A_0,b_1,\ldots,b_s)$ has size $s+1=A$. The rest of the coordinates always contain $Y,Q_1,\ldots,Q_s$, and may
also contain $Z$, the two coordinates for each $I\in\cP$,
and $W_J$ for each $J\in\cB$.
Its largest size is
$$
A+1+2\left(\binom{s}{2}-s\right)+\binom{s}{3}-1
=\binom{A+1}{3}-2A+3.
$$
Omitting the additional coordinates independently gives every size
from $A$ to this upper bound.
For $A=4$, regroup the nine coordinates as
$(b_1,b_2,b_3)$, $(A_0,Q_1,Q_2)$, and $(Y,Q_3,Z)$.
The polynomial has degrees at most $(3,2,2)$
in this division, proving \textup{(2)}.
\item[step 3:] In this step, we treat $(A,B)=(3,3)$.
Let $K=k(a,x,y)$ and $\alpha=[a;x,y]$.
Let $E=K(\theta)$, with $\theta^2+\theta=a$, and let
$h=\langle1,x,y,xy\rangle_E$ be the diagonal Hermitian
form, also viewed as its associated quadratic form.
Define
$w_0=xe_0+\bar\theta(1+x)e_1+\bar\theta xe_2$ and
$w_1=tye_1+\bar\theta e_{1,2}$. Let $B=a+t^2y$ and $L=xyB$, and define
$$
\begin{aligned}
    G_t^+=&X^2+X(x+R)+aR^2+a^2x+a^2x^2(x+y)\\
       &+L(P^2+PS+aS^2)+taxy(1+x)S.
\end{aligned}
$$
We prove that $\{G_t^+=0\}$ and its restrictions obtained
by letting either or both of $R,S$ be zero are reference
varieties. The $S=0$ case will
give the desired equation for our purpose. We check the conditions for Proposition~\ref{prop:charts} are satisfied. 

For Condition~(2) of Proposition~\ref{prop:charts}, it is convenient to consider the
homogeneous quadratic form
$$
q^+=\Nm\bigl(X+\theta(xQ+R)\bigr)
+h\bigl((W+\theta Q)w_0+(P+\theta S)w_1\bigr)+xRW.
$$
Its restriction at $W=0,Q=1$ is $G_t^+$.
We now show that $q^+$ embeds locally into $h$.
Fix a henselian discretely valued field $F$ extending
$K$, with $\mu(t)=0$, and with a $C_2$ residue field, such that
$\alpha_F\ne0$.
Then $h_F$ is anisotropic and $E_F/F$ is a field.
Write $\Nm_b(U,V)=U^2+UV+bV^2$ and let
$D=a(1+x)^2+L$ and $H=t^2(1+x)^2+xB$.
We use the convention $\mu(0)=+\infty$.
Here $B\ne0$, because otherwise $y=a/t^2=\Nm_{E_F/F}(\theta/t)$ is a norm and the vector $\frac{\theta}{t}e_0+e_2\neq 0$ makes $h_F$ isotropic, a contradiction.
Let $M=axD/B$. The vectors $w_0-xe_0,w_1$ are independent,
as their $e_2$- and $e_{1,2}$-coefficients show.
Their Hermitian pairings form the matrix
$$
\begin{pmatrix}
ax\bigl((1+x)^2+xy\bigr)&txy(1+x)\bar\theta\\
txy(1+x)\theta&L
\end{pmatrix}.
$$
Since $\theta\bar\theta=a$ and $B+t^2y=a$, its determinant is
$$
\begin{aligned}
&ax\bigl((1+x)^2+xy\bigr)L+at^2x^2y^2(1+x)^2\\
&\qquad=ax^2y\bigl(a(1+x)^2+xyB\bigr)=ax^2yD=LM.
\end{aligned}
$$
If this determinant were zero, a nonzero vector in their span
would be orthogonal to both vectors and hence to itself,
contradicting anisotropy. Thus $D\ne0$. Let $\psi=xB/D$.

Replace $X,R$ by $X+xW,R+xQ$.
The terms involving $e_0$ then reduce to
$\Nm_a(X,R)+ax^2Q^2$.
Completing the Hermitian square by
$C'=P+\theta S+\frac{txy(1+x)\bar\theta}{L}(W+\theta Q)$
gives
$$
q^+\simeq\Nm_a\perp L\Nm_a\perp M\Nm_{a+\psi}.
$$
The same orthogonalization of $h$ gives
$$
h\simeq\cR\perp M\cR,\qquad
\cR=\Nm_a\perp L\Nm_a.
$$
Indeed, the first three Hermitian coefficients are $1,L,M$ and
the fourth is $LM$ modulo norms, since the determinant of
$h$ is the square $x^2y^2$.

We claim that $[\psi;aL]=0$.
We use the norm criterion for Artin--Schreier symbols
\cite[Proposition~9.2.5 and Corollary~4.7.5]{GilleSzamuely}.
The norm calculation in step~1 of Lemma~\ref{lem:endpoints}
is valid for every Artin--Schreier parameter.
Also $\mu(a)\leq0$, and $\mu(\psi)>0$ already implies
$[\psi;aL]=0$ by Hensel's lemma.
Assume then that $\mu(\psi)\leq0$.
When $H\ne0$, let $\zeta=t^2D/(aH)$.
The identities
$$
aH+t^2D=xB^2,\qquad
\Nm_\psi(1+x,\psi^{-1})=yH/B
$$
give
$$
\frac{aL}{1+\zeta}=a^2yH/B,
\qquad
[\psi;aL]=[\psi;1+\zeta].
$$
The first quotient is a nonzero norm for the parameter $\psi$.

If $\mu(B)<0$, then
$\mu(D)\geq\mu(x)+\mu(B)>\mu(x)+2\mu(B)$.
Thus $H\ne0$ and $\mu(aH)=\mu(x)+2\mu(B)$, so
$
\mu(\zeta)+\mu(\psi)=-\mu(B)>0.
$
The near-one norm calculation gives $[\psi;1+\zeta]=0$.
If $\mu(B)>0$, the expression
$y=(a/t^2)(1+B/a)$ makes $y$ a norm for $a$ by the
same calculation, since $\mu(B/a)>-\mu(a)$.
This contradicts $\alpha_F\ne0$.

We are left with the case $\mu(B)=0$.
If $\mu(a)=0$, then $\mu(y)\geq0$.
For $\mu(x)\ne0$, the two terms of $D=a(1+x)^2+xyB$
have different valuations and
$\mu(D)=2\min(0,\mu(x))$, giving $\mu(\psi)>0$.
If $x,y$ are both units, the reduced nonsingular form
$h$ has dimension eight over a $C_2$ field.
It is isotropic and its nonzero zero lifts by Hensel's
lemma, again a contradiction. Hence $\mu(x)=0<\mu(y)$.
If $\mu(1+x)>0$, the near-one calculation makes $x$
a norm. Otherwise $1+x$ is a unit and
$$
\bar\psi=\frac{\bar x}{(1+\bar x)^2}
=\wp\left(\frac1{1+\bar x}\right).
$$
Hensel's lemma makes $\psi$ Artin--Schreier, proving the claim.

If $\mu(a)<0$, then $\mu(y)=\mu(a)$.
When $\mu(x)\ne0$, the first summand of $D$ has strictly
smaller valuation and again $\mu(\psi)>0$.
Suppose $\mu(x)=0$, then $H$ is integral.
If it is a unit, we have 
$\mu(D)=\mu(a)<0$, hence $\mu(\psi)>0$.
And if $\mu(H)>0$, reduction gives
$$
\bar x\bar B=t^2(1+\bar x)^2,\qquad
\bar x=\bar B
 \left(\frac{\bar x}{t(1+\bar x)}\right)^2.
$$
Here $1+\bar x\ne0$.
This is a nonzero solution of the norm equation for
the unit parameter $B$.
Its derivative with respect to the first norm coordinate
is nonzero, so Hensel's lemma lifts it and gives $[B;x]=0$.
Since $a=B+t^2y$ and $[t^2y;y]=0$, it follows that
$[a;x,y]=0$, a contradiction.
This proves $[\psi;aL]=0$ in every case.

To get the embedding, use $[\psi;aL]=0$ and
\cite[Chapter~1, Exercise~4]{GilleSzamuely} to choose a nonzero
solution of $U^2+UV+\psi V^2+aLZ^2=0$.
The form $\cR$ is anisotropic, so $aL$ is not a square.
Indeed, if $aL=u^2$, the vector $(u,0;0,1)$ is isotropic for $\cR$.
Thus $V\ne0$, and division by $V^2$ gives
$\psi=f^2+f+aLr^2$, with $f=U/V$ and $r=Z/V$.
In the quaternion algebra $[a;L]$, use generators $i,j$ with
$i^2+i=a$, $j^2=L$ and $ji=(i+1)j$.
Let $v=ij$ and $i'=i+f+rv$.
Then $v^2=aL$, $(i')^2+i'=a+\psi$, and
$i'v+vi'=v$.
Taking the reduced norm gives
$$
\cR\simeq\Nm_{a+\psi}\perp aL\Nm_{a+\psi}.
$$
As a result, we have $q^+\subset\cR\perp M\cR\simeq h$.

For each nontrivial discrete valuation of $K/k$, the
residue field has transcendence degree at most two and
is $C_2$. The above argument is valid over
its completion whenever $\alpha$ survives.
An injective embedding persists after any field extension,
so a point of the quadric splits $\alpha$.
This proves condition~\textup{(2)} of
Proposition~\ref{prop:charts}.
For $W=0$, the boundary $Q=0$ embeds into $h$ by
$$
(X,R,P,S)\longmapsto
(X+\theta R)e_0+(P+\theta S)w_1.
$$
This map is injective and preserves the form.
It also works after omitting $R$ or $S$, this proves condition~\textup{(3)}.

For integrality and condition~\textup{(1)}, first let $S=0$,
let $\tau=t^{-1}$, and replace $P=\tau Z$.
At $\tau=0$, the change of coordinate $T=yZ+a(1+x)$ gives
$$
G_0=B_0+a^2x^2y,\qquad
B_0=\Nm(X+\theta R)+xX+xT^2.
$$
After localizing at $a,x,y$, this equation is linear in $y$
with unit coefficient $a^2x^2$.
It is integral, with function field
$\Omega(a,x,X,R,T)$, also with $R$ omitted.
The symbol there is $[a;x,B_0]$.
The slots $a,x,B_0$ are units at the generic divisors $T=0$ and
$R=0$. Thus unramified specialization at these divisors gives
$[a;x,X(X+x)]$. If $R$ is absent, only the specialization at $T=0$ is needed.
Its residues at $X=0$ and then $x=0$ give $[a]\ne0$.
Hence the original symbol is nonzero. Since the coefficient of $X^2$ is $1$,
no factor of the special polynomial becomes a unit after localizing at $a,x,y$.
Thus the special polynomial is irreducible. Both the transformed polynomial and its specialization
have total degree five.
Lemma~\ref{lem:constant} then gives integrality and shows $\alpha\ne 0$ over $k=\overline{\Omega(\tau)}$.
Then with Remark~\ref{rmk:quadric-changes}\textup{(2)}, 
Proposition~\ref{prop:charts} proves the conclusion for all choices of $R,S$.

Finally, take $S=0$ and keep $R$. Let $A_0=a^{-1}$, $b=x/a$ and $c=y/a$.
Substitute $X=a^2Y$, $R=aZ$ and $P=aQ$.
Multiplication by $A_0^5$ gives
$$
A_0Y^2+A_0^2\bigl(bY+YZ+Z^2+b\bigr)
+b^2(b+c)+bc(1+t^2c)Q^2.
$$
It has total degree at most three in $(A_0,b,c)$
and at most two in $(Y,Z,Q)$.
It is primitive, since its $Y^2$ coefficient is $A_0$
and its constant coefficient is not divisible by $A_0$.
The changes preserve the reference variety property by
Remark~\ref{rmk:quadric-changes}\textup{(1)}, which proves
the case $(A,B)=(3,3)$.
\end{itemize}

\end{proof}

The next lemma provides us more allowed modifications for our constructed reference varieties.

\begin{lemma}\label{lem:cycle-charts}
After the following changes to their defining equations, the
constructed affine reference varieties are still reference varieties:
\begin{enumerate}
\item In characteristic different from two, replace
$2tY_se_{1,s}$ in $\iota_t$ by $2tz_1^{-1}Y_se_{1,s}$,
and use either affine chart of
Construction~\ref{constr:gram-cycle}.
\item In characteristic two, use the second or third choice of
Definition~\ref{def:endpoint-choices} and the chart $Q_0=1$ of
$$
X^2+h\left(\theta Q_0w_0+\sum_{j=1}^sP_jw_j\right),
$$
including its restriction $X=aQ_0$, where $h$ is the form of
Definition~\ref{def:hermitian-data}.
\end{enumerate}
\end{lemma}

\begin{proof}

\begin{itemize}
\item[step 1:] In this step, we prove \textup{(1)}. We use the complementary vectors $r_0,\ldots,r_{s-1}$ in the proof of
Lemma~\ref{lem:gram-local}, and let
$r_s=\frac{(1-\Sigma_s)e_0+2\sum_{i=1}^se_i-t^{-1}e_{1,s}}{1+\Sigma_s}$.
Their quadratic values are
$\Pi_s(r_j)=1+(t^2D_j)^{-1}$, where $D_0,\ldots,D_{s-1}$ are the same as in the proof of
Lemma~\ref{lem:gram-local}, and 
$D_s=\frac{(1+\Sigma_s)^2}{z_1z_s}$.
Let $\mu$ be a nontrivial valuation of $K/k$ at which the
coordinate symbol is nonzero over the henselization, and let
$\rho_i=\mu(z_i)$. If every $\mu(D_j)\geq0$, the inequalities
for $1\leq j<s$, together with
Lemma~\ref{lem:laurent-values}, give
$$
\rho_s\leq\cdots\leq\rho_2\leq-|\rho_1|,
\qquad
m_j=\min(0,\rho_1,\ldots,\rho_j)=\rho_j.
$$
The latter equality holds for $j\geq2$.
Then $\mu(D_s)\geq 0$ gives $\rho_s-\rho_1\geq0$.
Hence $\rho_1\leq0$ and
$\rho_2=\cdots=\rho_s=\rho_1$.
Now $\mu(D_0)=\rho_1\geq0$ forces every $\rho_i=0$.
The same discussion as in Lemma~\ref{lem:laurent-values} then implies that
$\mu$ is trivial on $K$, a contradiction.
Thus some $\mu(D_j)<0$, and Hensel's lemma gives a vector
of value one orthogonal to the image of the changed map,
as in Lemma~\ref{lem:gram-local}.
This proves condition~\textup{(2)} of
Proposition~\ref{prop:charts}.

The coefficients at the distinct basis indices used to prove
injectivity are still nonzero. Thus the maps on the homogenizing
hyperplanes are injective, proving condition~\textup{(3)} of
Proposition~\ref{prop:charts}.
At $t=0$, the changed term vanishes, so the specialization is
the same as in Lemma~\ref{lem:gram-local}.
Its proof gives integrality and condition~\textup{(1)}  of
Proposition~\ref{prop:charts}.
The complementary vectors have support in $\cT_s$, so the same
argument permits the terms $4z_IR_I^2$ with $I\notin\cT_s$,
using Remark~\ref{rmk:quadric-changes}\textup{(2)}.
Taking primitive numerators and applying
Proposition~\ref{prop:charts} then proves our case \textup{(1)}.

\item[step 2:] In this step, we prove \textup{(2)}.
Lemma~\ref{lem:endpoints} gives the local embeddings for both choices.
On $Q_0=0$, the scalar forms are restrictions of $h$ under the injective
maps $(X,\boldsymbol P)\mapsto Xe_0+\sum_jP_jw_j$ and
$\boldsymbol P\mapsto\sum_jP_jw_j$. Thus conditions~\textup{(2)} and
\textup{(3)} of Proposition~\ref{prop:charts} hold.

At $t=0$, the equations specialize to the same scalar forms as
$H_s^1$ and $H_s^0$ in Proposition~\ref{prop:nonvanishing}.
For every finite extension $D'/\Omega(t)$, use the DVR
$R_{D'}[a,z_1,\ldots,z_s]_{(\pi)}$ from that proof.
The factors $z_1/a$, $z_1^{-1}$ and $a^{-1}$ are units there.
Both the equation and its reduction have degree two in the fiber
coordinates. The irreducible special fiber and its anisotropic Pfister
form show that the generic quadric over $D'(a,\boldsymbol z)$ is integral
and that $\alpha$ remains nonzero in its function field.
Passing to $k$ gives the same conclusions, as in that proposition.

The added norm coordinates have supports outside $\cT$, so they are
orthogonal to both the local complementary vectors and the boundary
images. Their coordinate maps remain injective.
Remark~\ref{rmk:quadric-changes}\textup{(2)} preserves integrality and
symbol nonvanishing. Taking primitive numerators and applying
Proposition~\ref{prop:charts} proves \textup{(2)}.
\end{itemize}
\end{proof}
By virtue of Lemma~\ref{lem:cycle-charts}, we give another construction of reference varieties that is used in Section~\ref{sec:bounds} over characteristic two.  
\begin{construction}\label{constr:cycle-coordinates}
For the modified admissible choice of Definition~\ref{def:endpoint-choices},
let $u=z_1/a$, $Z=z_1$ and $x_1=u$.
Let $x_i=z_i/a$ for $i\geq2$. Start with the equation of Lemma~\ref{lem:cycle-charts}\textup{(2)},
replace $X$ by $a^2uX$ and $P_s$ by $auP_s$. After division by $a^3$, the equation becomes
$$
\begin{aligned}
F={}&ZuX^2+u^2+u+Z^2\sum_{i=1}^s x_i
       +t^2Z^2ux_2x_3\\
&+\sum_{j=1}^{s-1}
 x_jx_{j+1}(x_j+x_{j+1}+t^2)P_j^2
+x_s(u^2+t^2u)P_s^2=0.
\end{aligned}
$$
For choice $b$ of Definition~\ref{def:endpoint-choices}, replace $X$ by $a^2uX$ and $P_s$ by $aP_s$. After division by $a^3$, only the last term changes, to $x_s(1+t^2x_b)P_s^2$.

For $q$ in Definition~\ref{def:hermitian-equations},
Proposition~\ref{prop:nonvanishing}\textup{(2)} already proves
the reference variety property, including the permitted restrictions.
On the chart $Q_0=1$, first let $S_j=U_I=V_I=0$.
Make the same substitutions as above and replace $R$ by $aR$.
After division by $a^3$, the replacement of $X^2$ by
$\Nm(X+\theta(z_1+R))$ in the original equation adds
$
u^2X+uXR+u^2+R^2
$
to the corresponding scalar equation.
Replacing $R$ by $auR$ instead adds
$
u^2(X+XR+1+R^2).
$

For each $I\notin\cT$, let $d=|I|$ and
$J_I=I\setminus\{1\}$.
Consider the additional term
$
z_I\Nm(U_I+\theta V_I).
$
If $d$ is odd, replace $U_I$ by $a^{(3-d)/2}U_I$ and
$V_I$ by $a^{(3-d)/2}V_I$.
If $d$ is even, replace $U_I$ by $a^{(4-d)/2}V_I$ and
$V_I$ by $a^{(2-d)/2}U_I$, simultaneously.
After division by $a^3$, this term becomes
$
x_I(U_I^2+U_IV_I+aV_I^2).
$
If $1\notin I$, further replace $V_I$ by $uV_I$.
Since $a=Z/u$, the two expressions are
$$
\begin{cases}
x_{J_I}(U_I^2+uU_IV_I+ZuV_I^2),&1\notin I,\\
x_{J_I}(uU_I^2+uU_IV_I+ZV_I^2),&1\in I.
\end{cases}
$$
We may let either $U_I$ or $V_I$ be equal to zero.

All these fiber changes are invertible over the base function
field and preserve the homogenizing hyperplane.
Remark~\ref{rmk:quadric-changes}\textup{(1)} then shows these modifications still give reference varieties.
\end{construction}

\subsection{Additional Constructions}\label{subsec:small}
In this section, we give explicit constructions for the low dimensional cases beyond the previous uniform constructions. 

\begin{proposition}\label{prop:quartics}
In every characteristic different from two, there exists an algebraically closed field $k$ of that characteristic such that there are integral affine hypersurfaces $V$ over $k$ with
$\delta_V\ne0$ in $\CH_0(V_{k(V)},\ZZ/2)$ in each of the
following cases:
\begin{enumerate}
\item hypersurfaces in $\mathbb{A}^{n+1}_{k}$ of total degree at most four with $3\leq n\leq 8$.
\item hypersurfaces with profiles $(4,4,4;2,2,2)$ and $(5,5,5;2,2,2)$,
\item quintic hypersurfaces with profile $(s,s+1;3,2)$, where $s\geq3.$
\item quintic hypersurfaces with profile $(s,s+2;3,2)$, where $s\geq4.$
\end{enumerate}
\end{proposition}

\begin{proof}
Let $k=\overline{\Omega(t)}$ be the field in
Construction~\ref{constr:gram-cycle}, with $\Omega$ uncountable.
\begin{itemize}
\item[step 1:] In this step, we prove Cases~\textup{(3)} and~\textup{(4)},
and Case~\textup{(1)} in dimensions three through seven.

For Case~\textup{(3)}, use the chart $X=1$ of $\cG_s$ in
Construction~\ref{constr:gram-cycle}, denoted by $G_s$.
We then have $\deg G_s=5$.
Divide its coordinates into
$(z_1,\ldots,z_s)$ and $(Y_0,\ldots,Y_s)$, and let
$V=\{G_s=0\}\subset\bA_k^s\times\bA_k^{s+1}$.
The equation has degree at most three in the first group
and two in the second, so $V$ has profile $(s,s+1;3,2)$.
Taking $\cU=\varnothing$ in Lemma~\ref{lem:gram-local}
proves that $V$ is geometrically integral of dimension $2s$
and
$
0\ne\delta_V\in\CH_0(V_{k(V)},\ZZ/2).
$
This proves Case~\textup{(3)}.

For Case~\textup{(4)}, since $s\geq4$, the subset $\{1,3\}$
does not belong to $\cT_s$ in
Construction~\ref{constr:gram-cycle}.
We may then add an independent coordinate $Z$ and use
$V=\{G_s+4z_1z_3Z^2=0\}\subset\bA_k^s\times\bA_k^{s+2}$,
with coordinates divided into 
$(z_1,\ldots,z_s)$ and $(Y_0,\ldots,Y_s,Z)$.
The added term has degree two in each part and total degree
four. Thus the equation still has total degree five and
profile $(s,s+2;3,2)$.
Then Lemma~\ref{lem:gram-local}, with $\cU=\{\{1,3\}\}$, proves
that $V$ is geometrically integral of dimension $2s+1$ and
that $\delta_V\ne0$ with coefficients $\ZZ/2$.
This proves Case~\textup{(4)}.

For Case~\textup{(1)} in dimensions three and four, take $n'=m=2$
in \cite[Theorem~4.1]{LangeZhang2026}. That theorem gives a nonzero
generic-point class after deleting the locus $\partial_zf=0$ from the
affine hypersurface. Localization then implies that the class on
the full affine hypersurface is nonzero as well.
In dimension five, take the quartic $X'\subset\PP^6$ constructed in
\cite[Example~4.6]{LS} and its chart $x_0=1$. The cited obstruction holds
after deleting $\{x_0y_1z=0\}$, so localization again gives a nonzero
class on the larger affine chart $x_0=1$.

For dimensions six and seven, let $G=1+xu^2+yv^2+2xy$ and
$$
 W=\left\{w^2+(2G+x+y)w+xa^2+yb^2+xyc^2+
                    \frac{(x-y)^2}{4}=0\right\}.
$$
With independent coordinates $(x,y,u,v,w,a,b,c)$, this defines
$W_7\subset\bA_k^8$. Substituting $v=u$ defines $W_6\subset\bA_k^7$.
The latter comes from \cite[Proposition~3.2 and
Equation~(4.1)]{Ottem2026}, with $a(z)=b(z)=1+z^2$,
$A(z)=B(z)=z$, and $c_1=c_2=2$ in his notation, after rescaling
two coordinates. We verify the affine obstruction for both models
in every characteristic different from two.

Let $\rho=(x,y)_2$, $\Pi=\langle1,x,y,xy\rangle$, and
$g_1=G+x$. Let $g_2=G+y$ and $\alpha_i=\rho\cup(g_i)_2$.
The residue of $\alpha_1$ at $y=0$ is
$(x,1+x(1+u^2))_2$. Taking its residue at
$1+x(1+u^2)=0$ gives the nonsquare class
$(-1/(1+u^2))_2$ over $k(u,v)$, or over $k(u)$ for $W_6$.
It is nonsquare because $1+u^2$ has simple zeros.
Interchanging $x,u$ with $y,v$ proves $\alpha_2\ne0$.

Completing the square gives
$$
 \Pi\left(w+G+\frac{x+y}{2},a,b,c\right)=g_1g_2.
$$
Lemma~\ref{lem:common-symbol} gives geometric integrality of the
generic quadric. The affine equation is monic in $w$, so Gauss's
lemma gives integrality of $W$.

We check the boundary condition of Lemma~\ref{lem:pfister}.
Let $\mu$ be a nontrivial discrete valuation of the base field over
$k$, and work over its henselization. Suppose that both $\alpha_i$
are nonzero there. Then $\rho\ne0$.
Let $X=\mu(x)$, $Y=\mu(y)$, and $g=\mu(G)$.
Since $G=\Pi(1,u,v,\sqrt2)$, Lemma~\ref{lem:norm-values} gives
$$
 g=\min\{0,X+2\mu(u),Y+2\mu(v),X+Y\}.
$$
Also $\rho\cup(G)_2=0$. If $X>g$, then $g_1/G$ belongs to
$1+\fm_\mu$ and is a square, contradicting $\alpha_1\ne0$.
Thus $X\leq g$, and similarly $Y\leq g$.
It follows that $g\leq X+Y\leq2g$ and $g\leq0$, so
$g=X=Y=0$ and $\mu(u),\mu(v)\geq0$.

Lemma~\ref{lem:laurent-values} gives $\mu(y-x)=0$.
If $\mu(g_1)>0$, then $g_2/(y-x)$ belongs to $1+\fm_\mu$,
and the Steinberg relation $(x,y,y-x)_2=0$ gives $\alpha_2=0$.
The other case is symmetric. Hence both $g_i$ are units.
The center of $\mu$ belongs to the affine open
$B^\circ=\{xyg_1g_2\ne0\}$. This verifies the boundary condition.
The completed-square change is an isomorphism over this open and
preserves the hyperplane at infinity. Taking the closure of the graph
of the projection in projective compactifications gives the proper
morphism in Lemma~\ref{lem:pfister}. That lemma now gives
$\delta_{W_i}\ne0$ modulo two for $i=6,7$.

\item[step 2:] In this step, we prove Case~\textup{(2)}.
Let $a,p,c$ be independent base variables, let $b=p-a$,
and choose $\varepsilon\in k$ with $\varepsilon^2=2$.
Let $G_0=1+c(pu^2+2ab)$ and define
$$
\begin{aligned}
    F_0={}&w^2+[2(1+c)G_0+cp]w
       +aY_a^2+bY_b^2+abY_{ab}^2\\
&+c\left[
R^2+(2a-p)RQ+\frac{p^2Q^2}{4}
+aT^2+2aTS+pS^2
\right]
+\frac{c^2(2a-p)^2}{4},
\end{aligned}
$$
$$H=ae^2+2acue+c(1+c)v^2+abf^2
       +2\varepsilon abcf.$$
Let $F_1:=F_0+2wH$. The polynomial $F_0$ has profile $(4,4,4;2,2,2)$.
Its coordinate groups are $(c,Y_a,Y_b,Y_{ab})$, $(a,u,Q,S)$ and $(p,w,R,T)$.
The polynomial $F_1$ has profile $(5,5,5;2,2,2)$ with coordinate
groups $(c,Y_b,Y_{ab},e,f)$, $(a,u,Q,S,v)$ and $(p,w,R,T,Y_a)$.

We have to prove $\delta_V\ne0$ for $V=\{F_i=0\}$. To do so, we apply Lemma~\ref{lem:pfister} over the affine base
with coordinates $(a,b,c,u)$ for $F_0$, and
$(a,b,c,u,e,f,v)$ for $F_1$.
Let $\rho=(a,b,c)$. Let $h=\frac{c}{1+c}$, $x=\frac{e}{1+c}$,
and $y=\frac{f}{1+c}$. Consider the 3-fold Pfister form, see Definition~\ref{def:pfister},
$\Pi=\langle\!\langle a,b,c\rangle\!\rangle=\langle1,a,b,ab,c,ac,bc,abc\rangle$.
Define
$$
\begin{aligned}
    G:=&1+ax^2+ac(u+x)^2+bcu^2+cv^2+aby^2+abc(\varepsilon+y)^2\\
    =&\Pi(1,x,0,y,v,u+x,u,\varepsilon+y)
\end{aligned}
$$
Let $g_1=G+ah$, $g_2=G+bh$, and $\alpha_i=\rho\cup(g_i)_2$.
 Because $
\left(c,\frac{c}{1+c}\right)_2=0$, we have 
$
\rho\cup(h)_2=0.
$

For $F_0$, we simply take $e=f=v=0$ in the above construction.  Hence $x=y=0$, and 
$G=1+acu^2+bcu^2+2abc.$
Hence, with $D=1+(1+c)u^2$, we have 
$$\alpha_1=(a,b,c,1+\frac{acD}{1+c}+bc(u^2+2a)).$$
By \cite[Lemmas~3.1--3.2 and Remark~3.3]{Schreieder}, taking residues of $\alpha_1$ first at $b=0$, then at
$1+acD/(1+c)=0$, and then at $D=0$ gives
$$
0\neq\left(-\frac{1+u^2}{u^2}\right)_2\in H^1(k(u),\mu_2) .
$$
This proves $\alpha_1\ne0$. Interchanging $a$ and $b$ proves
$\alpha_2\ne0$. For $F_1$, the symbol entries are units at the generic
coordinate divisors $e=0$, then $f=0$, then $v=0$. Thus the classes are
unramified there and specialize in this order to the classes for
$F_0$. These final specializations are nonzero by the residue computation
above. Hence the original classes for $F_1$ are nonzero as well.
This proves $\alpha_1,\alpha_2\ne0$ for both constructions.

We verify the third hypothesis of Lemma~\ref{lem:pfister}.
Let $\mu$ be a nontrivial discrete valuation of the base field over
$k$, and work over its henselization. Suppose that both $\alpha_i$
are nonzero there. Then $\rho\ne0$.
Let $A=\mu(a)$, $B=\mu(b)$, and $C=\mu(c)$.
Let $g=\mu(G)$.
If $C<0$, then $h=(1+c^{-1})^{-1}$ is a square. Replacing the
zero coordinate in the vector representing $G$ by $\sqrt h$
shows that $g_2=G+bh$ is represented by $\Pi$, a contradiction.
Thus $C\geq0$.
Lemma~\ref{lem:laurent-values} gives $\mu(1+c)=0$ and $\mu(h)=C$.

The form $\Pi$ represents $G$, so $\rho\cup(G)_2=0$.
If $A+C>g$, then $g_1/G\in1+\fm_\mu$ is a square,
contradicting $\alpha_1\ne0$.
Similarly, $\alpha_2\ne0$ gives $B+C\leq g$.
Lemma~\ref{lem:norm-values} computes $g$ as the minimum of the
term values in the expression for $G$. In particular, $g\leq0$.
Since $y$ and $\varepsilon+y$ cannot both have positive valuation,
the terms $aby^2$ and $abc(\varepsilon+y)^2$ give
$g\leq A+B+C$. Thus
$$
 g+C\leq A+B+2C\leq2g,\qquad C\geq0,\qquad g\leq0.
$$
It follows that $g=C=0$. The inequalities
$A,B\leq0$ and $A+B\geq0$ then give $A=B=0$.
The same minimum formula gives $\mu(u),\mu(x),\mu(y),\mu(v)\geq0$,
and hence $\mu(e),\mu(f)\geq0$.

The difference $g_2-g_1=h(b-a)$ is a unit by
Lemma~\ref{lem:laurent-values}. If $\mu(g_1)>0$, then
$g_2/(h(b-a))\in1+\fm_\mu$ is a square. The relations
$\rho\cup(h)_2=0$ and $(a,b,b-a)_2=0$ give $\alpha_2=0$,
a contradiction. The second relation is the Steinberg relation
for $a/b$. The argument with $g_1,g_2$ interchanged is identical.
Thus both $g_i$ are units and the center belongs to
$B^\circ=\{abc(1+c)g_1g_2\ne0\}$ in the affine base.
Every valuation centered outside this open kills at
least one $\alpha_i$, as needed.

We then verify the coordinate change in Lemma~\ref{lem:pfister}.
Expansion gives $G=G_0+H/(1+c)$, with $H=0$ for $F_0$.
After homogenizing $F_i$ in the eight fiber coordinates
$w,Y_a,Y_b,Y_{ab},R,Q,T,S$
with the additional coordinate $\zeta$, direct expansion gives
$$
\begin{aligned}
\Pi\biggl(
&w+\frac{(1+c)(g_1+g_2)}2\zeta,
Y_a,Y_b,Y_{ab},\\
&R+(a-p/2)Q,T+S,S,Q
\biggr)
=(1+c)^2g_1g_2\zeta^2.
\end{aligned}
$$
Divide the eight arguments of $\Pi$ by $1+c$, leaving
$\zeta$ fixed.
Since $\Pi$ is homogeneous of degree two, the equation
becomes
$$
\Pi(\boldsymbol Z)-g_1g_2\zeta^2=0,
$$
which is the quadric in Lemma~\ref{lem:common-symbol}.

Before division by $1+c$, the change is a composition
of a permutation of coordinates and linear changes
with determinant one.
After division, its determinant is
$\pm(1+c)^{-8}$, a unit on $B^\circ$.
All its coefficients are regular there.
Thus it gives an isomorphism over $B^\circ$ and
preserves $\zeta=0$.
Restricting to $\zeta=1$ verifies the second hypothesis
of Lemma~\ref{lem:pfister}.
The first hypothesis follows from the definition of
$B^\circ$.

Lemma~\ref{lem:common-symbol} proves that the generic
projective quadric is smooth and geometrically integral.
Its affine complement is then integral.
Since $F_i$ is monic in $w$, it is primitive over the
polynomial ring of the base coordinates.
Gauss's lemma then proves that $V=\{F_i=0\}$ is integral.

Choose projective compactifications of $V$ and of its
affine base, and take the closure of the graph of the
projection.
This gives the proper compactifications and the morphism
between them in Lemma~\ref{lem:pfister}.
Applying that lemma, we obtain
$$
0\ne\delta_V\in\CH_0(V_{k(V)},\ZZ/2).
$$
Together with the coordinate groups above, this proves
Case~\textup{(2)}.

\item[step 3:] In this step, we construct the quartic eightfold $V_8$.
Let
$K=k(a,b,c)$, $P=abc$, and $L=a+1$. Let $\alpha=(a,b,c)$
and $\Pi=\langle1,a,b,c,ab,ac,bc,P\rangle$.
Let $\lambda=y_1+y_3+Ly_4+y_5+(c-a-b-1)y_6$ and define
$$
\sigma_6=ay_1^2+by_2^2+cy_3^2+aby_4^2+acy_5^2
             +4bcy_6^2+\lambda^2.
$$
Let $q_6=\sigma_6-PL\zeta^2$ and
$V_8=\{\sigma_6-PL=0\}\subset\bA_k^9$.
Its defining equation has total degree four.
We verify the three conditions of
Proposition~\ref{prop:charts}.

Coordinate residues give $\alpha\ne0$, so $\Pi$ is
anisotropic.
The vector
$(\lambda,y_1,y_2,y_3,y_4,y_5,2y_6,0)$
defines an injection $\sigma_6\subset\Pi$.
Let
$D=aL^2+2L(a+2)(b+c)+a(b+c)^2+4bc$, $\beta=(L,D)$.
Computing the determinant gives $\det\sigma_6=P^2D$.
The vector $e_P$ is orthogonal to $\sigma_6$ and has
norm $P$, so
$$
\Pi\simeq\sigma_6\perp P\langle1,D\rangle.
$$
In particular, $\sigma_6$ is nondegenerate, and
$Q_6=\{q_6=0\}$ is a smooth geometrically integral quadric.

We prove $\beta=(L,D)_2\ne0$ directly.
Let $F=k(a,c)$, $T=c^2+(a+2)c+L^2$, and
$B_0=(a+2)(c+L)$.
The homogeneous quadratic defining $T$ has determinant $-1/4$,
so $T$ is irreducible. It is coprime to $L=a+1$.
The valuations at $L$ and $T$ show that their square classes in
$F^\times/F^{\times2}$ are independent.
The identity
$$
 aD=(ab+B_0)^2-4LT
$$
then shows that $D$ is irreducible in $F[b]$ and that its residue
field is $F(\sqrt{LT})$. The class of $L$ remains nonsquare in
this quadratic extension, since $L$ and $T$ have independent
square classes. We then have
$$
 \partial_D\beta=(L)_2\ne0.
$$

We then prove condition~\textup{(1)} for $Q_6$.
Suppose $\alpha|_{K(Q_6)}=0$.
If $q_6$ were isotropic, $Q_6$ would be rational over $K$,
contradicting $\alpha\ne0$.
The Pfister subform theorem
\cite[Theorem~2.3]{SchSmall} and roundness
\cite[Corollary~9.9]{EKM} then give $q_6\subset\Pi$.
Here both forms represent $L$, with $y_1=1$ and the
other coordinates zero for $q_6$.

Witt cancellation \cite[Theorem~8.4]{EKM}, applied to $\sigma_6$ in
$$
q_6=\sigma_6\perp\langle-PL\rangle,\qquad
\Pi\simeq\sigma_6\perp P\langle1,D\rangle
$$
shows that $\langle1,D\rangle$ represents $-L$.
Thus $\beta=(L,D)=0$, contrary to $\beta\ne0$.
This proves condition~\textup{(1)}.

For condition~\textup{(2)}, let $\mu$ be a nontrivial
discrete valuation of $K/k$ with
$\alpha|_{K_\mu^h}\ne0$, and work over $K_\mu^h$.
Let
$(A,B,C)=(\mu(a),\mu(b),\mu(c))$.
By Lemma~\ref{lem:laurent-values}, these three integers
cannot all be zero, since otherwise $\mu$ would be
trivial on $K$.
We prove that $\beta=0$ over $K_\mu^h$.

If $A>0$, Hensel's lemma makes $L=1+a$ a square.
Suppose $A<0$, and let $r=\min(A,B,C)$.
In
$$
D=a(a+b+c)^2+2a^2+6a(b+c)+a+4bc+4(b+c),
$$
the first term has valuation $A+2r$, strictly smaller
than that of every other term.
Hence
$\frac{D}{a(a+b+c)^2}$, $\frac{L}{a}$
belong to $1+\fm_\mu$ and are squares.
It follows that $\beta=(a,a)=0$.

Suppose $A=0$ and $r=\min(B,C)<0$.
Let $D_2=a(b+c)^2+4bc$.
Lemma~\ref{lem:laurent-values} gives
$$
\mu(D_2)=2r,\qquad
\mu(D-D_2)\geq r>2r.
$$
Thus $D/D_2$ is a square.
The identity
$$
aD_2=(ab+(a+2)c)^2-4Lc^2
$$
and the Steinberg relation $(L,a)=0$ give
$(L,D_2)=0$, and then $\beta=0$.

Finally, suppose $A=0$ and $B,C\geq0$.
After interchanging $b,c$, we may assume $B>0$.
Let
$D_0=ac^2+2L(a+2)c+aL^2$.
Lemma~\ref{lem:laurent-values} gives $\mu(D_0)=0$,
and $D-D_0\in b\cO_\mu$.
Hence $D/D_0$ is a square.
The identity
$$
aD_0=(ac+L(a+2))^2-4L^3
$$
gives $(L,D_0)=0$, so again $\beta=0$.

It follows that $P\langle1,D\rangle$ represents $-PL$.
The orthogonal decomposition above gives
$$
q_6\subset\Pi
$$
over $K_\mu^h$, and hence over the completion.
The last paragraph of
Remark~\ref{rmk:quadric-changes} proves
condition~\textup{(2)}.
Condition~\textup{(3)} is the injection
$q_6|_{\zeta=0}=\sigma_6\subset\Pi$.

The coefficients of $y_1^2$ and $y_3^2$ in the affine
equation are $a+1$ and $c+1$, which are coprime.
Thus the equation is primitive over $k[a,b,c]$.
Integrality of $Q_6$ and Gauss's lemma prove that $V_8$
is integral and is an affine quadric model of
$Q_6\setminus\{\zeta=0\}$.
Proposition~\ref{prop:charts} gives
$$
0\ne\delta_{V_8}\in\CH_0((V_8)_{k(V_8)},\ZZ/2).
$$

This completes Case~\textup{(1)} and the proof.
\end{itemize}
\end{proof}

\begin{definition}\label{def:binary-quartics}
Assume $\cha k=2$. In projective coordinates $[\ell:a:b:d:U:V:W]$ on $\PP^6_k$, define
$$
 X=\{\ell aU^2+\ell bV^2+(\ell d+ab)(W^2+ab)
                 +(\ell(a+b+\ell)+d^2)^2=0\}\subset\PP^6.
$$
Write $R_5=X_{\ell=1}$, $R_4=(R_5)_{W=0}$ and
$R_3=(R_5)_{W=V=0}$. On this affine chart let
$q=W^2+ab$ and $h=a+b+d^2+1$.
Thus $R_i$ has $i+1$ affine coordinates before adding its one
equation. On $\ell=1$, $a,b,d,U,V,W$ denote the corresponding
affine coordinates, with these omissions.
\end{definition}

\begin{proposition}\label{prop:binary-models}
The hypersurfaces $R_3,R_4,R_5$ are geometrically integral quartics.
For $i\in\{3,4,5\}$, every extension $L/k$, and every closed
point $p\in(R_i)_L$, let
$a_p=\{a(p),q(p)\}\in\mathrm k_2^M(\kappa(p))$ if $a(p)q(p)\ne0$,
and let $a_p=0$ otherwise.
These classes satisfy Condition~\textup{(C)} and hence define
$$
 \CH_0((R_i)_L,\Ftwo)\longrightarrow\mathrm k_2^M(L),
 \qquad [p]\longmapsto\cores_{\kappa(p)/L}(a_p).
$$
The value on $\delta_{R_i}$ over $k(R_i)$ is nonzero.
In particular, $\Tortwo(R_i)=2$ for $i=3,4,5$.
\end{proposition}
\begin{proof}
We prove integrality first, and then check the curve classes at
all points of their normal projective compactifications.
\begin{itemize}
\item[step 1:] In this step, we prove integrality and nonvanishing of the
symbol over the function field. On $\ell=1$ write the defining equation as
$$
 F=aU^2+bV^2+(d+ab)q+h^2=0,
$$
where $q=W^2+ab$ and $h=a+b+d^2+1$.
Regarding $F=aU^2+G$ as a polynomial in $U$, the identity
$\partial_d(G/a)=q/a\ne0$ shows that $G/a$ is not a square in
the rational function field of the other variables. Also
$$
 G\bmod a=bV^2+dW^2+(b+d^2+1)^2\ne0.
$$
Thus $F$ is primitive and irreducible. The same proof works
after the fixed omissions, with $q=ab$. Since $k$ is algebraically closed, the varieties are
geometrically integral. All equations keep total degree four. Over $k(R_i)$, we have $F_d=q\ne0$, so the variables other than
$d$ give a separating transcendence basis.

In particular
$$
 d\log\{a,q\}=\frac{da\wedge db}{q}\ne0.
$$
Thus the Milnor symbol over the function field $\beta=\{a,q\}$ is nonzero.

We next give the symbol identities that will be used in the proof. Let $s=W^2+d$, $P=V^2+as$, and $Q=U^2+bs$. Expansion modulo $F$ gives
$$
 a(h+q)^2+(aU+WV)^2+qP=0.
$$
Using $\{x,r^2+xs^2\}=0$ and its symmetric versions gives
$$
 \beta=\{a,q\}=\{b,q\}=\{a,P\}=\{b,Q\}
              =\{s,P\}=\{s,Q\},
$$
whenever the slots are nonzero. We use a presentation suited to each valuation, and none is used when one of its slots is identically zero.

\item[step 2:] In this step, we prove the valuation formula used below.
Let $(E,\mu)$ be a discretely valued field of characteristic two.
Write $\partial_\mu$ for the Milnor residue and
$s_\mu^\pi(\gamma)=\partial_\mu(\{-\pi\}\cdot\gamma)$ for
specialization with respect to a uniformizer $\pi$, as in
\cite[Section~1, especially rule~R3f]{Rost1996}.
For $x=\pi^rx_0$ and $y=\pi^sy_0$, with $x_0,y_0$ units,
we have $\partial_\mu\{x,y\}=r\{\bar y_0\}+s\{\bar x_0\}$ and
$s_\mu^\pi\{x,y\}=\{\bar x_0,\bar y_0\}$.
When the residue vanishes, specialization is independent of $\pi$.
For $x,y\in E^\times$, consider
$q_{x,y}=Z_0^2+xZ_1^2+yZ_2^2+xyZ_3^2$.
Its value is the minimum of its four term values if either
$\partial_\mu\{x,y\}\ne0$, or both $\mu(x),\mu(y)$ are even
and their normalized residue symbol is nonzero. To see this, absorb even powers of a uniformizer into the
coordinates and write $x=\pi^rx_0$, $y=\pi^sy_0$ with units
$x_0,y_0$ and $r,s\in\{0,1\}$.

For $(r,s)=(1,0)$, nonzero residue means $\bar y_0$ is nonsquare, so the two parity groups of terms have anisotropic residual binary
forms proportional to $\langle1,\bar y_0\rangle$.

The case $(0,1)$ is symmetric.

For $(1,1)$, the nonzero residue means $\bar x_0\bar y_0$ is
nonsquare, so the residual forms
$\langle1,\bar x_0\bar y_0\rangle$ and
$\langle\bar x_0,\bar y_0\rangle$ are anisotropic.

For $(0,0)$, a nontrivial zero of the residual four dimensional
form would either make $\bar x_0$ square or express $\bar y_0$
as a quotient of two nonzero binary values $A^2+\bar x_0B^2$.
Both alternatives contradict its nonzero Milnor symbol.

Reducing the terms of least value then proves the claimed
exact minimum, with $\mu(0)=\infty$.

We then use two consequences of these residue and specialization formulas. A symbol with one unit slot
whose residue is square has both zero residue and zero
specialization. 

If the residue is zero and at least one
slot has odd valuation, its specialization is zero: in parity
$(1,0)$ or $(0,1)$ the opposite normalized unit is square, and
in parity $(1,1)$ their square classes coincide. In the last case
use $\{u,u\}=0$. These claims concern residue classes and
do not assert that a unit with square residue is itself square.

\item[step 3:] In this step, we define the classes on curves and check
their specializations at points on the affine variety.
Fix $L/k$ and an integral closed curve $C\subset(R_i)_L$.
If $C\subset\{aq=0\}$, take its class to be zero. Otherwise take the
restriction of $\beta$ to $L(C)$ and use the
normal projective compactification of the normalization of $C$. Off $aq=0$, unit slots give the fixed interior specialization. At a center on $aq=0$ the following residue calculations give
a square unit slot in an available presentation, and all equalities
in this paragraph are in the residue field. If $a=0$ and $q\ne0$, then $q=W^2$. If $q=0$ and $a\ne0$, then $ab=W^2$ and
$(aU+WV)^2=ah^2$. Thus $a$ is square when $h\ne0$.

If $h=0$ but $a+b\ne0$, then
$$
 a=\left(\frac{a+W}{d+1}\right)^2.
$$

If $q=h=a+b=0$ and $a\ne0$, then $b=a$, $d=1$, $W=a$, and $U=V$.
Use the unit $a$ if $a=1$, and otherwise use the unit
$s=(a+1)^2$ in $\{s,P\}$ or $\{s,Q\}$.

If $a=q=0$ and $b\ne0$, the relation
$bV^2=(b+d^2+1)^2$ makes $b$ square, and when $V=0$ this follows
from $b=(d+1)^2$ directly. If $a=b=q=0$, then $d=1$ and $s=1$. Step~2 gives zero residue and the fixed
zero specialization in each case.

For clarity, the identically zero entries in these symbol presentations cause no
exception. If $b=0$ in $L(C)$, then $q=W^2$ and $\beta=0$. If $s\ne0$ and $P=Q=0$, then $ab=(UV/s)^2$, again making $q$
square. 

If $s=P=Q=0$, then $U=V=0$, $d=W^2$, and the equation
gives $h=q$. For $a\ne1$ we have
$$
 q=\frac{(W+a)^2+a(W^2+1)^2}{1+a}.
$$
Its numerator and denominator are nonzero values of
$A^2+aB^2$, so they both annihilate $\{a,-\}$ and $\beta=0$. For $a=1$ the symbol is already zero.

\item[step 4:] In this step, we prove vanishing at points outside the
affine curve whose projective centers lie in the chart $a\ne0$.
At such a point, let
$(T,B,D,u,v,w)=(\ell,b,d,U,V,W)/a$.
Then $\mu(T)>0$ and all other chart coordinates are integral. Let
$\cQ=w^2+B$, $C_0=v^2+D$, and $L_0=w^2+TC_0$. Let
$\cP=u^2+Bv^2+D\cQ$, $u'=u+wv$, and
$H'=T(1+B+T)+D^2+\cQ$.
The homogeneous equation gives
$$
 \cQ L_0=H'^2+Tu'^2,\qquad
 T\cP=H'^2+\cQ w^2,
$$
and hence, whenever defined,
$$
 \beta=\{T,\cQ\}=\{T,L_0\}=\{C_0,L_0\}=\{\cP,\cQ\}.
$$
The residue equation is
$\bar B(\bar w^2+\bar B)+\bar D^4=0$. If $\bar w\ne0$, use the square unit residue
$\bar L_0=\bar w^2$. If $\bar w=0$ and $\bar B\ne0$, use
$\bar{\cQ}=\bar B=\bar D^2$. If these both vanish and $\bar v\ne0$, use
$\bar C_0=\bar v^2$, and if they and $\bar v$ vanish but $\bar u\ne0$,
use $\bar{\cP}=\bar u^2$. Each gives zero residue and specialization by Step~2. If $C_0=0$ but $L_0\ne0$, then $L_0=w^2$ and
$\beta=0$. If $L_0=0$ but $C_0\ne0$, use
$T=w^2/C_0$ and $\{C_0,\cQ\}$ in the case involving the unit $C_0$.

If $C_0=L_0=0$, then $w=0$, $D=v^2$, and $H'^2=Tu^2$. For $u\ne0$ this makes $T$ square. Otherwise $H'=0$ and
$$
 \cQ=T+\frac{v^4}{1+T}
       =\frac{(T+v^2)^2+T}{1+T}.
$$
Both numerator and denominator are binary values for $T$, so
$\{T,\cQ\}=0$ here as well.

We then treat the case where $w,B,v,u$ all have positive
valuation or are zero. The residue equation also gives $\bar D=0$. Assume the curve symbol is not identically zero, and let
$r=\mu(T)>0$, $s=\mu(\cQ)>0$. We have two more identities:
$$
 \begin{aligned}
 T\cQ D&=H'^2+Tu'^2+\cQ w^2+T\cQ v^2,\\
 H'&=(T+D)^2+T(1+w)^2+\cQ+T\cQ.
 \end{aligned}
$$
Suppose its residue were nonzero, or that $r,s$ were both even
with nonzero normalized specialization. The exact minimum for
$q_{T,\cQ}$ works for both identities. If $D\ne0$, the first gives
$$
 2\mu(H')\geq r+s+\mu(D)>r+s,
$$
whereas the second gives $\mu(H')\leq\min(r,s)$, since $1+w$
is a unit. This is impossible.

If $D=0$, anisotropy in the first identity gives
$H'=u'=w=v=0$. The second then gives
$$
 0=(1+T)(T+\cQ).
$$
Since $1+T$ is a unit, $\cQ=T$, and
$\beta=\{T,T\}=0$, again a contradiction. Thus the residue vanishes. If either parity is odd its
specialization vanishes by Step~2, and if both are
even, nonzero specialization has just been excluded.

\item[step 5:] In this step, we prove vanishing at points outside the
affine curve whose projective centers satisfy $a=0$.
The chart $b\ne0$ follows by interchanging $(a,U)$ with $(b,V)$.

At a center with $\ell=a=b=0$, the homogeneous equation
forces $d=0$. At least one of $U,V,W$ is a unit. For unit $W$, the second slot $W^2+ab$ has square unit residue.

For unit $U$, let $I=U^2+bd$, also a unit with square residue. If $b\ne0$ and $bW^2+\ell I\ne0$, the binary norm identities give
$$
 \beta=\left\{I,\frac{bW^2+\ell I}{b}\right\}.
$$
If the second slot is zero, use $b/\ell=I/W^2$ and
$\{I,W^2+ab\}$. If $b=0$, the original second slot is square.

The case of unit $V$ follows by interchanging the coordinates. Thus every missing point has
zero residue and zero specialization.

\item[step 6:] In this step, we apply reciprocity on the proper curve.
The curve classes are unramified everywhere, have precisely
the fixed interior values, and have zero boundary values.
Reciprocity on proper curves for Milnor theory modulo two,
including inseparable field extensions,
\cite[Proposition~2.2~\textup{(RC)}, Remark~2.4 and rules R2d, R3f]{Rost1996},
gives Condition~\textup{(C)} and hence the homomorphism explained
after that condition. Its generic value is the nonzero
symbol computed at the start. Hence $\Tortwo(R_i)=2$ for
$i=3,4,5$.
\end{itemize}
\end{proof}

\begin{definition}\label{def:ternary-quartics}
Let $\Omega$ be algebraically closed of characteristic two,
let $\lambda$ be transcendental over $\Omega$, and let
$k=\overline{\Omega(\lambda)}$.
In independent affine coordinates $a,b,W,Y,Z$ define
$$
 \Psi_4=W^3+aY^3+bZ^3+abYZ+a^2b^2
         +ab(Y+Z+\lambda W+a+b+1).
$$
Let $C_4=\{\Psi_4=0\}$ and $C_3=\{\Psi_4|_{Y=0}=0\}$.
The letter $Z$ here is one affine coordinate, not a subset indexed
family. The associated class is the cyclic symbol $(a,b)_3$.
\end{definition}

\begin{proposition}\label{prop:ternary-models}
The hypersurfaces $C_3,C_4$ and their projective closures are
geometrically integral. For $i=3,4$ and every field extension $L/k$,
assigning $(a(p),b(p))_3$ at a closed point with $a(p)b(p)\ne0$,
and zero otherwise, gives homomorphisms
$$
 \CH_0((C_i)_L,\ZZ/3)\longrightarrow\Br(L)[3]
$$
by corestriction. For $L=k(C_i)$, the value on $\delta_{C_i}$ is nonzero.
Thus $\Torthree(C_i)=3$.
For $i=3,4$, their profiles include
$$
 \begin{gathered}
 C_i:(2,i-1;3,3),\quad(i-1,1,1;3,2,2),\\
 C_3:(1,3;2,3),\qquad C_4:(1,4;2,4).
 \end{gathered}
$$
Together with the modulo two models, they give divisibility
by six for the torsion orders of very general bidegree $(4,3)$
hypersurfaces in $\PP^2\times\PP^2$, $\PP^2\times\PP^3$ and
bidegree $(2,4)$ hypersurfaces in
$\PP^1\times\PP^3$, $\PP^1\times\PP^4$ in characteristic two.
\end{proposition}
\begin{proof}
We first express the equation as a reduced norm.
Choose a primitive cube root $\zeta\in k$ and use the cyclic
algebra of Definition~\ref{def:norm-forms}, with
$i^3=a$, $j^3=b$ and $ij=\zeta ji$.
Expanding its reduced norm gives
$$
 \operatorname{Nrd}(W+Yi+\zeta^2Zj+\zeta i^2j+ij^2+ij+(ij)^2)
 =\Psi_4-\lambda abW.
$$
The coordinate $W$ has weight $(0,0)$, while the monomial on
the right has exponent $(1,1)$.
Thus the point in Theorem~\ref{thm:norm-criterion} is $(3/2,3/2)$.
It lies strictly inside the convex hull of
$(2,1),(1,2),(1,1),(2,2)$, the four fixed weights.

All other free variables are norm coordinates, so the integrality
condition for those variables follows directly. The same argument works after letting $Y=0$. Theorem~\ref{thm:norm-criterion} then proves integrality,
nonvanishing of the symbol over the function field and the
evaluation statement on the affine variety. The projective closures are geometrically integral because they are
closures of geometrically integral affine hypersurfaces.

Finally, we verify the profiles by specifying the coordinate groups. For $C_3$, use $(a,Z)$ and $(b,W)$, or $(b)$ and $(a,W,Z)$,
for the profiles with two groups. For its profile with three groups, use $(W,Z)$, $(a)$ and $(b)$.

For $C_4$, the two groups are $(a,Z)$ and $(b,W,Y)$,
or $(a)$ and $(b,W,Y,Z)$.
For its profile with three groups, use $(W,Y,Z)$, $(a)$ and $(b)$.
The last divisibility statements follow from
Corollary~\ref{cor:embedding}, using these modulo three profiles and the
modulo two profiles in Lemma~\ref{lem:catalogue}\textup{(2)}.
\end{proof}
\section{Quartic hypersurfaces in characteristic two}
\label{sec:char2-quartics}

In this section we prove that very general quartic hypersurfaces of every
positive dimension admit no decomposition of the diagonal in
characteristic two. We use the affine criterion of
Proposition~\ref{prop:charts}. The quadratic forms occurring below can be
singular. Their associated Pfister forms are nonsingular, and the affine
criterion does not need a resolution of the reference quadric.

\subsection{A norm estimate}

We begin with a valuation estimate for norms. It will give the local
vanishing needed for the affine quadric construction.

Throughout this section, $k$ is an algebraically closed field of
characteristic two. Let $K=k(a,z_1,\ldots,z_s)$ and
$E=K(\theta)$, where $\theta^2+\theta=a$. Let $z_I=\prod_{i\in I}z_i$.
We use the Hermitian form with diagonal coefficients $z_I$, whose
underlying quadratic form is
$h=\mathop{\perp}_{I\subseteq\{1,\ldots,s\}}z_I[1,a]$.
Its associated symbol is $\alpha=[a;z_1,\ldots,z_s]$. We write $e_I$
for its Hermitian coordinate vectors. Thus
$h(ce_I)=z_I N_{E/K}(c)$, and a Hermitian line has dimension two as a
quadratic space over $K$.

\begin{lemma}\label{lem:char2-norm-estimate}
Let $v$ be a nontrivial discrete valuation of $K/k$, and suppose that
$h$ is anisotropic over $F=\widehat K_v$. Let
$C=\min_{c\in F}\max\{0,-v(a+c^2+c)\}$.
Then $C\le -v(a)$, and every $1+\eta$ with $v(\eta)>C$ is a norm
from $E_F/F$. For nonzero mutually Hermitian orthogonal vectors $v_i$,
$$
 \min_i v(h(v_i))\le v\Bigl(\sum_i h(v_i)\Bigr)
 \le\min_i v(h(v_i))+C.
$$
\end{lemma}
\begin{proof}
The possible values form a nonempty set of nonnegative integers,
so they have a minimum. No representative $a+c^2+c$ has positive valuation, since
Hensel's lemma would then split $E_F/F$, contradicting anisotropy.
Choose $c$ attaining the minimum, and let $a_0=a+c^2+c$ and
$\theta_0=\theta+c$, so that $\theta_0^2+\theta_0=a_0$. If $v(\eta)>C$, Hensel's
lemma applied to $a_0\eta T^2+T+1$ gives $q\in\mathcal O_F$ satisfying
$q+a_0\eta q^2=1$. Then
$$
 N(1+\eta q\theta_0)=1+\eta.
$$
The first inequality in the conclusion is the valuation inequality.
For the second, let $S=\sum_i h(v_i)$ and choose a summand of least
valuation, say $h(v_j)$. If $v(S)>v(h(v_j))+C$, the first conclusion
gives $\lambda\in E_F^\times$ with $N(\lambda)=1+S/h(v_j)$. Then
$$
 h\Bigl(\lambda v_j+\sum_{i\ne j}v_i\Bigr)
 =\bigl(1+S/h(v_j)\bigr)h(v_j)+\sum_{i\ne j}h(v_i)=0.
$$
The vector is nonzero because orthogonal vectors of nonzero Hermitian
norm are $E_F$-linearly independent. This contradicts anisotropy.
\end{proof}

We will use two choices of coefficients. For the first, let $b\ge2$.
Write the variables
$z_i$ in the form
$u_1,\ldots,u_u,x_1,y_1,\ldots,x_b,y_b$.
Let
$U=\prod_i u_i$, $X=\prod_{j=2}^b x_j$ and
$Y=\prod_{j=2}^b y_j$. Let $D=x_1Y+y_1X$.
Let $G=\prod_i(1+u_i)\prod_j(1+x_j+y_j)$ and $a'=a+G^2+G$.
Define $f_1=a^3a'x_1XD$ and $f_2=U\prod_jx_j^2y_j$.
For the second, let $b=1$. Let $U=\prod_i u_i$,
$G=\prod_i(1+u_i)(1+x+y)$ and $a'=a+G^2+G$.
Define $f_1=a^2a'Ux(1+y)$ and $f_2=(Ux)^2y$.
Empty products are one.

\begin{lemma}\label{lem:char2-coefficients}
For either choice above, every nontrivial discrete valuation $v$ of $K/k$
for which $h_{\widehat K_v}$ is anisotropic satisfies
$$
 v(f_2/f_1)>C,
$$
where $C$ is as in Lemma~\ref{lem:char2-norm-estimate}.
\end{lemma}
\begin{proof}
Write $A=-v(a)\ge C\ge0$ and $A'=-v(a')\ge C$.
\begin{itemize}
\item[step 1:] In this step, we prove the estimate for the first choice,
apart from the case where $a$ and every $z_i$ are units.

Let $c_i=v(u_i)$, $p_j=v(x_j)$ and $q_j=v(y_j)$.
Let $-r=\sum_i\min(0,c_i)+\sum_j\min(0,p_j,q_j)$.
Then $r\ge0$. The expansion of $G$ is a single orthogonal sum of
squarefree monomials.
Let $g=v(G)$. Lemma~\ref{lem:char2-norm-estimate} then gives
$$
 -r\le g\le-r+C.
$$
In particular, the error is $C$, independently of the number of factors.
Let $d_0$ be the minimum of $p_1+\sum_{j\ge2}q_j$ and $q_1+\sum_{j\ge2}p_j$.
Applying the same lemma to the two terms of $D$ gives $v(D)\le d_0+C$.
Since
$$
 \sum_i c_i+\sum_j(p_j+q_j)-d_0
 =\sum_i c_i+\max\left(p_1+\sum_{j\ge2}q_j, q_1+\sum_{j\ge2}p_j\right)
 \ge-r,
$$
we obtain
$$
 v(f_2/f_1)-C\ge3A+A'-r-2C.
$$
If $r\le2A$, the right side is at least $A-C\ge0$. Equality requires
$A=C$, $r=2A$, and $A'=C$. When $A>0$, the bound $g\le-A$ implies
$v(a+G^2+G)=2g$ and hence $A'\ge2A$, a contradiction.
If $r>2A$, then $g<-A$, so $A'=-2g\ge2r-2C$ and the right side
of the above inequality is at least $r-A>0$.
It remains to consider $A=r=0$. Then $C=0$, all coordinate valuations
are nonnegative, and $v(D)=d_0$. The exact expression for $v(f_2/f_1)$
is positive if any $u_i,x_j,y_j$ has positive valuation.

\item[step 2:] In this step, we prove the same estimate for the second
choice, with the same exception.

Write $g=v(G)$, $d=v(1+y)$, and
$H=v(Uxy)-d$.
The identity
$$
 G(1+y)=\prod_i(1+u_i)\bigl((1+y)^2+x+xy\bigr)
$$
is an orthogonal sum of the values of
$(1+y)e_{u_I},e_{xu_I},e_{xyu_I}$. One summand is $Uxy$.
Thus Lemma~\ref{lem:char2-norm-estimate} gives $g+d\le v(Uxy)+C$,
and hence
$$
 v(f_2/f_1)=H+2A+A'\ge g-C+2A+A'.
$$
For $A>0$, the inequality $g+2A+A'>2C$ follows in three cases.
If $g\ge0$, its left side is at least $2A+C$. If
$-A/2\le g<0$, it is at least $3A/2+C$. If $g<-A/2$, then
$A'=-2g$, so it equals $2A-g$.
For $A=0$, we have $C=0$, and the same conclusion holds when $g\ne0$.
If $g=0$, applying the orthogonal estimate to the full expansion of
$G$ shows that all $u_i,x,y$ have nonnegative valuation. The estimate
for $1+y$ then gives $d=0$. Then $H>0$ if any of these coordinate
valuations is positive.

\item[step 3:] In this step, we exclude the case where $a$ and all $z_i$
are units.

Since the valuation of $K/k$ is nontrivial, its residue field
has transcendence degree at most $s$ over $k$. This follows, for example,
by lifting a residue transcendence basis and adjoining an element of
positive valuation. The residue field is a $C_s$-field by
\cite[Theorem~2a]{Nagata1957}. This conclusion does not require finite
generation: the coefficients of any homogeneous form generate a field
of transcendence degree at most $s$ over $k$, and Nagata's theorem
gives the $C_s$ property for that field. The reduction of $h$ is nonsingular and has
$2^{s+1}>2^s$ variables, so it is isotropic.
A nonzero isotropic residue vector lifts by Hensel's lemma, contradicting
anisotropy. This treats the last case.
\end{itemize}
\end{proof}

\subsection{Affine reference quadrics}

We now construct affine quadric models with nonzero generic-point class.
The norm estimate verifies the local hypothesis of
Proposition~\ref{prop:charts}.

\begin{proposition}\label{prop:char2-quadric-reference}
Choose either set of coefficients in Lemma~\ref{lem:char2-coefficients}.
Let $\rho$ be the vector defined below. Let $\sigma$ be a quadratic
subform of the underlying quadratic form on $(E\rho)^\perp$, where
orthogonality is Hermitian. Assume that $\sigma$
contains the full norm block representing $f_2$. Let $t$ be transcendental
and let $\Omega=\overline{k(t)}$. By an affine quadric model we mean an
integral variety mapping to an affine model of the rational base, with
generic fiber identified with the indicated affine quadric, as in
Definition~\ref{def:affine-quadric}. Then every such model over
$\Omega$ of
$Q_t\setminus\{V=0\}$, where $Q_t=\{\sigma+(f_1+t^2f_2)V^2=0\}$,
has nonzero generic-point class modulo two.
\end{proposition}
\begin{proof}
We give orthogonal representations of the
two coefficients.
For the first construction, use the two distinct indices
$A_0=x_1y_1$, $B_0=XY$, and let
$\rho=\theta^3(G+\theta)(Xe_{A_0}+x_1e_{B_0})$.
For the second, let
$\rho=a(G+\theta)(e_{Ux}+e_{Uxy})$.
In both cases $h(\rho)=f_1$. The coefficient $f_2$ is represented by
$x_1Xe_{Uy_1Y}$ in the first construction and by $Uxe_y$ in the second.

We prove nonvanishing over the generic
quadric after algebraic extension of the constant field.
We check the first hypothesis of Proposition~\ref{prop:charts} over
$\Omega$. The symbol $[a;z_1,\ldots,z_s]$ is nonzero over $K$:
iterated residues at the variables $z_i$ give the nonzero
Artin--Schreier class $[a]$ over $k(a)$. Thus $h$ is anisotropic.
At $t=\infty$, after replacing $V$ by $t^{-1}V$, the quadric specializes
to $\sigma+f_2V^2=0$. The representing vector for $f_2$, together with
$V=1$, gives a smooth $K$-point, because that vector lies in a full
nonsingular norm block. Projection from this smooth point makes the special quadric
$K$-rational \cite[Section~70, before Lemma~70.3]{EKM}.
It is geometrically integral: its polar form has rank at least two,
and its radical contains $V$ with nonzero coefficient $f_2$.
The generic quadric is geometrically integral for the same reason.

For any finite extension of $k(t)$, choose a valuation above infinity.
Its residue field is $k$. The generic point of the rational special
quadric provides a valuation of the generic quadric function field
with residue field $K(Q_\infty)$. If $h$ became isotropic over the
generic field, properness of its projective quadric would make it
isotropic over this residue field. This contradicts anisotropy over a
purely transcendental extension of $K$ \cite[Lemma~7.15]{EKM}. Since an isotropic vector is
defined over a finite constant extension, $h$ remains anisotropic over
$\Omega(Q_t)$. This proves the first hypothesis.

We check the local vanishing and the
hyperplane condition.
For a nontrivial base valuation over $\Omega$, either $h$ already splits
over its completion, or we can use Lemma~\ref{lem:char2-coefficients}
with $\Omega$ in place of $k$. The constant $t$ has valuation zero. Hence
$1+t^2f_2/f_1$ is a norm, and multiplying $\rho$ by an element of this
norm realizes $f_1+t^2f_2$ on the same Hermitian line. We get an injective
quadratic-space map
$$
 \sigma\perp\langle f_1+t^2f_2\rangle\hookrightarrow h.
$$
A point on $Q_t$ over any extension of the completion then
splits $h$ and kills the symbol. Finally, the hyperplane $V=0$ is the
subform $\sigma\subset h$. Proposition~\ref{prop:charts} gives the conclusion.
\end{proof}

We next explain the specialization of the affine obstruction. Suppose a
polynomial family of quartics has, after a polynomial change of affine
coordinates and a positive weighted rescaling by $\tau$, an integral
affine special fiber as in Proposition~\ref{prop:char2-quadric-reference}.
Divide the rescaled equation by its least power of $\tau$. The new
affine family is flat over $\Omega[\tau]_{(\tau)}$, because the equation
is not divisible by $\tau$. Extend the constant field to
$L=\Omega(c_1,\ldots,c_M)$, where the $c_i$ are algebraically
independent, and add the quartic with these coefficients multiplied by
a strictly higher power of $\tau$ before rescaling. Work over
$L[\tau]_{(\tau)}$. Over $L(\tau)$, translation by the
fixed coefficients gives coefficients algebraically independent over $\Omega(\tau)$ for a
universal generic quartic. Its generic fiber is smooth, and its special
affine fiber is the scalar extension of the original reference variety.
By \cite[Lemma~3.6(e)]{LS}, its geometric obstruction persists under
this extension. Apply \cite[Lemma~3.8]{LS} with empty boundary after
every finite extension of $L(\tau)$, using a valuation above $\tau=0$.
This gives a nonzero generic-point class on the geometric generic affine
chart. A decomposition of the diagonal on the smooth projective geometric
generic quartic would make its zero-cycle degree map an isomorphism after
every field extension. Over $\overline{L(\tau)}$, choose a point on its
nonempty hyperplane boundary.
The generic-point class would equal this point, and localization to the
affine chart would make it zero, a contradiction.

\subsection{The Segre construction}

We next realize the affine quadrics as weighted limits of quartic
equations. We first use products of projective lines and planes.

\begin{proposition}\label{prop:char2-segre-family}
The following Segre constructions give quartic reference degenerations
in characteristic two.
\begin{enumerate}
\item Let $a_0\ge1$ and $b\ge2$. The embedding of
$(\PP^1)^{a_0}\times(\PP^2)^b$ and the first coefficients of
Lemma~\ref{lem:char2-coefficients} give every dimension in the interval
$$
 2^{a_0}3^b-2\le n\le2^{a_0+2b}+a_0+2b-3.
$$
\item The second coefficients and the embedding of
$(\PP^1)^{u+1}\times\PP^2$ give dimensions $5,6$ for $u=0$, and
for $u\ge1$ give every dimension in the interval
$$
 2^{u+1}3-2\le n\le2^{u+3}+u-2.
$$
In particular, $u=1$ gives dimensions $10,\ldots,15$.
\end{enumerate}
Very general quartic hypersurfaces in all these dimensions admit no
decomposition of the diagonal.
\end{proposition}
\begin{proof}
Let $a_0$ and $b$ be the numbers of line and plane factors.
Distinguish one line parameter $a$ as the Artin--Schreier parameter,
and write $u_i$ for the rest of the line parameters. Let $u=a_0-1$
and $m=a_0+2b$.
Let $\mathcal Q$ be the Segre monomial set, and let $\mathcal P$ be its
complement after deleting the top monomial $aU\prod_jx_j$ and its $m$
neighbors differing in one factor. Thus
$r_0=|\mathcal P|=2^{a_0}3^b-m-1$.
Let $\mathcal N$ be the Segre monomials of ordinary degree at least
two. The bijection $g:\mathcal P\to\mathcal N$ and quadratic
minors $D_{1,\mu},D_{2,\mu}$ of \cite[Lemma~4.3 and equations~(4.5)--(4.6)]{BohningBothmerPirutkaZhang2026} give normal coordinates $\eta_{g(\mu)}$ and positive weights
$w_\mu$. The terms of least weight in $D_{1,\mu}$ and $D_{2,\mu}$ are
$\eta_{g(\mu)}$ and $\mu\eta_{g(\mu)}$, respectively.
This follows by expanding the minors in the normal coordinates,
so the calculation holds in characteristic two. Every other linear
normal term has larger weight, as does every quadratic normal term.
The same proof of the bijection and the weight inequalities works for
$b=1$: it needs at least two Segre factors, which are present because of
the distinguished line factor and the plane. Let $F_\mu=D_{1,\mu}D_{2,\mu}$.

\begin{itemize}
\item[step 1:] In this step, we treat $b\ge2$.

Neither reserved index $A_0=x_1y_1$ nor $B_0=XY$ is a Segre monomial,
since $b\ge2$. The additional Hermitian line generated by
$\xi=Ye_{A_0}+y_1e_{B_0}$
is orthogonal to $\rho$ and has coefficient
$\gamma=y_1Y(x_1Y+y_1X)$.
Both $\gamma$ and $a\gamma$ have degree at most two on every Segre
factor. The orthogonal complement of $E\rho$ is the sum of every
coordinate norm block except $A_0,B_0$, together with $\gamma[1,a]$.
Its dimension over $K$ is $2^m-2$.

Under pairing $\mu\leftrightarrow a\mu$, the Segre set has
$2^{a_0-1}3^b$ pairs. The deleted neighborhood removes one full pair
and only the upper member from $m-1$ other pairs. Thus $\mathcal P$
contains complete pairs and lower singletons, with no upper singleton.
The full pair with coefficient $Uy_1Y$ belongs to $\mathcal P$ and
represents $f_2$. For every $r_0\le r\le2^m-2$, choose an
$r$-dimensional scalar coordinate subspace $\sigma_r$ in the above
orthogonal complement, containing all forced coordinates. Adjoin the
rest of the scalar basis vectors one at a time, and include the norm cross
term whenever both members of a pair occur.

We construct quartic lifts of the quadratic
terms.
Let $s_i$ be the extra
ambient linear variables, with positive weights $v_i$, and choose $L$
at least every sum of two weights. A forced diagonal is lifted by
$$
 \tau^{L-2w_\mu}F_\mu.
$$
For a complete forced pair, add the cross term
$$
 \tau^{L-w_\mu-w_{a\mu}}D_{2,\mu}D_{1,a\mu}.
$$
Its initial form is $\mu\eta_{g(\mu)}\eta_{g(a\mu)}$.
If a forced lower coordinate is paired with an extra upper coordinate
$s_i$, use the cross term
$$
 \tau^{L-w_\mu-v_i}D_{2,\mu}s_iW_1.
$$
Its upper diagonal term is
$$
 \tau^{L-2v_i}Q_{a\mu}s_i^2,
$$
where $Q_{a\mu}$ is a quadratic Segre lift of $a\mu$.
For a pair of extra coordinates, choose quadratic Segre lifts
$Q_c,Q_{ac}$ of its coefficients and use
$$
 \tau^{L-2v_i}Q_cs_i^2+\tau^{L-v_i-v_j}Q_cs_is_j
       +\tau^{L-2v_j}Q_{ac}s_j^2.
$$
This includes the added pair with coefficient $c=\gamma$. An unpaired
extra coordinate requires its diagonal term alone. All these polynomials
are homogeneous quartics. Each of these initial forms follows by multiplying
the unique least weighted terms of its factors. For coefficient lifts,
restriction to the Segre base is the weight-zero term.

We lift the coefficients and specialize the
quartic equation.
The degrees of $f_1$ on the line and plane factors are respectively at
most four and four. Those of $f_2$ are at most one and three. Hence they
admit quartic Segre lifts $F_1,F_2$. Add
$\tau^L(F_1+t^2F_2)$ and a general quartic in higher $\tau$-degree.
Normal rescaling gives the special affine equation
$\sigma_r+f_1+t^2f_2=0$. Its coefficient indexed by $1\in\mathcal P$
is one. It is then primitive over the base polynomial ring, and
its geometrically integral generic quadric makes the full affine model
integral. Proposition~\ref{prop:char2-quadric-reference} and the affine
specialization argument apply. Finally, the ambient projective dimension
is $m+r$, so $n=m+r-1$, which gives the above interval.

\item[step 2:] In this step, we treat $b=1$.

Here $a_0=u+1$ and $b=1$. Use the above normal minors and
norm cross terms, and omit the full norm blocks $Ux,Uxy$.
The first is the deleted top lower entry and the second is not a Segre
monomial. Thus every forced coefficient lies in the subform formed by the rest of the norm blocks,
whose rank is $2^{u+3}-4$. For $u\ge1$, the full $y[1,a]$ pair is forced:
both $y$ and $ay$ differ from the top entry in at least two factors.
The number of forced coordinates is $r_0=2^{u+1}3-u-4$, and
$n=u+2+r$, giving the above range. In particular, $u=1$ gives
$n=10,\ldots,15$.
When $u=0$, the forced entries are $1,y$. Adjoin the upper $ay$ entry
so that the full $y[1,a]$ block is present. Ranks three and four give
dimensions five and six. The degree bounds for $f_1$ are three in each
line parameter and four on the plane, so all coefficient lifts are
quartic. The coefficient indexed by $1$ is again one, so the same
primitivity and specialization argument proves the conclusion.
\end{itemize}
\end{proof}

\subsection{Three smaller ambient spaces}

The one-plane families of Proposition~\ref{prop:char2-segre-family}(2)
give dimensions five and six, and ten through fifteen.
We now construct the three dimensions between these ranges.

\setcounter{theorem}{5}
\begin{lemma}\label{lem:char2-small-lifts}
There are quartic reference degenerations in dimensions seven, eight,
and nine in characteristic two.
\end{lemma}
\begin{proof}
We use the second coefficient pair of Lemma~\ref{lem:char2-coefficients}.
Let $t$ be transcendental over $k$, and let $\Omega=\overline{k(t)}$.
In each case we give a polynomial of degree at most four whose
weighted specialization is an affine quadric of
Proposition~\ref{prop:char2-quadric-reference}.

\begin{itemize}
\item[step 1:] In this step, we construct the degeneration in dimension seven.
Use the one-plane construction with $u=0$.
Let $K=k(a,x,y)$ and $G=1+x+y$.
The vectors $\rho=a(G+\theta)(e_x+e_{xy})$ and $\xi=ye_x+e_{xy}$
are Hermitian orthogonal, and $h(\xi)=\gamma=xy(1+y)$.
Thus $\sigma=[1,a]\perp y[1,a]\perp\langle\gamma\rangle$
has rank five and lies in $(E\rho)^\perp$.
Its full $y[1,a]$ block represents $f_2=x^2y$.

Use coordinates $a,x,y,A,B,C,D,s$ on an affine chart of $\PP^8_\Omega$.
Let $r_A=A+ax$, $r_B=B+ay$, and $r_C=C+xy$.
Let $r_D=D+aC$.
The base is the graph $r_A=r_B=r_C=r_D=s=0$.
Let $R_A=A+ax$, $S_A=xB+yA$, and $T_A=xD+AC$.
Let $R_B=B+ay$, $R_C=C+xy$, and $S_C=aC+xB$.
Let $R_D=D+aC$ and $S_D=yD+BC$.
Two useful expansions are
$$
 T_A=xyr_A+xr_D+r_Ar_C,\qquad
 S_D=yr_D+xyr_B+r_Br_C.
$$
Let $F_1=a(A+D)(a+G^2+G)$ and $F_2=xC$.
On the graph, $a(A+D)=a^2x(1+y)$, so these polynomials restrict to
$f_1=a^2(a+G^2+G)x(1+y)$ and $f_2=x^2y$.
Their degrees are at most four and two, respectively.
Consider
$$
\begin{aligned}
 P_\tau={}&\tau^6(R_A+S_A)T_A+R_B^2
   +\tau^4R_CS_C+\tau^2R_BR_C\\
 &+\tau^2R_DS_D+\tau^{13}S_Ds+\tau^{24}Bs^2
   +\tau^{26}(F_1+t^2F_2).
\end{aligned}
$$
Every term has degree at most four in the ambient coordinates.
Substitute $r_A=\tau^{10}v_A$, $r_B=\tau^{13}v_B$, and $r_C=\tau^{11}v_C$.
Also substitute $r_D=\tau^{12}v_D$ and $s=\tau v_s$.
Every term of $P_\tau$ is divisible by $\tau^{26}$.
After division and setting $\tau=0$, the first product gives
$\gamma v_A^2$, the next three terms give
$v_B^2+v_Bv_C+av_C^2$, and the next three give
$y(v_D^2+v_Dv_s+av_s^2)$.
The special equation is
$$
 \gamma v_A^2+(v_B^2+v_Bv_C+av_C^2)
       +y(v_D^2+v_Dv_s+av_s^2)+f_1+t^2f_2=0.
$$
It is the affine quadric associated with $\sigma$.

\item[step 2:] In this step, we construct the degeneration in dimension eight.
Use one ordinary line parameter $u$.
Let $G=(1+u)(1+x+y)$,
$f_1=a^2(a+G^2+G)ux(1+y)$, and $f_2=u^2x^2y$.
Use coordinates $a,u,x,y,A,B,C,D,E_0$ on an affine chart of $\PP^9_\Omega$.
Let $\eta_A=A+au$, $\eta_B=B+ax$, and $\eta_C=C+ux$.
Let $\eta_D=D+uy$ and $\eta_E=E_0+auy$.
Their common zero locus is the graph of the base.
Let $R_A=A+au$, $S_A=aD+yA$, and $R_B=B+ax$.
Let $S_B=xE_0+BD$, $S_C=aC+xA$, and $T_C=xD+yC$.
Let $R_D=D+uy$, $R_E=E_0+aD$, and $S_E=uE_0+AD$.

Substitute $\eta_A=\tau^{13}v_A$, $\eta_B=\tau^{10}v_B$, and $\eta_C=\tau^{11}v_C$.
Also substitute $\eta_D=\tau^{14}v_D$ and $\eta_E=\tau^{12}v_E$.
Expansion gives the following congruences modulo $\tau$ after division
by the indicated powers:
$$
\begin{aligned}
 \tau^{-13}(R_A,S_A)&\equiv(v_A,yv_A),\\
 \tau^{-10}(R_B,S_B)&\equiv(v_B,uyv_B),\\
 \tau^{-11}(S_C,T_C)&\equiv(av_C,yv_C),\\
 \tau^{-14}R_D&=v_D,\qquad
 \tau^{-12}(R_E,S_E)\equiv(v_E,uv_E).
\end{aligned}
$$
For example, $S_B=uy\eta_B+x\eta_E+ax\eta_D+\eta_B\eta_D$
and $R_E=\eta_E+a\eta_D$, which explain the two different powers
used for these expressions.
It follows that the quartic polynomial
$$
 H_\tau=\tau^2R_AS_A+\tau^6S_CT_C+\tau^4R_AT_C
          +\tau^8R_BS_B+R_D^2+\tau^4R_ES_E
$$
is divisible by $\tau^{28}$ after substitution, and its quotient at
$\tau=0$ is
$$
 \sigma=y[1,a](v_A,v_C)+uyv_B^2+v_D^2+uv_E^2.
$$
The four summands lie in the norm blocks indexed by $y,uy,1,u$.
They avoid the blocks indexed by $ux,uxy$ used for $\rho$, and the full
$y[1,a]$ block represents $f_2$.

Let $\widetilde G=1+u+x+y+C+D$.
On the graph, $\widetilde G=G$ and $B(A+E_0)=a^2ux(1+y)$.
Hence the polynomials $F_1=B(A+E_0)(a+\widetilde G^2+\widetilde G)$
and $F_2=xCD$ restrict to $f_1,f_2$ and have degrees at most four.
Use the quartic polynomial $H_\tau+\tau^{28}(F_1+t^2F_2)$.
The same substitution and division give the special equation
$\sigma+f_1+t^2f_2=0$.

\item[step 3:] In this step, we obtain dimension nine from dimension eight.
Adjoin an ambient coordinate $s$ and add $\tau^{26}xs^2$ to the
polynomial above. Substitute $s=\tau v_s$.
After division by $\tau^{28}$, the special quadratic form is
$\sigma\perp\langle x\rangle$.
The additional summand lies in the norm block indexed by $x$, so it is
orthogonal to the blocks used for $\rho$.
\end{itemize}

In each case, the change to the normal coordinates is a polynomial
automorphism before the weighted rescaling.
The coefficient of $v_B^2$ in dimension seven, and of $v_D^2$ in
dimensions eight and nine, is one.
The special equations are primitive over the base polynomial ring.
Their generic quadrics are geometrically integral by
Proposition~\ref{prop:char2-quadric-reference}, so the affine special
fibers are integral.
Homogenize each equation to degree four and add a general quartic
multiplied by $\tau^M$, with $M>26$ in dimension seven and $M>28$ in
dimensions eight and nine.
The affine specialization argument following
Proposition~\ref{prop:char2-quadric-reference} proves the three conclusions.
\end{proof}

\begin{theorem}\label{thm:char2-all-quartics}
Let $k$ have characteristic two and let $X\subset\PP^{n+1}_k$ be a very general quartic hypersurface
of positive dimension. Then $X_{\bar k}$ admits no integral
decomposition of the diagonal.
In particular, $X$ is not geometrically stably rational.
\end{theorem}
\begin{proof}
For $n=1,2$, adjunction gives a nonzero top differential form, and the
differential-form obstruction follows from Lemma~\ref{lem:forms}. For $n=3,4,5$, use the
binary quartic reference varieties of Proposition~\ref{prop:binary-models}
and Corollary~\ref{cor:embedding}. Proposition~\ref{prop:char2-segre-family}(2)
covers dimension six and dimensions ten through fifteen, while
Lemma~\ref{lem:char2-small-lifts} covers seven through nine.

For the rest of the dimensions, Proposition~\ref{prop:char2-segre-family}
with $(a_0,b)=(1,2)$, $(2,2)$, $(1,3)$, $(2,3)$ and $(1,4)$ gives respectively
$[16,34]$, $[34,67]$, $[52,132]$, $[106,261]$, $[160,518]$.
For every $b\ge4$, the intervals with $(a_0,b)=(1,b)$ and $(1,b+1)$
overlap, since
$$
 6\cdot3^b\le2\cdot4^b+2b-1.
$$
The inequality holds at $b=4$ and remains true after increasing $b$.
These intervals then cover every $n\ge16$. The affine specialization argument proves the claim for the
universal geometric generic quartic. Definition~\ref{def:general} and
\cite[Lemmas~2.5 and~3.6(e)]{LS} give the formulation above.
\end{proof}

The ambient coordinate-degree bounds needed for products of projective
spaces are proved in Lemma~\ref{lem:quartic-singleton-profiles}.
\section{Coordinate bounds and consequences}\label{sec:bounds}
We use the affine varieties of Sections~\ref{sec:varieties}
and~\ref{sec:char2-quartics} to prove degree bounds for complete
intersections in products of projective spaces. The main task is to
divide the coordinates of a reference variety into groups so that its
equations fit the multidegrees of the complete intersection.
We first treat coordinate placement and quintic rows, then quadratic
rows, and finally cubic and quartic rows.

\begin{setup}\label{setup:coordinates}
Unless otherwise indicated, let $k$ be a field and let
$X\subset\prod_{j=1}^r\PP^{N_j}_k$ be a positive-dimensional
geometrically integral very general complete intersection of $c$ ample
hypersurfaces. Let $L=\sum_jN_j$, $M=\max_jN_j$, and
$N_{\min}=\min_jN_j$. Thus $\dim X=L-c$.
We write $x_+=\max\{x,0\}$.
The affine reference varieties are defined over an algebraically
closed field $\Omega$ of the same characteristic as $k$.
This field may be enlarged during a construction.
\end{setup}

\begin{definition}\label{def:placement}
A \emph{row} is the multidegree of one defining equation of $X$.
A quadratic, cubic or quartic row has every entry at least $2$, $3$
or $4$, respectively. A row \emph{suffices} if, under the other
hypotheses of the statement, its presence implies that $X_{\bar k}$
admits no decomposition of the diagonal.
The row used for the main equation of a reference variety is the
\emph{distinguished row}.

Each factor $\PP^{N_j}$ provides $N_j$ affine coordinate positions.
We call $N_j$ the \emph{capacity} of that group.
A \emph{placement} assigns the coordinates of an affine reference
variety to distinct positions in these groups. If $a_j$ coordinates
are assigned to the $j$th group, then $a_j\leq N_j$.
For an affine hypersurface, a \emph{profile}
$(a_1,\ldots,a_r;e_1,\ldots,e_r)$ means that its coordinate groups
have sizes $a_j$ and its defining equation has degree at most $e_j$
in the $j$th group. The notation $d^{[r]}$ means $r$ repeated entries
equal to $d$. A tuple of coordinates in a proof denotes one group,
and the groups are listed in the order of the profile.

An affine hypersurface of dimension $L-c$ uses $L-c+1$ coordinates.
Its placement leaves $c-1$ unused positions for the other rows.
Joining groups combines their coordinates. Splitting a group
partitions its coordinates and cannot increase the degree in any
of the smaller groups.
\end{definition}

More generally, let $V$ be a nonempty integral affine complete intersection
of codimension $q\leq c$, with a nonzero generic-point class for the
coefficients of the construction.
Place its coordinates in distinct ambient positions, using $a_j$
positions in the $j$th group.
To use Corollary~\ref{cor:embedding}, we check three conditions.
First, $\sum_ja_j-q=L-c$, so the dimensions coincide.
Second, match its $q$ defining equations to distinct rows so that
each equation has multidegree at most its row.
Third, the other $c-q$ rows can be assigned to distinct unused
positions, with positive degree in the factor containing that position.
The equations for these rows are the corresponding coordinate
functions. Their elimination recovers the reference variety.

For ample rows every unused position is permitted.
For other rows we use Hall's theorem: every set of rows must have
at least as many permitted positions as its number of rows.
This is the only matching condition used below.

The constructions may use different coefficient rings.
An integral decomposition of the diagonal would make the generic-point
class of a smooth affine chart vanish integrally, by
Lemma~\ref{lem:comparison}, and hence with any of these coefficients.
Corollary~\ref{cor:embedding} then contradicts the nonzero class on
the reference variety. This also allows the rational coefficients
in Lemma~\ref{lem:forms}.

\subsection{Coordinate placement and quintic rows}\label{subsec:placement}
We first explain how to assign the coordinates of the reference models
to groups with given sizes. We then prove the cubic and quintic bounds
used below. In the combinatorial proofs, each color represents one
coordinate group.

\begin{lemma}\label{lem:coloring-matching}
Let $n\geq2$ and let positive integers $a_j$ sum to $n$.
If $\max_j a_j\leq\floor{ n/2}$, a cyclic list of $n$ positions can be properly colored
with exactly $a_j$ occurrences of color $j$.
Moreover, when each of $n$ positions forbids at most two colors,
a fixed color quota can be matched to the positions if
Hall's inequalities hold for every set of one or two colors.
\end{lemma}
\begin{proof}
Let $a_{\max}$ be the largest multiplicity. Place the $a_{\max}$
occurrences of that color around the cycle, leaving a gap between each
pair. Fill these gaps cyclically with the other colors, completing one
color before moving to the next and continuing from the last gap used.
Each other color occurs at most once in a gap, since its multiplicity
is at most $a_{\max}$. Every gap is nonempty, since the other
multiplicities sum to at least $a_{\max}$. Entries within a gap have
distinct colors, and the color on either side of the gap is different.
Thus the cyclic coloring is proper.

For the matching, replace a color by as many copies as its quota and
use Hall's theorem \cite[Theorem~2.1.2]{Diestel}.
Every position permits at least one color from any set of three colors.
Such a set has all $n$ positions as neighbors, so its Hall inequality
follows from the total quota. Only sets of one or two colors need to be
checked.
\end{proof}

\begin{definition}\label{def:catalogue}
The symbols $V_5,V_8,W_6,W_7,R_3,R_4,R_5$ refer to the
affine hypersurfaces of Section~\ref{subsec:small}.
The $x_0=1$ equation of $V_5$ in
\cite[Example~4.6]{LS} is
$$
 z^2y_1^2+x_1x_2y_2^2+x_2y_3^2
 +x_1(1+x_1^2+x_2^2-2x_1-2x_2-2x_1x_2)+z-ty_1=0.
$$
The symbol $q_4$ denotes the polynomial
$$
 x_1z_1^2+x_2z_2^2+x_1x_2z_3^2
       +1+x_1^2+x_2^2-2x_1-2x_2-2x_1x_2.
$$
We also write $\{q_4=0\}$ for its affine hypersurface.
The polynomial $\Phi_5$, in independent coordinates
$u,v,z,y_1,y_2,y_3$, is
$$
 \Phi_5=z^2y_1^2+(v^2-uv)y_2^2+(v-u)y_3^2
             +v(u+1)^2-4v^2+z-ty_1.
$$
Here $t$ is the constant in the cited fivefold construction.
The equations $T_i$ and $F_i$ are defined in the proof below from the
models $R_i$.
\end{definition}

\begin{lemma}\label{lem:catalogue}
The affine models of Definition~\ref{def:catalogue} admit the
following profiles. Each model $V$ is integral and satisfies
$\delta_V\ne0$ in $\CH_0(V_{\Omega(V)},\ZZ/2)$.
We use the profile notation of Definition~\ref{def:placement}.

\smallskip
\noindent\textup{(1)} Suppose $\cha k\ne2$.
The conic constructed in the proof has profiles
$$
(3,1;4,2),\qquad (1,3;3,3),\qquad
(2,2;3,2),\qquad (2,1,1;2^{[3]}).
$$
The profiles of $q_4$ and the models $V_5,W_6,W_7$ include
$$
\begin{aligned}
q_4:\quad & (2,3;2,2),\quad (1,4;2,3),\\[3pt]
V_5,W_6,W_7:\quad & (r,m+1-r;3,2),\quad (1,m;2,4),
\end{aligned}
$$
where $m=5,6,7$, respectively, and $3\leq r<m$.
In addition, $V_5$ has profiles
$$
(2,2,2;2^{[3]}),\qquad (3,2,1;2^{[3]}),
$$
and $V_8$ has profiles
$$
(r,9-r;4,2),\qquad (2,7;4,3),
$$
where $3\leq r\leq8$ in the first profile.
For three coordinate groups, we have
$$
\begin{aligned}
V_8:\quad &(1+h,2+\eta,6-h-\eta;2,2,2),\\[3pt]
W_7:\quad &(1+h,2+\eta,5-h-\eta;2,2,2),
\end{aligned}
$$
where $0\leq h\leq2$ for $V_8$ and $0\leq h\leq1$ for $W_7$,
with $\eta\in\{0,1\}$ in both cases.

\smallskip
\noindent\textup{(2)} Suppose $\cha k=2$.
The profiles of $R_i$ and $T_i$ include
$$
\begin{aligned}
R_3:\quad &(3,1;4,2),\\[3pt]
R_m:\quad &(r,m+1-r;4,2),\\[3pt]
T_3:\quad &(1,3;3,3),\\[3pt]
T_4:\quad &(2,3;3,3),\\[3pt]
T_5:\quad &(4,2;3,3).
\end{aligned}
$$
For $R_m$, the ranges are $m=4,5$ and $3\leq r\leq m$.
The affine conic bundle constructed in the proof has profile
$(2,2;2,2)$.
Proposition~\ref{prop:norm-extension} gives $(3,3;3,2)$,
and Lemma~\ref{lem:binary-profiles} below gives the stronger profile
$(3,3,3;2,2,2)$.
The other affine charts have profiles
$$
\begin{gathered}
(2,1,1;2^{[3]}),\qquad (2,2,1;2^{[3]}),\\[3pt]
(3,1,1,1;2^{[4]}),\qquad (2,2,1,1;2^{[4]}),\\[3pt]
(2,2,2;3,2,2).
\end{gathered}
$$
\end{lemma}
\begin{proof}
We give the coordinate partitions for the equations above
and construct the two conic bundles.

\begin{itemize}
\item[step 1:] In this step, we check the profiles of the explicit equations
in characteristic different from two. The nonzero class for $q_4$ is
proved in \cite[Proposition~2.14]{FiammengoLuders}. Its first group is
$(x_1,x_2)$ for $(2,3;2,2)$ and $(x_1)$ for $(1,4;2,3)$.
The other coordinates form the second group.

For $V_5$, \cite[Example~4.6]{LS} gives a nonzero class after removing a
larger boundary. Restriction from our affine chart to that smaller open
subset sends $\delta_{V_5}$ to the class in that example, so
$\delta_{V_5}\ne0$. Its first group can be $(x_1,x_2,z)$,
$(x_1,x_2,z,y_3)$ or $(y_2)$, giving the profiles
$(3,3;3,2)$, $(4,2;3,2)$ or $(1,5;2,4)$, respectively.

For $W_6$, use $(x,y,u)$, $(x,y,u,a)$, $(x,y,u,b)$ or
$(x,y,u,a,b)$ as the first group. For $W_7$, use $(x,y,w)$, $(x,y,u,v)$, $(x,y,u,v,a)$,
$(x,y,u,v,b)$ or $(x,y,u,v,a,b)$. Again the other coordinates form the second group. For the
profiles $(1,m;2,4)$ of $W_6,W_7$, use $(x)$ as the first group.
The profile $(4,3;3,2)$ also follows from
\cite[Theorem~1.2]{FiammengoLuders}.

For $V_8$, use $(a,b,c)$ together with fiber coordinates,
or use $(a,y_6)$, for its profiles with two groups. For the profiles with three groups, start with $(a)$, $(b,c)$
and $(y_1,\ldots,y_6)$. Move $h$ coordinates from $\{y_2,y_3\}$ to the first group,
and move $y_1$ to the second group when $\eta=1$.
For $W_7$, take $(x)$, $(y,w)$ and $(u,v,a,b,c)$ when
$(h,\eta)=(0,0)$. For $(h,\eta)=(0,1)$, move $a$ from the third group
to the second. For $(h,\eta)=(1,0)$, use $(x,y)$, $(u,v)$ and
$(a,b,c,w)$. For $(h,\eta)=(1,1)$, move $a$ from the third group
to the second in this last partition.
The groups $(x,v)$, $(y,u)$, $(w,c)$ and $(a,b)$ also give
degree at most two in each group. These bounds follow from the
three quartic monomials $xu^2w,yv^2w,xyc^2$ and the lower terms.

The invertible linear substitution $x_1=v$, $x_2=v-u$ gives
$\Phi_5$ of Definition~\ref{def:catalogue}.
The groups $(v,z)$, $(u,y_1)$ and $(y_2,y_3)$ give the profile
$(2,2,2;2^{[3]})$. The groups $(y_2,y_3,z)$, $(v,y_1)$ and $(u)$
give $(3,2,1;2^{[3]})$. The degree in each group is read directly
from $\Phi_5$.

\item[step 2:] In this step, we construct the conic in characteristic
different from two. Work over $\Omega$ and choose $\lambda\in\Omega^\times$
with $16\lambda\ne1$. Let $D=\lambda(1+x)^2(1+y)^2-xy$.
For the equation $xU^2+xyV^2+D=0$, the first three profiles use the
two groups $(x,y,V)$ and $(U)$, then $(x)$ and $(y,U,V)$, then
$(x,V)$ and $(y,U)$. The last profile uses the three groups
$(U,V)$, $(x)$ and $(y)$.
Let $K=\Omega(x,y)$, $\rho=(x,y)_2$ and
$\Pi=\langle1,x,y,xy\rangle$.

We first show that $\rho$ remains nonzero over the function field
of the projective conic defined by $q=\langle x,xy,D\rangle$. The discriminant of $D$ as a quadratic in $x$ is $yB$,
where $B=y-4\lambda(1+y)^2$ has distinct nonzero roots. Thus $D$ is irreducible and
$\Omega(D)=\Omega(y)(\sqrt{yB})$. The element $y$ is not a square in this extension. Otherwise its square class over $\Omega(y)$ would be
$0$ or $[yB]$, making $y$ or $B$ a square. Hence $(y,D)_2$ has nonzero residue at $D=0$.

If restriction to the conic made $\rho$ zero, the Pfister
subform theorem and roundness, using the represented value $x$,
would give $\langle1,y,D/x\rangle\subset\Pi$. Cancelling $\langle1,y\rangle$ would make $D/x^2$ a value
of $\langle1,y\rangle$. This would imply $(y,D)_2=0$, a contradiction.

We next verify condition~\textup{(2)} of
Proposition~\ref{prop:charts}. Suppose that $\rho$ is nonzero over the henselization at a
nontrivial valuation of $K/\Omega$. We cannot have $\mu(x)=\mu(y)=0$, since the residue field
has transcendence degree at most one and $\rho$ would be
an unramified class of degree two. The Pfister minimum formula gives
$\mu(1+x)=\min(0,\mu(x))$ and the analogous equality for $y$. We then have
$$
 \mu(xy)-2\mu((1+x)(1+y))=|\mu(x)|+|\mu(y)|>0.
$$
Since $\lambda$ is a square, Hensel's lemma makes $D$ a square. Thus $q$ embeds into $\Pi$. Every further extension with a point on the quadric of $q$
then splits $\rho$, proving condition~\textup{(2)}.

On $H=\{\zeta=0\}$, the form is $xU^2+xyV^2$,
which embeds injectively into $\Pi$. The generic affine fiber is $Q\setminus H$, where
$Q=\{xU^2+xyV^2+D\zeta^2=0\}$. Its equation is primitive over $\Omega[x,y]$ because $x\nmid D$. Gauss's lemma and integrality of the generic conic give an
affine quadric model. Proposition~\ref{prop:charts} proves
$\delta_V\ne0$ modulo two.

\item[step 3:] In this step, we give the coordinate changes for the quartics
in characteristic two. For $R_i$, use $(a,b,d)$ together with fiber coordinates,
or $(a,d)$, as the first group. To define $T_i$ for $i=3,4,5$, substitute
$a=(A+d^2+b+1)/(b+1)$ in the equation of $R_i$,
replace $V$ by $V+W$, and multiply by $b+1$. Keep $W=0$ for $i<5$ and $V=0$ for $i=3$. For $T_3$, the groups are $(d)$ and $(A,b,U)$. For $T_4$, they are $(d,V)$ and $(A,b,U)$. The degree in each group is at most three.

At $b=1$, the new polynomial is $(A+d^2)(U^2+d)$ for
$i=3,4$, and $(A+d^2)(U^2+d+W^2)$ for $i=5$. These polynomials are nonzero, so no component is contained in $b=1$.
On $b\ne1$, the substitution is invertible, which proves integrality.
In the original model $R_i$, the evaluations on $b=1$ are
$\{a,1\}$ for $i<5$ and $\{a,W^2+a\}$ for $i=5$. They vanish by the Steinberg relation.
Remark~\ref{rmk:localization-criterion} then gives the nonzero affine class. For $T_5$, let $\omega^2+\omega+1=0$ and $D=d^2+d$. Replacing $(A,U,W)$ by $(A+\omega D,U+d,W+d)$ gives the
groups $(A,U,V,W)$ and $(b,d)$.

Define $F_5$ by substituting $\ell=d+1$ in the
homogeneous equation of Definition~\ref{def:binary-quartics}. Let $F_4=F_5|_{W=0}$ and $F_3=F_4|_{V=0}$. For $F_5$, use $(U,V,W)$, $(a)$, $(b)$, $(d)$,
or $(a,V)$, $(b,U)$, $(d)$, $(W)$. For $F_4$, use $(U,V)$, $(a)$, $(b)$, $(d)$,
or $(a,V)$, $(b,U)$, $(d)$. For $F_3$, use $(U)$, $(a)$, $(b)$, $(d)$,
or $(a)$, $(b,U)$, $(d)$.

To compare these charts, we check the divisor $d=1$ in $R_i$, where
$\ell=1$. Use the symbol presentations and the notation $q,s,P,Q$
of Proposition~\ref{prop:binary-models}. On this divisor, $s=(W+1)^2$. If $s\ne0$, the representatives $\{s,P\}$ and
$\{s,Q\}$ in Proposition~\ref{prop:binary-models} vanish whenever
the second slot is nonzero. If $P=Q=0$, the square expression for
$q$ gives zero. If $s=0$, then $W=1$, $P=V^2$ and $Q=U^2$.
A nonzero $P$ or $Q$ is a square slot. If both vanish, the old
equation becomes $((a+1)(b+1))^2=0$, so $a=1$ or $b=1$ and
the symbol is again zero. Points with $aq=0$ have value zero by
definition.

The two charts have the common open subset
$R_i\setminus\{d=1\}\simeq F_i\setminus\{\ell=0\}$.
To pass between them, divide the homogeneous coordinates
$[1:a:b:d:U:V:W]$ by $1+d$ to reach the chart $\ell+d=1$.
The pairing on $R_i$ vanishes on the removed divisor, so localization
gives a pairing on this open subset with nonzero generic value.
Restriction from $F_i$ to it proves $\delta_{F_i}\ne0$.
Each $F_i$ is integral because it is a nonempty chart of the
integral projective closure of $R_i$. Finally, substituting $a=A+W,b=B+W$ in $F_5$ gives
$(d,W)$, $(A,V)$ and $(B,U)$, with bounds $(3,2,2)$.

\item[step 4:] In this step, we construct the affine conic bundle in
characteristic two. Let $Y$ be the $(2,2)$ bundle of
\cite[Definition~6.1 and Theorem~6.2(b),(c)]{ABBBConic}
over an algebraically closed field $\Omega$. Its first projection is the conic bundle in that reference. For the other projection, the coefficients $z^2,x^2+xz,y^2+yz$
have no common projective zero. Its fiber at $(x:y:z)=(1:0:0)$ is the smooth conic
$uv+uw+w^2=0$. Hence both projections are flat and
generically smooth.

Choose general lines in the two bases which avoid the finitely
many images of singular points. Let $D_1,D_2$ be their inverse
images. Their generic fibers are smooth conics over the function
fields of lines. Tsen's theorem \cite[Theorem~59.67.10 (\href{https://stacks.math.columbia.edu/tag/03RD}{Tag~03RD})]{stacks-project} makes these conics rational,
so $D_i$ are rational surfaces. Flatness excludes additional surface components in individual
fibers.

On the smooth threefold, each projection has differential of
rank at least one, since a plane conic fiber has tangent space
of dimension at most two. At a point of rank one, only one line through its image can
have a singular pullback there. The incidence of these points and lines has dimension at most
two. Thus a general line has at most finitely many such points,
and each $D_i$ is integral with only isolated singularities.

A rational resolution of $D_i$ has universally trivial $\CH_0$. Its pushforward to $D_i$ is universally surjective. Every
exceptional image is rational over $\Omega$, and each exceptional
fiber contains a point rational over $\Omega$. The degree map then shows that $D_i$ is universally $\CH_0$
trivial, as in \cite[Proposition~6.13]{ABBBConic}.

Let $\widetilde Y\to Y$ be the universally $\CH_0$ trivial
resolution of \cite[Theorem~6.2(c)]{ABBBConic}. The nonzero Brauer class on $\widetilde Y$ evaluates to zero on
the inverse images of $D_i$. Indeed, these lie over the smooth locus and are isomorphic to
$D_i$. Localization of the pairing gives
$\delta_{\widetilde Y\setminus(D_1\cup D_2)}\ne0$.

For $V=Y\setminus(D_1\cup D_2)$, the restricted resolution
$\widetilde V\to V$ is universally $\CH_0$ trivial by
\cite[Proposition~6.13]{ABBBConic}. Its pushforward is an
isomorphism on $\CH_0$ after every field extension and sends
$\delta_{\widetilde V}$ to $\delta_V$. Hence $\delta_V\ne0$
modulo two. The pairing here is the characteristic-two Brauer
pairing of \cite[Proposition~3.1]{ABBBBrauer}.
\end{itemize}
\end{proof}

\begin{remark}\label{rmk:adjunction}
We will also use affine models obtained by adjunction.
Let $(e_{ij})$ be a positive degree matrix with $q$ rows, and let
$a_j$ be the dimensions of the projective factors. Bertini, adjunction and
Lemma~\ref{lem:forms} give an integral affine variety $V$ with
$\delta_V\ne0$ whenever
$0<\sum_ja_j-q$ and $\sum_i e_{ij}\geq a_j+1$,
with smooth coordinate boundary components and Chow coefficients $\QQ$ in
characteristic zero or $\FF_p$ in characteristic $p$.
For a quartic equation, choose a zero of its leading homogeneous part
as one coordinate axis. The coefficient of the fourth power of that
coordinate is then zero, so its degree in that coordinate is at most three.
For a quartic surface model quadratic in one coordinate, take a
general member of $4H-2E$ on the blow up of
$[0:0:1:0]\in\PP^3_{[x:y:z:w]}$. This system is very ample, and
its general member has trivial canonical bundle, boundary
$E\cup\{w=0\}$ with smooth irreducible components, and an affine equation quadratic in $z$.
Lemma~\ref{lem:forms} again gives a nonzero class.
We will also use the following consequence of the catalogue and
the quadratic dependence on fiber coordinates in the uniform quintic
models. For every $m\geq3$, there is an integral affine variety $V$
with $\delta_V\ne0$ modulo two, of profile
$$(\floor{ m/2}+1,\ceil{ m/2};5,2)$$
and total degree at most five.
\end{remark}

\begin{lemma}\label{lem:balanced-cubic}
Let $m\geq3$ and let positive integers $a_j$ satisfy
$\sum_j a_j=m+1$ and $\max_j a_j\leq\floor{(m+1)/2}+1$.
There is an integral affine hypersurface $V$ of dimension $m$.
Its profile is $(a_1,\ldots,a_r;3^{[r]})$, and
$\delta_V\ne0$ in $\CH_0(V_{\Omega(V)},\ZZ/2)$.
Thus, in Setup~\ref{setup:coordinates}, a cubic row
suffices whenever $2M\leq L+c+1$.
\end{lemma}
\begin{proof}
We first treat the dimensions covered by the uniform equations.
These are $m\geq6$ in characteristic two, and $m=6$ or $m\geq8$
outside characteristic two. Let $\ell=m+1$, and write
$m=2s+\varepsilon$ with $\varepsilon\in\{0,1\}$.
Suppose one group has size $\floor{\ell/2}+1$.
Outside characteristic two, place $Y_0,\ldots,Y_s$ in that group,
together with $Z$ when $\varepsilon=1$.
In characteristic two, place $a,P_1,\ldots,P_s$ there, together with
$X$ when $\varepsilon=1$.
The other groups contain only the bases $z_i$ and have degree at most three.

When $m$ is odd and there are two equal groups, use the chart
$G_s^\sharp$ outside characteristic two.
In characteristic two, use $H_s^1$ with groups
$(a,P_1,\ldots,P_s)$ and $(z_1,\ldots,z_s,X)$.

In every other case, Lemma~\ref{lem:coloring-matching} gives a proper
cyclic coloring with multiplicities $a_j$. For $m=2s$, list the
coordinates as $Y_0,z_1,Y_1,\ldots,z_s,Y_s$ outside characteristic two,
and as $a,z_1,P_1,\ldots,z_s,P_s$ in characteristic two.
For $m=2s+1$, append $Z$ to the first list and $X$ to the second.
Let $\gamma_i$ be the color of the $i$th entry, with indices modulo $\ell$.

We choose a starting point for which neither
$\gamma_0=\gamma_3=\gamma_5$ nor
$\gamma_1=\gamma_5=\gamma_{\ell-1}$ holds.
If the first equality holds at a starting point $i$, move the start to
$i+2$. The first equality then fails because
$\gamma_{i+2}\ne\gamma_{i+3}=\gamma_{i+5}$.
The second fails because $\gamma_{i+1}\ne\gamma_i=\gamma_{i+3}$.
If the first equality never holds but the second holds at every start,
then $\gamma_{i+1}=\gamma_{i-1}$ for every $i$. This is the period-two
coloring with two equal groups, which was treated above.

The proper coloring bounds the degrees of the linear forms in
Construction~\ref{constr:gram-cycle} and of the Hermitian edge terms.
The two additional inequalities exclude a degree-four contribution
from $z_1z_2z_3Y_0^2$, $z_1z_3Z^2$ or $a^2z_1z_2z_3$.
The terms $a^2z_1^2z_i$ and $az_1z_sP_s^2$ also have degree at most
three in each group, since the adjacent entries in the cyclic lists
have different colors.

For $m=3$, refine groups of sizes $(3,1)$ or $(2,2)$.
For $m=4$, refine $(3,2)$ or $(2,2,1)$.
For $m=5$, refine $(4,2)$, $(3,3)$ or $(2,2,2)$.
The equations in Lemma~\ref{lem:catalogue} have degree at most three in
each group for these partitions.
For $m=7$ outside characteristic two, use $W_7$.
Its only quartic monomials are $xu^2w,yv^2w,xyc^2$.
If there are at least three groups, place $x,y,w$ in distinct groups.
If there are two groups, use $(3,5;3,2)$ or $(4,4;3,2)$, and permute
the groups if necessary.

We now prove the consequence for complete intersections.
If $c\geq M-1$, every column sum is at least $c+2\geq N_j+1$,
so adjunction gives the conclusion.
Otherwise, let $m=L-c$, $p=c-1$ and
$v=\floor{(m+1)/2}+1$.
The hypothesis $2M\leq L+c+1$ gives $M\leq p+v$.
If two capacities exceed $v$, truncating both to $v$ already gives
$2v\geq m+1$ coordinates.
If at most one capacity exceeds $v$, truncation removes at most $p$
coordinates, and $L-p=m+1$.
Thus $\sum_j\min\{N_j,v\}\geq m+1$ in both cases.

Choose integers $0\leq a_j\leq\min\{N_j,v\}$ with sum $m+1$,
and omit the groups with $a_j=0$.
The affine construction above uses these coordinates.
Exactly $c-1$ coordinates are unused, and the other ample rows can be
assigned to them. Corollary~\ref{cor:embedding} proves the conclusion.
\end{proof}

\begin{lemma}\label{lem:quintic-placement}
Let $C_0,\ldots,C_r\geq1$, $C_0=\max_jC_j$ and
$m=\sum_jC_j-1>0$. There is a geometrically integral affine hypersurface $V$ with these
group sizes and degrees at most $(5,2,\ldots,2)$ such that
$\delta_V\ne0$ with the coefficients of Lemma~\ref{lem:forms}
for $m\leq2$, and with coefficients $\ZZ/2$ for $m\geq3$.
Outside characteristic two the first degree bound can be replaced by
three when
$m=2s+\varepsilon$, $s\geq3+\varepsilon$, $\varepsilon\in\{0,1\}$ and $C_0\leq s$.
\end{lemma}
\begin{proof}
Let $M=C_0$, and let $a_j=C_j$ for $j>0$.
We call the first group color zero and the other groups colors $j>0$.
We assign the base coordinates first and then the fiber coordinates.

\begin{itemize}
\item[step 1:] In this step, we prove the cubic bound outside characteristic two.
Use $G_s$ from Construction~\ref{constr:gram-cycle}, with the extra
coordinate $Z$ when $\varepsilon=1$.
Let $\ell=s-M$. Place $M$ bases in color zero.
Choose a color $j_0>0$ for $Y_0$, and place no bases in that color.
For $j>0$, let $b_j$ be the number of bases of color $j$.
We need $b_{j_0}=0$, $b_j\leq a_j$ and $\sum b_j=\ell$.
The sum of the capacities outside colors zero and $j_0$ is at least
$2\ell+1+\varepsilon$, so these choices are possible.

Start the base cycle at $z_1$ of color zero.
If $\ell\leq M$, use the color-zero bases to separate all other bases.
If $\ell>M$, Lemma~\ref{lem:coloring-matching} gives a proper coloring,
since every base multiplicity is at most $M<s/2$.
After placing $Y_0$, the fiber quota for color $j$ is
$f_j=a_j-b_j-\delta_{j=j_0}$.
These quotas sum to $s+\varepsilon$.
The coordinate $Y_i$ cannot use either color of $z_i,z_{i+1}$,
and $Y_s$ cannot use the color of $z_s$.
The coordinate $Z$, when present, cannot use the color of $z_3$.

For one color, at most $2b_j+\varepsilon$ positions are forbidden.
Since $b_j\leq\ell=s-M$, we have
$f_j\leq M-b_j\leq s-2b_j$.
The extra position when $\varepsilon=1$ compensates for the possible
extra forbidden position.
For two colors, no position forbids both when $\ell\leq M$.
When $\ell>M$, at most $b_i+b_j$ positions forbid both, and
$f_i+f_j+b_i+b_j\leq a_i+a_j\leq2M<s$.
Lemma~\ref{lem:coloring-matching} then gives the fiber assignment.
Color $j_0$ contains no bases, so the term involving $Y_0^2$ has
degree at most two there. The other degree bounds follow from the
forbidden colors and the formulas for $G_s$.

\item[step 2:] In this step, we treat the small cases in both characteristics
and the other cases outside characteristic two.
Adjunction gives the conclusion for $m\leq2$.
For one group, use $G_s,H_s^\varepsilon$ and the small quartics of
Section~\ref{subsec:small}.
If every group other than the first has size one, place $z_1$ in the
first group for a uniform equation, or place $d$ there for $R_i$.
These equations have total degree at most five and degree at most two
in each other individual coordinate.

For $3\leq m\leq5$ and $M\geq\floor{m/2}+1$, use the final
$(\floor{m/2}+1,\ceil{m/2};5,2)$ profile of Remark~\ref{rmk:adjunction},
whose generic-point class is nonzero modulo two. Move coordinates from
the second group to give the first group size $M$, and split the second group.
The total degree is at most five, and every part of the second group
has degree at most two.
Otherwise refine groups of sizes $(2,2,1)$ or $(2,2,2)$ from
Lemma~\ref{lem:catalogue}.
In characteristic two, the group whose degree may be three is made
the first group.
For $m=7$ outside characteristic two and $M\geq3$, start with the
$(3,5;3,2)$ profile of $W_7$. Move $M-3$ fiber coordinates into the
first group and split the other fiber coordinates among the other groups.
The first degree is at most four, hence at most five, and the other
coordinates still occur quadratically. For $M=2$, refine the four
groups of size two and degree at most two.

For the other dimensions outside characteristic two, if
$M\geq\floor{m/2}+1$, start with the $(5,2)$ profile in
Remark~\ref{rmk:adjunction}. Move enough fiber coordinates into the
first group to give it size $M$, and split the other group.
The total degree remains at most five, and every group of fiber
coordinates has degree at most two.
If $M<\floor{m/2}+1$, use the cubic construction in step~1.

\item[step 3:] In this step, we prove the quintic bound in even dimension
in characteristic two.
The dimensions at most five were treated in step~2, so we may assume $s\geq3$.
Use $H_s^0$, so $M+\sum_j a_j=2s+1$.
If $M\geq s+1$, place $a,z_1,\ldots,z_s$ in the first group,
fill its other $M-s-1$ positions with fiber coordinates, and distribute
the other fiber coordinates among the other groups.
The first degree is at most five and the other degrees are at most two.

Suppose now that $M\leq s$.
Place $a$ and $M-1$ of the bases $z_i$, including $z_1$, in the first group.
For $j>0$, let $b_j$ be the number of bases in group $j$.
We choose these numbers so that
$0\leq b_j\leq\min(a_j,s-a_j)$ and $\sum_jb_j=s-M+1$.
Such a choice is possible if $\sum_j(2a_j-s)_+\leq s$.
At most two capacities $a_j$ exceed $s/2$, because
$M\geq a_j$ and $M+\sum_j a_j=2s+1$.
If only one does, its contribution is at most $s$ since $a_j\leq M\leq s$.
If two do, let their sizes be $u,v$.
Then $u+v\leq\min(2M,2s+1-M)\leq\floor{3s/2}$,
so their two contributions also sum to at most $s$.

Start the base cycle with color zero.
If $2M\geq s+2$, separate the bases of the other colors by those
of color zero.
If $2M\leq s$, properly color the cycle.
If $2M=s+1$, alternate color zero with the other colors, placing the
extra base of a nonzero color last and making the final two nonzero
colors distinct. Each nonzero color occurs at most $M-1$ times,
so this can be done.

The fiber quotas are $f_j=a_j-b_j$, and their sum is $s$.
Each other color can occur at most once in
$(z_i,z_{i+1},P_i)$ and in $(z_s,P_s)$.
For a single color, the matching inequality is $a_j+b_j\leq s$.
For two colors, no position forbids both when their bases are separated.
In a proper coloring, at most $b_j+b_k$ positions forbid both,
and $a_j+a_k\leq2M\leq s$.
In the alternating case, at most one position forbids both.
If that happens, $b_j,b_k\geq1$, so
$f_j+f_k\leq a_j+a_k-2\leq s-1$.
Lemma~\ref{lem:coloring-matching} gives all fiber assignments.

\item[step 4:] In this step, we treat odd dimension in characteristic two.
The dimensions at most five were treated in step~2, so $m=2s+1\geq7$.
For $H_s^1$, place $X$ in a group other than the first and reduce that
group's capacity by one. Discard a group if its capacity becomes zero.
Step~3 assigns all the other coordinates.
The formulas for $H_s^\varepsilon$ give degree at most five in the first
group and at most two in every other group.
\end{itemize}

Only the grouping of coordinates has changed. The affine variety
and its nonzero generic-point class are those of the reference model.
\end{proof}

\begin{corollary}\label{cor:quintics}
Let $P=Q\times R$, where $Q$ is a nonempty product of projective spaces
and every factor of $R$ has dimension at most $\dim Q$, allowing $R$ to be a point.
A positive-dimensional geometrically integral very general complete
intersection $X\subset P$ has no geometric decomposition of the diagonal
if one equation has degree at least five in every factor of $Q$ and
at least two in every factor of $R$, and all other equations are
relatively ample over $Q$ (positive in every $R$ factor, with arbitrary
nonnegative degrees in the $Q$ factors).
\end{corollary}

\begin{proof}
Let $c$ be the codimension of $X$, $a=\dim Q$ and $N=\dim R$.
Let $p=c-1$ and $m=a+N-c=\dim X$.
We first assign the $p$ rows other than the distinguished row to
coordinate positions on which their degrees are positive.

For a set $I$ of these rows, let $P_I$ be the product of the factors
on which at least one row in $I$ has positive degree.
At each point of $P_I$, the equations give $|I|$ independent linear
conditions on their independent coefficient vectors.
If $|I|>\dim P_I$, the incidence variety has dimension smaller than
the coefficient space, so the very general common zero locus is empty.
Since $X$ is nonempty, $|I|\leq\dim P_I$.
The number of coordinate positions available to the rows in $I$ is
$\dim P_I$. Hall's theorem then gives distinct positions for all $p$ rows.

Every one of these rows has positive degree in every factor of $R$.
If a row has been assigned to a position in $Q$ and a position in $R$
is unused, move its assignment to that position in $R$.
After these moves, exactly $\min(p,N)$ positions in $R$ are used.
There are $a+N-p=m+1$ unused positions in total.

If $N\leq p$, all $m+1$ unused positions lie in $Q$.
Join them into one group and use the one-group case of
Lemma~\ref{lem:quintic-placement}.
If $N>p$, all $a$ positions in $Q$ are unused.
Join them into one group. Every nonempty group of unused positions
in $R$ has size at most $a$ by the hypothesis on its factors.
Lemma~\ref{lem:quintic-placement} gives an affine hypersurface of
dimension $m$, with degree at most five in the joined group and at most
two in each group coming from $R$.

Split the joined group according to the original factors of $Q$.
Its degree in each factor is still at most five.
The $p$ positions assigned at the start provide the coordinate equations
for the other rows. Corollary~\ref{cor:embedding} proves the conclusion.
The first case also includes $R$ a point.
\end{proof}

\subsection{Quadratic rows}\label{subsec:quadratic-rows}
We now construct affine models whose equations have degree at most two
in each coordinate group. We first treat $\PP^2\times\PP^2$ and the
smaller characteristic-two examples. We then count the coordinates
that can be added to the cycle constructions. This gives the intervals
needed for equal factors and for five factors of comparable dimensions.
The graph equations of Remark~\ref{rmk:graphs} give models for complete
intersections of higher codimension.

\begin{lemma}\label{lem:planes}
In every characteristic a very general $(d,e)$ hypersurface in
$\PP^2\times\PP^2$, $d,e\geq2$, admits no decomposition of the diagonal.
\end{lemma}
\begin{proof}
In characteristic zero, \cite[Proposition~2.8]{FiammengoLuders} gives
$\delta_V\ne0$ for the affine chart of a very general $(2,2)$ hypersurface.
In characteristic two, Lemma~\ref{lem:catalogue} gives an affine model
with the same bidegree and a nonzero generic-point class.
Corollary~\ref{cor:embedding} gives the conclusion in these two cases.

Suppose now that the characteristic is odd and positive.
We work over an algebraically closed field $\Omega$ and use the
construction of \cite{BB}. We check that it needs division by two only.
In $\PP^2_{[u:v:w]}$, let $a=-u^2-v^2-w^2$ and $A_0=\{a=0\}$.
Choose two disjoint reduced divisors of degree three on the smooth conic
$A_0\simeq\PP^1$. The squares of their binary cubic equations lift to
plane cubics $C_1,C_2$. Indeed, the restriction map is surjective and
its kernel is $aH^0(\PP^2,\cO(1))$.

Choose the lifts generally. Off $A_0$, the kernel separates first jets.
At a point of either contact divisor, it varies the normal derivative.
Thus the plane cubics are smooth and meet transversely off $A_0$.
They have no common point on $A_0$, since their contact divisors are disjoint.
Lift the product of the two binary cubics to a plane cubic $b$.
Then $C_1C_2=b^2-ac$ for a quartic $c$.
At a point with $a=b=0$, exactly one $C_i$ vanishes.
Its normal derivative and the value of the other cubic are nonzero.
Differentiating the identity in the normal direction gives $c\ne0$.
Hence $\{a=b=c=0\}=\varnothing$.

On each $C_i$, adjoining $\sqrt a$ gives an unramified double cover,
since the contact divisor occurs with multiplicity two.
The cover is nontrivial. Otherwise its degree-three contact divisor
would be linearly equivalent to a line section of $C_i$.
Every section of this line bundle comes from a line in $\PP^2$,
but three distinct points of a smooth conic are not collinear.

We next construct the deformation of \cite[Section~3]{BB}.
Let $x=(u,v,w)^{\mathrm t}$ and define $S$ by $Sx'=x\times x'$.
Choose a column $h$ of quadratics with $b=-x^{\mathrm t}h$.
Let $g=Sh$, $J=(\partial g_i/\partial x_j)$, and $f=c+h^{\mathrm t}h$.
Choose a symmetric quadratic matrix $B$ with $x^{\mathrm t}Bx=f$.
Such a matrix exists by distributing each quartic monomial among
products $x_ix_j$. Only the off-diagonal entries need division by two.
For a general symmetric quadratic matrix $C$, let
$$
 N(t)=\begin{pmatrix}
 -S^2-\frac t2(J+J^{\mathrm t})+t^2B+t^3C&g-tBx-t^2Cx\\
 (g-tBx-t^2Cx)^{\mathrm t}&f+t x^{\mathrm t}Cx
 \end{pmatrix}.
$$
Euler's identity and differentiation of $x^{\mathrm t}g=0$ give
$$
 Jx=3g,\qquad J^{\mathrm t}x=-g,
 \qquad N(t)(x,t)^{\mathrm t}=0.
$$
No division by three occurs, so the identities also hold in characteristic three.

The map $\cO(-1)\xrightarrow{(x,t)}\cO^3\oplus\cO(-1)$ is a subbundle
inclusion with quotient of rank three. The matrix $N(t)$ defines a
conic family on this quotient. For $t\ne0$, the quotient is isomorphic
to $\cO^3$, and the quadratic form is the upper block of $N(t)$.
Its dependence on $C$ has invertible coefficient $t^3$.
Taking the coefficients of $C$ to be independent makes the
geometric generic fiber a geometric generic $(2,2)$ hypersurface,
by the coefficient change used in Theorem~\ref{thm:transfer}.

We check the central fiber on $w=1$. Its quotient matrix is
$$
 \begin{pmatrix}
 v^2+1&-uv&vh_3-h_2\\
 -uv&u^2+1&h_1-uh_3\\
 vh_3-h_2&h_1-uh_3&f
 \end{pmatrix}.
$$
The determinant is $C_1C_2=b^2-ac$, and the upper minor is $-a$.
The other two charts are obtained by permuting the coordinates.
Where $a\ne0$, the matrix has rank at least two.
Where $a=0$ and $b\ne0$, it has rank three.
Where $a=b=0$, the condition $c\ne0$ makes its rank two.
Thus no central conic is a zero equation or a double line.
The nonzero fiber equations give flatness, after shrinking the
parameter curve about zero if necessary.

The central total space is smooth over every simple zero of the discriminant.
At a point of $C_1\cap C_2$, take local base parameters $r,s$.
After an \'etale change of coordinates, the conic has equation
$XY+rsZ^2=0$. Its unique singular point is the ordinary double point
$XY+rs=0$. Blowing up these finitely many points resolves the threefold.
The exceptional divisors are smooth rational quadrics, so the resolution
morphism is universally $\CH_0$ trivial.

The two discriminant components and their nontrivial unramified double
covers satisfy the conic bundle criterion following
\cite[Corollary~10]{HKT}. It gives a nonzero unramified Brauer class of
order two on the resolved central fiber. The specialization argument of
\cite[Propositions~6,~8 and Theorem~9]{HKT} works modulo two.
A decomposition on the geometric generic fiber would specialize through
the universally $\CH_0$ trivial resolution. The Brauer pairing vanishes
on constant cycles and is nonzero on the diagonal class over the function
field, giving a contradiction.

It remains to pass to an affine chart. Let $H$ be a very general smooth
$(2,2)$ hypersurface, and let $D=D_1\cup D_2$ be general coordinate
hyperplane sections, one from each factor. Each $D_i$ is a smooth
$(2,2)$ surface in $\PP^1\times\PP^2$. Its smooth generic conic over
$\Omega(\PP^1)$ has a point by Tsen's theorem
\cite[Theorem~59.67.10 (\href{https://stacks.math.columbia.edu/tag/03RD}{Tag~03RD})]{stacks-project}.
Thus $D_i$ is rational and universally $\CH_0$ trivial.

Choose $x_0\in D_1\cap D_2$ over $\Omega$.
After every field extension, the image of $\CH_0(D,\Ftwo)$ in
$\CH_0(H,\Ftwo)$ is generated by $[x_0]$.
If $\delta_{H\setminus D}=0$, localization expresses $\delta_H$ as a
multiple of this class. The degree map makes the multiple equal to one
modulo two. This contradicts the obstruction just proved. Hence
$$
 \delta_{H\setminus D}\ne0,
 \qquad \Tortwo(H\setminus D)=2.
$$
Corollary~\ref{cor:embedding} gives every $d,e\geq2$.
\end{proof}

\begin{lemma}\label{lem:binary-profiles}
In characteristic two the following profiles occur, with $N=3,4$,
$2\leq B\leq7$ and $1\leq j\leq3$:
$$
\begin{gathered}
 (N,N,N;2,2,2),\ 
 (2,2,B;2,2,2),\  
 (2,2,2,j;2,2,2,2),\\
 (1,2,2;2,2,2),\ (1,2,3;2,2,2),\ (1,3,3;2,2,2),\ (1,2,4;2,2,2),\\
 (2,3,3;2,2,2),\ (2,4,4;2,2,2),\ (3,3,4;2,2,2),\ (3,4,4;2,2,2).
\end{gathered}
$$
\end{lemma}
\begin{proof}
We give the coordinate groups for each profile.
Unless we say otherwise, all coordinates are included.
The base substitutions are birational. The changes of fiber coordinates
extend to invertible linear changes preserving the homogenizing hyperplane.
By Remark~\ref{rmk:quadric-changes}\textup{(1)} and
Proposition~\ref{prop:charts}, the primitive numerators define models with
$\delta_V\ne0$. Removing content does not increase any group degree.
It is thus enough to check the degree bounds for the equations below.

\begin{itemize}
\item[step 1:] In this step, we obtain the profiles from $G_t^+$ in
Proposition~\ref{prop:norm-extension}.
The coordinates $R,S$ may be omitted independently.
Let $u=x/a$, $Z=x$, and $v=y/a$.
Replace $X$ by $a^2uX$ and $R$ by $aR$, keep $P,S$, and divide by $a^3$.
The four groups $(u,R)$, $(v,X)$, $(Z,P)$ and $(S)$ have degree at most two.
They give $(2,2,2,1;2,2,2,2)$.
Omitting $S$ gives $(2,2,2;2,2,2)$.

For the next substitution, let
$(a,x,y)=(u^{-1},(uw)^{-1},(uv)^{-1})$.
Replace $(X,R,P,S)$ by $(X,R,uvP,u^2vS)/(u^2vw)$ and multiply the equation
by $u^5v^2w^3$.
Omitting $R,S$ gives the groups $(u)$, $(v,X)$ and $(w,P)$.
Omitting only $R$ gives $(u)$, $(v,X)$ and $(w,P,S)$.
With both coordinates included, use either $(u)$, $(v,X,R)$ and $(w,P,S)$,
or $(u,R)$, $(v,X)$ and $(w,P,S)$.
These are the profiles $(1,2,2)$, $(1,2,3)$, $(1,3,3)$ and $(2,2,3)$,
all with degrees at most two.

For the last substitution in this step, let
$(a,x,y)=((uv)^{-1},(uw)^{-1},w^{-1})$.
Replace $(X,R,P,S)$ by $(X,vR,uvwP,uvwS)/(u^2vw)$ and again multiply by
$u^5v^2w^3$. The groups $(u)$, $(v,w)$ and $(X,R,P,S)$ give
$(1,2,4;2,2,2)$.

\item[step 2:] In this step, we use equations with three base parameters.
For $s=3$, Construction~\ref{constr:cycle-coordinates} gives the four groups
$(u,P_2)$, $(x_2,X)$, $(x_3,P_1)$ and $(Z,P_3)$.
Each group has size two and degree at most two.
To increase the last size to three, use $q$ from
Definition~\ref{def:hermitian-equations}, replace $R$ by $auR$, and place
$R$ in the last group.

For the next profiles, use $Q=q|_{Q_0=1}$ from
Proposition~\ref{prop:nonvanishing}\textup{(2)} with
$(s,J,I,h_0,r,\epsilon)=(3,\{1,3\},\{1,2,3\},0,-1,\theta)$.
Write $S,T$ for $S_1,S_3$. The coordinates $R,S,T$ are independently
optional. Here $\ell=1$ and $d_J=2$, so $d_J+2h_0=2\ell$.
First set $S=T=0$.
Let $a=(uW)^{-1}$, $z_1=W^{-1}$, and $z_2=x_2/(uW)$.
Let $z_3=(uVW)^{-1}$.
Replace $(X,R)$ by $(X,R)/(uVW^2)$ and replace
$(P_1,P_2,P_3)$ by $(P_1/(VW),P_2/W,P_3/(uW))$.
After multiplication by $u^3V^2W^5$, the groups are
$(V,X,P_1)$, $(W,x_2,R)$ and $(u,P_2,P_3)$.
They give $(3,3,3;2,2,2)$.
Omitting $R$ gives $(3,2,3;2,2,2)$, which can be reordered as
$(2,3,3;2,2,2)$.

To obtain $(2,2,B;2,2,2)$ for $4\leq B\leq7$, let
$a=K^{-1}$, $z_1=N^{-1}$, and $z_2=U/(KN)$.
Let $z_3=L/(KN)$.
Replace $(X,R,P_1,P_2,P_3,S,T)$ by
$(X/(KN),R/(KN),P_1,P_2,P_3/K,S,T)$.
The polynomial $K^3N^3Q$ is quadratic in each of
$(K,L)$, $(N,U)$ and $(X,P_1,P_2,P_3,R,S,T)$.
Including zero, one, two or three of $R,S,T$ gives the four last-group sizes.
Step~1 already gives $B=2,3$.

For $(2,4,4;2,2,2)$, set $T=0$.
Let $a=(KM)^{-1}$, $z_1=N^{-1}$, and $z_2=(KN)^{-1}$.
Let $z_3=(KVN)^{-1}$.
Replace $X$ by $X/(KMNV)$, $R$ by $R/(KNV)$, and $P_1$ by $P_1/V$.
Replace $P_3$ by $P_3/(KM)$ and $S$ by $S/V$, and keep $P_2$.
Then $K^3V^2M^2N^3Q$ has degree at most two in each of
$(M,N)$, $(V,X,P_1,S)$ and $(K,R,P_2,P_3)$.

\item[step 3:] In this step, we obtain $(3,3,4)$, $(3,4,4)$ and $(4,4,4)$.
For the first profile, use
$X^2+h(\theta Q_0w_0+\sum_{j=1}^4P_jw_j)=0$
with $(s,J,I,h_0,r,\epsilon)=(4,\{1,4\},\{2,4\},0,-1,\theta)$.
At $t=0$, the degree-five equation and its symbol are those of
Proposition~\ref{prop:nonvanishing}\textup{(1)}.
Here $\ell=1$ and $d_J=3\geq2\ell$.
The functions $D_0=z_4/z_1$ and $D_4=z_2T_4/(a^2z_4)$ verify
condition~\textup{(2)} of Proposition~\ref{prop:charts} by
Lemma~\ref{lem:endpoints}.
The terms involving $t$ use distinct basis vectors, which gives
condition~\textup{(3)}.

Use the substitutions in Construction~\ref{constr:cycle-coordinates},
replacing $P_4$ by $aP_4$.
In its formula for $b=2$, this replaces $t^2Z^2ux_2x_3$ by $t^2Zx_4$.
The groups $(u,x_3,P_4)$, $(x_2,x_4,X)$ and $(Z,P_1,P_2,P_3)$
give $(3,3,4;2,2,2)$.
Adding $\theta e_{1,3}S$ adds the term $ux_3S^2$.
Place $S$ in the second group to obtain $(3,4,4;2,2,2)$.

Finally, use $(s,b)=(5,2)$ in Construction~\ref{constr:cycle-coordinates}.
Let $T=x_4$ and $H=x_4+x_5$.
The groups $(u,x_3,P_4,P_5)$, $(x_2,H,X,P_3)$ and $(Z,T,P_1,P_2)$
give $(4,4,4;2,2,2)$.
\end{itemize}
\end{proof}

\begin{definition}\label{def:truncated-counts}
For nonnegative integers $e,t$, define
$S_e(t)=\sum_{i=0}^{\min(e,t)}\binom ti$ and $S(t)=S_2(t)$.
Thus $S_e(t)$ counts subsets of a set with $t$ elements whose
cardinality is at most $e$. Products of these numbers count
independent subset choices in disjoint coordinate groups.
\end{definition}

\begin{proposition}\label{prop:arrangement}
The following profiles occur.

\begin{itemize}
\item[(i)] If $\cha k\ne2$, let $s\geq3$, $e=\floor{ s/2}$,
$o=\ceil{ s/2}$. Let $\varepsilon\in\{0,1\}$,
$0\leq u\leq\binom e2$ and $0\leq v\leq\binom o2-\delta_{2\nmid s}$.
Then the profiles are
$$
(o+u,e+\varepsilon+v,B;2,2,2),\qquad
s+1-\varepsilon\leq B\leq S(e)S(o)-s-1-\varepsilon-u-v.
$$

\item[(ii)] If $\cha k=2$, let $e\geq2$, $\varepsilon\in\{0,1\}$,
$0\leq\delta\leq\varepsilon$ be integers. Let
$s=2e+\varepsilon$, $F=S_2(e+\varepsilon)S_2(e)$ and
$f_A=\binom e2+e-2+\varepsilon$. Let
$f_C=\binom{e+\varepsilon}{2}-\varepsilon$, $a_0=f_A+\binom e2$ and $c_0=2f_C$.
For integers $u,v\geq0$ with
$u\leq a_0+1$, $v\leq c_0+1$ and
$(u-a_0)_++(v-c_0)_+\leq1$, the profiles are
$$
\begin{gathered}
(e+\varepsilon+u,e+1+v+\delta,B;2,2,2),\\
s+1-\delta\leq B\leq
F-s-(u-f_A)_+-(v-f_C)_+-\delta.
\end{gathered}
$$
\end{itemize}
\end{proposition}
\begin{proof}
Use Lemma~\ref{lem:cycle-charts}\textup{(1)} when $\cha k\ne2$, and
Construction~\ref{constr:cycle-coordinates} when $\cha k=2$.
We count separately the coordinates added to a first group and those
moved there from the last group. Only a move decreases the number of
coordinates available for the last group.
All added form coordinates are distinct. Set every unused optional
coordinate equal to zero.

\begin{itemize}
\item[step 1:] In this step, we prove part~\textup{(i)}.
Write $z_1=\xi$, $z_j=x_j$ for odd $j>1$, and $z_j=\xi x_j$ for even $j$.
Replace $X$ by $\xi X$ and $Y_s$ by $\xi Y_s$, and divide the equation by $\xi$.
The first group contains the $o$ odd base variables, and the second
contains the $e$ even base variables.
For odd $i>1$, the corresponding term from
Construction~\ref{constr:gram-cycle} is a multiple of
$$
 \xi x_i(1-x_{i-1}Y_{i-1}+x_{i+1}Y_i)^2.
$$
For even $i$, it is a multiple of
$$
 x_i(\xi-x_{i-1}Y_{i-1}+x_{i+1}Y_i)^2.
$$
Both expressions, and the term from the last vector, have degree at
most two in each base group.
The scalar term is $4\xi X^2$, so $X$ may stay in the last group or move
to the even group. These choices give $\varepsilon=0$ and $\varepsilon=1$,
respectively. The initial group sizes are $(o,e+\varepsilon,s+1-\varepsilon)$.

For an unused support $I$, let $h_I=\delta_{1\in I}$ and
$r_I=|I\cap\{2,4,\ldots\}|$.
Let $\delta_I\in\{0,1\}$ satisfy $\delta_I\equiv h_I+r_I-1\pmod2$.
Rescaling its coordinate by $\xi^{-\floor{(h_I+r_I-1)/2}}$ makes its
coefficient $4\xi^{\delta_I}x_{I\setminus\{1\}}$.
The coordinate can belong to the last group exactly when
$|I\cap\{3,5,\ldots\}|+\delta_I\leq2$ and $r_I\leq2$.
For each set of even indices, the two choices of $h_I$ give the two
values of $\delta_I$. The odd-index count is
$S_2(o-1)+S_1(o-1)=S_2(o)$.
Thus there are $S(e)S(o)$ supports, including the $2s+2$ supports
already used by the original equation.

A square coordinate can move into a base group when its coefficient
contains no variable from that group.
For the odd group, the unused supports are $\{1,i,j\}$ with $i,j$ even.
There are $\binom e2$ of them.
For the even group, they are $\{1\}\cup J$, where $J$ consists of one or
two odd indices other than $1$. Exclude $\{1,s\}$ when $s$ is odd.
There are $\binom o2-\delta_{2\nmid s}$ of these supports.
The two collections are disjoint.

Moving $u$ coordinates to the odd group and $v$ to the even group gives
first-group sizes $o+u,e+\varepsilon+v$.
The largest last-group size is
$$
 (s+1-\varepsilon)+S(e)S(o)-(2s+2)-u-v
 =S(e)S(o)-s-1-\varepsilon-u-v.
$$
Omitting optional last coordinates gives every integer down to
$s+1-\varepsilon$. This proves the interval in part~\textup{(i)}.

\item[step 2:] In this step, we prove part~\textup{(ii)}.
Let $A$ contain the odd base variables and $C$ the even base variables
and $X$. Let $D$ contain $Z,P_1,\ldots,P_s$.
If $s$ is odd and $\delta=1$, use $p_s=a^{-1}e_{1,s}$.
Its coefficient contains no variable from $C$, so $P_s$ may move to $C$.
The initial group sizes are $(e+\varepsilon,e+1+\delta,s+1-\delta)$.

A coordinate $U_I$ in the normalized norm polynomial has coefficient
degree at most two in each base group exactly when $I$ contains at most
two odd and two even indices.
There are $F=S_2(e+\varepsilon)S_2(e)$ such supports.
After removing the $2s+2$ supports used by the original equation,
the largest last-group size is $F-s-1-\delta$.

We first count the coordinates for $A$.
For a pair of even indices $i,j$, add $V_{1,i,j}$ to $A$.
We can also add $V_{1,j}$ for even $j$, except for the path support
$\{1,2\}$ and, when $s$ is even, the last-vector support $\{1,s\}$.
These additions do not remove a coordinate from $D$.
Their number is $f_A=\binom e2+e-2+\varepsilon$.
We can also move the $\binom e2$ coordinates $U_{i,j}$ with even $i,j$
from $D$ to $A$. Their supports differ from those of the added coordinates.
Thus $A$ can receive $a_0=f_A+\binom e2$ coordinates before using the
extra scalar coordinate below.

For $C$, take a pair of odd indices, excluding $\{1,s\}$ when $s$ is odd.
For a pair $\{1,j\}$, add $V_{1,j}$ and move $U_{1,j}$ from $D$.
For a pair $\{i,j\}$ with $1\notin\{i,j\}$, add $V_{1,i,j}$ and move
$U_{i,j}$ from $D$.
There are $f_C=\binom{e+\varepsilon}{2}-\varepsilon$ pairs.
Thus $f_C$ additions do not change the size of $D$, and $f_C$ moves
remove one last coordinate each, giving $c_0=2f_C$ possibilities.
Both coordinates of a binary norm can be used, since their mixed term
has degree two in $C$.

The collections for $A$ and $C$ are disjoint because their indices
other than $1$ have different parities.
Their degree bounds follow from the two expressions for additional norm
terms in Construction~\ref{constr:cycle-coordinates}.
When only $V_I$ is used, set $U_I=0$.
In particular, such an addition does not contribute a coordinate to $D$.

Use $q$ from Definition~\ref{def:hermitian-equations} in place of the
equation with $X^2$. After replacing $R$ by $aR$ in
Construction~\ref{constr:cycle-coordinates}, the extra coordinate $R$
can belong to any of $A,C,D$.
Initially place it in $D$, increasing the largest last-group size to
$F-s-\delta$.
Fill $A,C$ first with the $f_A,f_C$ additions that do not use a last coordinate.
Each further coordinate placed there costs one coordinate of $D$.
The coordinate $R$ can be one of these moves, or can give the single
coordinate beyond $a_0,c_0$ allowed by
$(u-a_0)_++(v-c_0)_+\leq1$.
Thus the largest last-group size is
$$
 F-s-(u-f_A)_+-(v-f_C)_+-\delta.
$$
The first two group sizes are $e+\varepsilon+u,e+1+v+\delta$.
Independent omissions of optional last coordinates give every integer
down to the initial size $s+1-\delta$.
\end{itemize}

We finally check nonvanishing in both characteristics.
Conditions~\textup{(2)} and~\textup{(3)} of Proposition~\ref{prop:charts}
hold for every coordinate choice above, including both coordinates of a
binary norm. Setting the added coordinates to zero gives an equation
from Lemma~\ref{lem:cycle-charts} or
Proposition~\ref{prop:nonvanishing}\textup{(2)}.
Remark~\ref{rmk:quadric-changes}\textup{(2)} preserves integrality and
condition~\textup{(1)} of Proposition~\ref{prop:charts}, also for the mixed terms.
Part~\textup{(1)} of that remark covers the coordinate scalings.
Taking primitive numerators and using Proposition~\ref{prop:charts}
gives $\delta_V\ne0$ modulo two.
\end{proof}

\begin{lemma}\label{lem:interval}
In every characteristic and for $N\geq3$, $N\leq B\leq4N$, there is
an integral affine hypersurface $V$ of profile $(N,N,B;2,2,2)$ with
$$
 0\ne\delta_V\in\CH_0(V_{\Omega(V)},\ZZ/2).
$$
For $N=2$ they exist for $2\leq B\leq7+\delta_{\cha k\ne2}$.
\end{lemma}
\begin{proof}
We use Proposition~\ref{prop:arrangement} to obtain two overlapping
ranges of last-group sizes. Lemma~\ref{lem:binary-profiles} and the
smaller models give the first few values.

\begin{itemize}
\item[step 1:] In this step, we treat characteristic different from two.
In part~\textup{(i)}, take $s=2N$, $\varepsilon=0$, and $u=v=0$.
The first two group sizes are $N,N$, and the last-group interval is
$[2N+1,S(N)^2-2N-1]$.
Its upper bound is at least $4N$ for $N\geq2$.

For $N=3,4,5$, use respectively
$(s,\varepsilon,u,v)=(4,1,1,0),(5,1,1,1),(6,1,2,1)$.
These give the intervals $[4,9]$, $[5,19]$ and $[6,38]$.
Each joins the interval starting at $2N+1$.
The value $B=N$ follows from $V_8$ for $N=3$ and
Proposition~\ref{prop:quartics} for $N=4,5$.

For $N\geq6$, take
$(s,\varepsilon,u,v)=(N,1,\floor{N/2},\ceil{N/2}-1)$.
All parameter bounds in part~\textup{(i)} hold.
This gives the interval
$[N,S(\floor{N/2})S(\ceil{N/2})-2N-1]$,
whose upper bound is at least $2N$.
It is adjacent to or overlaps the interval starting at $2N+1$.

For $N=2$, the equations $\Phi_5$, the cycle with $s=3$, and $W_7$
give $B=2,3,4$. The first interval gives every $5\leq B\leq8$.
This proves the range in every case.

\item[step 2:] In this step, we treat characteristic two.
For $N\geq5$, use part~\textup{(ii)} with
$e=\floor{(N-1)/2}$, $\varepsilon=1$, and $s=2e+1$.
Let $u=N-e-1$, $v=u-1$, and $\delta=1$.
The first two group sizes are $N,N$.
The last-group interval begins at $2e+1\leq N$ and reaches at least
$2N-2$. It thus covers $[N,2N-2]$.

For the upper range with $N\geq3$, take $s=2e$ and
$(e,\varepsilon,u,v,\delta)=(N-1,0,1,0,0)$.
The first two group sizes are again $N,N$, and
$$
 2N-1\leq B\leq S(N-1)^2-2N+2.
$$
The upper bound is at least $4N$ for $N\geq3$.
Thus this interval joins the lower range for $N\geq5$.

For $N=3,4$, Lemma~\ref{lem:binary-profiles} gives $B=N$.
It also gives $(3,3,4;2,2,2)$, after which the upper interval covers
$5\leq B\leq12$.
For $N=4$, the tuple $(e,\varepsilon,u,v,\delta)=(2,1,1,0,1)$ gives
$5\leq B\leq22$.
The same lemma gives every $2\leq B\leq7$ when $N=2$.\qedhere
\end{itemize}
\end{proof}

\begin{lemma}\label{lem:graphs}
All models in this lemma are integral affine complete intersections
$V$ with $\delta_V\ne0$ in $\CH_0(V_{\Omega(V)},\ZZ/2)$.
Outside characteristic two, there are such models with two equations
and coordinate groups of sizes
$(4,3),(4,4),(5,4),(5,5)$, of dimensions $5,6,7,8$, respectively.
In characteristic two, there is a fourfold with two equations and
coordinate groups of sizes $(3,3)$, and a fivefold with three equations
and coordinate groups of sizes $(4,4)$.
In every characteristic and for $s\geq4$, $\epsilon\in\{0,1\}$,
there is such an affine model of dimension $2s+\epsilon$ with $s+3$
equations in two groups of sizes differing by at most one.
In all these presentations, the main equation is quadratic in each
group and every other equation is linear in each group.
\end{lemma}
\begin{proof}
We use the monic graph equations of Remark~\ref{rmk:graphs}.
A copied coordinate may be used independently in different monomials,
since elimination of its graph equation restores the same polynomial.
Suppose a monomial has total degree $p+q+t\leq4$.
Here $p,q\leq2$ count the occurrences fixed in the two groups, and $t$
counts occurrences for which a copy is available in either group.
If $t_1$ of these occurrences use the first group, the degree bounds are
$p+t_1\leq2$ and $q+t-t_1\leq2$.
An integer $t_1$ exists because
$\max\{0,q+t-2\}\leq\min\{t,2-p\}$.
We use this choice monomial by monomial below.

\begin{itemize}
\item[step 1:] In this step, we construct the six models in the smaller dimensions.
For $\Phi_5$, add $u'=u$ and use the groups
$(y_2,y_3,z,u)$ and $(v,y_1,u')$.

For $W_6,W_7$, add $x'=x$.
Replace $x$ by $x'$ only in the terms $xu^2w$ and $xyc^2$.
The groups are $(x,b,u,c)$ and $(y,a,w,x')$ for $W_6$,
and $(x,b,u,v,c)$ and $(y,a,w,x')$ for $W_7$.

For $V_8$, first replace $c$ by $c-b$, $y_2$ by $y_2-y_3-2y_6$,
and $y_4$ by $y_4+2y_6-y_5$.
Add $a'=a$ and use $(b,a,y_3,y_5,y_6)$ and $(a',c,y_1,y_2,y_4)$.

For $F_4$, add $d'=d$ and use $(a,V,d)$ and $(b,U,d')$.
For $F_5$, add $d'=d$ and $W'=W$ and use
$(a,V,d,W)$ and $(b,U,d',W')$.

For $\Phi_5,V_8,F_4,F_5$, distribute the copied occurrences by the rule
above. Every copy equation is linear in each group.
Elimination gives the original affine variety and preserves its
nonzero generic-point class.

\item[step 2:] In this step, we construct the uniform models outside
characteristic two.
Use $G_s^\sharp$ when $\epsilon=0$ and
$G_s^\sharp+4z_1z_3Z^2$ when $\epsilon=1$.
Let $s_1=\floor{s/2}$ and $s_2=s-s_1$.
Place the first $s_1$ bases in $A$ and the other $s_2$ bases in $B$.
For $j\ne s_1$, place $Y_j$ opposite $z_j,z_{j+1}$.
For $j=s_1$, use one copy in each group. In each $\ell_i$, use the copy
opposite the base variable multiplying it.
Thus only $Y_{s_1}$ needs two copies.
In a term involving a square of this coordinate, the product of its
two copies may be used.
Let $Y_s$ belong to $A$.
Adjoin $T_i=\ell_i$ in the group opposite $z_i$ for each $i$, and
adjoin the copy $d=z_1$ in $B$.
There are $s+2$ graph equations: the $s$ equations for $T_i$, one for
the second copy of $Y_{s_1}$, and one for $d$.

If $\epsilon=1$, place $Z$ in $A$ when $s_1=2$ and in $B$ when $s_1\geq3$.
Replace $\ell_i$ by $T_i$ in the main equation.
Distribute the copied occurrences as above.
Each monomial then has total degree at most four, with at most two
occurrences fixed in either group.
The graph equations are linear in each group, and the main equation is
quadratic in each group.
Before placing $X$, the group sizes are
$s+s_2+1+\epsilon\delta_{s_1=2}$ and
$s+s_1+1+\epsilon\delta_{s_1\geq3}$.
Place $X$ in a smaller group.

\item[step 3:] In this step, we construct the uniform models in characteristic two.
Use the standard Hermitian equation, writing $\omega$ for its parameter $a$.
Let $s_1=\floor{s/2}+1$ and $s_2=s-s_1$.
Place the first $s_1$ bases in $A$ and the other $s_2$ bases in $B$.
Place $\omega$ in $B$.
Adjoin $\omega'=\omega$ in $A$ and $h=\omega z_1$ in $B$.

For each $j<s$ with $j\ne s_1$, the bases $z_j,z_{j+1}$ belong to the
same group. Place $P_j$ and $R_j=z_jP_j$ in the other group.
For the edge $j=s_1$, place $P_{s_1}$ in $B$ and adjoin $v,w$ in $A$
with equations $v=z_{s_1}+z_{s_1+1}+t^2\omega$ and $w=vP_{s_1}$.
Place $P_s$ in $A$.
Let $\Psi=\omega^2$ when $\epsilon=0$, and let $\Psi=X^2$ when $\epsilon=1$.
The main equation is
$$
\begin{aligned}
&\Psi+\omega z_1^2+\omega^2z_1+h^2\sum_{i\leq s_1}z_i
 +\omega'hz_1\sum_{i>s_1}z_i+t^2\omega hz_2z_3\\
&+\sum_{j<s_1}z_{j+1}(z_j+z_{j+1}+t^2\omega')R_jP_j\\
&+\sum_{s_1<j<s}z_{j+1}(z_j+z_{j+1}+t^2\omega)R_jP_j
 +z_{s_1}z_{s_1+1}P_{s_1}w+z_s(1+t^2h)P_s^2=0.
\end{aligned}
$$
The two sums use the copy of $\omega$ in the same group as their base
variables. Each term has degree at most two in each group.
The graph equations are linear in each group.
They consist of the two equations for $\omega',h$, the $s-2$ equations
for $R_j$, and the two equations for $v,w$.
Before the optional coordinate $X$, the group sizes are
$s+s_2+2$ and $s+s_1+1$.
When $\epsilon=1$, place $X$ in a smaller group.
\end{itemize}

In both uniform constructions, the final group sizes differ by at most one.
There are $3s+3+\epsilon$ coordinates and $s+3$ equations, so the dimension
is $2s+\epsilon$.
The $s+2$ graph equations can be eliminated one after another.
This gives the original affine variety, preserving integrality and
$\delta_V\ne0$.
In characteristic two, elimination also gives $h^2=\omega'hz_1$.
\end{proof}

\begin{proposition}\label{prop:equal-profiles}
Let $N\geq2$, $r\geq2$ and $1\leq c<N$. In every characteristic there
are an algebraically closed field $\Omega$ and an integral affine
complete intersection
$V=\{F_1=\cdots=F_q=0\}\subset\prod_{j=1}^r\bA^{C_j}_{\Omega}$,
where $1\leq C_j\leq N$, $1\leq q\leq c$ and $\sum_jC_j-q=rN-c$,
with $\Tortwo(V)=2$, such that $\deg_jF_i\leq1$ for $i\geq2$ and
$$
 (\deg_1F_1,\ldots,\deg_rF_1)\leq
 \begin{cases}(3,2),&r=2,\\(2,\ldots,2),&r\geq3.\end{cases}
$$
One may take $q=1$ whenever $r\geq3$ or $c=1$.
\end{proposition}
\begin{proof}
We construct a model of dimension $rN-c$ in each case. For a
hypersurface, this means using $rN-c+1$ coordinates. If graph equations
are needed, each new equation comes with one new coordinate.

\begin{itemize}
\item[step 1:] In this step, we treat $r\geq4$ and the case $r=3,c=1$.
For $r=3,c=1$, Lemma~\ref{lem:interval} with $B=N$ gives
the profile $(N,N,N;2,2,2)$, including $N=2$.

Suppose $r\geq4$. Let $\ell=\floor{(r-1)/3}$, $t=\ell N$ and
$B=(r-2\ell)N-c+1$. Then $t<B\leq(\ell+3)N\leq4t$.
Except when $t=2,B=8$, Lemma~\ref{lem:interval} gives a hypersurface
with profile $(t,t,B;2,2,2)$. Split each of the first two groups into
$\ell$ groups of size $N$, and split the last group into $r-2\ell$
nonempty groups of sizes at most $N$. This is possible because
$B=(r-2\ell)N-c+1$ and $c<N$.
The sum of their sizes is $rN-c+1$, so $q=1$.

The exception $t=2,B=8$ occurs only for $(N,r,c)=(2,6,1)$.
Use the profile $(4,4,4;2,2,2)$ from Lemma~\ref{lem:interval}
and split each group into two groups of size two.

\item[step 2:] In this step, we treat $r=3$ and $2\leq c<N$.
We construct three groups directly from
Proposition~\ref{prop:arrangement}.
Suppose first that the characteristic is different from two.
Let $e=\floor{N/2}$ and $o=\ceil{N/2}$. Let $s=B=N$ and
$\varepsilon=1$. Choose nonnegative integers $u,v$ such that
$u+v=N-c$, $u\leq\min\{\tbinom e2,N-o\}$ and
$v\leq\min\{\tbinom o2-\delta_{2\nmid N},N-e-1\}$.
For $N\geq4$, the two upper bounds for $u,v$ sum to at least $N-2$,
so this choice is possible. Also $S(e)S(o)-2N-1\geq N$.
Proposition~\ref{prop:arrangement}\textup{(i)} then gives group sizes
$(o+u,e+1+v,N)$. Each is at most $N$, and their sum is $3N-c+1$.
For $(N,c)=(3,2)$, use the profile $(2,3,3;2,2,2)$ of $W_7$.

Now suppose the characteristic is two and $N\geq5$.
Let $s=N-1$, $e=\floor{(N-1)/2}$ and
$\delta=\varepsilon=\delta_{2\mid N}$. Let $B=N$.
Choose $u,v$ by $u+v=N-c+1-\delta$ and $u=\min\{u+v,e+1\}$.
Then $0\leq v\leq e-1$, and the inequalities of
Proposition~\ref{prop:arrangement}\textup{(ii)} hold.
For $e\geq3$, its constants satisfy $f_A\geq e+1$ and $f_C\geq e-1$,
so neither loss term in the upper bound for the last group occurs.
Since $F\geq(e+1)(e+\varepsilon+1)\geq4e+1+3\varepsilon$,
that upper bound is at least $N$.
For $e=2$, the upper bound is at least ten when $\varepsilon=0$
and at least twenty-one when $\varepsilon=1$.
If $e=2$, $\varepsilon=0$ and $u=3$, use the additional scalar coordinate.
The first two group sizes are $e+\varepsilon+u$ and
$e+1+v+\delta$. They are at most $N$, and all three sizes sum to
$3N-c+1$.
For $(N,c)=(3,2),(4,2),(4,3)$, use respectively the profiles
$(2,3,3;2,2,2)$, $(3,4,4;2,2,2)$ and $(2,4,4;2,2,2)$
of Lemma~\ref{lem:binary-profiles}. Thus $q=1$ in every case with $r\geq3$.

\item[step 3:] In this step, we treat $r=2,c=1$.
For $N=2$, use the conic of Lemma~\ref{lem:catalogue}.
Outside characteristic two, use $V_5,W_7$ for $N=3,4$.
For $N\geq5$, use the odd-dimensional chart of
Construction~\ref{constr:gram-cycle} with $s=N-1$.
In characteristic two, use Proposition~\ref{prop:norm-extension}.
These hypersurfaces have dimension $2N-1$, and each equation has
bidegree at most $(3,2)$.

\item[step 4:] In this step, we treat $r=2$ and $c\geq2$.
Let $m=2N-c$, and write $m=2s+\varepsilon$ with
$\varepsilon\in\{0,1\}$. If $m\leq5$, the only cases are
$(N,c)=(3,2),(4,3)$.
Outside characteristic two, use $q_4$ and the $\Phi_5$ graph.
In characteristic two, use the $F_4,F_5$ graphs of
Lemma~\ref{lem:graphs}, with $q=2,3$, respectively.

Suppose $m\geq6$ and the characteristic is different from two.
The chart of Construction~\ref{constr:gram-cycle} has group sizes
$(s+1,s+\varepsilon)$ and degrees $(3,2)$.
For $m=7$, use instead the profile $(4,4;3,2)$ of $W_7$.
The group sizes are at most $N$ because $c\geq2$.

In characteristic two, adjoin $h=a$ to $H_s^{\varepsilon}$.
Replace $a$ by $h$ only in the coefficients of the edges and the last vector.
For $\varepsilon=0$, use the groups
$(z_1,\ldots,z_s,h)$ and $(a,P_1,\ldots,P_s)$.
For $\varepsilon=1$, use
$(z_1,\ldots,z_s,h,X)$ and $(a,P_1,\ldots,P_s)$.
Their sizes are $(s+1+\varepsilon,s+1)\leq(N,N)$, and the
main equation has degrees $(3,2)$.
The equation $h-a=0$ is linear in each group. Eliminating $h$
gives the original affine variety and preserves its nonzero class.
Here $q=2\leq c$. Outside characteristic two, the construction above
uses $q=1$. In both cases the number of coordinates minus $q$ is
$2N-c$.\qedhere
\end{itemize}
\end{proof}

\begin{lemma}\label{lem:five-factors}
Let integers $n_1,\ldots,n_5$ satisfy
$2\leq N_{\min}=n_1\leq n_2\leq\cdots\leq n_5\leq N_{\min}^2$.
There is an integral affine hypersurface $V$ with profile
$(n_1,n_2,n_3,n_4,n_5;2^{[5]})$ and
$\delta_V\ne0$ in $\CH_0(V_{\Omega(V)},\Ftwo)$.
\end{lemma}
\begin{proof}
Let $b=N_{\min}$. We first describe the coordinate changes and the
degree bounds used in the three cases below.
For $s\geq5$, replace $z_s$ by $h-\sum_{i<s}z_i$ in $G_s$ of
Construction~\ref{constr:gram-cycle}.
Let $B=(z_i)_{i\leq s-2,\,2\nmid i}$ and
$C=((z_i)_{i\leq s-2,\,2\mid i},h)$.
Let $A=(z_{s-1})$ and $D=(Y_0,\ldots,Y_s)$.
The equation has degree at most two in each group.
For the last adjacent pair, let $\xi=z_{s-1}$ and
$\eta=h-\sum_{i<s-1}z_i$. The possible cubic terms in $\xi$ cancel because
$$
 \xi(\eta-\xi)^2+(\eta-\xi)\xi^2=\xi\eta(\eta-\xi).
$$
Any of $Y_0,\ldots,Y_{s-3}$ may be moved into $A$.
When $s$ is odd, $Y_{s-2}$ may also be moved into $C$.
Their square coefficients do not involve the destination group.
A mixed term involving that group contains at most one moved fiber
coordinate, so its degree in the group stays at most two.
The unused term $4z_2z_{s-1}R^2$ permits a new coordinate $R$ in $B$.

In characteristic two, write $a_0$ for the parameter $a$ of
Definition~\ref{def:hermitian-data}.
The scalar Hermitian equations have no mixed $P_iP_j$ terms.
Their adjacent coefficients satisfy the same identity, and
$a_0^2z_1^2\Sigma_s$ becomes $a_0^2z_1^2h$.
Thus the same coordinate groups and moves satisfy the degree bounds.
Write $(F_j,U)=(Y_j,Y_0)$ outside characteristic two and
$(F_j,U)=(P_j,a_0)$ in characteristic two.
Write $V$ for the optional odd coordinate.

We add coordinates by using pure square terms $z_IR_I^2$.
To place $R_I$ in a group, its coefficient must avoid that group.
If $s\notin I$, allow at most two indices of each parity among
$1,\ldots,s-2$. If $s\in I$, allow at most one of each parity there.
The substitution for $z_s$ then adds at most one variable from either
parity group. The index $s-1$ may also occur.
These conditions keep every coefficient degree at most two.

We order the five groups as $(B,C,A,D,Q)$, with final sizes
$(n_1,n_2,n_3,n_4,n_5)$.
We count the available supports before adding coordinates, so the
coordinates used to enlarge different groups can be distinct.

\begin{itemize}
\item[step 1:] In this step, we treat $n_3=b\geq3$.
Then $n_1=n_2=n_3=b$.
Let $s=2b-1$, $\tau=s-1$ and $O=\{1,3,\ldots,s-2\}$.
Let $E=\{2,4,\ldots,s-3\}$.
Let $B=(z_i)_{i\in O}$, $C=((z_i)_{i\in E},h,F_{s-2})$ and
$A=(z_\tau,F_1,\ldots,F_{b-1})$.
Let $D=(V,(F_j)_{b\leq j\leq s,\ j\ne s-2})$ and $Q=(U)$.
The five sizes are $(b-1,b,b,b,1)$.

There are $2S_2(b-1)S_2(b-2)$ allowed supports without $s$ and
$2b(b-1)$ with $s$.
The factor two counts whether $\tau$ occurs.
Remove the $4b$ supports used in the initial equation and the support
$\{1,3\}$. The latter is used by the odd-dimensional coordinate in
Construction~\ref{constr:gram-cycle}.
In characteristic two we leave it unused.
The number of supports available for new coordinates is
$$
 \begin{aligned}
 W_b&=2S_2(b-1)S_2(b-2)+2b(b-1)-4b-1,\\
 W_b-(2b^2-b)&=\frac{(b-3)(b^3-b^2+6b-2)}2\geq0.
 \end{aligned}
$$
Use $\{2,\tau\}$ to add one coordinate to $B$.
Its coefficient avoids $B$, so this brings $B$ to size $b$.
The first three groups now have their final sizes.
The combined deficit in $D,Q$ is at most
$(b^2-b)+(b^2-1)=2b^2-b-1$.
At least $W_b-1$ supports are still available, and their coefficients
avoid both $D$ and $Q$. Use distinct coordinates from these supports
to fill the two deficits.

\item[step 2:] In this step, we treat $n_3=b=2$.
Use $s=4$ and the substitution $z_4=h-z_1-z_2-z_3$.
Take the five groups
$(z_1,F_2)$, $(h,z_2)$, $(z_3,F_1)$, $(F_3,F_4)$ and $(U,V)$.
Direct substitution gives degree at most two in each group, by the
same adjacent-pair identity above.
Each group has size two, so the first three already have their final sizes.
The unused supports $\{2,4\}$, $\{1,2,4\}$, $\{1,3,4\}$,
$\{2,3,4\}$ and $\{1,2,3,4\}$ have coefficient degree at most two
in the first three groups and avoid the last two.
Each of the last two groups has final size at most four.
Their combined deficit is at most four, so these five supports suffice.

\item[step 3:] In this step, we treat $n_3>b$.
Let $s=2b+1$, $\tau=s-1$ and $O=\{1,3,\ldots,s-2\}$.
Let $E=\{2,4,\ldots,s-3\}$, so $|O|=b$ and $|E|=b-1$.
In the even-dimensional equation of Construction~\ref{constr:gram-cycle},
adjoin the coordinate with support $I_*=\{1,2,\tau\}$.
In characteristic two, take the odd scalar equation and leave $I_*$ unused.
Let $\nu=\min(n_3-1,2b-2)$.
Let $B=(z_i)_{i\in O}$, $C=((z_i)_{i\in E},h)$ and
$A=(z_\tau,F_1,\ldots,F_\nu)$.
Let $D=(F_{\nu+1},\ldots,F_s)$ and $Q=(U,V)$.
The five sizes are $(b,b,1+\nu,s-\nu,2)$.
They do not exceed the corresponding $n_i$, since
$1+\nu\leq n_3$ and $s-\nu\leq n_3\leq n_4$.
The group $B$ already has its final size.

To enlarge $A$, use supports containing one or two indices from each
of $O,E$, and neither $s$ nor $\tau$.
Their coefficients avoid the only base variable in $A$.
For $b\geq3$, remove the path pairs and $\{1,2,3\}$.
The number of supports left satisfies
$$
 \frac{b^2(b^2-1)}4-2b+1\geq(b-1)^2\geq n_3-1-\nu.
$$
Here $n_3-1-\nu=(n_3-2b+1)_+\leq(b-1)^2$.
For $b=2$, the deficit is at most one. If it is one, use $\{1,3\}$.

To enlarge $C$, use two indices from $O$, with optional $\tau$, and
the support $\{1,\tau\}$.
There are $2\binom b2+1=b(b-1)+1$ such supports.
They contain no index from $E$ and do not contain $s$, so their
coefficients avoid $C$.
They are disjoint from the supports just used for $A$, except for
$\{1,3\}$ when $b=2$.
Discard that support if it was used for $A$.
At least $b^2-b$ choices are left, which is enough for the deficit $n_2-b$.

Finally enlarge $D,Q$.
Use supports containing at most one index from each of $O,E$, with
$s,\tau$ optional. Before exclusions there are $4b(b+1)$ choices.
Of the $4b+4$ reserved supports, only $\{1,2,3\}$
is outside this collection. Thus $4b+3$ must be excluded.
Also exclude $I_*$, at most one support used for $C$, and at most
$b(b-1)$ used for $A$.
At least
$$
 4b(b+1)-(4b+3)-1-1-b(b-1)=3b^2+b-5
$$
choices are left.
For $b\geq2$, this is at least $2b^2+3b-5$.
The combined deficit of $D,Q$ is at most
$$
 (b^2-(s-\nu))+(b^2-2)\leq2b^2-5.
$$
Thus there are enough distinct supports to fill both groups.
Their coefficients avoid $D,Q$ and satisfy the degree bounds above.
This also covers $b=2$.
\end{itemize}

Omit all unused optional coordinates.
Lemma~\ref{lem:gram-local} proves integrality and nonvanishing after
these additions and omissions outside characteristic two.
Proposition~\ref{prop:nonvanishing} gives the same conclusion in
characteristic two.
The coordinate groupings do not change the affine variety, so
$\delta_V\ne0$ modulo two for every profile above.
\end{proof}

\begin{corollary}\label{cor:quadratic-cases}
The following hold in every characteristic.
\begin{enumerate}
\item For $r\geq3$ and $2\leq c<N$, a very general complete intersection
in $(\PP^N_k)^r$ with one row $(2,\ldots,2)$ and $c-1$ rows
$(1,\ldots,1)$ is Fano and geometrically stably irrational.
\item For a positive-dimensional geometrically integral very general
intersection of $c$ ample hypersurfaces in $(\PP^N_k)^2$, one row
at least $(2,2)$ suffices for no geometric decomposition if at least one
of the following holds: $N=2$, $c\geq N$, $3c\geq2N+5$, or
$c\geq2$ and $2N-c\leq8-3\epsilon$, where
$\epsilon=\delta_{\cha k=2}$.
\item A quadratic row suffices if $r\geq5$ and $M\leq N_{\min}^2$.
\end{enumerate}
\end{corollary}
\begin{proof}
For \textup{(1)}, the very general complete intersection is smooth.
Adjunction gives $\omega_X=\cO_X(c-N,\ldots,c-N)$.
Since $c<N$, its inverse is ample, so $X$ is Fano.
Theorem~\ref{thm:quadratic} excludes a decomposition of the diagonal
and hence stable rationality.

For \textup{(2)}, adjunction treats $c\geq N$, and
Lemma~\ref{lem:planes} treats $N=2$.
Suppose $c<N$, $c\geq2$ and $m=2N-c\leq8-3\epsilon$.
Use $q_4$ and the small profiles of Lemma~\ref{lem:graphs}.
Only the characteristic-two fivefold needs three equations.
In that case $m=5$ forces $c$ to be odd, so $c\geq3$.
For $(N,c)=(3,2)$ in characteristic zero, the intersection with rows
$(2,2)$ and $(1,1)$ is the K\"uchle fourfold of
\cite[Section~4]{FiammengoLuders}.

Now suppose $c<N$ and $3c\geq2N+5$.
These inequalities imply $N\geq8$ and $m=2N-c\geq N+1$.
Write $m=2s+\varepsilon$ with $\varepsilon\in\{0,1\}$.
Then $s\geq4$, and $3c\geq2N+5$ gives $q=s+3\leq c$.
The uniform model of Lemma~\ref{lem:graphs} has $m+q\leq2N$
coordinates, divided into two groups whose sizes differ by at most one.
Both group sizes are at most $N$.
Match the other $c-q$ rows to distinct unused coordinates.
Their positive degrees allow the corresponding coordinate equations,
and monic elimination gives the affine model.
Corollary~\ref{cor:embedding} proves the conclusion.

For \textup{(3)}, adjunction treats $N_{\min}=1$.
Suppose $N_{\min}\geq2$.
We first combine the ambient factors into five groups.
List the factors in decreasing order of dimension and assign each
to a group whose current total dimension is smallest.
Let $b$ be the smallest final total, so $b\geq N_{\min}\geq2$.
Consider the last factor added to a group with at least two factors.
Before this addition, the group's total was at most $b$.
The added factor also has dimension at most $b$, since each group
already contained a factor at least as large.
Thus every such group has total at most $2b\leq b^2$.
A group with only one factor has total at most
$M\leq N_{\min}^2\leq b^2$.

If $c$ is at least the largest group total, every original factor
has dimension at most $c$, and adjunction gives the obstruction.
Otherwise let $S=L-c+1$.
If $S\geq5b$, choose five coordinate-group sizes summing to $S$,
each between $b$ and its available total.
This is possible because the five available totals sum to $L$.
The largest size is at most $b^2$, and the smallest is at least $b$.

If $S<5b$, the other four groups have total at least $b$ and $c$
is smaller than the largest total. Hence $S\geq4b+2$.
Choose five sizes differing by at most one and summing to $S$.
Each is at most $b$ and at least two.
Their largest is at most the square of their smallest.

In both cases Lemma~\ref{lem:five-factors} gives the affine model.
Split the five coordinate groups among the original factors.
This does not increase the degree in any factor.
There are $L-S=c-1$ unused coordinates.
Match them to the other ample rows and use
Corollary~\ref{cor:embedding}.
\end{proof}

\subsection{Cubic and quartic bounds}\label{subsec:degree-bounds}
We first obtain quartic equations from the Gram and Hermitian constructions.
We then count the coordinates in the binary and cubic norm constructions
with given degree bounds in the first $r$ groups.
For a hypersurface, $B$ denotes the number of coordinates in the last group.
When graph equations are added, we subtract their number from the size of
that group, as explained below. A norm coordinate is indexed by its
weight vector.

\begin{proposition}\label{prop:quartic-profiles}
Let $a_j\geq1$, $e_j\geq2$, and $a=\sum_ja_j\geq3$.
Let $P=\prod_jS_{e_j}(a_j)$ and $\epsilon=\delta_{\cha k=2}$.
There are geometrically integral affine hypersurfaces $V$ with
$\delta_V\ne0$ in $\CH_0(V_{\Omega(V)},\ZZ/2)$ and profile
$(a_1,\ldots,a_r,B;e_1,\ldots,e_r,4)$ whenever
$3a+1+\epsilon\leq B\leq P S_2(a+1)-(1+4\epsilon)a-1$.
\end{proposition}
\begin{proof}
In both cases, the first $r$ groups contain the base parameters indicated
below. Each optional form coordinate is placed in the last group.
We check the degree bounds before counting how many such coordinates
can be added.
\begin{itemize}
\item[step 1:] In this step, we treat characteristic different from two.
Use the charts $Y_0=1$ of $G_{2a}^{\sharp}$ and $G_{2a+1}^{\sharp}$.
Place the even indexed base parameters in groups of sizes $a_j$.
Place all other coordinates in the last group.
If $i$ is odd, both neighbors of $i$ in $L_i|_{Y_0=1}$ are even.
Thus $z_i(L_i|_{Y_0=1})^2$ has total degree at most two in the first
$r$ groups together and degree at most three in the last group.
If $i$ is even, these bounds are one and four.
The terms from the last vector, the scalar term and the perturbation
have degrees at most $(2,4)$ as well. Since $e_j\geq2$, these bounds
give the degree bound in each first group.

For $s=2a$, the last group has $a$ odd base parameters, $2a$ fiber
coordinates and $X$, hence $3a+1$ coordinates before any additions.
There are $PS_2(a)$ supports compatible with the degree bounds, and
$4a+2$ are already used by the construction.
For $s=2a+1$, the initial size is $3a+3$, and $4a+4$ of the
$PS_2(a+1)$ compatible supports are already used.
Add the coordinates $R_I$ of the other supports one at a time.
This gives every integer in the intervals
$$
 [3a+1,PS_2(a)-a-1],\qquad
 [3a+3,PS_2(a+1)-a-1].
$$
Every subset of the first $a$ indices of size at most two obeys
the bounds in each group. Hence $P\geq S_2(a)$.
The inequality $S_2(a)^2\geq4a+3$ for $a\geq3$ makes the first
upper bound at least the second lower bound. The intervals overlap.
Lemma~\ref{lem:gram-local} proves nonvanishing for the initial charts
and for every subset of the added coordinates.

\item[step 2:] In this step, we treat characteristic two.
Use $s=2a$ in Construction~\ref{constr:cycle-coordinates}, placing
odd $x_i$ in the first groups and the even bases, $Z,X,P_1,\ldots,P_s$
in the last group. Its initial size is $3a+2$, and the degrees are
$(2,4)$. The added coordinates $U_I$ allow at most two last base
indices, and $V_I$ allow at most one. In the first group containing
$u$, their support count is
$$
 S_{e_j}(a_j-1)+S_{e_j-1}(a_j-1)=S_{e_j}(a_j).
$$
The counts in the other first groups multiply independently.
For each of the two kinds of coordinates, $4a+2$ supports are already
used. The largest size of the last group is thus
$$
 (3a+2)+PS_2(a)+P(a+1)-2(4a+2)
 =PS_2(a+1)-5a-2.
$$
The form $q$ in Definition~\ref{def:hermitian-equations} permits one
more coordinate $R$ without increasing the degrees.
Keeping any subset of the optional coordinates gives every size
from $3a+2$ to $PS_2(a+1)-5a-1$.

Proposition~\ref{prop:nonvanishing} and Lemma~\ref{lem:cycle-charts}
give nonvanishing before these additions. The coordinate changes
are covered by Remark~\ref{rmk:quadric-changes}\textup{(1)}.
The added coordinates preserve conditions~\textup{(2)} and~\textup{(3)}
of Proposition~\ref{prop:charts}, and
Remark~\ref{rmk:quadric-changes}\textup{(2)} proves condition~\textup{(1)}
for every subset. Proposition~\ref{prop:charts} proves the conclusion.\qedhere
\end{itemize}
\end{proof}

\begin{definition}\label{def:paired-graphs}
Let $g\in\{0,1\}$ and $h=1+g$ in
Construction~\ref{constr:binary-data}.
For $g=0$, use its affine hypersurface without adding an equation.
For $g=1$, adjoin $V_1$ with equation $V_1-u_1u_2=0$.
Replace $u_1u_2$ by $V_1$ in the coefficients and fixed products,
and only then take the squares in the defining equation.
For $J\subseteq[h]$, let $\widetilde u_J$ denote the polynomial
obtained from $u_J$ in this way. It has degree at most one in the
last group. Eliminating $V_1$ when $g=1$ recovers the original
affine hypersurface, so the two affine varieties are isomorphic.
\end{definition}

\begin{definition}\label{def:norm-counting}
Let $a_1,\ldots,a_r$ be positive integers, and let $e_j\geq a_j+1$ be integers.
Let $a=\sum_j a_j$, $\boldsymbol e=(e_1,\ldots,e_r)$, and $p=\cha k$.
Let $q\geq0$ be an integer.
For the use of these numbers in Theorem~\ref{thm:uniform}, $e_j$
is the degree of the distinguished defining row in $\PP^{a_j}_k$.
The number $q$ counts the other rows of degree at least two in every
factor of $R$. The degree $d$ is the minimum of the distinguished
row's degrees in the factors of $R$, with $d=3$ when $R$ is the
empty product.

For $f\in\ZZ[z_1,\ldots,z_r]$, let
$[z_1^{\nu_1}\cdots z_r^{\nu_r}]f$ denote the coefficient of
$z_1^{\nu_1}\cdots z_r^{\nu_r}$ in $f$.

We first define the integers $T_3$ and $T_4$.
If $a=1$, let $(T_3,T_4)=(4,11)$.
If $a=2$, let $(T_3,T_4)=(C+1,3C+1)$, where
$C=\sum_{\nu_1=0}^{e_1}\cdots\sum_{\nu_r=0}^{e_r}[z_1^{\nu_1}\cdots z_r^{\nu_r}]\prod_j(1+z_j+z_j^2)^{a_j}$.
If $a\geq3$, let
$$
 D=\max
       \sum_{\nu_1=0}^{e_1}\cdots\sum_{\nu_r=0}^{e_r}
       [z_1^{\nu_1}\cdots z_r^{\nu_r}]
       \prod_j(1+z_j+z_j^2)^{t_j}(1+z_j)^{a_j-t_j},
$$
where the maximum is over integers $t_1,\ldots,t_r$ such that
$0\leq t_j\leq a_j$ and $\sum_jt_j=3$.
In this case, let $(T_3,T_4)=(D,3D+1)$.

For each integer $d\geq3$, we define the bound $B_d$ used in
Theorem~\ref{thm:uniform} as follows.
Let $g_3=\min(q,1)$ and
$B_3=\max\bigl\{
       \delta_{p\ne2}(2^{a+1+g_3}+g_3),\,
       \delta_{p\ne3}T_3\bigr\}$.
Let $B_d=\infty$ for $d\geq4$.
\end{definition}

To compare the bound with the dimension in Theorem~\ref{thm:uniform},
let $c$ be the codimension of $X$ and $N=\sum_\ell N_\ell$.
Let $B=N-c+1$. Then $\dim X=a+B-1$, so the dimension condition is
$B\leq B_d$. Thus $B$ is determined by $X$ and its ambient dimensions.
The upper bound $B_d$ is computed from $(a_j)$, $(e_j)$, $p$, $q$, and $d$.
When $R$ is empty, $N=0$ and $B=1-c\leq0$, and the proof of the
theorem uses adjunction.

Proposition~\ref{prop:binary-bound} gives the first term of $B_3$
when $p\ne2$. Lemma~\ref{lem:exact-norm} gives the second term when
$p\ne3$. For $p\ne3$, it also gives a degree-four hypersurface with
$T_4$ last coordinates.
The integer $T_4$ is a finite bound for that norm construction.
It is used in the characteristic-two proof of
Proposition~\ref{prop:quartic-all-capacities}, which allows every
$B\geq4$ and gives the final bound $B_d=\infty$ for $d\geq4$.

Relative ampleness in Theorem~\ref{thm:uniform} requires a positive
degree in every factor of $R$. The degrees of the other rows in the factors
of $Q$ may be zero. The case $q=0$ includes hypersurfaces.
Every choice of distinguished row and splitting into $Q$ and $R$ may be used.

For the constructions below, choose a partition
$[a]=A_1\sqcup\cdots\sqcup A_r$ with $|A_j|=a_j$.
The base coordinates $x_i$ with $i\in A_j$ belong to the $j$th group.
The coefficient sums in Definition~\ref{def:norm-counting} depend on
the sizes $a_j$ and the bounds $e_j$. The sets $A_j$ are used below to
check the degree in each group.

Suppose an affine model is the complete intersection of one
distinguished equation and $g$ graph equations, with $b$ coordinates
in the last group. Let $B=b-g$. Its dimension is $a+B-1$, so it uses
$B+g$ last coordinates. The number of graph equations is denoted by $g$.
To place the model in Theorem~\ref{thm:uniform}, these equations must
fit distinct rows among the $q$ rows counted in
Definition~\ref{def:norm-counting}. Thus we need $g\leq q$.

\begin{proposition}\label{prop:binary-bound}
Suppose $p\ne2$. Let $g\in\{0,1\}$ and $h=1+g$.
For every $h+1\leq B\leq2^{a+1+g}+g$ there is an integral affine
complete intersection $V$ with $\Tortwo(V)=2$, first group sizes
$a_j$, and $B+g$ last coordinates. One equation has degrees at most
$(e_1,\ldots,e_r;3)$, and the other $g$ equations have degrees
$(0,\ldots,0;2)$. Allowing at most $g$ additional equations gives
every $2\leq B\leq2^{a+1+g}+g$ with coefficients $\ZZ/2$.
\end{proposition}
\begin{proof}
Start with the equation $G=0$ in Construction~\ref{constr:binary-data}.
It is primitive in the base parameters. Indeed, the coefficient of
$W^2$ is $C=x_{[a]}u_{[h]}$, and the coefficient independent of the
free form coordinates has constant term one in the base parameters.
No base parameter dividing $C$ divides all the coefficients.
Lemma~\ref{lem:binary-detection} gives an integral hypersurface
with nonzero generic-point class.

Make the substitution of Definition~\ref{def:paired-graphs}.
For $g=1$, eliminating $V_1$ recovers $G=0$.
The graph equation and the hypersurface equation form a regular sequence.
Indeed, the quotient by the graph equation is a polynomial ring and $G$
is nonzero in that ring. The affine variety and its generic-point class
are preserved by elimination.

If $S=\varnothing$ in Construction~\ref{constr:binary-data}, one
of the $2^{a+h}$ norm coordinates is replaced by $P_0$.
If $S\ne\varnothing$, two norm coordinates are expressed in terms
of $P_0$ and one new coordinate $Y$. Both cases leave
$2^{a+h}-1$ free norm coordinates.
Together with the $h$ last base parameters and $g$ graph coordinates,
they give $2^{a+h}+h+g-1$ last coordinates. Subtracting $g$ gives
the largest value $B=2^{a+1+g}+g$.

After substitution, every $\widetilde u_J$ and the fixed product
$\widetilde P_0$ have last degree at most one. Every term consisting
of a coefficient times a squared free coordinate then has last
degree at most three. The same holds for the fixed and coupled
terms. In the $j$th first group, the coupled terms obey
$$
 2k_j+|S\cap A_j|=a_j+1\leq e_j,\qquad
 |S\cap A_j|\leq k_j.
$$
These inequalities also bound the single fixed square when
$S=\varnothing$. All other norm terms have first degree at most
$a_j$. The right side has first degree $a_j$ and last degree at
most two. The graph equation has degrees $(0,\ldots,0;2)$.

Keep the base parameters, the fixed polynomials, the graph
coordinates and $W$. Set any subset of the other free coordinates
equal to zero. Lemma~\ref{lem:binary-detection} permits all these
omissions. The smallest model has $h+g+1$ last coordinates,
so $B=h+1$. Omitting one coordinate at a time gives every integer
up to $2^{a+1+g}+g$.
For $g=1$, the interval starts at three, and the construction with
$g=0$ also gives $B=2$. The intervals overlap, proving the final
claim about at most $g$ graph equations.
\end{proof}

For the cubic norms, let $|w|_j=\sum_{i\in A_j}w_i$ for each weight $w$.
This is the degree of the pure-power coefficient in the $j$th group.
For a polynomial $f$, its \emph{Newton polytope} is
$\operatorname{Newt}(f)=\operatorname{conv}(\operatorname{supp}f)$.
A linear function on a convex polytope has a minimum at a vertex.
Thus the valuation bounds from the fixed monomials can be checked on
their Newton polytopes.
A \emph{scalar coordinate} is the coefficient of the identity in a
norm algebra. In a cubic cyclic basis, the notation
$X(i)=X_0+X_1i+X_2i^2$, and likewise $Y(i),Z(i)$, uses three independent
coefficient coordinates. 

\begin{lemma}\label{lem:exact-norm}
Assume $p\ne3$.
Integral affine hypersurfaces $V$ with $\Torthree(V)=3$ and degrees at most
$(e_1,\ldots,e_r;3)$ and $(e_1,\ldots,e_r;4)$ reach last capacities
$T_3$ and $T_4$, respectively. Omissions give every capacity
$2\leq B\leq T_3$ in the cubic case and $3\leq B\leq T_4$ in the quartic case.
The lower capacities are given by adjunction with the
coefficients of Lemma~\ref{lem:forms}, not necessarily $\ZZ/3$.
\end{lemma}
\begin{proof}
We construct equations to which Theorem~\ref{thm:norm-criterion}
can be applied. The norms are those of Remark~\ref{rmk:specified-forms}.
For each equation, we count its free coordinates and exhibit the
interior point needed in that theorem. Lemmas~\ref{lem:norm-values}
and~\ref{lem:laurent-values} give the valuation bounds from the fixed
monomials.

\begin{itemize}
\item[step 1:] In this step, we treat $a\geq3$.
We first calculate the number of norm coordinates and then construct
the cubic and quartic equations.
For $e_j=a_j+1$ in every group, the count in
Definition~\ref{def:norm-counting} is
$$
 D=27\cdot2^{a-3}-\begin{cases}
 a+1,&r=1,\\
 3\cdot2^{\min(a_1,a_2)-1},&r=2,\\
 0,&r\geq3.
 \end{cases}
$$
A group of size $b$ containing zero, one, two or three Tits
parameters contributes, respectively,
$2^b$, $3\cdot2^{b-1}$, $9\cdot2^{b-2}-1$, or
$27\cdot2^{b-3}-b-1$ allowed weights, whenever it has enough
parameters for that case. These are the coefficient sums of
$(1+z+z^2)^t(1+z)^{b-t}$ through degree $b+1$.

If $r\geq3$, place the three Tits parameters in different groups.
No weight is then excluded. If $r=2$, placing two parameters in a
larger group and one in the other group excludes
$3\cdot2^{\min(a_1,a_2)-1}$ weights. Placing all three parameters
in a group of size $b$ excludes $(b+1)2^{a-b}$ weights, which is
no smaller. If $r=1$, use the contribution with three Tits parameters.
This proves the formula for $D$.

In particular, $D\geq2^{a+1}$: in the one block case use
$11\cdot2^{a-3}\geq a+1$, and in the two block case use
$3\cdot2^{\min(a_1,a_2)-1}\leq11\cdot2^{a-3}$.
Larger degrees cannot decrease $D$.

Choose integers $(t_j)$ for which the maximum defining $D$ holds.
Place the three Tits parameters in the groups with these multiplicities.
Suspend by each of the other base parameters $x_i$, using $x_i^2$
if $x_i$ is the only parameter in its group.
Let $\Psi$ be the coordinate subform consisting of the weights $w$
with $|w|_j\leq e_j$ for every $j$, and let $\cW$ be its weight set.
Its number of coordinates is $D$.
For a group of size one containing no Tits parameter, the suspension
has weights $0,2$ in place of $0,1$. Both weights are allowed because
$e_j\geq2$, so this substitution does not change the count.

We next choose the coordinate that will remain free and the one
that will be replaced by a fixed polynomial.
Let $\eta$ have entry two in each group of size one.
In a group of size at least two, let its entries be one, except that
one Tits entry is increased to two if that group contains a Tits parameter.
These choices give a weight of $\Psi$.

Let $\rho_j\in\{0,1,2\}$ be the residue of $a_j+1$ modulo three.
For $a_j\geq2$, choose $\omega_j$ with zero or one entries and with
exactly $\rho_j$ entries equal to one. If $\eta$ has an entry two in
that group, take the corresponding entry of $\omega_j$ to be zero.
For $a_j=1$, let $\omega_j=0$.
The vectors $\omega_j$ together give a weight $\omega$ of $\Psi$
different from $\eta$.

Let $\bx_j=(x_i)_{i\in A_j}$. For $w\in\cW$, let
$$
 \begin{aligned}
 P_0&=\prod_{a_j>1}
         \sum_{\substack{\nu\in\ZZ_{\geq0}^{a_j}\\
                        3|\nu|\leq a_j+1-\rho_j}}
                          \bx_j^\nu,\\
 \widehat P_w&=\prod_j
         \sum_{\substack{\nu\in\ZZ_{\geq0}^{a_j}\\
                        3|\nu|\leq e_j-|w|_j}}
                          \bx_j^\nu.
 \end{aligned}
$$
Let $H_0=\omega+3\operatorname{Newt}(P_0)$,
$H_1=\operatorname{conv}\bigcup_{w\in\cW}(w+3\operatorname{supp}\widehat P_w)$,
and $\tau_0=(3\one-\eta)/2$.
The polynomials $P_0$ and $\widehat P_w$ have fixed coefficients.
They introduce no free norm coordinates.
For a group of size $b\geq2$, let
$\Delta_b=\{v\in\RR_{\geq0}^b:\sum_i v_i\leq1\}$.
The corresponding factor of $H_0$ is
$\omega_j+(b+1-\rho_j)\Delta_b$.
The projection of $H_1$ to this group contains $b\Delta_b$.
Indeed, on each coordinate axis, a weight entry zero or one together
with a monomial in $\widehat P_w$ gives an exponent at least $b$
and at most $b+1$. The zero vector also occurs.
The choices in different groups are independent, so $H_1$ contains
the product of these simplexes.

For a group of size one, $H_0$ has coordinate zero and $H_1$
contains the interval $[0,2]$ in that coordinate.
Let $1/4<\alpha<1$. In each group of size $b\geq2$,
the polytope $(1-\alpha)H_0+\alpha H_1$ contains
$$
 (1-\alpha)\omega_j+
       \bigl(b+(1-\alpha)(1-\rho_j)\bigr)\Delta_b.
$$
Subtract $(1-\alpha)\omega_j$ from the $j$th group of $\tau_0$.
Every coordinate is positive. Its coordinate sum is at most
$b-(1-\alpha)\rho_j$, which is smaller than the radius
$b+(1-\alpha)(1-\rho_j)$ by at least $1-\alpha$.
For a group of size one, the coordinate of $\tau_0$ is $1/2$,
and $0<1/2<2\alpha$.
The product of the group interiors then gives
$$
 \tau_0\in\operatorname{int}\bigl((1-\alpha)H_0+\alpha H_1\bigr).
$$
Increasing the degrees only enlarges $H_1$.

For the cubic equation, adjoin a last base parameter $u$.
Let $Y_\omega=P_0$ and $W=Y_\eta$.
Leave the other coordinates of $\boldsymbol Y$ free and use
$$
 \Psi(\boldsymbol Y)-u^2\Psi((\widehat P_w)_w)
                    =\lambda x_{[a]}uW.
$$
The left side is a binary suspension by $u^2$.
Its fixed monomials give $H_0\times\{0\}$ and
$H_1\times\{2\}$.
At last coordinate $3/2$, their convex hull has slice
$\tfrac14H_0+\tfrac34H_1$.
The interior calculation above with $\alpha=3/4$ shows that
$(\tau_0,3/2)$ is an interior point of the full convex hull.
The degree bounds on pure and mixed terms give degrees at most
$(\boldsymbol e,3)$.
There are $D-1$ free norm coordinates and one last base parameter,
so the last group has size $D=T_3$.

For the quartic equation, adjoin last base parameters $u,v$.
Let $Y^{(0)}_\omega=P_0(1+u+v)$ and $W=Y^{(1)}_\eta$.
Leave the other coordinates of the three tuples
$\boldsymbol Y^{(0)},\boldsymbol Y^{(1)},\boldsymbol Y^{(2)}$ free and use
$$
 \Psi(\boldsymbol Y^{(0)})-u\Psi(\boldsymbol Y^{(1)})
 -v\Psi(\boldsymbol Y^{(2)})+uv\Psi((\widehat P_w)_w)
                      =\lambda x_{[a]}uvW.
$$
The fixed bodies are $H_0\times3\Delta_2$ and
$H_1\times\{(1,1)\}$, while the interior point is $(\tau_0,1,3/2)$. Use $\alpha=1/3$ and
$$
 (1,3/2)=\tfrac23(1,7/4)+\tfrac13(1,1),\qquad
                         (1,7/4)\in\operatorname{int}(3\Delta_2).
$$
Together with the first-group interior calculation for $\alpha=1/3$,
this gives an interior point of the full convex hull of fixed weights.
The three tuples have $3D-1$ free norm coordinates.
Adding $u,v$ gives last-group size $3D+1=T_4$.
The last degree is at most four, including the mixed terms.

For both equations, Remark~\ref{rmk:degree-control} bounds the mixed
terms by the same group degrees as the pure terms.
The variables other than the bases are free norm coordinates, so the
integrality condition in Theorem~\ref{thm:norm-criterion} follows from
integrality of those coordinates. Lemmas~\ref{lem:norm-values}
and~\ref{lem:laurent-values}, together with the two interior calculations,
verify the other conditions of that theorem.

\item[step 2:] In this step, we treat $a=2$.
There is either one group of size two or two groups of size one.
The cyclic norm on $x_1,x_2$ has $C$ coordinates allowed by the
group-degree bounds.
For degree three, adjoin a last base parameter $u$ and use the first Tits equation
$$
 N_u(E,P,P)=\lambda x_1x_2uW,
$$
where all allowed coordinates of $E$ are free and $W$ is its
coordinate of weight $(1,1)$.
For one group of size two, let the scalar coordinate of $P$ be
$1+x_1+x_2$ and its other allowed coordinates be one.
For two groups of size one, let all allowed coordinates of $P$ be one.
The convex hull of fixed weights contains
$3\Delta_2\times[1,2]$ in the first case and
$[0,2]^2\times[1,2]$ in the second.
Both have $(1,1,3/2)$ as an interior point.
The $C$ free coordinates of $E$ and the parameter $u$ give
$T_3=C+1$.

For degree four, use the quartic equation of step~1 with the cyclic
norm in place of $\Psi$.
For one group of size two, let $\omega=0$, $P_0=1+x_1+x_2$,
and $\eta=(2,1)$.
Then $H_0=3\Delta_2$ and $H_1$ contains $3\Delta_2$, so
$\tau_0=(1/2,1)$ is an interior point of their convex combinations.
For two groups of size one, let $\omega=0$, $P_0=1$, and $\eta=(2,2)$.
Here $H_0=\{0\}$ and $H_1$ contains $[0,2]^2$.
The point $\tau_0=(1/2,1/2)$ is in the interior of
$\tfrac23H_0+\tfrac13H_1$.
The same last-coordinate calculation as in step~1 gives the full
interior point. Three copies of the cyclic norm have $3C-1$ free
coordinates after the one fixed substitution. The two last base
parameters give $T_4=3C+1$.

\item[step 3:] In this step, we treat $a=1$.
Let $x_1$ be the first-group parameter, and let $u$ be a last-group parameter.
Take cyclic generators satisfying $i^3=x_1$, $j^3=u$ and $ij=\zeta ji$.
For degree three, use
$$
 \operatorname{Nrd}\bigl(X_0+X_1i+X_2i^2+(1+i^2)(j+j^2)\bigr)
                            =\lambda x_1uX_1.
$$
The fixed coordinates have weights $(0,1)$, $(2,1)$, $(0,2)$ and $(2,2)$.
Their convex hull contains $(1,3/2)$ in its interior.
The coordinates $X_0,X_1,X_2$ and the parameter $u$ give $T_3=4$.

For quartics, let $u,v$ belong to the last group and use
$$
 N_v\bigl(X(i)+Y(i)j+j^2,Z(i),1+i^2+i^2j^2\bigr)=\lambda x_1uvY_1.
$$
The fixed weights $(0,0,2)$, $(2,0,2)$, $(0,2,0)$ and $(2,2,2)$
are affinely independent. Their barycenter is $(1,1,3/2)$,
so this point lies in the interior of their convex hull.
The three tuples $X(i),Y(i),Z(i)$ have nine free coordinates.
Together with $u,v$, they give $T_4=11$.
Each pure term has degree at most two in $x_1$ and last degree at most four.
Remark~\ref{rmk:degree-control} gives the same bounds for the mixed terms.

\end{itemize}

In each case, keep the fixed polynomials, the last base parameters
and $W$. Set the other free norm coordinates equal to zero one at a time.
The fixed convex hull and its interior point do not change.
Theorem~\ref{thm:norm-criterion} continues to give a nonzero class
for every such omission. This proves the full intervals in the lemma.
Adjunction and Lemma~\ref{lem:forms} give the smaller last-group sizes,
with the coefficients of that lemma.
\end{proof}

\subsection{Quartic bounds in all dimensions}\label{subsec:quartic-products}

We now prove that quartic rows have no upper dimension bound.
We first construct quartic hypersurfaces in large affine spaces while
keeping degree at most two in certain groups of one or two coordinates.
We then use the finite norm constructions to cover the smaller dimensions.

\setcounter{theorem}{24}
\begin{lemma}\label{lem:quartic-singleton-profiles}
Fix $s$ disjoint groups of one coordinate and $t$ disjoint groups of
two coordinates. Let $b_0=\max\{2,t+\floor{s/2}\}$,
$u_0=\max\{0,s-2(b_0-t)\}$ and $D=2^{u_0}3^{b_0}$.
Then $u_0\in\{0,1\}$. Let $\epsilon=\delta_{\cha k=2}$.
For every $N\ge2^\epsilon D-1$ there is, over an algebraically
closed extension $\Omega/k$, a geometrically integral affine
hypersurface $V\subset\bA^N_\Omega$ such that
$$
 0\ne\delta_V\in\CH_0(V_{\Omega(V)},\ZZ/2),\qquad \deg V\le4,
$$
and its equation has degree at most two in each of these groups.
\end{lemma}
\begin{proof}
We may replace $k$ by its algebraic closure. Let $T$ be transcendental
over $k$ and let $\ell=\overline{k(T)}$.
For $u\ge0$ and $b\ge2$, use the ordinary Segre embedding of
$(\PP^1)^u\times(\PP^2)^b$. The affine parameters are $u_i$ on the
line factors and $(x_j,y_j)$ on the plane factors.
Let $C=\prod_i u_i$, $P=\prod_jx_j$, and $Q=\prod_jy_j$.
Let $A=x_1\prod_{j>1}y_j$ and $B=y_1\prod_{j>1}x_j$.
Write $W_q$ for the coordinate indexed by a Segre monomial $q$.
Let $L=\sum_qW_q$ and
$G=\prod_i(1+u_i)\prod_j(1+x_j+y_j)$.
On the Segre variety, $L$ restricts to $G$.

We work on the affine chart $W_1=1$.
For the coordinate placement, take $b\ge t$.
Use one plane for each of the $t$ two-coordinate groups, with coordinates
$W_{x_j},W_{y_j}$. The other $b-t$ planes and the $u$ lines provide
$2(b-t)+u$ singleton coordinates for the $s$ one-coordinate groups.
All additional ambient coordinates will lie outside these groups.

\begin{itemize}
\item[step 1:] In this step, we construct the affine reference and
coefficient lift in characteristic different from two.
Use the reference quadric of
\cite[Propositions~3.1 and~3.3]{BohningBothmerPirutkaZhang2026}.
Let $m=u+2b$ and $r_0=2^u3^b-1-m$.
Let $f_1=G^2P(A+B)$ and $f_2=CP^2Q$.
For each $r_0\le r\le2^m-2$, the affine reference hypersurface is
$$
 V_0=\{\sigma_r(U)+f_1-T^2f_2=0\}.
$$
Here $\sigma_r$ has $r$ coordinates, so $V_0$ has $m+r$ ambient
coordinates. One coefficient of $\sigma_r$ is $1$, since the deleted
monomials all have positive degree. Thus the equation is primitive over
the polynomial ring in the base parameters. Its generic quadric is
geometrically integral, and so is $V_0$.
The cited proposition gives a nonzero coordinate symbol $\alpha$ over
$\ell(V_0)$.

We verify the affine obstruction before constructing the degeneration.
Let $(F,\nu)/\ell$ be a henselian discretely valued extension with
$\nu|_{\ell^\times}=0$. Consider an $F$-point with nonzero base
coordinates at which $\alpha$ is nonzero. The corresponding Pfister form
is anisotropic. The calculation in the proof of
\cite[Lemma~3.6]{BohningBothmerPirutkaZhang2026}, applied to the coordinates
of this point, gives
$$
 \nu(f_2/f_1)\ge\sum_i|\nu(u_i)|+
 \frac12\sum_j\bigl(|\nu(x_j)|+|\nu(y_j)|\bigr).
$$
Every factor in this quotient is nonzero by anisotropy of the
coordinate subforms. If some base coordinate is not a unit, the right
side is positive. Hensel's lemma then makes $1-T^2f_2/f_1$ a square.
It follows that $\sigma_r\perp\langle f_1-T^2f_2\rangle$ is isometric
to a subform of the anisotropic Pfister form. Its point $(U,1)$ gives
a contradiction. Hence all base coordinates are units.

The reduction of the Pfister form is also anisotropic. Indeed, an
isotropic residual vector would lift by Hensel's lemma, because two
and every diagonal coefficient are units. Thus
$\nu(\sigma_r(U))=2\min_i\nu(U_i)$ for $U\ne0$.
The equation of $V_0$ now makes all $U_i$ integral. This is immediate
when $U=0$. Corollary~\ref{cor:affine-detection} gives
$\delta_{V_0}\ne0$ with coefficients $\ZZ/2$.

Use the weighted normal coordinates of
\cite[Section~4]{BohningBothmerPirutkaZhang2026}. A quartic lift of
$f_1-T^2f_2$ is
$$
 L^2W_P(W_A+W_B)-T^2W_1W_{CP}W_PW_Q.
$$
Since $b\ge2$, each of $P,Q,A,B,CP$ has support of size at least
two. Among the factors of this lift, only $L^2$ involves the designated
singleton coordinates. Its degree in each designated group is at most
two. After the weighted coordinate change, the special equation is
that of $V_0$. As $r$ varies, the ambient coordinate count $N=m+r$
runs through every integer from $2^u3^b-1$ to
$2^{u+2b}+u+2b-2$.

\item[step 2:] In this step, we construct the coefficient lifts
in characteristic two.
Use Proposition~\ref{prop:char2-segre-family} with one additional
line factor. Write $\mathfrak a$ for its Artin--Schreier parameter,
so that the ordinary parameters are still $u_i,x_j,y_j$.
Let $L_{\mathfrak a}=\sum_qW_{\mathfrak aq}$, where $q$ runs over
the ordinary Segre monomials. On the Segre variety this polynomial
restricts to $\mathfrak aG$.
The coefficients $f_1=\mathfrak a^3(\mathfrak a+G^2+G)P(A+B)$ and
$f_2=CP^2Q$ have the quartic lifts
$$
\begin{aligned}
 F_1={}&L_{\mathfrak a}^2W_{\mathfrak aP}(W_A+W_B)
 +W_{\mathfrak aP}(W_{\mathfrak aA}+W_{\mathfrak aB})W_{\mathfrak a}^2\\
 &+L_{\mathfrak a}W_{\mathfrak aP}
                  (W_{\mathfrak aA}+W_{\mathfrak aB})W_1,\\
 F_2={}&W_{CP}W_PW_QW_1.
\end{aligned}
$$
The additional Artin--Schreier coordinate lies outside the designated
groups. Neither lift contains a designated ordinary singleton coordinate.
\end{itemize}

We next check the degree bounds for the normal terms in both
characteristics. A flattening minor has degree at most one in each
designated group. For a two-coordinate group, its coordinates occupy
one row when the corresponding plane is the flattening factor, and
one column otherwise. A determinant term uses at most one entry from
that row or column.

A product of two minors has degree at most two in each group.
This also bounds the characteristic-two cross terms between two forced
normal directions. A cross term between a forced direction and an extra
direction is a minor times the extra coordinate and $W_1$.
For two extra directions, the diagonal and cross coefficients have
quadratic Segre lifts. The extra coordinates lie outside the designated
groups, so these terms satisfy the same bounds.
In characteristic two, the additional Hermitian line has coefficient
$\gamma=Q(A+B)$. Its two coefficients have lifts $W_Q(W_A+W_B)$ and
$W_{\mathfrak aQ}(W_A+W_B)$, which also have group degree at most two.

We now choose the generic fiber to be smooth and geometrically integral.
In both constructions use the quartic perturbation
$$
 J=W_1^3\left(\lambda_0W_1+\sum_i\lambda_iW_i\right),
$$
where the sum runs over all other ambient coordinates.
On $W_1=1$, this is a general affine-linear polynomial, so it preserves
all group-degree bounds. Over $\ell((\tau))$, the linear system spanned
by the unperturbed equation, $1$, and the affine coordinates defines
a graph immersion.
The Bertini theorems \cite[Lemma~33.47.3 (\href{https://stacks.math.columbia.edu/tag/0FD6}{Tag~0FD6})]{stacks-project}
and \cite[Lemma~37.32.3 (\href{https://stacks.math.columbia.edu/tag/0G4F}{Tag~0G4F})]{stacks-project} give
a smooth geometrically integral generic member. For the
codimension-two hypothesis of the latter lemma, use the sections $1,w_1,w_2$.
The intersection $\{w_1=w_2=0\}$ has codimension two, and the base
locus is empty. The open coefficient chart with coefficient one
on the unperturbed equation is precisely this perturbation family.

Let the $\lambda_i$ be independent over $\ell$ and
$R=\ell(\boldsymbol\lambda)[[\tau]]$. Comparing Laurent coefficients
gives an injection
$\ell((\tau))(\boldsymbol\lambda)\longrightarrow\operatorname{Frac}R$,
so geometric integrality persists over this generic field.
Multiply the perturbation by a power of $\tau$ larger than the
weighted degree divided out in the normal rescaling. Its reduction
modulo $\tau$ is then zero, so the special equation is the affine
reference equation. This equation is nonzero. Hence $\tau$ is a
nonzerodivisor in the hypersurface coordinate ring, and the family
is flat over $R$. By \cite[Lemma~3.6(e)]{LS}, the geometric special
fiber has nonzero generic-point class. Applying \cite[Lemma~3.8]{LS}
with empty boundary gives nonvanishing on the geometric generic
fiber. Take $\Omega=\overline{\operatorname{Frac}R}$.

It remains to check which ambient dimensions occur while the
designated groups fit. Use $t$ planes for the two-coordinate groups.
The inequality $2(b_0-t)+u_0\ge s$ places the one-coordinate groups
in the other planes and lines. The ordinary Segre variety then has
$D=2^{u_0}3^{b_0}$ homogeneous coordinates. In characteristic two,
the extra Artin--Schreier line doubles this number.

For $\cha k\ne2$, the choices $u=0$ and $u=1$ give the intervals
$$
 [3^b-1,4^b+2b-2],\qquad
 [2\cdot3^b-1,2\cdot4^b+2b-1].
$$
The first interval joins the second because $4^b+2b\ge2\cdot3^b$.
The second joins the first interval for $b+1$ because
$2\cdot4^b+2b+1\ge3^{b+1}$.
For $\cha k=2$, Proposition~\ref{prop:char2-segre-family} gives
$$
 [2\cdot3^b-1,2\cdot4^b+2b-1],\qquad
 [4\cdot3^b-1,4\cdot4^b+2b].
$$
The corresponding inequalities are
$2\cdot4^b+2b+1\ge4\cdot3^b$ and
$4\cdot4^b+2b+2\ge6\cdot3^b$.
All four inequalities hold at $b=2$ and follow by induction for $b\ge2$.

Consider the sequence of pairs
$(0,2),(1,2),(0,3),(1,3),\ldots$, beginning at $(u_0,b_0)$.
The intervals for consecutive pairs overlap or are adjacent, and
their upper limits tend to infinity. Their union contains every
$N\ge2^\epsilon D-1$.
The number $u+2(b-t)$ of singleton coordinates available after the
$t$ planes are used increases by one at every pair in this sequence.
Thus all designated groups continue to fit. Place the other coordinates
outside these groups. This proves the lemma.
\end{proof}

\setcounter{theorem}{26}
\begin{proposition}\label{prop:quartic-all-capacities}
Let $a_j\ge1$, $e_j\ge a_j+1$, and $a=\sum_j a_j$.
For every $B\ge4$, in every characteristic, there are an algebraically
closed extension $\Omega/k$, a prime $\ell\in\{2,3\}$, and a
geometrically integral affine hypersurface $V$ with group sizes
$(a_1,\ldots,a_r;B)$ and degrees at most $(e_1,\ldots,e_r;4)$ such
that
$$
 0\ne\delta_V\in\CH_0(V_{\Omega(V)},\ZZ/\ell).
$$
\end{proposition}
\begin{proof}
Only groups with $a_j=1,e_j<4$ or $a_j=2,e_j=3$ need degree less
than four. Let $s$ and $t$ be their respective numbers, and use
$u_0,b_0,D$ from Lemma~\ref{lem:quartic-singleton-profiles}.
The $s$ groups have one coordinate each, and the $t$ groups have two.
Every other group permits degree four.

We have $D\le\max\{9,2^a\}$.
If $s+2t\le3$, then $b_0=2$, $u_0=0$, and $D=9$.
If $s+2t\ge4$, then $2b_0+u_0=s+2t\le a$, so
$D=2^{u_0}3^{b_0}\le2^{s+2t}\le2^a$.
The affine ambient space has $a+B$ coordinates. Thus the lemma gives
every $B\ge2^{\delta_{\cha k=2}}D-1-a$.
We prove that the finite constructions join this range.

\begin{itemize}
\item[step 1:] In this step, we fill the smaller capacities in
characteristic different from two.
Use Construction~\ref{constr:binary-data} with $h=2$, before adding
the graph equation. Write its last bases as $u,v$, so
$P_0=Q(1+u)(1+v)$. Lemma~\ref{lem:binary-detection} gives a nonzero
affine class for this equation and every permitted omission.
If $S\ne\varnothing$, its coupled terms satisfy
$$
 (P_0-x_SY)^2+x_S(P_0+Y)^2
 =(1+x_S)P_0^2+(x_S^2+x_S)Y^2.
$$
Their degree in the $j$th first group is at most $a_j+1$, since
$2k_j+|S\cap A_j|=a_j+1$ and $2|S\cap A_j|\le2\le a_j+1$.
Their last degree is at most four. If $S=\varnothing$, the only
fixed term is $P_0^2$, with the same bounds. Every other norm term
has first degree at most $a_j$ and last degree at most four.
The right side $\lambda x_{[a]}uvW$ has last degree three.

The single or coupled assignment leaves $2^{a+2}-1$ free norm
coordinates. Together with $u,v$, these give $2^{a+2}+1$ last
coordinates. Keep $u,v,W$ and omit the other free coordinates one at
a time. This gives every $3\le B\le2^{a+2}+1$.
The large-dimensional construction starts at most at
$\max\{9,2^a\}-1-a$.
For $a=1$, this number is seven and the finite range reaches nine.
For $a\ge2$, it is at most $2^{a+2}+1$. Hence the two ranges overlap.

\item[step 2:] In this step, we fill the smaller capacities in
characteristic two.
Lemma~\ref{lem:exact-norm} gives every $3\le B\le T_4$ with
coefficients $\ZZ/3$.
For $a\ge3$, the coefficient sum in Definition~\ref{def:norm-counting}
includes all $2^a$ binary weights. Hence $T_4\ge3\cdot2^a+1$.
For $a\ge4$, the large-dimensional range starts at most at
$2^{a+1}-1-a$, which is smaller than this finite bound.
For $a=3$, it starts at most at fourteen, while $T_4\ge25$.
For $a=2$, the cyclic count is at least eight for a plane and nine
for two lines. Thus $T_4\ge25$, and the large-dimensional range
starts at most at fifteen. These ranges overlap in every case $a\ge2$.

For $a=1$, the two ranges give $3\le B\le11$ and $B\ge16$.
We fill the integers between them using the one-plane construction of
Proposition~\ref{prop:char2-segre-family}\textup{(2)} with one ordinary
line parameter $u$. It has $11\le N\le16$ ambient coordinates.
Assign the coordinate $W_u$ to the first group, and all other
coordinates to the last group.
Let $G=(1+u)(1+x+y)$, $L_{\mathfrak a}=\sum_qW_{\mathfrak aq}$,
and $D_0=W_1+W_y$.
Its coefficients $f_1=\mathfrak a^2(\mathfrak a+G^2+G)ux(1+y)$
and $f_2=u^2x^2y$ have quartic lifts
$$
\begin{aligned}
 F_1={}&L_{\mathfrak a}^2W_{ux}D_0
       +W_{\mathfrak aux}W_{\mathfrak a}^2D_0
       +L_{\mathfrak a}W_{\mathfrak aux}D_0W_1,\\
 F_2={}&W_{ux}^2W_yW_1.
\end{aligned}
$$
Neither lift involves $W_u$.
The minor calculation and the affine-linear perturbation in
Lemma~\ref{lem:quartic-singleton-profiles} give degree at most two
in $W_u$ for all other terms and preserve the affine obstruction.
There are $B=N-1$ last coordinates, so this gives every
$10\le B\le15$.
Together with the finite norm range and the large-dimensional range,
these intervals cover every $B\ge4$.
\end{itemize}
\end{proof}

\begin{theorem}\label{thm:quartic-small-first-degrees}
Use the notation $Q$ and $R$ of Theorem~\ref{thm:uniform}. Let
$X\subset Q\times R$ be a positive-dimensional smooth geometrically
integral very general complete intersection of hypersurfaces of
nonnegative multidegrees, of codimension $c$.

Suppose the distinguished row is at least $(e_1,\ldots,e_r;4,\ldots,4)$,
where $e_j\ge2$ for $a_j\le2$ and $e_j\ge4$ for $a_j\ge3$.
Suppose every other row is relatively ample over $Q$. Let
$s_1=\#\{j:a_j=1, e_j<4\}$, $s_2=\#\{j:a_j=2, e_j<4\}$,
$b_0=\max\{2,s_2+\floor{s_1/2}\}$.
Let $u_0=\max\{0,s_1-2(b_0-s_2)\}$, $D=2^{u_0}3^{b_0}$, and
$\epsilon=\delta_{\cha k=2}$.
If $B=\sum_\ell N_\ell-c+1\ge1$ and $\dim X\ge2^\epsilon D-2$,
then $X_{\bar k}$ admits no integral decomposition of the diagonal.
\end{theorem}
\begin{proof}
Let $a=\sum_ja_j$ and $N=\sum_\ell N_\ell$.
Use the $s_1$ one-coordinate groups and $s_2$ two-coordinate groups
in Lemma~\ref{lem:quartic-singleton-profiles}.
The hypothesis gives $a+B=\dim X+1\ge2^\epsilon D-1$, so that lemma
gives an affine hypersurface in $a+B$ coordinates.
Its equation has degree at most two in each of the designated groups
and total degree at most four. Distribute its other coordinates among
the other first groups and the last group. Its group degrees are then
at most $(e_1,\ldots,e_r;4)$.

Place the $a_j$ coordinates of the $j$th first group in the corresponding
factor of $Q$. Since $1\le B=N-c+1\le N$, the $B$ last coordinates fit
among the affine coordinates of $R$.
Splitting the last group among these factors does not increase its
degree in any factor. Exactly $N-B=c-1$ coordinates of $R$ are unused.
Every other row has positive degree in every factor of $R$, so assign
one unused coordinate to each of these rows. On the affine chart, take
its equation to be that coordinate. Eliminating these equations
recovers the reference hypersurface. Corollary~\ref{cor:embedding}
proves the conclusion.
\end{proof}

\begin{remark}
The count $2^{u_0}3^{b_0}$ is minimal among the ordinary Segre
embeddings used in Lemma~\ref{lem:quartic-singleton-profiles}.
For the integers $s,t$ in that lemma, the constraints on the factors are
$b\ge\max\{2,t\}$ and
$u+2(b-t)\ge s$. Replacing two line factors by one plane preserves
the number of available singleton coordinates and lowers their
contribution to the Segre count from four to three.
With $b\ge2$, the least count occurs at $b=b_0$ and $u=u_0$.
The conclusion $B_d=\infty$ for $d\ge4$ removes the dimension
restriction in Theorem~\ref{thm:uniform} under its degree hypotheses.
It does not give optimal degree thresholds for arbitrary products.
\end{remark}

\setcounter{theorem}{30}
\begin{theorem}\label{thm:unequal}
A very general hypersurface of bidegree at least $(3,2)$ in
$\PP^A\times\PP^D$ has no geometric decomposition if
$4\leq A\leq D\leq\binom{A+1}{3}-2A\delta_{\cha k=2}
                    +3\delta_{\cha k=2}$.
Outside characteristic two $(A,D)=(3,4),(3,5)$ are also allowed.
\end{theorem}
\begin{proof}
In characteristic two, Proposition~\ref{prop:norm-extension} gives
every pair in the indicated range. Applying
Corollary~\ref{cor:embedding} proves the conclusion in this case.

Suppose now that the characteristic is different from two.
Place the $A$ base coordinates of $G_A$ in the first group and its
$A+1$ fiber coordinates in the second.
For each unused support of size at most three, adjoin its free square
coordinate. The new term has degree at most three in the first group
and two in the second. There are $S_3(A)-(2A+2)$ such supports.
Adding their coordinates one at a time gives every second-group size
from $A+1$ to
$$
 A+1+S_3(A)-(2A+2)=\binom{A+1}{3}.
$$
Proposition~\ref{prop:equal-profiles} with $r=2$, $N=A$, and $c=1$
gives the pair $D=A$.
For $A=3$, the equations $W_6,W_7$ of
Proposition~\ref{prop:quartics} give $D=4,5$.
Lemma~\ref{lem:gram-local} preserves the affine obstruction for the
intermediate coordinate counts, and every equation has bidegree at
most $(3,2)$. Corollary~\ref{cor:embedding} proves the result.
\end{proof}

\begin{remark}\label{rmk:quartic-torsion}
In characteristic two, a very general quartic threefold or fourfold
also satisfies
$$6\mid\Tor(X_{\bar k})\mid24.$$
Apply the modulo-two models $R_i$ and the modulo-three models
$C_i$ of Propositions~\ref{prop:binary-models} and
\ref{prop:ternary-models} as separate degenerations of the same
very general quartic.
Corollary~\ref{cor:embedding} and Lemma~\ref{lem:comparison} give
divisibility by both primes, hence by six. The upper bound $24=4!$ is
\cite[Proposition~5.2]{CL}.
\end{remark}

\section{Proofs of the main results}\label{sec:proofs}

We use Corollary~\ref{cor:embedding} to place the affine constructions
in the ambient products. A complete intersection of dimension $n$ given
by $q$ equations uses $n+q$ coordinates. We match its equations to distinct
defining rows and its unused coordinate slots to the other rows.
Monic elimination of the graph equations preserves the affine scheme
and its generic-point class. A nonzero class then excludes a geometric
decomposition of the diagonal by Lemma~\ref{lem:comparison}, and hence
stable and retract rationality.

\begin{proof}[Proofs of Theorems~\ref{thm:quadratic} and~\ref{thm:cubic}
and Corollaries~\ref{cor:bidegrees} and~\ref{cor:rationality}]
Let $c$ be the codimension of $X$. If every column sum is at least
$N+1$, adjunction gives a globally generated canonical bundle with a
nonzero section. Then \cite[Lemma~2.2]{Totaro} excludes an integral
decomposition of the diagonal. This includes $c\geq N$, since ampleness
and the distinguished row make every column sum at least $c+1$.

Otherwise $c<N$. Proposition~\ref{prop:equal-profiles} gives an integral
affine variety $V$ of dimension $rN-c$ with $\delta_V\ne0$ and $q\leq c$
equations. Its distinguished equation fits the distinguished row, and
its other $q-1$ equations have degree at most one in each group.
The $rN-c+q$ coordinates leave $c-q$ unused slots. Match the unused
ample rows to these slots and use Corollary~\ref{cor:embedding}.
This proves both theorems. If $d,e\geq2$ and $d+e\geq5$, permuting the
two factors gives $(d,e)\geq(3,2)$. Theorem~\ref{thm:cubic} with $c=1$
then gives Corollary~\ref{cor:bidegrees}.

For Corollary~\ref{cor:rationality}, suppose $d_j=1$ and project away
from the $j$th factor. Over the rational function field of the product
of the other factors, the defining equation is a nonzero linear equation
in $\PP^N$, so its zero locus is $\PP^{N-1}$. Thus $\bar k(X_{\bar k})$ is
purely transcendental over $\bar k$, and $X_{\bar k}$ is rational.
Rationality implies stable rationality and a decomposition of the
diagonal. If every $d_j\geq2$, Theorem~\ref{thm:quadratic} with $c=1$
excludes a decomposition of the diagonal, proving the converse.
\end{proof}

\begin{proof}[Proof of Theorem~\ref{thm:uniform}]
Let $c$ be the codimension, $N=\sum_\ell N_\ell$, and $B=N-c+1$.
Then $\dim X=a+B-1$, where $a=\sum_j a_j$.
If $B\leq d-1$, every $Q$ column sum is at least $e_j\geq a_j+1$,
and every $R$ column sum is at least
$$
 d+c-1=d+N-B\geq N+1\geq N_\ell+1.
$$
Adjunction makes the canonical bundle globally generated, so it has a
nonzero section. By \cite[Lemma~2.2]{Totaro}, there is no decomposition
of the diagonal. This includes $B\leq0$ and the empty product $R=\operatorname{Spec}k$,
for which $B=1-c\leq0$.

Suppose $B\geq d$. We first treat $d=3$. If $p\ne2$,
Proposition~\ref{prop:binary-bound} gives every
$3\leq B\leq2^{a+1+g}+g$ with $g=\min(q,1)$.
For $p\ne3$, Lemma~\ref{lem:exact-norm} gives every
$3\leq B\leq T_3$ with $g=0$.
Since $B\leq B_3$, at least one of these constructions gives an affine
model of dimension $a+B-1$ with $B+g$ last coordinates and $g+1$
equations. Its distinguished equation has degrees at most
$(e_1,\ldots,e_r;3)$, and its $g$ other equations have degrees at most
$(0,\ldots,0;2)$, with $g\leq q$.
For $d\geq4$, Proposition~\ref{prop:quartic-all-capacities} gives such a
model with $g=0$ and last degree at most four.
In each case the model is integral over an algebraically closed field
of characteristic $p$, and its affine generic-point class is nonzero
with the coefficients of the construction.

Place the first $a_j$ coordinates in the corresponding factor of $Q$.
Since $g\leq q\leq c-1$, the $B+g$ last coordinates fit among the $N$
affine slots of $R$. Splitting this group does not increase its degree
in any factor, so the distinguished equation fits its row.
The $g$ additional equations fit distinct rows counted by $q$.
There remain $N-(B+g)=c-g-1$ unused slots, all in $R$.
Relative ampleness allows the other rows to be matched to these slots.
On the affine chart their equations are the corresponding coordinate
functions, whose elimination recovers the affine model.
Corollary~\ref{cor:embedding} proves the conclusion. 
\end{proof}

\begin{proof}[Proofs of Corollaries~\ref{cor:hypersurfaces}
and~\ref{cor:quartic-products}]
For every $n\geq1$, there is an integral affine hypersurface
$V\subset\bA^{n+1}$ of degree at most four over an algebraically closed
field of the ground characteristic, with nonzero generic-point class
for suitable coefficients. For $n=1,2$, use adjunction and
Lemma~\ref{lem:forms}. In characteristic two and $n\geq3$, use the affine
models and generic affine quartics in the proof of
Theorem~\ref{thm:char2-all-quartics}. Outside characteristic two,
Proposition~\ref{prop:quartics} gives $3\leq n\leq8$, and
Lemma~\ref{lem:quartic-singleton-profiles} with $s=t=0$ gives $n\geq7$.
For an ordinary hypersurface of degree $d\geq4$, homogenizing its affine
equation to degree $d$ and using Corollary~\ref{cor:embedding} proves
Corollary~\ref{cor:hypersurfaces}.

For Corollary~\ref{cor:quartic-products}, let $c$ be the codimension and
$n=\sum_jN_j-c$. Distribute the $n+1$ coordinates of $V$ among the affine
slots of the product. Its equation has degree at most four in each group
and fits the distinguished row. Match the $c-1$ unused slots to the
other ample rows. Corollary~\ref{cor:embedding} proves the result.
\end{proof}

\printbibliography
\end{document}